\documentclass[12pt]{article}
\usepackage{amsmath}
\usepackage{graphicx}
\usepackage{enumerate}
\usepackage{natbib}
\usepackage{url} 

\usepackage{mathrsfs}
\usepackage{amsmath,amssymb}
\usepackage{bm}
\usepackage{natbib}
\usepackage[usenames]{color}
\usepackage{amsthm}
\usepackage{arydshln}

\usepackage{multirow} 
\usepackage{enumitem}
\usepackage{caption}
\usepackage{subcaption}
\usepackage{enumitem}
\usepackage{dsfont}
\usepackage{mathtools}
\usepackage{caption}
\usepackage{booktabs}

\usepackage[colorlinks,
linkcolor=red,
anchorcolor=blue,
citecolor=blue
]{hyperref}

\usepackage{mylatexstyle}

\newcommand{\blind}{1}

\usepackage{xcolor}

\ifdefined\final
\usepackage[disable]{todonotes}
\else
\usepackage[textsize=tiny]{todonotes}
\fi
\allowdisplaybreaks

\begin{document}

\def\spacingset#1{\renewcommand{\baselinestretch}%
{#1}\small\normalsize} \spacingset{1}


\if1\blind
{
  \title{\bf Estimation of Out-of-Sample Sharpe Ratio for High Dimensional Portfolio Optimization}
  \author{Xuran Meng
    \\
    Department of Biostatistics, University of Michigan\\   
    Yuan Cao \\
    Department of Statistics and Actuarial Science, University of Hong Kong\\
    and\\
    Weichen Wang\\
   Faculty of Business and Economics, University of Hong Kong}
  \maketitle
} \fi

\if0\blind
{
  \bigskip
  \bigskip
  \bigskip
  \begin{center}
    {\LARGE\bf Estimation of Out-of-Sample Sharpe Ratio for High Dimensional Portfolio Optimization}
\end{center}
  \medskip
} \fi

\bigskip
\begin{abstract}
Portfolio optimization aims at constructing a realistic portfolio with significant out-of-sample performance, which is typically measured by the out-of-sample Sharpe ratio. However, due to in-sample optimism, it is inappropriate to use the in-sample estimated covariance to evaluate the out-of-sample Sharpe, especially in the high dimensional settings. In this paper, we propose a novel method to estimate the out-of-sample Sharpe ratio using only in-sample data, based on random matrix theory. Furthermore, portfolio managers can use the estimated out-of-sample Sharpe as a criterion to decide the best tuning for constructing their portfolios. Specifically, we consider the classical framework of Markowits mean-variance portfolio optimization under high dimensional regime of $p/n \to c \in (0,\infty)$, where $p$ is the portfolio dimension and $n$ is the number of samples or time points. We propose to correct the sample covariance by a regularization matrix and provide a consistent estimator of its Sharpe ratio. The new estimator works well under either of the following conditions: (1) bounded covariance spectrum, (2) arbitrary number of diverging spikes when $c < 1$, and (3) fixed number of diverging spikes with weak requirement on their diverging speed when $c \ge 1$. We can also extend the results to construct global minimum variance portfolio and correct out-of-sample efficient frontier.
We demonstrate the effectiveness of our approach through comprehensive simulations and real data experiments. Our results highlight the potential of this methodology as a useful tool for portfolio optimization in high dimensional settings.
\end{abstract}

\noindent%
\vfill

\newpage
\spacingset{1.9} 

\section{Introduction}
The mean-variance (MV) portfolio, first introduced by  \citet{markowitz1952harry}, is a fundamental cornerstone in modern portfolio allocation theory. Suppose that the excessive returns of $p$ risky assets have population excess mean vector $\bmu\in\RR^p$ and population covariance matrix $\bSigma\in\RR^{p\times p}$.
In mean-variance portfolio optimization, one solves for the allocation vector $\wb$ that maximizes $\wb^\top \bmu - \lambda \wb^\top \bSigma \wb$, where $\lambda$ is the risk aversion. The solution $\wb^*$ of the MV portfolio optimization is proportional to $\bSigma^{-1}\bmu$. It is well known that this solution achieves the largest possible Sharpe ratio.

Despite the clean framework of the mean-variance portfolio, the key challenge in its practical application is that, the population mean $\bmu$ and population covariance matrix $\bSigma$ are unknown and have to be estimated using historical data. This leads to the phenomenon known as ``empirical optimism'', which states if we use the in-sample data both for estimation of $\bSigma,\mu$ and for out-of-sample assessment, we are overly optimistic about the true out-of-sample performance of the estimated portfolio \citep{Karoui2010}.  For example, if we have $n$ iid $p$-dimensional return vectors $\Rb_i \sim N(\bmu,\bSigma), i=1,\dots,n$, assuming $\bmu$ is known and $p/n \to c < 1$, the optimized MV portfolio based on the sample covariance $\hat\bSigma = \frac1n \sum_{i=1}^n (\Rb_i - \bmu)(\Rb_i - \bmu)^\top$, i.e. $\wb^*=C \hat\bSigma^{-1}\bmu$  for some  $C>0$, has the in-sample Sharpe ratio of $\sqrt{\bmu^\top \hat\bSigma^{-1}\bmu}$ while its out-of-sample Sharpe ratio is $\bmu^\top \hat\bSigma^{-1}\bmu / \sqrt{\bmu^\top \hat\bSigma^{-1} \bSigma \hat\bSigma^{-1}\bmu}$ (see \eqref{eq:def_SRQ} for Sharpe ratio calculation).  According to our Theorem \ref{thm:main_theorem} (with $\Qb=\0$), we will see that the latter is approximately $(1-c) \sqrt{\bmu^\top \hat\bSigma^{-1}\bmu}$, so that the in-sample Sharpe ratio is $1/(1-c)$ times larger than the actual Sharpe ratio. When $c$ is close to $1$, the portfolio performance will be significantly exaggerated. 
Through this example, we see the importance of developing a valid and accurate estimator to genuinely evaluate the out-of-sample Sharpe ratio using historical data. In this paper, we will propose such an estimator to help portfolio managers to achieve the desired portfolio performance.

 The example above  highlights the importance of addressing estimation errors in the population parameters. Although abundant literature studied the estimation of $\bmu$ \citep{green2013supraview, fischer2018deep,baba2020machine, petropoulos2022forecasting}, which is challenging as the signal level of financial data is weak, we will mainly focus on the issue caused by estimating the high-dimensional population covariance matrix $\bSigma$, i.e. the curse of high dimensionality.
When $\bSigma$ is an identity matrix, it is widely known that the sample covariance matrix $\hat\bSigma$ is a poor estimator of $\bSigma$ when the dimension of the matrix and the sample size increase to infinity proportionally---in this setting, the spectral distribution of $\hat\bSigma$ follows the MP-Law \citep{marchenko1967distribution, bai2009enhancement}, which is drastically different from that of true $\bSigma$.  Therefore, in high dimensional mean-variance portfolio optimization, we would not directly plug in the unreliable sample covariance matrix. As a remedy, in this paper, we consider adding a ridge-type of penalty to the MV optimization to alleviate the big estimation error from the sample covariance. 
Note that we care about covariance estimation via the specific function of Sharpe ratio, which is different from the majority of the literature on estimating covariance matrix itself, whose performance is typically measured using various matrix norms.

\subsection{Literature review}
Estimation errors in $\bmu$ and $\bSigma$ result in a substantial deterioration in out-of-sample portfolio performance, which is typically measured by  Sharpe ratio. 
Effort invested in improving the out-of-sample Sharpe ratio associated with Markowitz portfolio has a long history \citep{brown1976optimal,klein1976effect,jobson1981performance,craig2000asset}. More recently,  \citet{kan2007optimal} introduced a three-fund rule which can improve the out-of-sample portfolio performance  in the case where a risk-free asset
is available. 
\citet{demiguel2009optimal} empirically demonstrated that the ``$1/N$ rule'' consistently outperforms various optimally estimated portfolios across a wide range of datasets. \citet{tu2011markowitz} enhanced the three-fund rule by incorporating the ``$1/N$ rule''.

A second approach to enhance the out-of-sample Sharpe ratio involves integrating various regularization techniques.  \citet{bai2009enhancement,Karoui2010} found that using in-sample data for both creating and evaluating portfolios leads to an overly optimistic assessment of their true out-of-sample performance, particularly as the number of assets grows.  Recognizing the need for regularization, several studies \citep{jagannathan2003risk,fan2012vast,bruder2013regularization, ao2019approaching} have highlighted its significance in preventing the deterioration of portfolio performance with increasing number of assets.  \citet{jagannathan2003risk} empirically showed the advantage of no-short-sale constraints when estimating optimal portfolios, with \citet{fan2012vast} providing theoretical support for this strategy. Additionally, \citet{bruder2013regularization} reviewed various regularization methods that help stabilize mean-variance allocations, including weight constraints, resampling methods, and shrinkage techniques. 
\citet{demiguel2009generalized} considered a generalized approach to portfolio optimization by constraining different portfolio norms. More recently, 
\citet{li2022spectrally} gave a spectral correction to the sample covariance matrix to more accurately estimate the out-of-sample Sharpe ratio when  $p<n$. Additionally, \citet{kan2024sample} explored the asymptotic distribution of both in-sample and out-of-sample Sharpe ratios under the same condition with \citet{li2022spectrally}.
\citet{ao2019approaching} applied a Lasso penalty in portfolio optimization and estimated the maximal Sharpe ratio considering scenarios where the ratio $p/n$ approaches a constant $c \in (0,1)$.

Another line of research focuses on directly estimating the population covariance matrix. Current methods for high-dimensional population covariance estimation primarily utilize either truncation or shrinkage approaches. Truncation-based methods generally rely on a specific sparsity structure within the population covariance or precision matrix, aiming for more precise estimations by exploiting this sparsity \citep{bickel2008regularized, lam2009sparsistency, cai2011adaptive,fan2011high,fan2016overview,cai2020high}. These methods tend to achieve faster convergence rates or enhanced empirical performance, particularly when the sparsity assumption is valid.
In contrast, shrinkage-based methods adjust the spectrum of the sample covariance matrix to improve the estimation of the population covariance without relying on sparsity \citep{ledoit2004well, ledoit2012nonlinear, bodnar2022recent,bodnar2024reviving}. These techniques are largely driven by the bias-variance trade-off, aiming to minimize estimation variance and overall error by introducing a slight bias. For example, \citet{ledoit2004well} and \citet{ledoit2012nonlinear} developed a linear and nonlinear shrinkage estimator on the empirical eigenvalue spectrum that offers improved conditioning and accuracy compared to the traditional sample covariance matrix. In addition to correcting sample covariance matrix, factor models and its low-rank plus sparse covariance structure are also extensively studied in the literature.
Readers may refer to \citet{fan2021robust} to see an overview of factor models and their applications. 

\subsection{Our contributions}


In this paper, we investigate the MV portfolio optimization with ridge regularization. Specifically, under the setting where the unknown population covariance matrix is estimated by the sample covariance, we develop a novel method for estimating the out-of-sample Sharpe ratio for the ridge penalized MV portfolio. Our method can be applied to both the settings with and without the risk-free asset: 
with the risk-free asset, our approach directly proposes an estimator for the out-of-sample Sharpe ratio of the regularized MV portfolio; without the risk-free asset, our approach is given in the form of estimating the efficient frontier, that is, for any given target return, we estimate the variance of the regularized MV portfolio. 
Remarkably, the estimators given in our paper are consistent---while they can be calculated without the knowledge of the population covariance matrix, they still share the same asymptotic behavior as if they are calculated with the knowledge of the population covariance. 
We summarize our contributions as follows:
\begin{enumerate}[leftmargin = *]
    \item We introduce an innovative in-sample method that estimates the true out-of-sample variance and thus the true out-of-sample Sharpe ratio of a portfolio, considering both scenarios with and without the risk-free asset. With no risk-free asset, we actually provide a way to accurately estimate the efficient frontier.  
    Our method is backed by theoretical guarantees of almost sure convergence.
    \item Our method assists in the selection of optimal regularization parameters from a predetermined set.
    By providing valid estimations for the future performances of portfolios under different regularization levels, our method allows data-driven decision making for the regularization parameters, which can be particularly useful in situations where market conditions change rapidly and investment decisions must be made timely.
    \item  Our theoretical convergence results allow much relaxed assumptions on the population covariance matrix $\bSigma$. Traditional assumptions in random matrix theory require $\|\bSigma\|_{\op}$ to be finite, while our framework allows $\|\bSigma\|_{\op}$ to be unbounded with diverging eigenvalues. This makes our proposed estimators generally applicable to high-dimensional factor models with either strong or weak factors. This is achieved by rigorously handling the Stieltjes transform as $z \to 0$ and proving the interchangeability of the limits $z \to 0$ and $n \to +\infty$. Our method remains effective even with diverging $\|\bSigma\|_{\op}$, provided that either $c < 1$ or the number of diverging spikes is fixed with weak requirement on the diverging speed.
\end{enumerate}

We hope to clarify that this work is neither meant to compare different portfolio optimization settings (e.g., ridge penalty versus Lasso penalty, or Value-at-risk loss versus MV objective), nor to study the optimal design of the ridge regularization (although we will comment on this in Section~\ref{sec:mainresult}). Instead, we study the exact asymptotic behavior of adding the Ridge penalty to the classical mean-variance optimization so that our asymptotic results can be used for tuning the parameters in the Ridge penalty. This is in contrast to non-asymptotic study of concentration bounds, where typically only the proper rates for the tuning parameters are derived leaving the effect of their constants undetermined. Therefore, those results cannot be used to tune the tuning parameters. 
For example, \cite{ao2019approaching} applied Lasso regularization in the portfolio optimization. Our work differs from theirs in a few salient aspects. Firstly, we consider ridge regularization under more general and relaxed conditions allowing diverging spikes in the population covariance matrix with arbitrary orders. 
Moreover, they consider holding factor portfolios besides individual stocks. In contrast, since factor portfolios may not be always accessible, we only allocate portfolio among stocks. 
Most importantly, they only consider $p<n$ and give the estimation for the maximal Sharpe ratio, while we provide valid estimations for any regularization level, so that we can select the optimal tuning. It is unclear how the Lasso tuning affects Sharpe ratio in \citet{ao2019approaching}.


\subsection{Paper organization and notations}
\label{subsec:notation}
The remainder of the paper is organized as follows. We first introduce some notations below. In Section~\ref{sec:sharperatio} and Section~\ref{sec:frontier}, we give the problem settings for the regularized MV portfolio optimization with and without risk free assets, and propose our methodology to estimate the Sharpe ratio respectively. We then present simulation results in Section~\ref{sec:simulation} to support our theoretical findings. Section~\ref{sec:realdata} contains our demonstration about the effectiveness of our method for optimizing the MV portfolio using real financial return data. Finally, we conclude the paper in Section~\ref{sec:conclusion} and discuss potential future research. Additional experiments and proof details are provided  in the supplementary material.

 We employ lower case letters (e.g. $a$), Boldface letters (e.g. $\ab$ and $\Ab$) for scalar values, vectors and matrices. For a set $I$ of indices, we define $\one(I)$ such that $\one(i) = 1$ if $i \in I$, and $\one(i) = 0$ if $i \notin I$. Hence, $\ab(i)$ specifically refers to the $i$-th element in the vector $\ab$. Moreover, the vector $\1\in\RR^p$ ($\1_n\in\RR^n$) denotes the vector with all coordinates equal to 1. 
For a square matrix $\Ab$, we use $\lVert\Ab \rVert_{\tr}$ and  $\lVert\Ab\rVert_{\op}$ to denote its trace norm and operator norm, respectively, and use $\tr(\Ab)$ to denote its trace.  $\Ab\leq (\geq) \Bb$ means $\Bb-\Ab$ is non-negative(positive) definite. 
The sets of natural, real and complex numbers are denoted by $\NN$,  $\RR$ and $\mathbb{C}$, respectively. 
The set of integers from  $n_1$ to $n_2$ is denoted by $[n_1:n_2]=\{n_1,\ldots,n_2\}$ and $[n]=[1:n]$.

\section{High Dimensional Estimation of Sharpe Ratio}
\label{sec:sharperatio}
\subsection{Mean-variance portfolio and Sharpe ratio}
\label{sec:problemsetup}
In this section, we introduce the problem setting we consider in this paper. We first review the MV portfolio optimization under the existence of the risk-free rate $r_0 \ge 0$.

\begin{definition}[Mean-variance portfolio optimization]
\label{def:meanvarianceoptimization}
Given $p$ risky assets with mean $\rb\in\RR^p$ and covariance matrix $\bSigma\in\RR^{p\times p}$, the MV optimization optimizes the allocation vector $\wb$, which satisfies
\begin{align*}
    \wb^*=\argmin\limits_{\wb\in\RR^p} \wb^\top \bSigma \wb,\quad \st \wb^\top \bmu =\mu_0.
\end{align*}
Here, we denote by $\bmu=\rb-r_0\1 $ the excess return of the risky assets, and $\mu_0>0$ is the target excess return of the portfolio.
\end{definition}
The motivation behind Definition~\ref{def:meanvarianceoptimization} is clear. Each allocation vector $\wb$  represents a unique portfolio, and the optimal $\wb$ in Definition~\ref{def:meanvarianceoptimization} minimizes the variance of the portfolio with an expected excess return. It is worth noting that in the MV optimization, it is not required to have $\wb^\top \one=1$, because the existence of a risk-free asset makes investing or borrowing risk-less bond possible. To be exact, we could allocate $w_0 = 1 - \wb^\top \one$ to the risk-free asset, so that the total mean excess return is indeed $(\wb^\top \rb + w_0 r_0) - r_0 = \wb^\top \bmu$.

The solution in Definition~\ref{def:meanvarianceoptimization}    satisfies $\wb^* \propto \bSigma^{-1} \bmu$.  Here, vector $\ab\propto\bbb$ means that $\ab=C\bbb$ for some $C>0$. 
In practice, a portfolio manager needs to estimate $\bSigma$ in order to find out the optimal mean-variance portfolio $\wb^*$. However, the estimation of $\bSigma$ is very challenging due to its huge number of unknown parameters. To further investigate the effect of estimation error of $\bSigma$, we assume the following data generation process.
\begin{assumption}
    \label{def:X}
Given the mean return vector $\rb\in\RR^p$ and the covariance matrix $\bSigma\in\RR^{p\times p}$, we consider the observed data matrix $\Rb\in\RR^{n\times p}$ with the form
\begin{align*}
    \Rb=\one_n \rb^\top+\Xb,
\end{align*}
where $\Xb=\Zb\bSigma^{\frac{1}{2}}\in\RR^{n\times p}$. Here, the elements in the matrix $\Zb\in\RR^{n\times p}$ are i.i.d random variables with zero mean, variance 1 and finite $(8+\varepsilon)$-order moment for some $\varepsilon>0$.
\end{assumption}
%
For now, we assume the knowledge of $\bmu$ (or equivalently $\rb$) \footnote{ We first consider a known $\bmu$ because (1) it makes the theoretical analysis more transparent as it separates the randomness from estimating $\bmu$; (2) it reflects that hedge fund portfolio managers usually use many independent alternative datasets and complicated machine learning techniques to estimate $\bmu$ independent of the estimation of $\bSigma$. However, we also consider the out-of-sample SR estimation when $\bmu$ is unknown and estimated by the historical sample mean; see our Remark~\ref{remark:unknown_mu} and Appendix~\ref{sec:unknown_mu} for details.}, and focus only on the effect of estimation error from estimating $\bSigma$ using the sample covariance $\hat \bSigma = \Xb^\top \Xb / n$,  where $\Xb$ is defined in Assumption~\ref{def:X}. 
It has been widely known that when dimension $p$ is large with $p/n \to c>0$, that is $p$ grows with the sample size linearly, the sample covariance matrix is not a good estimator for the true covariance matrix.
To mitigate the bad condition number of the sample covariance, it is natural to consider adding regularization in the portfolio construction.
We add a regularization term $ \Qb$ to $\hat\bSigma$ in the MV portfolio. Therefore, we consider the following optimization:
\begin{align*}
    \wb^*=\argmin\limits_{\wb\in\RR^p} \wb^\top (\hbSigma+\Qb) \wb,\quad \st \wb^\top \bmu =\mu_0.
\end{align*}
The optimal $\wb^*$ satisfies  $\wb^* \propto(\hat\bSigma + \Qb)^{-1} \bmu$. 
Here, $\Qb$ is a positive definite matrix. The Sharpe Ratio (SR) of $\wb^*$ is defined as the mean of the excessive return of the portfolio over its standard deviation. Mathematically,
\begin{align}
SR(\Qb) = \frac{\EE_{\tilde\Rb}[\wb^{*\top} (\tilde\Rb-r_0\1)]}{\sqrt{\Var_{\tilde\Rb}[\wb^{*\top} (\tilde\Rb-r_0\1)]}} = \frac{\bmu^\top (\hat\bSigma + \Qb)^{-1} \bmu}{\sqrt{\bmu^\top (\hat\bSigma + \Qb)^{-1} \bSigma (\hat\bSigma + \Qb)^{-1} \bmu}}.\label{eq:def_SRQ}
\end{align}
Here, $\tilde\Rb$ is an out-of-sample data point with mean return $\rb=\bmu+r_0\1$ and covariance $\bSigma$. Note that the randomness from historical data in $\wb^{*}$ is treated as given in the Sharpe computation in \eqref{eq:def_SRQ}.
Define the following quantities:
\begin{align}
    T_{n,1}(\Qb)&=\tr \bigg[ \bigg(\frac{\Xb^\T\Xb}{n}+\Qb\bigg)^{-1}\Ab \bigg],\label{eq:Tn1}\\
    T_{n,2}(\Qb)&=\tr \bigg[ \bigg(\frac{\Xb^\T\Xb}{n}+\Qb\bigg)^{-1}\bSigma\bigg(\frac{\Xb^\T\Xb}{n}+\Qb\bigg)^{-1}\Ab \bigg],\label{eq:Tn2}
\end{align}
for some deterministic matrix $\Ab \in \RR^{p \times p}$.
Hence, by setting $\Ab=\bmu \bmu^\top$, it is easy to see that the Sharpe ratio in \eqref{eq:def_SRQ} could be expressed as 
\begin{align}
\label{eq:definition_SR}
    SR(\Qb)=\frac{T_{n,1}(\Qb)}{\sqrt{|T_{n,2}(\Qb)|}}.
\end{align}

\begin{remark}
We introduce a deterministic matrix $\Ab \in \RR^{p \times p}$ to obtain results more general than  $\Ab=\bmu \bmu^\top$. Here we have the underlying assumption that the mean vector of financial returns keeps unchanged over the time period under consideration. However, 
if daily excess returns have time-varying means, denoted as $\bmu_i$, where $i$ represents the trading date, and if the daily returns are independent across dates, the annual out-of-sample Sharpe ratio of the daily-balanced portfolio 
can be expressed via setting $\Ab = \sum_i\bmu_i\bmu_i^\top$ over all trading dates.
\end{remark}

To obtain the largest out-of-sample Sharpe ratio, one might select the regularization matrix $\Qb$ by a cross-validation procedure, but this is on the one hand computationally expensive, on the other hand not straightforward to do for a time series setting in practice. Note that although for theoretical analysis, we assume data are independent from time to time, but in practice, it is not recommended to shuffle time order in the computation of Sharpe ratio.
In this work, we hope to propose a method to directly estimate Sharpe ratio using in-sample data, which can be used as a good criterion to select the best regularization $\Qb$ to guarantee good out-of-sample Sharpe ratio of the regulated MV portfolio. 

\subsection{Out-of-sample estimation of Sharpe ratio}
\label{sec:mainresult}
In this section, without any information of $\bSigma$, we propose a way to effectively estimate the Sharpe ratio $SR(\Qb)$ in \eqref{eq:definition_SR} with in-sample data. We need the following assumptions. 
\begin{assumption}
\label{ass:assump1}
The portfolio dimension $p$ and the sample size $n$ both tend to infinity, with the regime of $p/n\to c$ where $c>0$ is a constant.
\end{assumption}

Assumption \ref{ass:assump1} provides a classic framework for high dimensional analysis, commonly used in random matrix theory \citep{yao2015sample}. Next, we make assumptions about the matrix $\Ab$ and the structure of the regularization matrix $\Qb$.
\begin{assumption}
\label{ass:assump2}  
The matrices $\Ab$ is deterministic.  The  positive definite regularization matrix $\Qb \in \mathcal{Q}$, independent of $\Xb$, satisfies the following: there exist constants $c_1, c_2 \geq 0$, with $c_1^2 + c_2^2 > 0$, and $\Qb_1$ and $\Qb_2$ such that $\Qb = c_1 \Qb_1 + c_2 \Qb_2$. The matrices $\Qb_1$ and $\Qb_2$ satisfy $c'\Ib \leq \Qb_1 \leq C'\Ib$ and $c'\Ib \leq \bSigma^{-1/2}\Qb_2\bSigma^{-1/2} \leq C'\Ib$ for some constants $c', C' > 0$, for any sequence $(n, p)$. In addition,  we allow the  case $\Qb=\zero$ when $c<1$.
\end{assumption}
Assumption~\ref{ass:assump2} says that it makes no sense to apply an ill-conditioned regularization matrix $\Qb$. For instance, we can choose $\Qb$ with bounded eigenvalues by setting $c_1 = 1$, $c_2 = 0$, $\Qb_1 = \Qb$. $\Qb$ can also be proportional to $\bSigma$, even when $\|\bSigma\|_{\op}$ is unbounded, by setting $c_1 = 0$, $c_2 = 1$, and $\Qb_2 = \bSigma$. We also allow no regularization $\Qb=\zero$ when $c<1$.  
\begin{assumption}
\label{ass:assump3}
$\bSigma$ is well scaled as $\|\bSigma/p\|_{\tr}\leq C$ for some constant $C>0$. Denote by $\lambda_1\geq \cdots\geq \lambda_{p}$ the eigenvalues of $\bSigma$.  It is assumed that $\lambda_p\geq c_1$ for some constant $c_1>0$. Moreover,  one of the following cases must hold:
\begin{enumerate}
    \item (Bounded spectrum) There exists  $C_1>0$ such that $ \lambda_1\leq C_1$. 
    \item (Arbitrary number of diverging spikes when $c<1$) $p/n \to c<1$ and we allow arbitrary number of top eigenvalues to go to infinity. 
    \item (Fixed number of diverging spikes when $c \ge 1$) $p/n \to c \ge 1$ and we let the number of diverging spikes be $K$, $K$ is fixed and $\lambda_1\leq C\lambda_K^2$ for some constant $C>0$. 
\end{enumerate}
\end{assumption}
Assumption~\ref{ass:assump3} is a very mild assumption and encompasses a wide range of structures of the covariance  $\bSigma$. Case 1 assumes that $\bSigma$ has bounded operator norm, which was typically adopted by RMT, such as \citet{li2022spectrally} and \citet{bodnar2024two}. Cases 2 and 3 correspond to the factor model structure, where $\bSigma$ has spiked eigenvalues that tend to infinity. 
Specifically, \cite{bai2003inferential} and \cite{fan2013large} considered a factor model with fixed $K$ number of strong or pervasive factors with $\lambda_1 \asymp \lambda_K \asymp p$, which falls in Case 2 or 3 depending on the dimensionality regime $p<n$ or $p > n$. \cite{onatski2012asymptotics} and \cite{wang2017asymptotics} considered semi-strong or weak factors, that is, $\lambda_k \asymp p^{\theta_k}$ for some $\theta_k \in (0,1)$ and $\theta_k$'s may not be the same. Our Case 2 can allow arbitrary $\theta_k$ when $c<1$ and our Case 3 can allow flexibility in setting $\theta_k$ to some degree ($\theta_K \le \theta_1 \le 2\theta_K$). 
In addition, Assumption~\ref{ass:assump3} actually allows unbounded diagonal elements of $\bSigma$. 
\citet{ao2019approaching} and the above factor model literature assume that the diagonal elements $\bSigma_{i,i}$ are uniformly bounded.  Here, we only require that the average $\sum_{i=1}^p \bSigma_{i,i}/p$ remains bounded.


\begin{theorem}
    \label{thm:main_theorem}
    Suppose Assumptions~\ref{def:X}, \ref{ass:assump1}, \ref{ass:assump2} and \ref{ass:assump3} hold. For any $\Qb\in\cQ$, a good estimator $\hat{SR}(\Qb)$ for $SR(\Qb)$ which is defined in \eqref{eq:definition_SR}, is given as follows. 
    \begin{align*}
        \hat{SR}(\Qb)=\frac{T_{n,1}(\Qb)}{\sqrt{\big|\hat{T}_{n,2}(\Qb)\big|}},\quad \text{where}\quad \hat{T}_{n,2}(\Qb)=\frac{\tr(\hat{\bSigma}+\Qb)^{-1}\hbSigma(\hat{\bSigma}+\Qb)^{-1}\Ab}{\big(1-\frac{c}{p} \tr\hbSigma(\hbSigma+\Qb)^{-1}\big)^2}.
    \end{align*}
If $\Ab$ is semi-positive definite, it holds that 
\begin{align*}
    \hat{SR}(\Qb)/SR(\Qb)~\cvas~1.
\end{align*}
If additionally  $\|\Ab\|_{\tr} $ is bounded\footnote{The condition can be further relaxed to bounded $\tr\big(\frac{\bSigma}{1+s_0}+\Qb\big)^{-1}\Ab $, where $s_0>0$ is defined by $s_0=\frac{c}{p}\tr~\bSigma\big( \frac{\bSigma}{1+s_0}+\Qb \big)^{-1}$ in Assumption~\ref{assump:assump5} below.}, then $SR(\Qb)$ is almost surely bounded and 
\begin{align*}
  \hat{SR}(\Qb)-SR(\Qb)~\cvas~0 .
\end{align*}
\end{theorem}
The proof of Theorem~\ref{thm:main_theorem} is given in Appendix~\ref{sec:appendix_proof_thm_bound} in the supplementary material.  The key in our proof is to extend a classical result (see Lemma C.1) in random matrix theory.
The first challenge in the extension arises from the relaxation of bounded $\|\bSigma\|_{\op}$ in Lemma~\ref{lemma:keylemma}. When spikes exist and $\|\bSigma\|_{\op}$ is no longer bounded, we cannot apply Lemma~\ref{lemma:keylemma} and have to extend Lemma~\ref{lemma:keylemma} to Lemmas~\ref{lemma:z=0converge} and \ref{lemma:z=0converge_c>1}, and prove that $z\to0$ and $ n\to+\infty$ can be interchanged. The second challenge arises from obtaining ratio consistency as well as difference consistency. As we see, ratio consistency has the advantage of circumventing additional scaling condition (i.e. bounded $\|\Ab\|_{\tr} $). 
To obtain ratio consistency, we establish Lemma~\ref{lemma:stronger_type_convergence} even when $\|\bSigma\|_{\op}$ is unbounded.

Theorem~\ref{thm:main_theorem} establishes an estimator for $SR(\Qb)$ with almost sure convergence guarantee in a wide range of scenarios. More importantly, it provides a  tractable framework for portfolio managers to determine the best regularization matrix $\Qb$ over some given candidate set $\cQ$.
Consider a candidate set $\cQ$ from which we want to select $\Qb\in\cQ$ in order to generate the maximal out-of-sample Sharpe ratio. It holds from Theorem~\ref{thm:main_theorem} that $\hat{SR}(\Qb)/SR(\Qb)~\cvas~1$, hence we can simply maximize $\hat{SR}(\Qb)$ over $\cQ$ to obtain the optimal $\Qb$. 
However, it is also important to note that  if we define $\hat{\Qb}=\arg\max_{\Qb}\hat{SR}(\Qb)$ without the constraint $\Qb\in\cQ$, $\hat{SR}(\hat{\Qb})$ may not align well with $SR(\hat{\Qb})$. This misalignment arises because optimizing $\hat{SR}(\Qb)$ over the entire $\Qb$ can overfit to the in-sample observed data, potentially inflating $\hat{SR}(\hat{\Qb})$ and  creating discrepancy between $\hat{SR}(\hat{\Qb})$ and $SR(\hat{\Qb})$.  Readers may refer to Appendix~\ref{sec:comparsion_optimal} for numerical results demonstrating this discrepancy. Nonetheless, optimizing \(\Qb\) over a finite deterministic candidate set \(\cQ\), or over \(\cQ\) parameterized by a finite number of flexible parameters (see Proposition~\ref{prop:optimal_Q_simple} below), can be practically useful.
\begin{proposition}\label{prop:optimal_Q_simple}
Given the conditions of Theorem~\ref{thm:main_theorem}, define
$\hat{\Qb}=\arg\max_{\Qb\in \cQ}\hat{SR}(\Qb)$. 
Suppose $\cQ$ has finite degree of freedom, under  assumptions  in Appendix~\ref{sec:optimalQ}, it holds  that
\[
\hat{SR}(\hat{\Qb})/SR(\hat{\Qb}) ~\cvas~ 1 \quad \text{and} \quad \hat{SR}(\hat{\Qb}) - SR(\hat{\Qb}) ~\cvas~ 0.
\]
\end{proposition}

Following Theorem~\ref{thm:main_theorem}, a natural question  is whether $SR(\Qb)$ can achieve or approximate the theoretical maximal Sharpe ratio $SR_{\max}$, where $SR_{\max} = \sqrt{\bmu^\top \bSigma^{-1} \bmu}$. 
In Appendix~\ref{sec:discussion_largest_SR}, we show that for any given $\varepsilon>0$, there always exists  $\Qb$ such that $1 - \varepsilon \leq SR(\Qb)/SR_{\max} \leq 1$ with probability one as $n\to+\infty$. 
The existence of $\Qb$ is proved by setting $\Qb = q \bSigma$ for some sufficiently large $q > 0$. When the covariance follows from a factor model, that is $\bSigma = \Bb\Bb^\top + \Db$ with the factor loading matrix $\Bb \in \mathbb{R}^{p \times K}$ and the residual covariance $\Db$, which is typically assumed as a diagonal matrix with eigenvalues bounded away from $0$ and $+\infty$, 
then the existence of $\Qb$ can also be proved by setting $\Qb$ proportional to $\Db$. 
An experienced and skillful portfolio manager may set an initial regularization matrix $\Qb_0$ as closely aligned (or proportional) to the diagonal matrix of residual variances as possible, say using some extra data to  estimate the residual variances. 
But sometimes it may be inevitable to suffer a loss in $SR(\Qb)$ from the maximal $SR_{\max}$ practically, as we only have partial information on $\bSigma$. 
Nonetheless, with any given $\Qb_0$, by Theorem~\ref{thm:main_theorem}, we can further explore $\Qb = q \cdot \Qb_0$ and fine-tune the scale (or shrinkage level) to achieve the best out-of-sample SR within this class of $\Qb$'s.



\begin{remark}
\label{remark:unknown_mu}
Theorem~\ref{thm:main_theorem} analyzes the out-of-sample Sharpe ratio under the assumption that $\bmu$ is known. When $\bmu$ is unknown and estimated by the sample mean $\hat{\bmu}$, the out-of-sample Sharpe ratio is given by $\frac{\hbmu^\top(\hbSigma+\Qb)\bmu}{\sqrt{\hbmu^\top(\hbSigma+\Qb)^{-1}\bSigma(\hbSigma+\Qb)^{-1}\hbmu}}$.  Under a slightly stronger assumption, where the data follow the Gaussian distribution, we also propose a corresponding estimator for the out-of-sample Sharpe ratio  $\hat{SR}(\Qb)=\frac{\hbmu^\top(\hbSigma+\Qb)^{-1}\hbmu-\frac{\tr(\hbSigma+\Qb)^{-1}\hbSigma}{n-\tr(\hbSigma+\Qb)^{-1}\hbSigma}}{\sqrt{\hbmu^\top(\hat{\bSigma}+\Qb)^{-1}\hbSigma(\hat{\bSigma}+\Qb)^{-1}\hbmu}}\cdot\big(1-\frac{c}{p} \tr\hbSigma(\hbSigma+\Qb)^{-1}\big)$. Detailed simulations, real data results, and theoretical justification on this estimator can be found in Appendix~\ref{sec:appendixA}, \ref{sec:appendixB} and  \ref{sec:unknown_mu}, respectively.
\end{remark}

\section{High Dimensional Estimation of Efficient Frontier}
\label{sec:frontier}
\subsection{Efficient frontier with no risk-free asset}
\label{sec:introefficientfrontier}
In this section, we give a brief introduction of the efficient frontier with no risk-free asset.
Suppose we have $p$ risky assets with mean $\rb\in\RR^{p}$ and covariance $\bSigma\in\RR^{p\times p}$, but no risk-free asset is available. Given target return $\mu_0>0$,  the portfolio optimization is given by
\begin{align}
\label{eq:frontier1}
    \wb^*=\arg\min_{\wb} \wb^\top \bSigma \wb, \quad \text{s.t. } \wb^\top \rb =\mu_0 \text{ and } \wb^\top \one=1.
\end{align}
Note that  $\mu_0$ is the target excess return in Section~\ref{sec:sharperatio},  while here it is the target raw return, with a slight abuse of notation. 
The solution $\wb^*$ from \eqref{eq:frontier1} is straightforward. 
However, in practice, $\bSigma$ is not observable and its sample estimation $\hbSigma=\Xb^\top\Xb/n$ can deviate significantly from the true $\bSigma$ in high dimensions. Similar to Section~\ref{sec:problemsetup},
the observed sample data $\Rb$ following Assumption~\ref{def:X} can be expressed as $    \Rb=\1_n\rb^\top+\Xb $.
Again we consider adding a ridge term $\Qb>0$ in the portfolio optimization problem:
\begin{align}
\label{eq:frontier2}
    \wb^*=\arg\min_{\wb} \wb^\top (\hbSigma+\Qb) \wb, \quad \text{s.t. } \wb^\top \rb =\mu_0 \text{ and } \wb^\top \one=1.
\end{align}
From \citet{merton1972analytic}, the optimal $\wb^*$ in \eqref{eq:frontier2} is given by 
$\wb^*=\gb+\mu_0\cdot\hb$, where
\begin{align}
\label{eq:def_gh}
\begin{array}{l@{~}l@{~}l@{~}l@{~}l}
&\gb=\frac{B}{D}(\hbSigma+\Qb)^{-1}\one-\frac{A}{D}(\hbSigma+\Qb)^{-1}\rb,~\hb=\frac{C}{D}(\hbSigma+\Qb)^{-1}\rb-\frac{A}{D}(\hbSigma+\Qb)^{-1}\one,\\
&A=\rb^\top(\hbSigma+\Qb)^{-1}\one,~B=\rb^\top(\hbSigma+\Qb)^{-1}\rb,~C=\one^\top(\hbSigma+\Qb)^{-1}\one,~D=BC-A^2.
\end{array}
\end{align}
Using the solution $\wb^*$, the efficient frontier can be represented as
\begin{align}
    \label{eq:frontierequation}
     \sigma_0^2= \wb^{*\top}\bSigma\wb^*=(\gb+\mu_0\cdot\hb)^\top \bSigma (\gb+\mu_0\cdot\hb),
\end{align}
which should be viewed as a curve of $(\sigma_0,\mu_0)$ as we change the target return $\mu_0$. Our objective is to estimate this curve represented by \eqref{eq:frontierequation} for a given regularization matrix $\Qb$.

\subsection{Out-of-sample estimation of frontier variances}
In this section, we present the main results on estimating the efficient frontier. 
We give an additional technical assumption on $\rb$, assuming that the vectors $\rb$ and $\one$ do not coincide. Otherwise, we have $\mu_0=1$ and we do not have a well-posed efficient frontier.
\begin{assumption} 
\label{assump:assump5}
Let $s_0>0$ to be the unique solution of the equation\footnote{The uniqueness and existence can be found in Lemma~\ref{lemma:constant_level}.}
\begin{align*}
    s_0=\frac{c}{p}\tr~\bSigma\bigg( \frac{\bSigma}{1+s_0}+\Qb \bigg)^{-1}.
\end{align*}
Define
\begin{align*}
    \cA_{rr}=\rb^\top \bigg(\frac{\bSigma}{1+s_0}+\Qb\bigg)^{-1}\rb,\quad \cA_{r1}=\rb^\top \bigg(\frac{\bSigma}{1+s_0}+\Qb\bigg)^{-1}\one,\quad \cA_{11}=\one^\top \bigg(\frac{\bSigma}{1+s_0}+\Qb\bigg)^{-1}\one.
\end{align*}
It is assumed that there exists a constant $\rho<1$ such that
$
   \cA_{r1}^2/(\cA_{11}\cA_{rr})\leq \rho<1.
$
\end{assumption}
This technical assumption states that the vectors $\rb$ and $\one$ do not coincide under the inner product $\langle \ab, \bb \rangle = \ab \big(\frac{\bSigma}{1+s_0}+\Qb\big)^{-1} \bb$, which is a mild assumption. 
In fact, Cauchy-Schwartz inequality gives  $\frac{\cA_{r1}^2}{\cA_{rr} \cA_{11}}\leq1$. The assumption $\frac{\cA_{r1}^2}{\cA_{rr} \cA_{11}}\leq \rho<1$ is a direct relaxation of the upper bound $1$.
Given this additional assumption, we have the following theorem.
\begin{theorem}
\label{thm:frontier}
Suppose that Assumptions~\ref{def:X}, \ref{ass:assump1}-\ref{ass:assump3} and \ref{assump:assump5} hold. 
Recall $\sigma_0^2$ is defined in \eqref{eq:frontier2}. Define 
\[
\hat{\sigma}^2=\frac{(\gb+\mu_0\hb)^\top\hbSigma(\gb+\mu_0\hb)}{(1-c/p\cdot\tr \hbSigma(\hbSigma+\Qb)^{-1})^2},
\]
where $\gb$ and $\hb$ are defined in \eqref{eq:def_gh}, 
it holds that
\begin{align*}
    \hat{\sigma}^2/\sigma_0^2~\cvas~1.
\end{align*}
Moreover, the following properties hold:
\begin{enumerate}
    \item If $\cA_{rr}$ is bounded, then for any $r_0=O(\mu_0)$ it holds that  $\frac{\mu_0-r_0}{\sigma_0}-\frac{\mu_0-r_0}{\hat{\sigma}}~\cvas~0$.
    \item If $\mu_0\leq C\sqrt{\cA_{rr}}$ for some large constant $C>0$, then $ \hat{\sigma}^2-\sigma_0^2~\cvas~0$.
\end{enumerate}
\end{theorem}
The proof of Theorem~\ref{thm:frontier} can be found in Appendix~\ref{sec:proofthmfrontier} in the supplementary material. Notice that when the target return $\mu_0$ is given, estimating the volatility naturally leads to an estimator of Sharpe ratio. 
The optimal true $SR(\Qb)$ without the risk-free asset equals $(\mu_0 - r_0)/\sigma_0$,
and Theorem~\ref{thm:frontier} implies that the estimator $\hat{SR}(\Qb)=(\mu_0-r_0)/\hat{\sigma}$ is ratio-consistent without additional scaling condition, and is also 
difference-consistent when $\cA_{rr}$ is bounded.
In Theorem~\ref{thm:main_theorem}, if we set $\Ab=\rb\rb^\top$, the relaxed condition that $\tr\Big(\frac{\bSigma}{1+s_0}+\Qb\Big)^{-1}\Ab $ is bounded (see the footnote for Theorem~\ref{thm:main_theorem}) for the difference consistency of Sharpe ratio is exactly the same as the condition that $\cA_{rr}$ is bounded in Theorem~\ref{thm:frontier}. 

Albeit the conditions are similar, we hope to point out new challenges in proving Theorem~\ref{thm:frontier}.
The first challenge arises from the presence of distant diverging spikes in the covariance matrix $\bSigma$ in Case 2 and Case 3 of Assumption~\ref{ass:assump3}. It is possible that the vectors $\rb$ and $\one$ lie within the space spanned by the eigenvectors of these spikes. Note that in Section~\ref{sec:sharperatio}, due to the existence of the risk-free asset, we directly work with the exccess return $\bmu$ without necessity to discuss the relationship between $\rb$ and $\one$. 
But in the proof of Theorem~\ref{thm:frontier}, we need to handle general $\cA_{rr}, \cA_{11}, \mu_0$ without specifying their statistical rates and without specifying a detailed relationship between $\rb$ and $\one$. 
Readers may refer to Appendix~\ref{sec:examplecase_frontier} for an example illustrating the values of $\cA_{rr}$, $\cA_{11}$, $\cA_{r1}$, the possible range of $\mu_0$, and how $\cA_{r1}^2 / (\cA_{11} \cA_{rr}) \leq \rho < 1$ can be satisfied under mild conditions. 
Moreover, we need to divide the proof into the cases $\mu_0\leq C\sqrt{\cA_{rr}/\cA_{11}}$ and $\mu_0\geq C\sqrt{\cA_{rr}/\cA_{11}}$, where distinct analyses are required. 
In addition to the above order-related challenge, the expression of $\sigma_0^2$ presents another obstacle. While $\sigma_0^2$ can be expressed as $\bxi^\top (\hbSigma+\Qb)^{-1}\bSigma(\hbSigma+\Qb)^{-1}\bxi$ for some vector $\bxi$, similar to \eqref{eq:Tn2} in Section~\ref{sec:sharperatio}, it is important to note that $\bxi$ depends on $\hbSigma$. To tackle this nontrivial dependency, we need to carefully substitute $\bxi$ with a deterministic counterpart and make sure this does not affect the almost sure convergence. More proof details can be found in Appendix~\ref{sec:proofthmfrontier} in the supplementary material. 
\section{Numerical Experiments}
\label{sec:simulation}
In this section, we perform numerical experiments to validate our theoretical findings. 
In Section~\ref{subsec:basic_simulation}, we provide the basic settings for SR estimation and verify Theorem~\ref{thm:main_theorem}, especially the asymptotic behavior as $n$ increases. 
In Section~\ref{sec:simu_frontier}, we verify our estimation of the efficient frontier in Theorem~\ref{thm:frontier}. 
Some additional simulations  are given in Appendix~\ref{sec:appendixA}. Appendix~\ref{subsec:simu_alternative} provides comparisons on different choices of the population covariance $\bSigma$, the mean vector of returns $\bmu$ and the regularization $\Qb$,  while Appendix~\ref{subsec:simu_unknown} gives simulations on unknown $\bmu$.  In all the simulations, we set $r_0=0$ for simplicity so that $\rb=\bmu$.

\subsection{Sharpe ratio estimation}
\label{subsec:basic_simulation}
We present the simulation results under the basic settings. We first fix the sample size as $n=1500$ and consider two dimensions: $p=750 $ ($c=1/2$) and $p=2250 $ ($c=3/2$). These cases correspond to values of $c$ less than and greater than 1. To simulate financial returns, we construct the population covariance $\bSigma_0$, mean vector $\bmu_0$, and the matrix $\Qb_0$ as follows:

\begin{enumerate}[leftmargin=*]
    \item The population covariance $\bSigma = \bSigma_0=\diag(\lambda_1, \dots, \lambda_p) + 2\one\one^\top$, where $\{\lambda_i\}_{i=1}^p$ are generated from a truncated $\Gamma^{-1}(1,1)$ distribution, truncated with the interval $[0.01,9]$, and then ranked in decreasing order. This particular construction of $\bSigma$ follows from the factor model.
    Specifically, the covariance  consists of the market
    systematic risk, where the vector $\one$ represents a market factor introducing a spike in covariance spectrum, and the specific risk, where $\lambda_i$ represents the individual residual risk of each asset.
    \item Assume the risk-free return is $r_0=0$. The mean vector of returns is given by $\bmu_0=\sqrt{5/p}\cdot (\one(S_+)-\one(S_-)) \in \RR^p$. Here, $\one(S)$ is defined in Subsection~\ref{subsec:notation}. The subsets $S_+$ and $S_-$ are randomly selected subsets of $[p]$ with $|S_+|=|S_-|=p/10$ and $S_+\cup S_-=\emptyset$. This formulation of $\bmu_0$ is motivated by the idea that a portion of the assets may yield low positive or negative returns, known as the ``alpha'' signal in asset pricing, while the majority assets simply have zero mean, ensuring no arbitrage opportunities.
    \item  The regularization matrix is $\Qb=q\cdot \Qb_0$, where $\Qb_0 = \diag(3,...,3,1,...,1)$ with $p/2$ numbers of $3$ and $1$ entries. By structuring $\Qb_0$ in this manner, we categorize stocks into high and low volatility regimes, capturing some weak knowledge of the covariance structure. We will vary the values of $q$ in order to verify Theorem~\ref{thm:main_theorem}.
\end{enumerate}
In the simulation, for each pair of $(n,p)$ considered, the parameters $\bSigma_0\in\RR^{p\times p}$, $\bmu_0\in\RR^{p}$, and $\Qb_0\in\RR^{p\times p}$ are generated, and they will remain fixed throughout the simulation procedure.
Then we generate iid Gaussian random vector $\Rb_i$ ($i\in[n]$)  with mean $\bmu=\bmu_0$ and  variance $\bSigma=\bSigma_0$. The observed random matrix is  $\Rb=(\Rb_1,\Rb_2,...,\Rb_n)^\top\in\RR^{n\times p}$. 
 In this simulation,  we assume the mean vector $\bmu$ is known, and thus $\Xb=\Rb-\one_n\bmu_0^\top$ is also well observed. 
Recall that $\hbSigma=\Xb^\top\Xb/n$, and $SR(\Qb)=\frac{\bmu^\top(\hbSigma+\Qb)^{-1}\bmu}{\sqrt{\bmu^\top(\hbSigma+\Qb)^{-1}\bSigma(\hbSigma+\Qb)^{-1}\bmu}}$.
We compare the true $SR(\Qb)$ and the estimator $\hat{SR}(\Qb)$ in Theorem~\ref{thm:main_theorem} which is
\begin{align*}
    \hat{SR}(\Qb)=\bigg(1-\frac{c}{p} \tr\hbSigma(\hbSigma+\Qb)^{-1}\bigg)\cdot \frac{\bmu^\top(\hbSigma+\Qb)^{-1}\bmu}{\sqrt{\bmu^\top(\hbSigma+\Qb)^{-1}\hbSigma(\hbSigma+\Qb)^{-1}\bmu}}.
\end{align*}
The Sharpe computation is then carried out as follows for each pair of $(n,p)$. We repeat the generation of Gaussian random matrices $\Rb\in\RR^{n\times p}$ for  $1000$  independent times. For each repetition, we vary the value of $q$ to compute and compare $SR(q\cdot\Qb_0)$ and $\hat{SR}(q\cdot\Qb_0)$. 
If $c<1$, we vary $q$ from $(1:30)/5$, while if $c>1$, we vary $q$ from $(1:30)/1.5$. 
For each given $q$, we will obtain  1000 different values of $SR(q\cdot\Qb_0)$ and $\hat{SR}(q\cdot\Qb_0)$.

\begin{figure}[t!]
    \centering
    \subfloat[$(n,p)=(1500,750)$]{\includegraphics[width=0.37\textwidth]{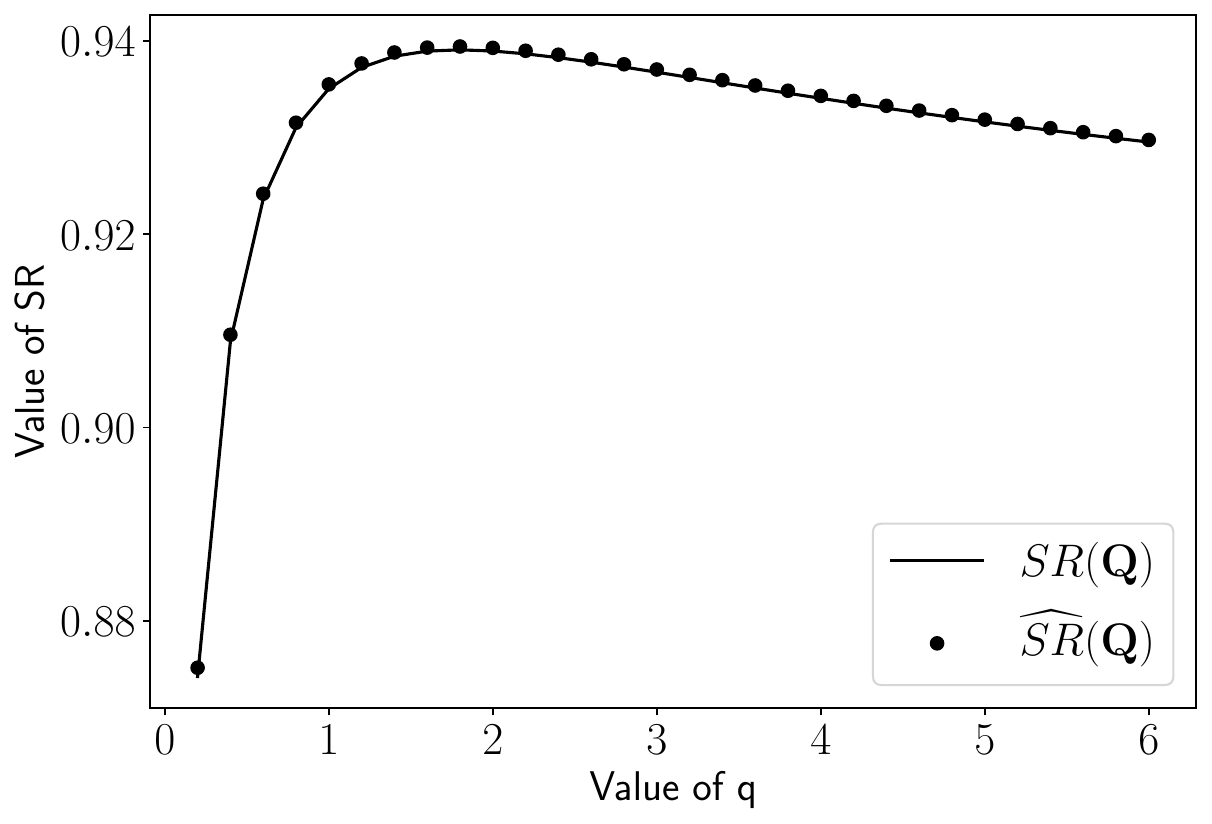}\label{fig_simu:Basen1500_c<1}}\qquad
    \subfloat[$(n,p)=(1500,2250)$]{\includegraphics[width=0.37\textwidth]{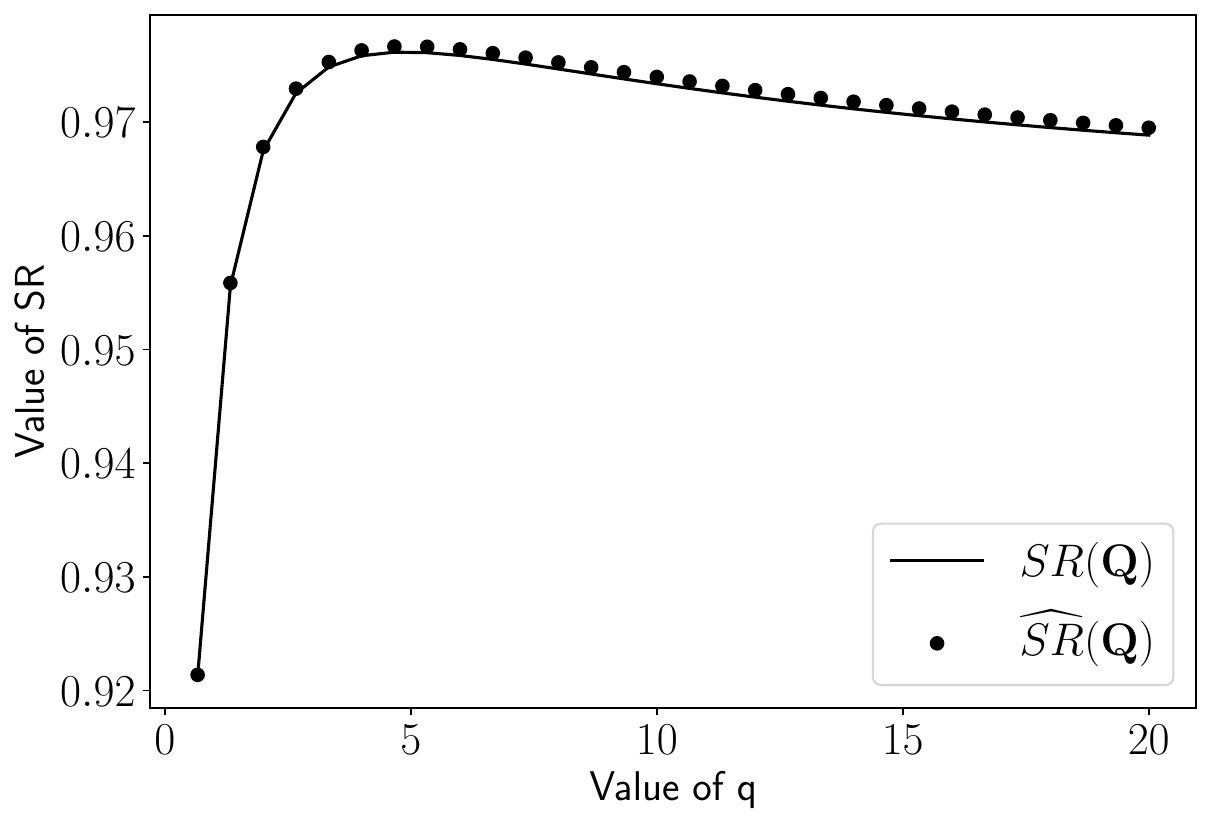}\label{fig_simu:Basen1500_c>1}}
    \caption{Simulation results in the basic settings. Figure~\ref{fig_simu:Basen1500_c<1} shows the case when $c=1/2$ and Figure~\ref{fig_simu:Basen1500_c>1} depicts the case when $c=3/2$. The x-axis corresponds to different $q$ values, and the y-axis is the value of $SR$. The black solid line connects the values of $SR(q\cdot\Qb_0)$, while the points represent the proposed statistics $\hat{SR}(q\cdot\Qb_0)$ in Theorem~\ref{thm:main_theorem}.}
    \label{fig_simu:Basen1500}
\end{figure}

In Figure~\ref{fig_simu:Basen1500}, the $x$-axis represents the different values of $q$, while the $y$-axis gives the Sharpe ratio. The simulation results are presented as the average of $SR(q\cdot\Qb_0)$ and $\hat{SR}(q\cdot\Qb_0)$ over  $1000$ independent trials.
It is evident that the value of $\hat{SR}(q\cdot\Qb_0)$  matches the value of $SR(q\cdot\Qb_0)$ very well. As $q$ increases, $SR(q\cdot\Qb_0)$ initially rises and then falls, indicating the presence of an optimal $q$ that maximizes $SR$. With our proposed statistics in Theorem~\ref{thm:main_theorem}, we can replicate this trend without requiring knowledge of the true population matrix $\bSigma$. Furthermore, the close match between $\hat{SR}(q\cdot\Qb_0)$ and $SR(q\cdot\Qb_0)$ suggests that our estimator can be used to determine the best regularization matrix $\Qb$ from a deterministic candidate set, such that $SR(\Qb)$ achieves the maximum over the candidate set too.

\subsubsection{Asymptotics with increasing sample size}
\label{subsubsec:simu_asymp_SR}
\begin{figure}[t!]
\centering
\subfloat[MSE of difference]{\includegraphics[width=0.235\textwidth]{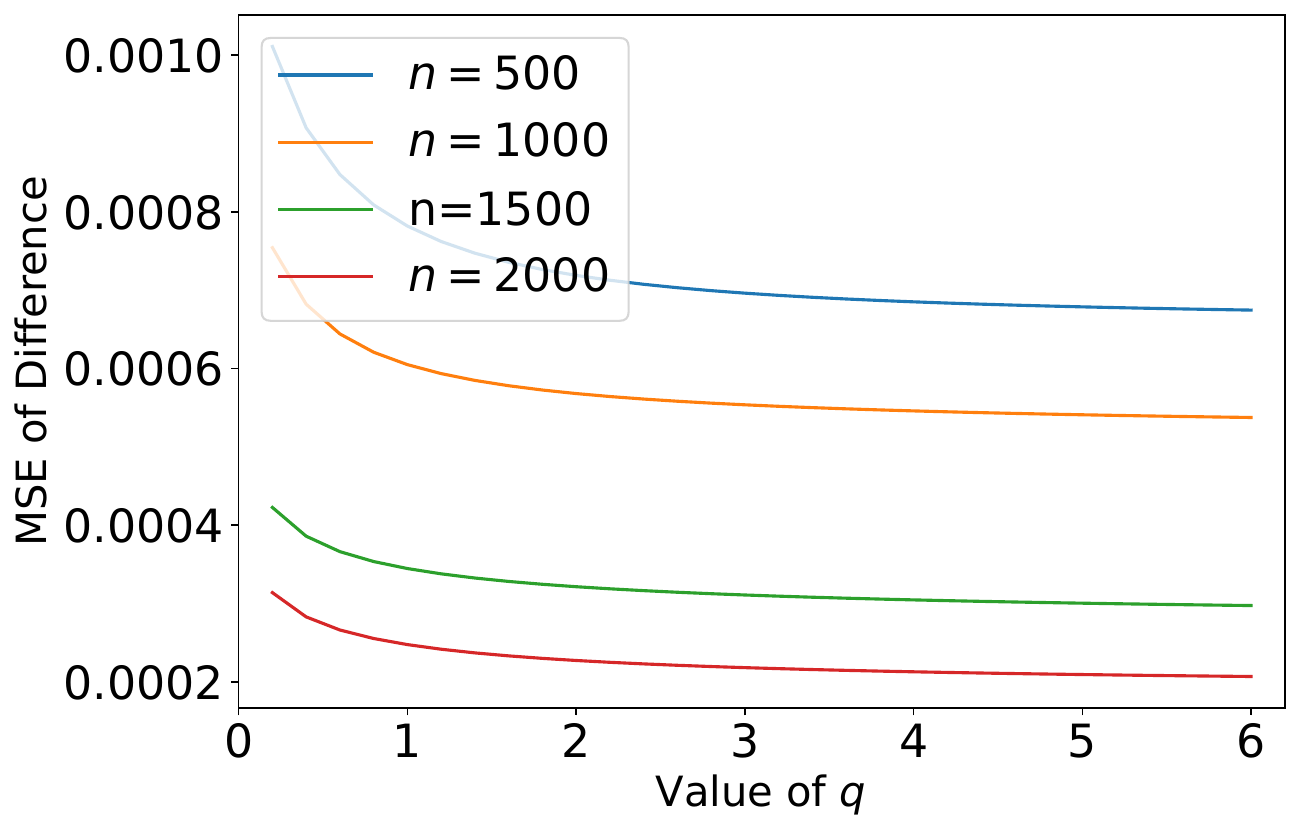}\label{fig_simu:asymp_msesmall}}
\subfloat[MSE of ratio]{\includegraphics[width=0.24\textwidth]{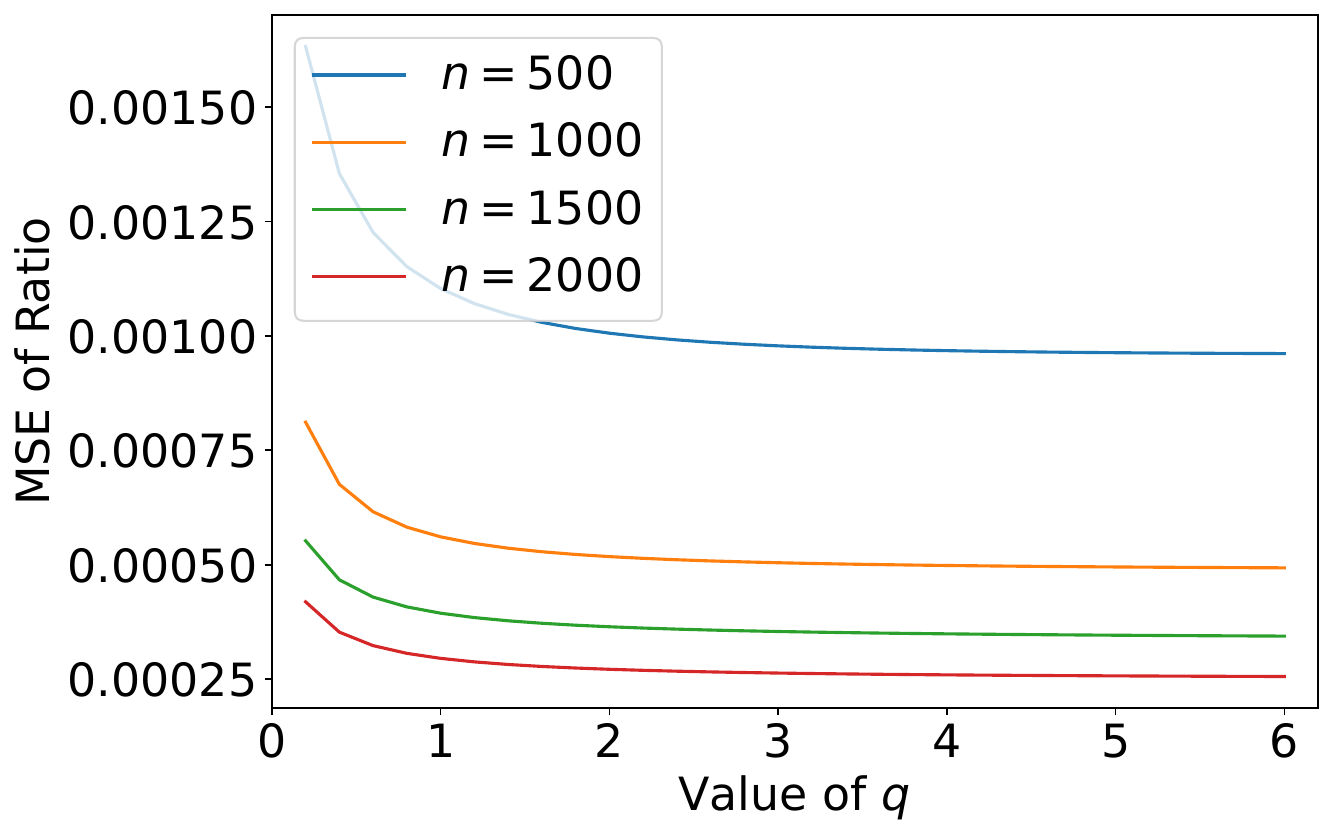}\label{fig_simu:asymp_mseratiosmall}}
\subfloat[Difference on $q^*$]{\includegraphics[width=0.22\textwidth]{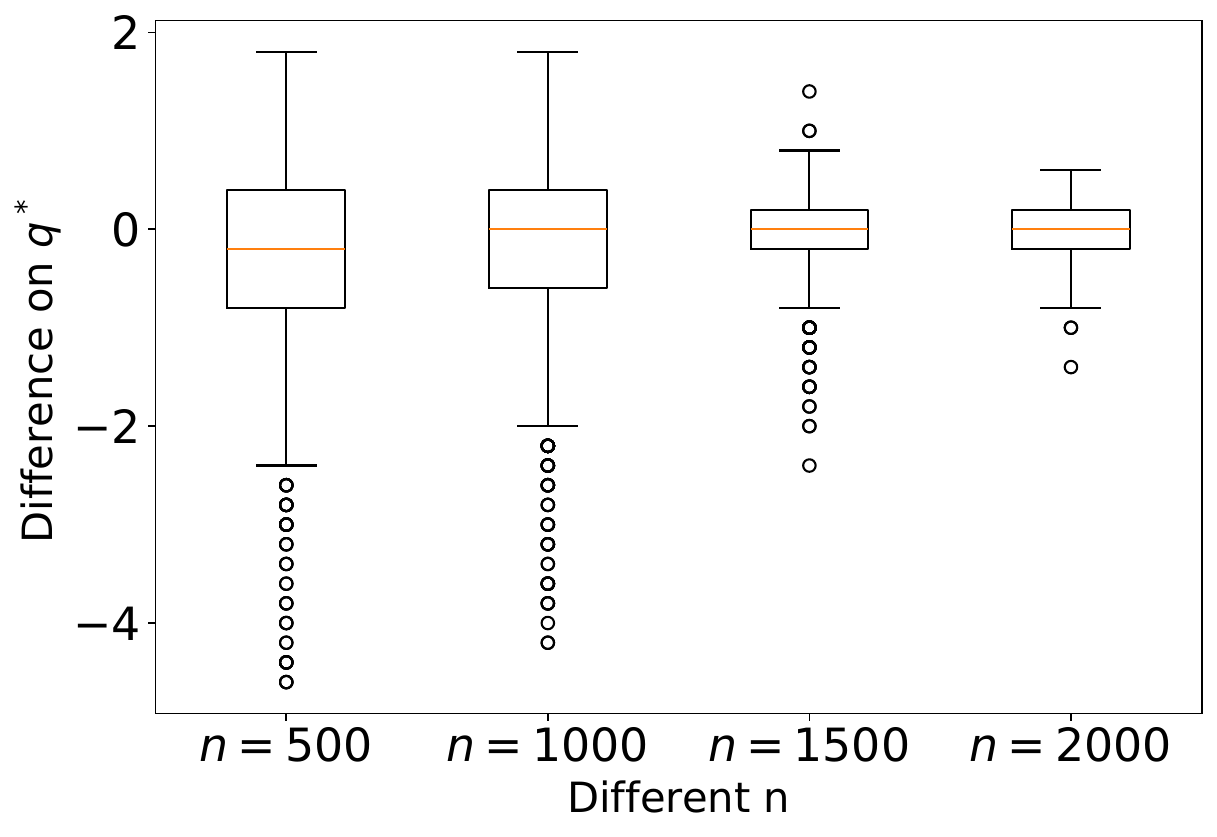}\label{fig_simu:asymp_qsmall}}
\subfloat[Difference on $\text{SR}^*$]{\includegraphics[width=0.24\textwidth]{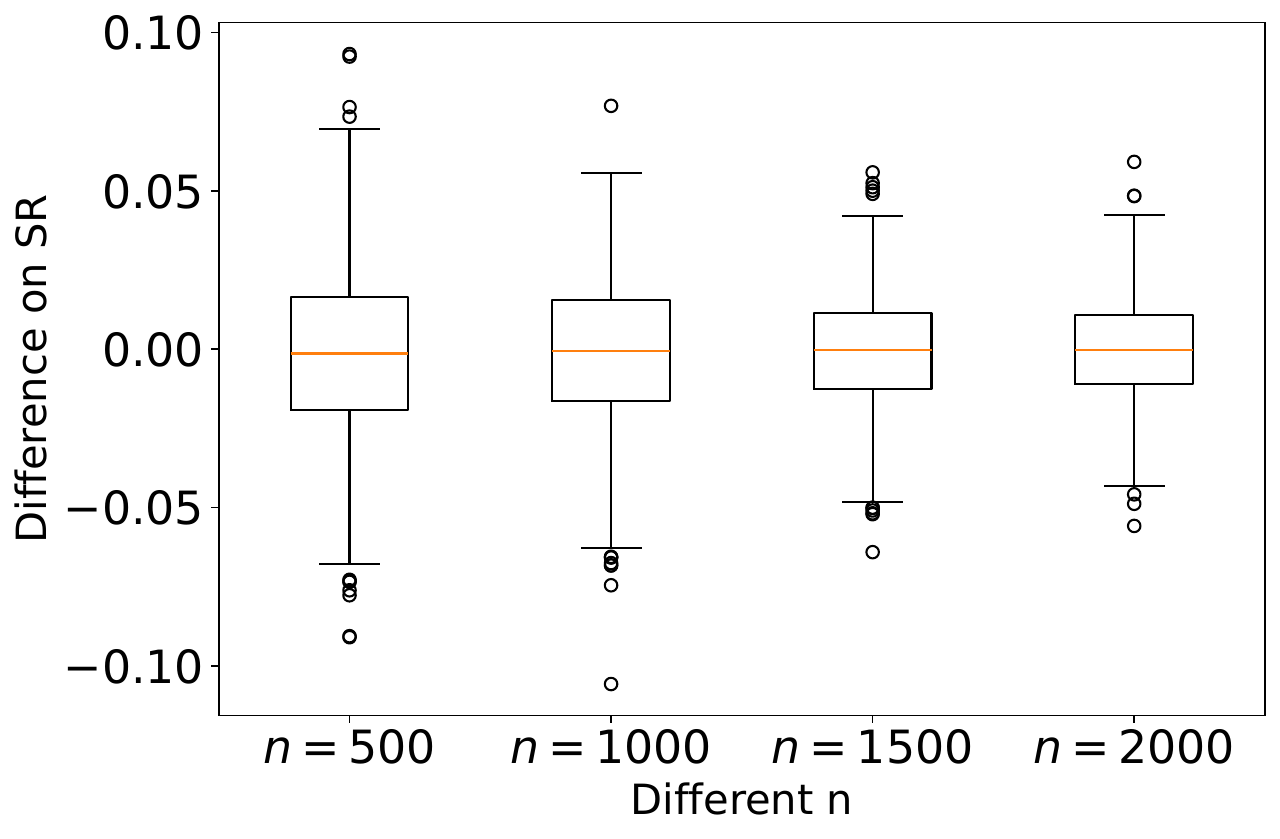}\label{fig_simu:asymp_SRsmall}}

\subfloat[MSE of  difference]{\includegraphics[width=0.235\textwidth]{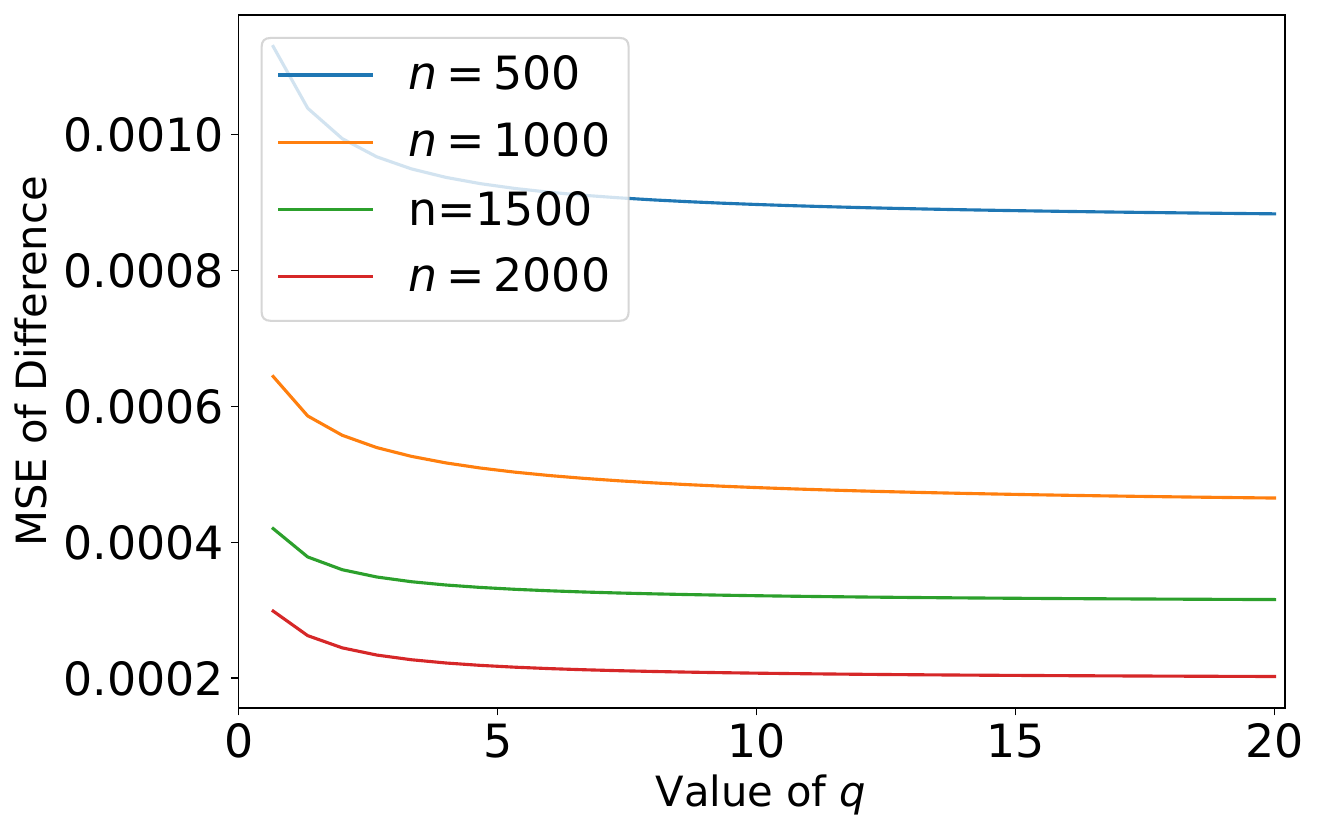}\label{fig_simu:asymp_msebig}}
\subfloat[MSE of ratio]{\includegraphics[width=0.235\textwidth]{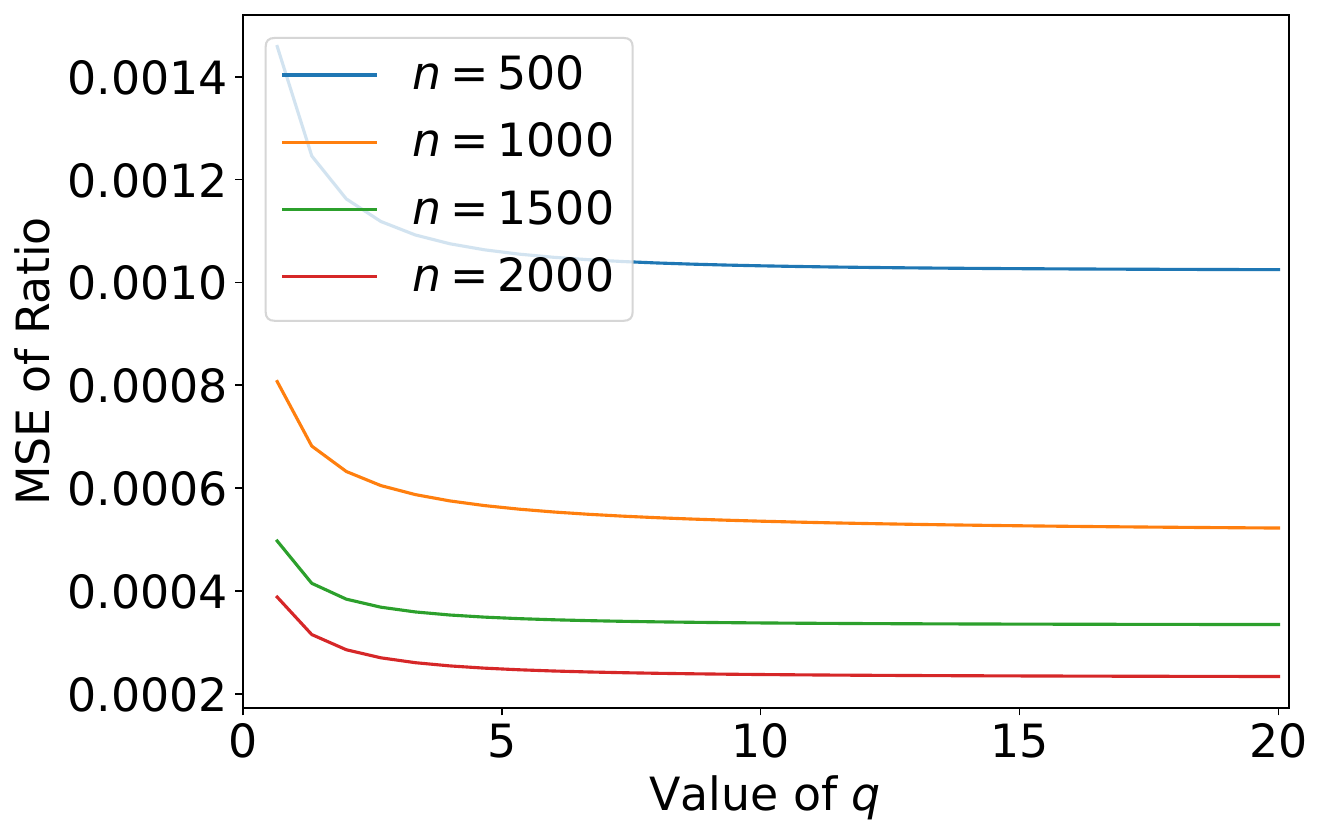}\label{fig_simu:asymp_mseratiobig}}
\subfloat[Difference on $q^*$]{\includegraphics[width=0.225\textwidth]{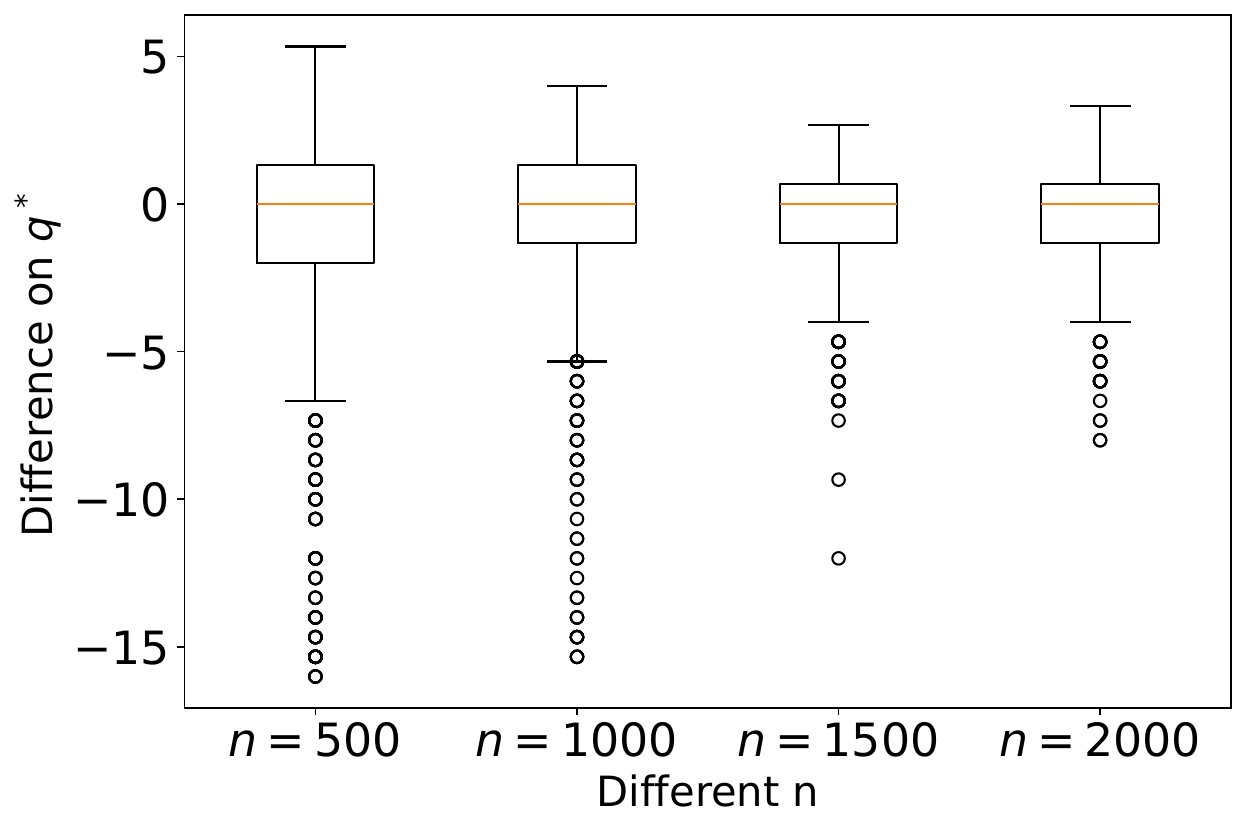}\label{fig_simu:asymp_qbig}}
\subfloat[Difference on $\text{SR}^*$]{\includegraphics[width=0.24\textwidth]{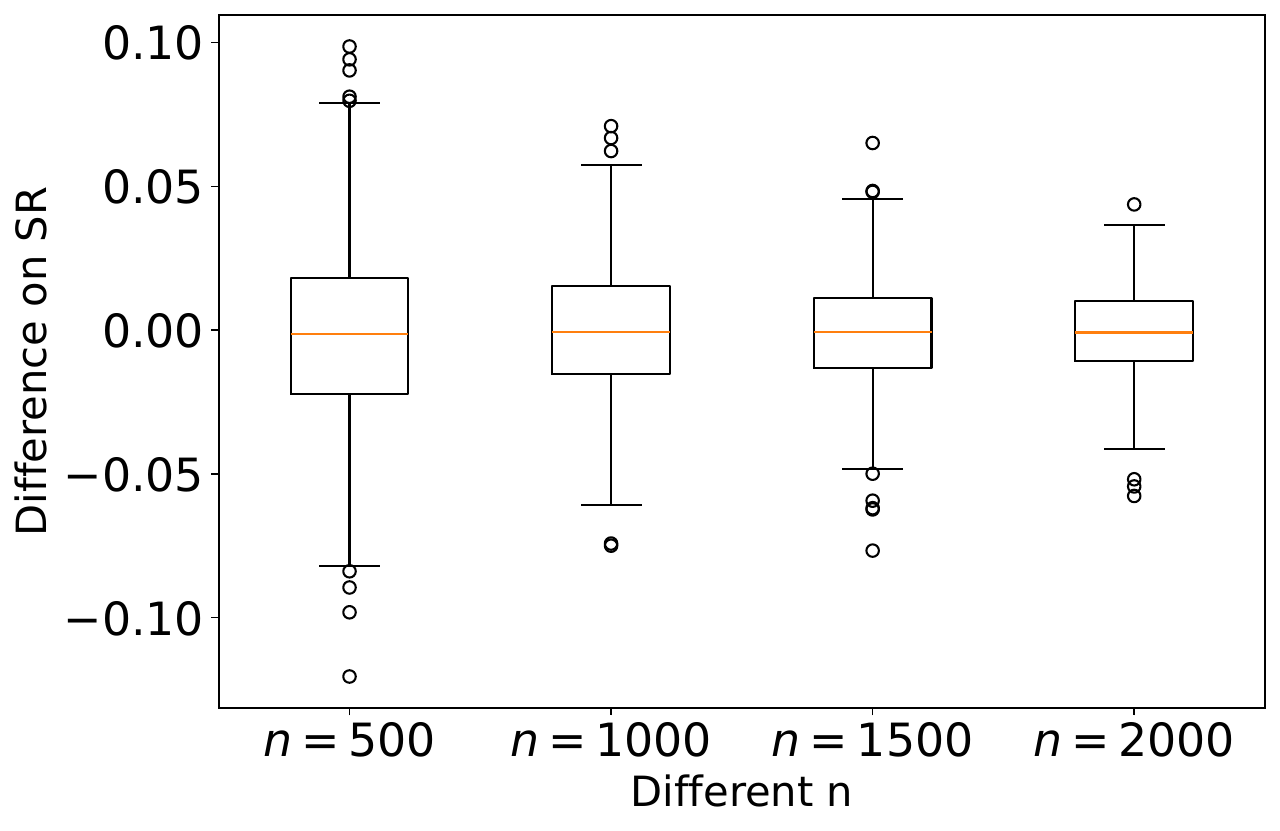}\label{fig_simu:asymp_SRbig}}
\caption{Simulation results with increasing $n$. Figures~\ref{fig_simu:asymp_msesmall}-\ref{fig_simu:asymp_SRsmall} and Figures~\ref{fig_simu:asymp_msebig}-\ref{fig_simu:asymp_SRbig} respectively correspond to $c=p/n=1/2$ and $c=p/n=3/2$. In Figures~\ref{fig_simu:asymp_msesmall} and \ref{fig_simu:asymp_msebig},  $y$-axis shows  $\sum_{b=1}^{1000}(SR_b(q\cdot\Qb_0;n,p)-\hat{SR}_b(q\cdot\Qb_0;n,p))^2/1000$ for different  $q$'s. In Figures~\ref{fig_simu:asymp_mseratiosmall} and \ref{fig_simu:asymp_mseratiobig},  $y$-axis shows  $\sum_{b=1}^{1000}(SR_b(q\cdot\Qb_0;n,p)/\hat{SR}_b(q\cdot\Qb_0;n,p)-1)^2/1000$  for different $q$'s. In Figures~\ref{fig_simu:asymp_qsmall} and \ref{fig_simu:asymp_qbig},  $y$-axis gives the boxplot of $\argmax_q SR_b(q\cdot\Qb_0;n,p)-\argmax_q\hat{SR}_b(q\cdot\Qb_0;n,p)$ across the  1000 trials for different $n$'s. In Figures~\ref{fig_simu:asymp_SRsmall} and \ref{fig_simu:asymp_SRbig}, $y$-axis displays the boxplot of $\max_q SR_b(q\cdot\Qb_0;n,p)-\max_q\hat{SR}_b(q\cdot\Qb_0;n,p)$ across the  1000 trials for different $n$'s.}
    \label{fig_simu:asymp}
\end{figure}

In this section, we aim to investigate the effect of varying $n$ and $p$ on our SR predictions, while keeping all other experimental settings consistent with previous settings. We continue to set the parameter $c=p/n$ at either $1/2$ or $3/2$, and focus on observing changes across different $n$ values, specifically $500, 1000, 1500$, and $2000$. 
Following the same procedure as described above in Section~\ref{subsec:basic_simulation}, for each $(n,p)$ and each given $q$, we obtain  $1000$ values of $SR(q\cdot\Qb_0)$ and $\hat{SR}(q\cdot\Qb_0)$. 
For clear illustration, we let $SR_b(q\cdot\Qb_0;n,p)$ and $\hat{SR}_b(q\cdot\Qb_0;n,p)$ denote the values obtained in the $b$-th independent trial ($b\in [1000]$) under the pair $(n,p)$. The outcomes are subsequently illustrated in Figure~\ref{fig_simu:asymp}. 

Figure~\ref{fig_simu:asymp} shows that as we increase the sample size $n$, the mean-squared-errors (MSE) of $SR(q\cdot\Qb_0;n,p)-\hat{SR}(q\cdot\Qb_0;n,p)$ (difference) and  $SR(q\cdot\Qb_0;n,p)/\hat{SR}(q\cdot\Qb_0;n,p)-1$ (ratio)  get smaller for each $q$ value. So we verified our theory  that when $n\to\infty$, the estimated value $\hat{SR}(q\cdot\Qb_0;n,p)$ converges to the value $SR(q\cdot\Qb_0;n,p)$ almost surely, and they are ratio consistent. Additionally, we examine how the optimal $q$ values and the maximal Sharpe ratios of $SR(q\cdot\Qb_0;n,p)$ and $\hat{SR}(q\cdot\Qb_0;n,p)$ can differ. Figure~\ref{fig_simu:asymp} clearly demonstrates that as $n$ grows, this discrepancy diminishes, lending further support to our theory.

\subsection{Efficient frontiers}
\label{sec:simu_frontier}
In this section, we conduct simulations to validate Theorem~\ref{thm:frontier} for estimating the out-of-sample efficient frontiers. We still apply $(n,p)=(1500,750)$ with $c=1/2$ and $(n,p)=(1500,2250)$ with $c=3/2$. Recall the optimization problem in \eqref{eq:frontier2}:
\begin{align}
    \wb^*=\arg\min_{\wb} \wb^\top (\hbSigma+\Qb) \wb, \quad \text{s.t. } \wb^\top \rb =\mu_0 \text{ and } \wb^\top \one=1.\label{eq:simu_frontier}
\end{align}

We will use the regularization matrix $\Qb = 0.2\Qb_0$, where $\Qb_0$ is defined in Section~\ref{subsec:basic_simulation}. The population covariance matrix $\bSigma$ and the mean $\rb = \bmu$ (since the risk-free rate $r_0$ is assumed to be zero) we use here are slightly different from those used in Section~\ref{subsec:basic_simulation}. We first generate a vector $\bxi_3\in\RR^p$ \footnote{Notations $\bxi_1,\bxi_2,\bmu_1,\bmu_2,\bSigma_1,\bSigma_2,\Qb_1,\Qb_2$ are used in Appendix \ref{subsec:simu_alternative}.} whose elements are independently distributed from $\Gamma(1,1)$. Then we let $\bSigma = \bSigma_3 = \diag(\lambda_1,\dots,\lambda_p)+2\cdot\one\one^\top+\bxi_3\bxi_3^\top$.  The matrix $\bSigma_3$ simulates a covariance matrix with two factors. In addition, the mean vector $\bmu$ here has two choices: $\bmu_3=p^{\frac{1}{4}}\bmu_0+2\cdot\one$ and $\bmu_4=\bmu_0+2\cdot\one+\bxi_3$, where $\bmu_0$ is again defined in Section~\ref{subsec:basic_simulation}. We design the vectors $\bmu_3$ and $\bmu_4$ such that $\cA_{rr}$ becomes unbounded when $\bmu = \bmu_3$, while it remains bounded when $\bmu = \bmu_4$. 

\begin{figure}[t]
    \centering
    \subfloat[$\bmu=\bmu_3,c=1/2$]{\includegraphics[width=0.24\textwidth]{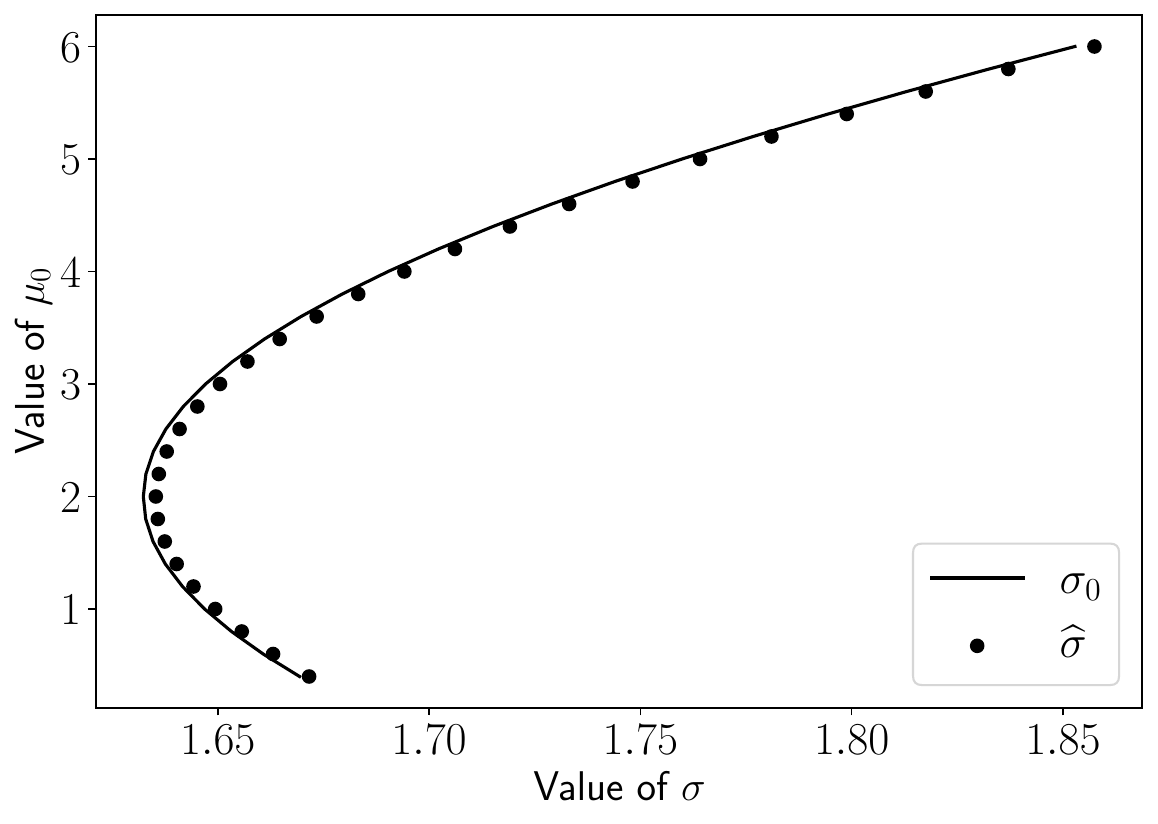}\label{fig_simu:mu1frontierc<1q=0.2}}
    \subfloat[$\bmu=\bmu_4,c=1/2$]{\includegraphics[width=0.24\textwidth]{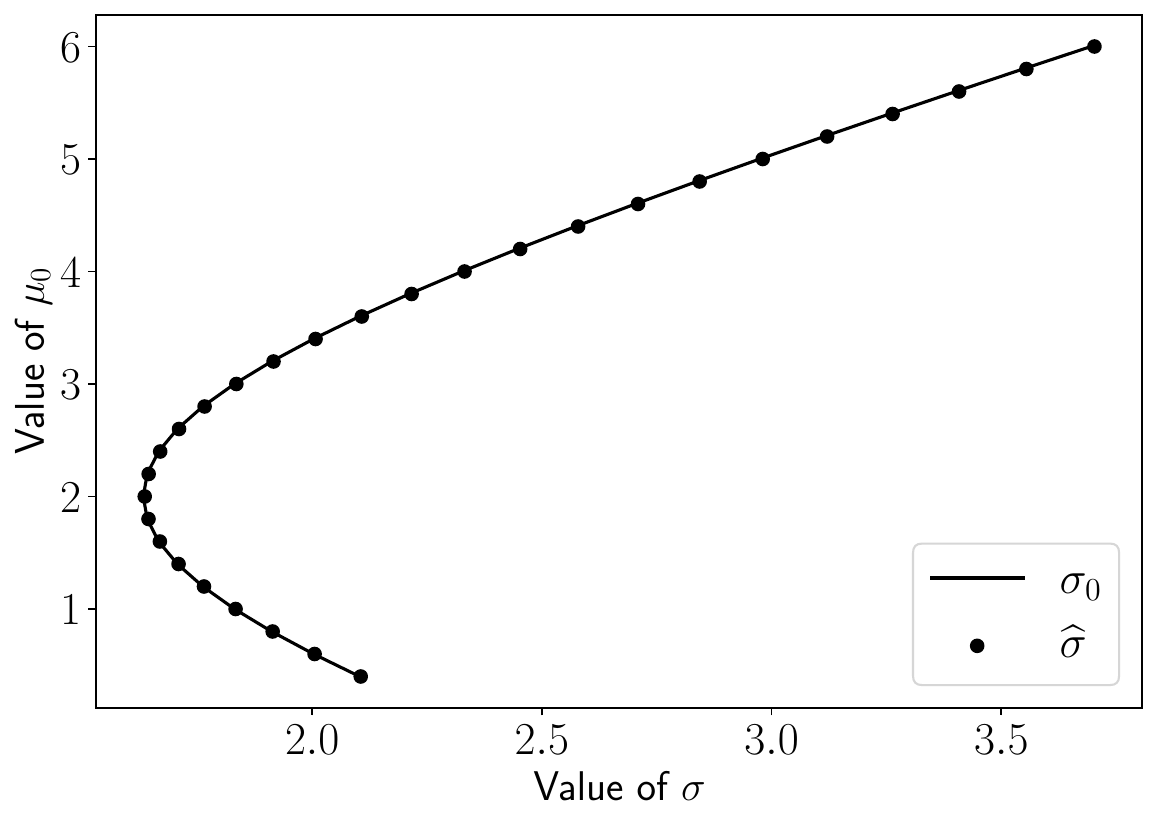}\label{fig_simu:basefrontierc<1q=0.2}}
    \subfloat[$\bmu=\bmu_3,c=3/2$]{\includegraphics[width=0.245\textwidth]{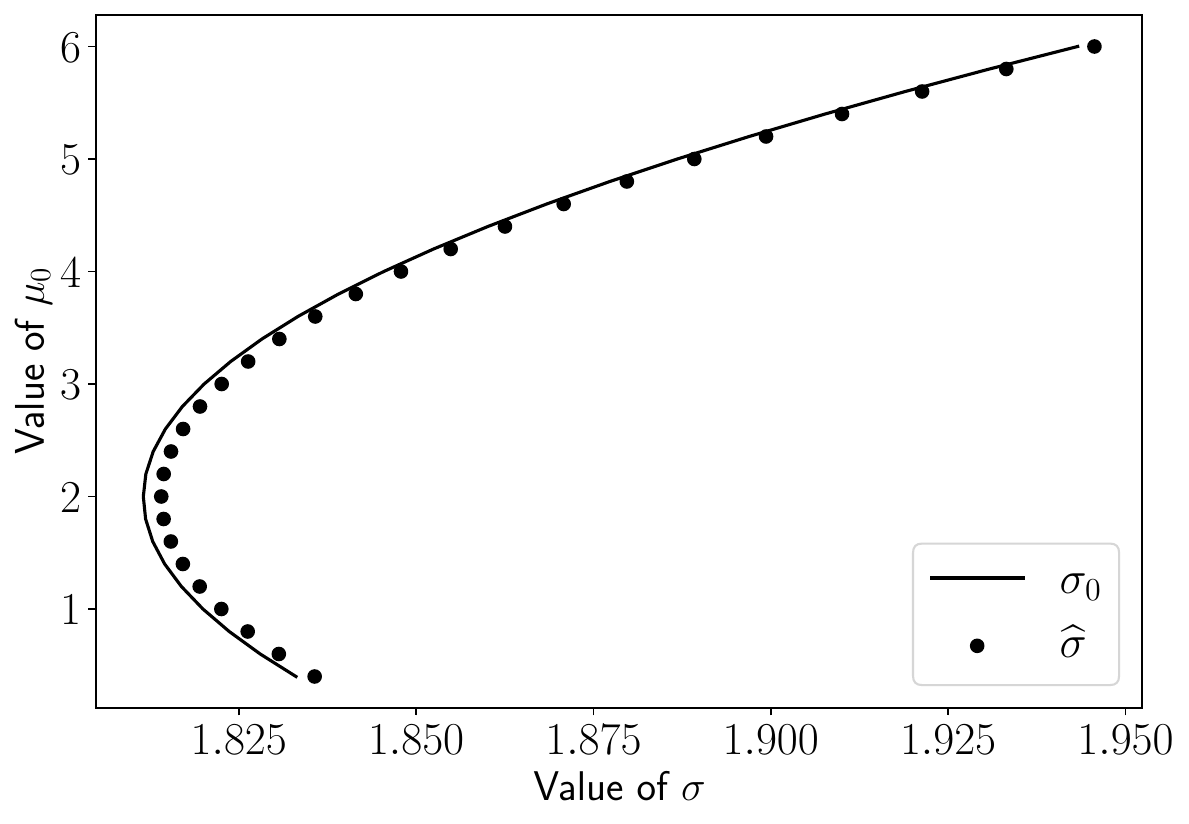}\label{fig_simu:mu1frontierc>1q=0.2}}
    \subfloat[$\bmu=\bmu_4,c=3/2$]{\includegraphics[width=0.24\textwidth]{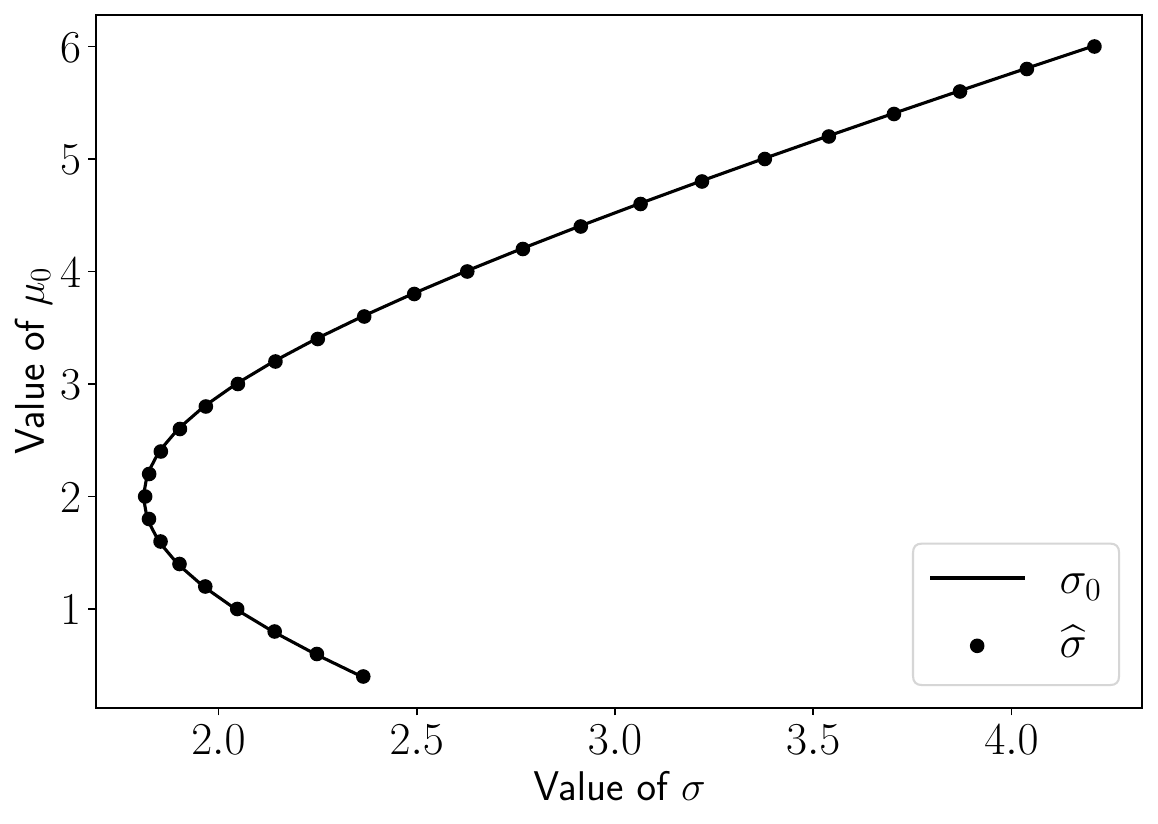}\label{fig_simu:mu2frontierc>1q=0.2}}
    \caption{Efficient frontiers of $\wb^*$ from \eqref{eq:simu_frontier}. The x-axis is the volatility level, and the y-axis is the target return $\mu_0$. The solid line characterizes the curve of $(\mu_0,\sigma_0)$, while the points represent $(\mu_0,\hat{
    \sigma})$, where $\hat{
     \sigma}$ is given in Theorem~\ref{thm:frontier}.}
    \label{fig_simu:basefrontier}
\end{figure}

The following procedure is similar to Section~\ref{subsec:basic_simulation}. For each $(n,p)$, we (1) specify $\Qb=0.2\Qb_0$, $\bSigma = \bSigma_3$, and $\bmu = \bmu_3$ or $\bmu_4$; (2) generate  1000 independent random matrices $\Rb\in\RR^{n\times p}$ with i.i.d rows following $\cN(\bmu, \bSigma)$; (3) for each given $\Rb$, compute the sample covariance matrix $\hbSigma =\frac{(\Rb-\one_n\bmu^\top)^\top(\Rb-\one_n\bmu^\top)}{n}$; (4)
with each $\hbSigma$, solve the optimization problem in \eqref{eq:simu_frontier}, for $\mu_0$ ranging from $0.2$ to $6$ with the increment of $0.2$; (5) for each $\mu_0$, calculate the corresponding values of $\sigma_0$ and $\hat{\sigma}$, whose formula is provided in Theorem~\ref{thm:frontier}.
According to Theorem~\ref{thm:frontier}, the  volatility of the optimal portfolio $\sigma_0$ is defined as $\sigma_0^2=\wb^*\bSigma\wb^*$, and the corresponding estimator is $\hat{\sigma}^2=\wb^*\hbSigma\wb^*/(1-c/p\cdot\tr \hbSigma(\hbSigma+\Qb)^{-1})^2$.
We conduct simulations to confirm that $\hat{\sigma}_0$ matches $\sigma_0$ well for each given $\mu_0$. The results displayed in Figure~\ref{fig_simu:basefrontier} are the average over the  1000 values of $\sigma_0$ and $\hat{\sigma}$ for different $\mu_0$'s.


As shown in Figure~\ref{fig_simu:basefrontier}, we observe that $\hat{\sigma}$ estimates $\sigma_0$ very accurately for various $\mu_0$ and different $c$. As $\mu_0$ increases, we notice that the volatility $\sigma_0$ initially decreases and then increases, actually forming a hyperbola. Remarkably, $\hat{\sigma}$  closely depicts this curve. 

\subsubsection{Asymptotics with increasing sample size}
\label{subsubsec:simu_asymp_frontier}

\begin{figure}[t!]
\centering
\subfloat[$\bmu=\bmu_3,c=1/2$]{\includegraphics[width=0.235\textwidth]{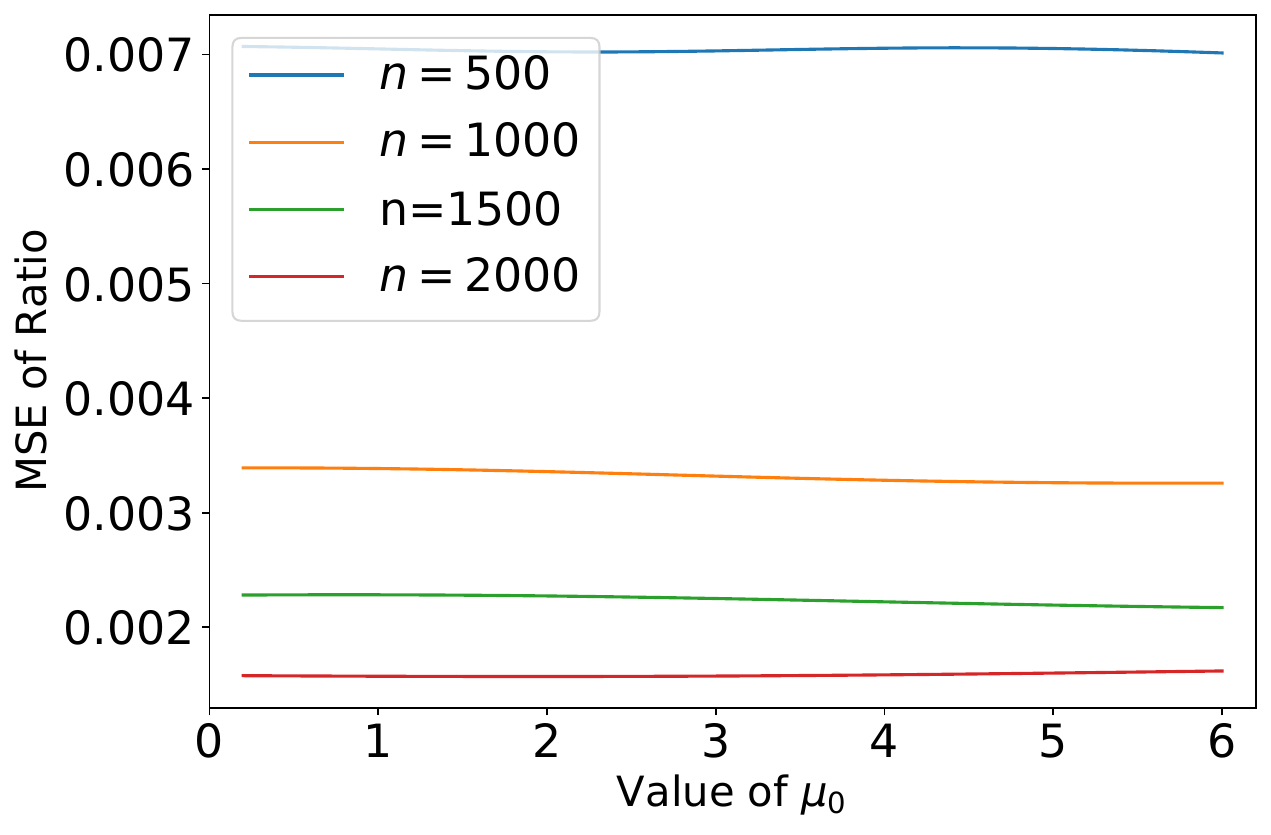}\label{fig_simu:frontier_asymp_mseratiosmallmu3}}
\subfloat[$\bmu=\bmu_4,c=1/2$]{\includegraphics[width=0.235\textwidth]{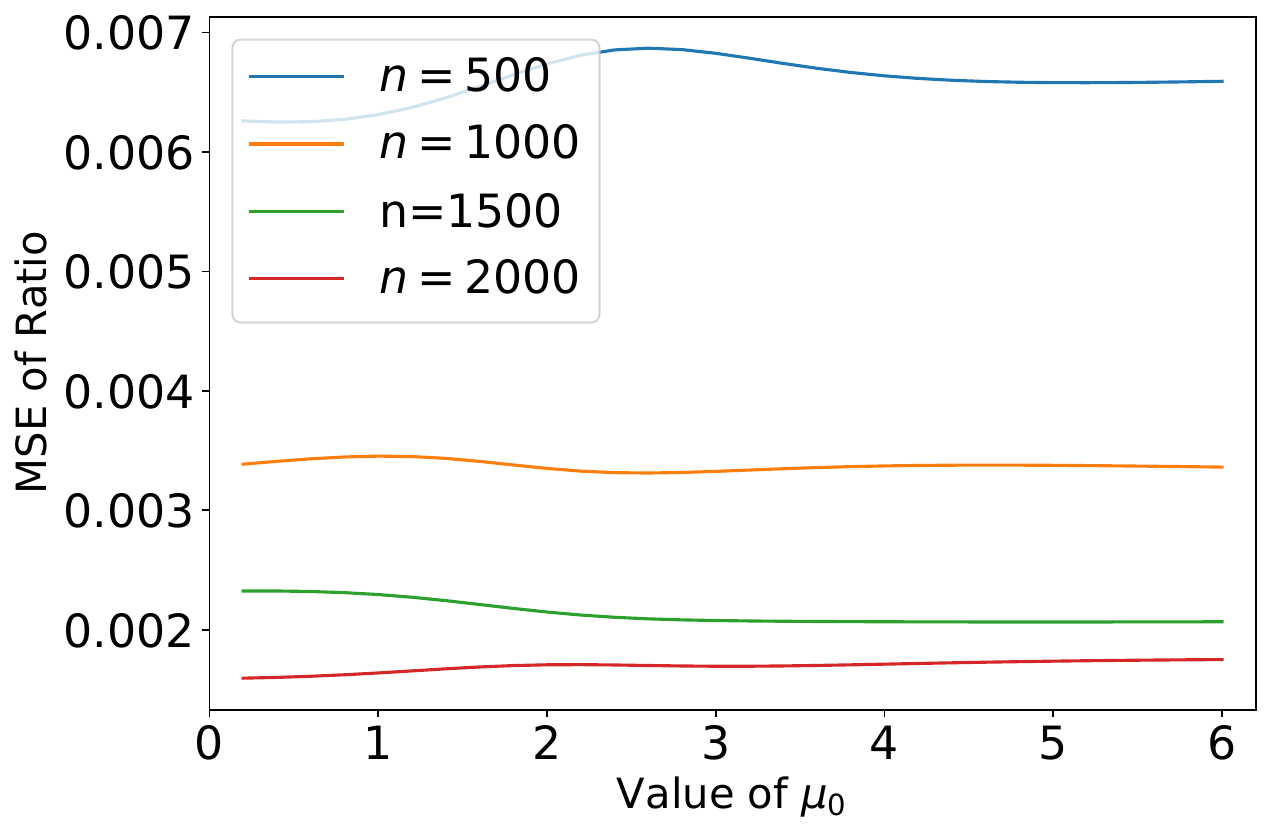}\label{fig_simu:frontier_asymp_mseratiosmallmu4}}
\subfloat[$\bmu=\bmu_3,c=3/2$]{\includegraphics[width=0.235\textwidth]{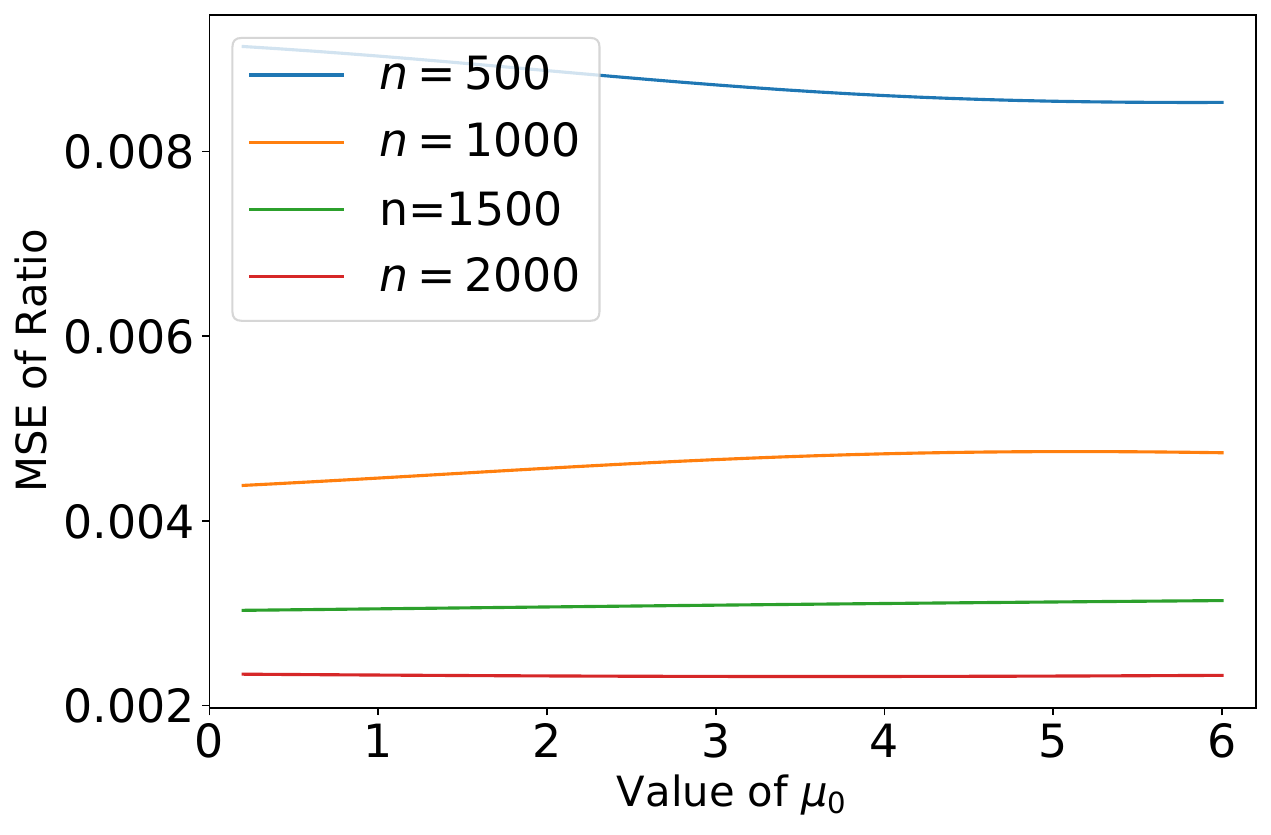}\label{fig_simu:frontier_asymp_mseratiobigmu3}}
\subfloat[$\bmu=\bmu_4,c=3/2$]{\includegraphics[width=0.235\textwidth]{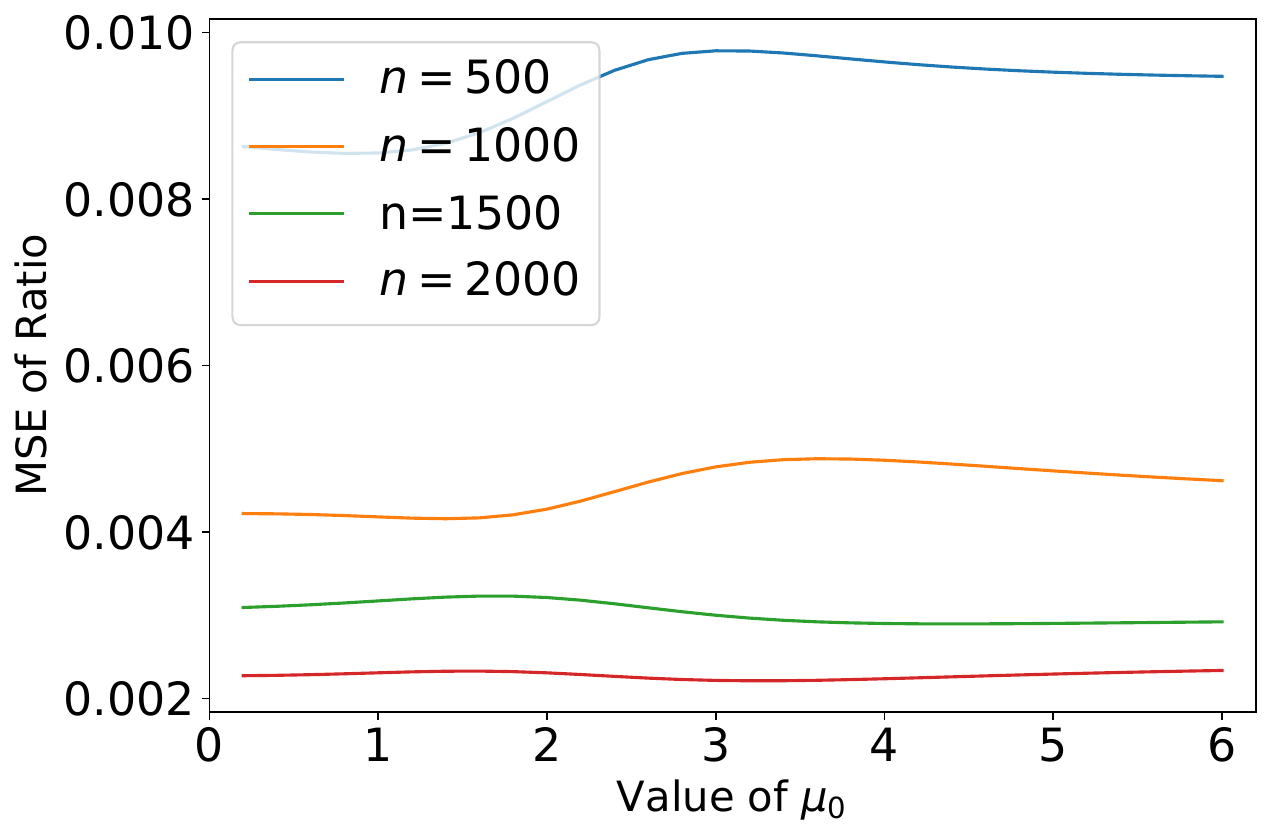}\label{fig_simu:frontier_asymp_mseratiobigmu4}}

\subfloat[$\bmu=\bmu_3,c=1/2$]{\includegraphics[width=0.235\textwidth]{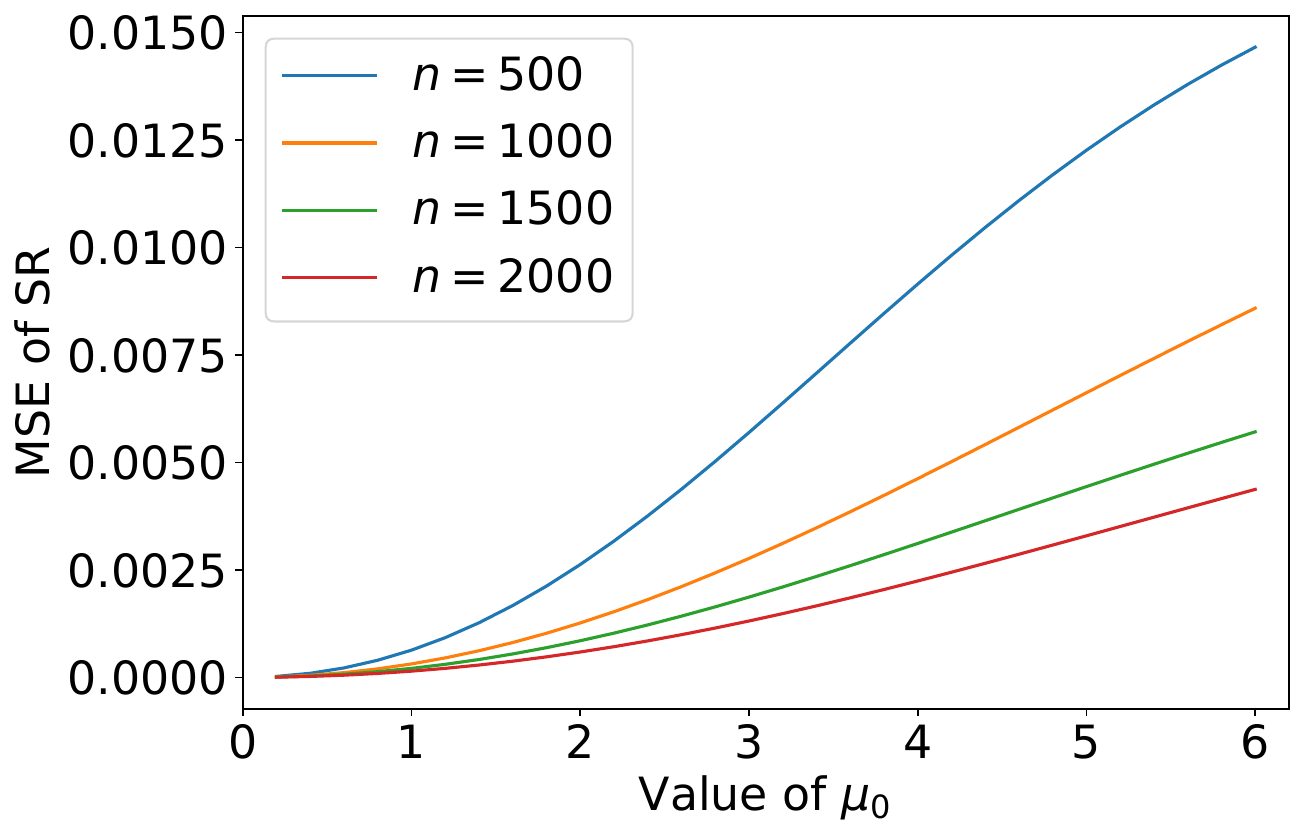}\label{fig_simu:frontier_asymp_mseSRsmallmu3}}
\subfloat[$\bmu=\bmu_4,c=1/2$]{\includegraphics[width=0.235\textwidth]{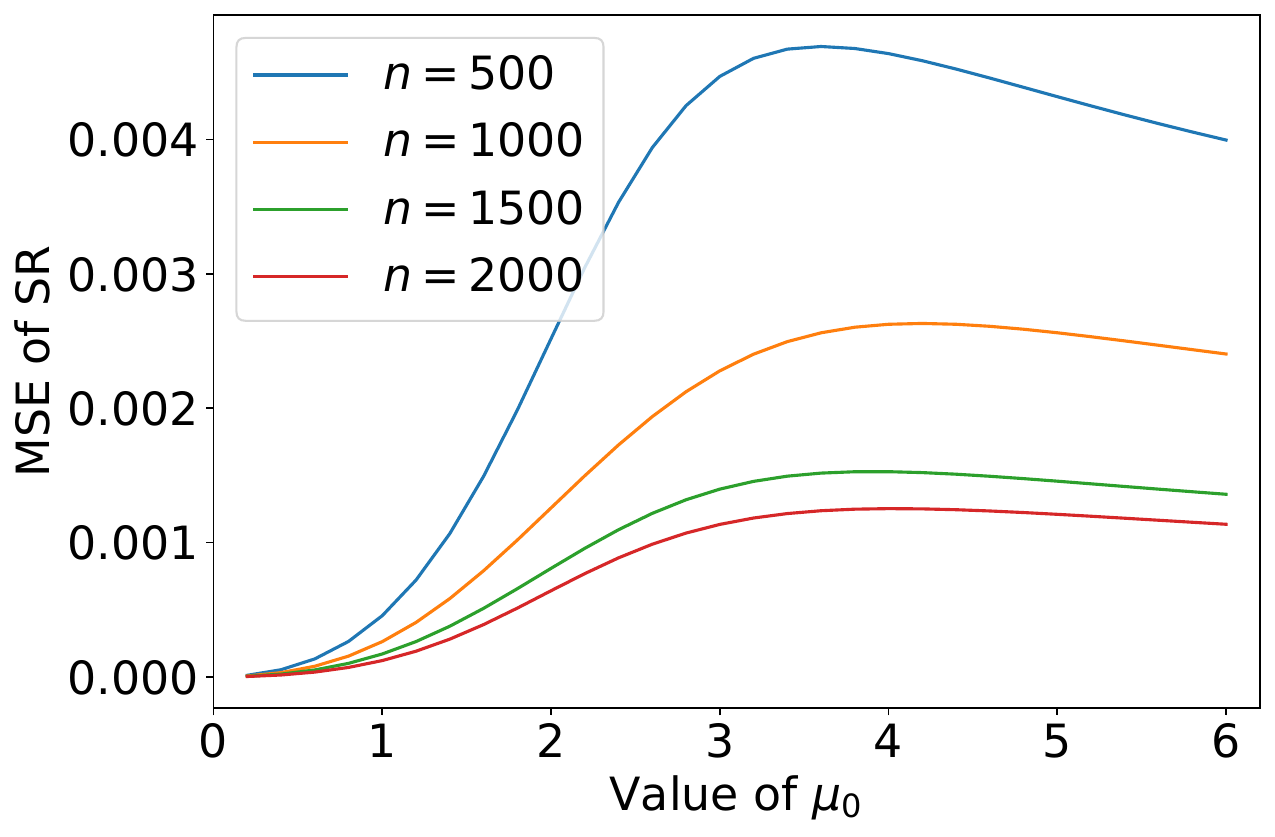}\label{fig_simu:frontier_asymp_mseSRsmallmu4}}
\subfloat[$\bmu=\bmu_3,c=3/2$]{\includegraphics[width=0.235\textwidth]{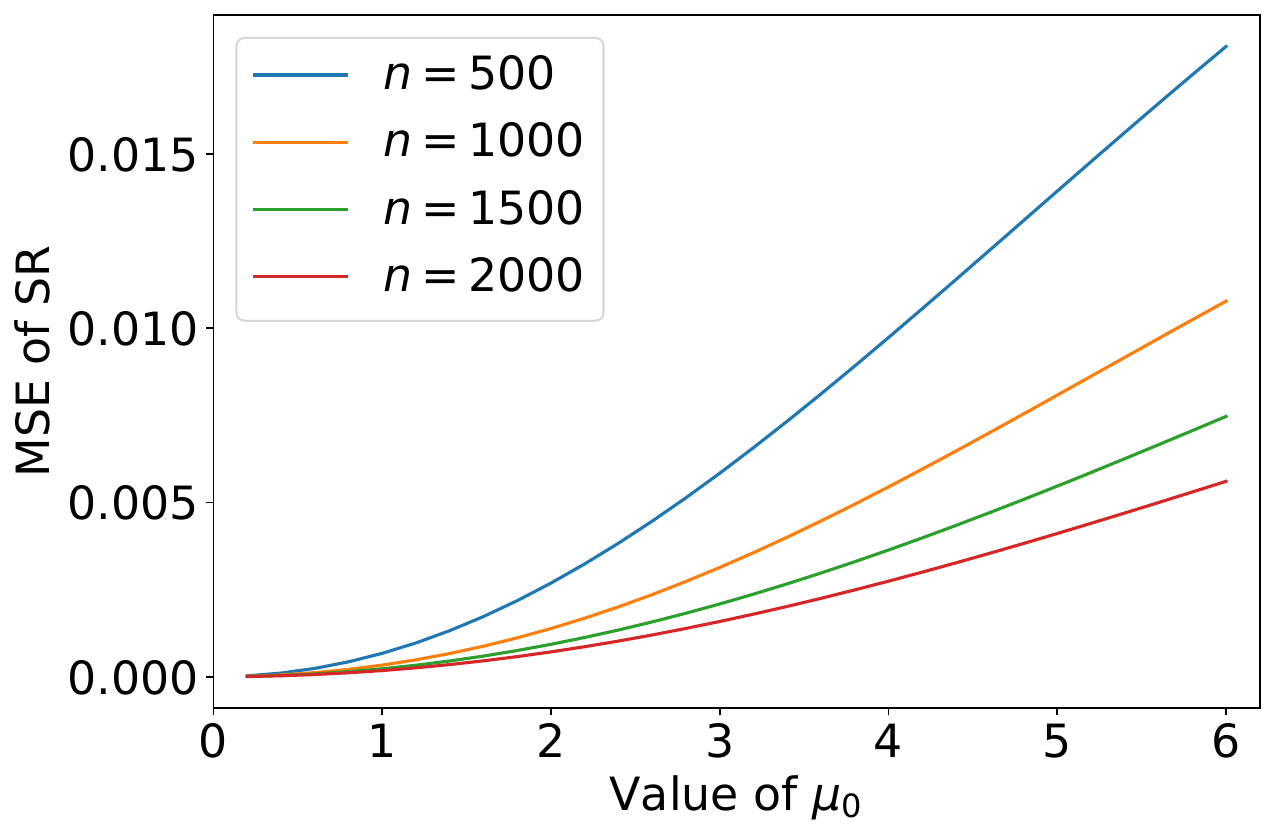}\label{fig_simu:frontier_asymp_mseSRbigmu3}}
\subfloat[$\bmu=\bmu_4,c=3/2$]{\includegraphics[width=0.23\textwidth]{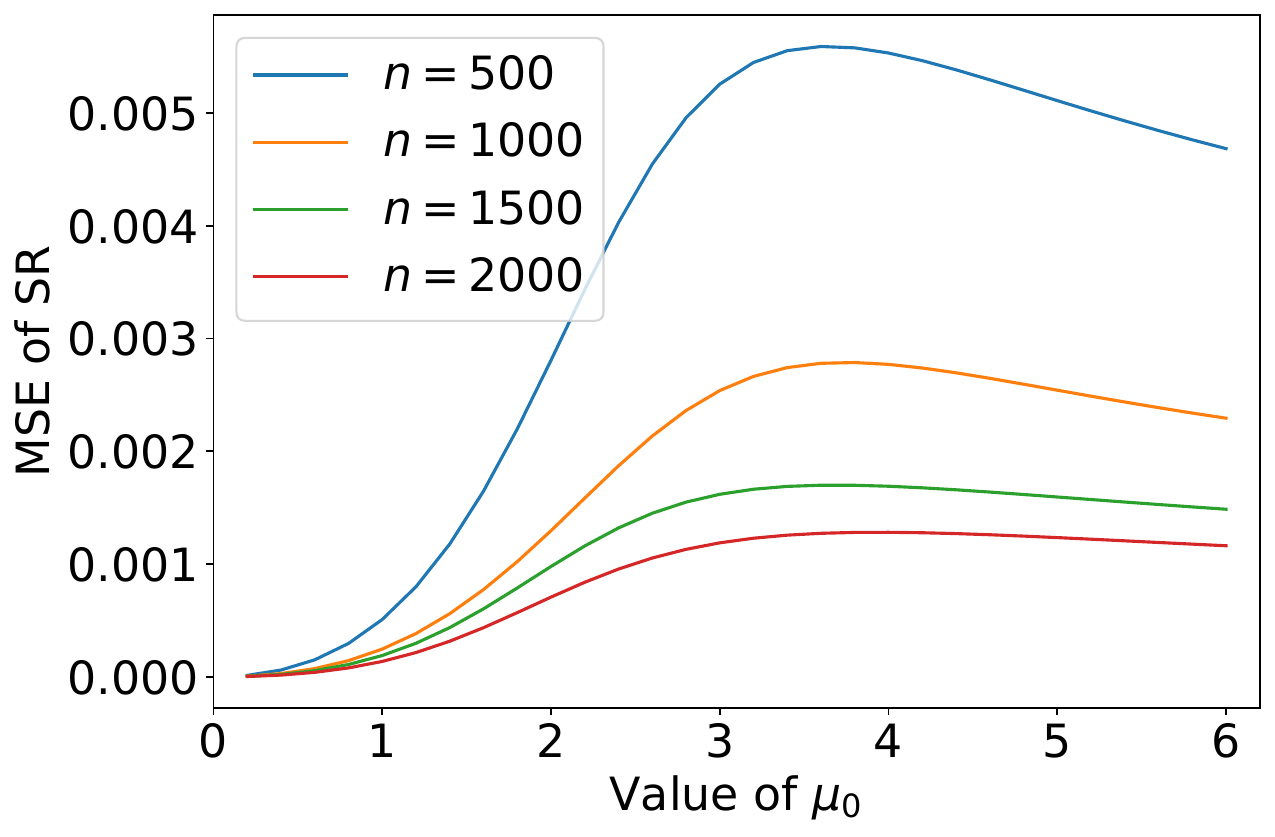}\label{fig_simu:frontier_asymp_mseSRbigmu4}}

\subfloat[$\bmu=\bmu_3,c=1/2$]{\includegraphics[width=0.235\textwidth]{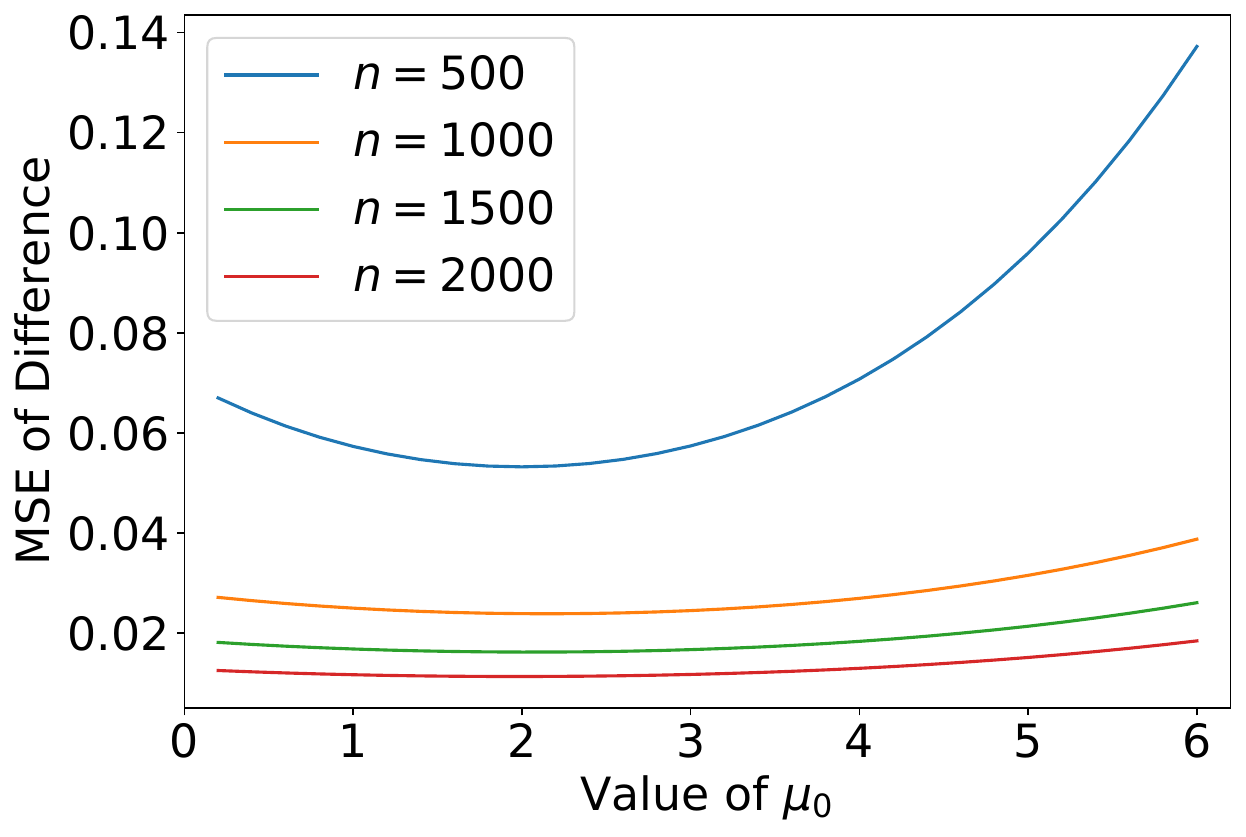}\label{fig_simu:frontier_asymp_msesmallmu3}}
\subfloat[$\bmu=\bmu_4,c=1/2$]{\includegraphics[width=0.235\textwidth]{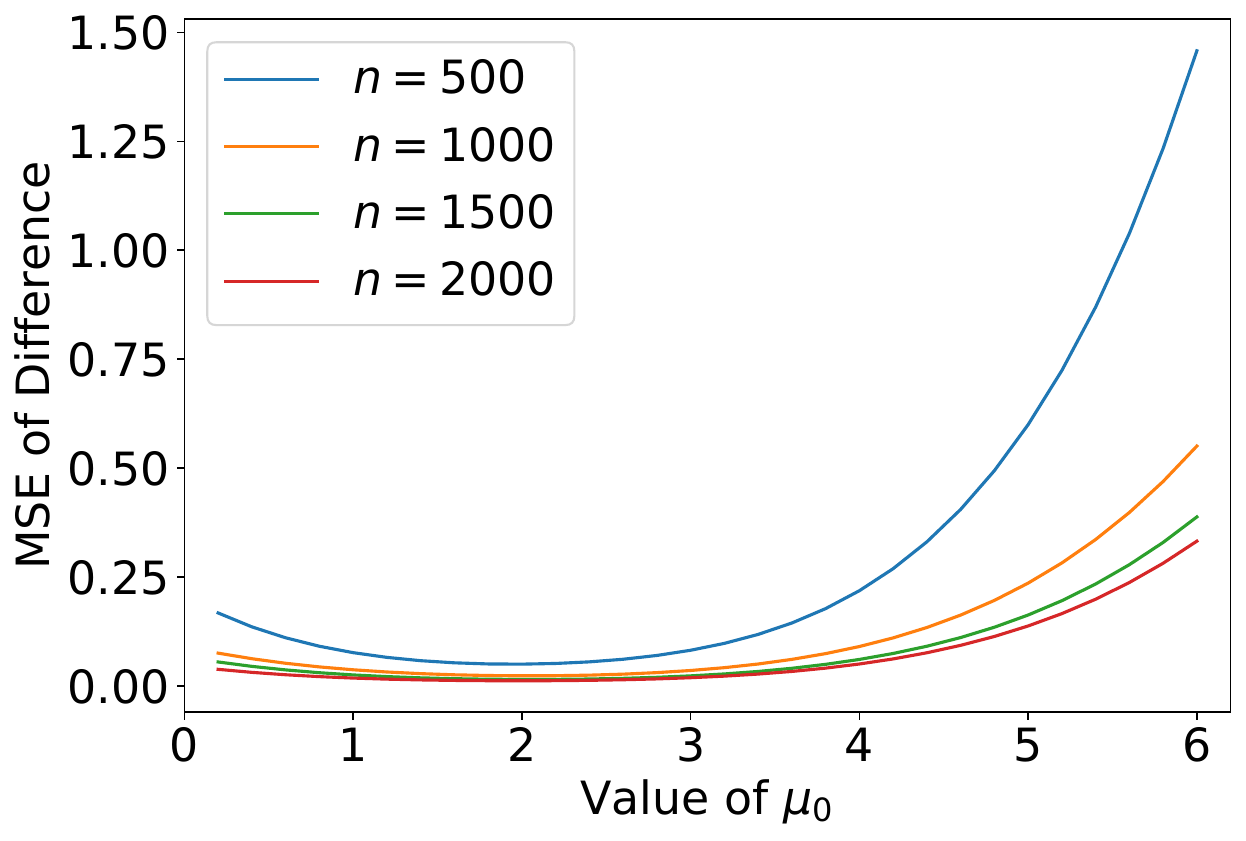}\label{fig_simu:frontier_asymp_msesmallmu4}}
\subfloat[$\bmu=\bmu_3,c=3/2$]{\includegraphics[width=0.235\textwidth]{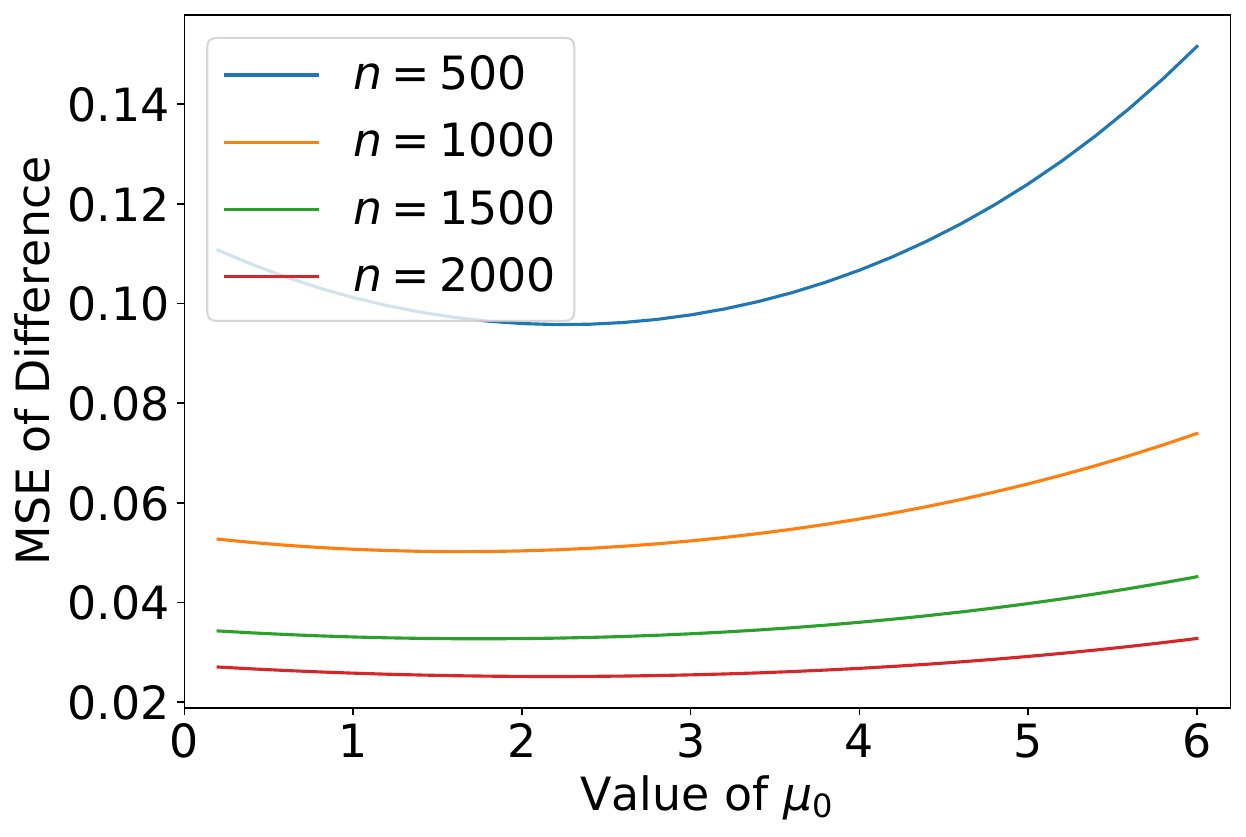}\label{fig_simu:frontier_asymp_msebigmu3}}
\subfloat[$\bmu=\bmu_4,c=3/2$]{\includegraphics[width=0.23\textwidth]{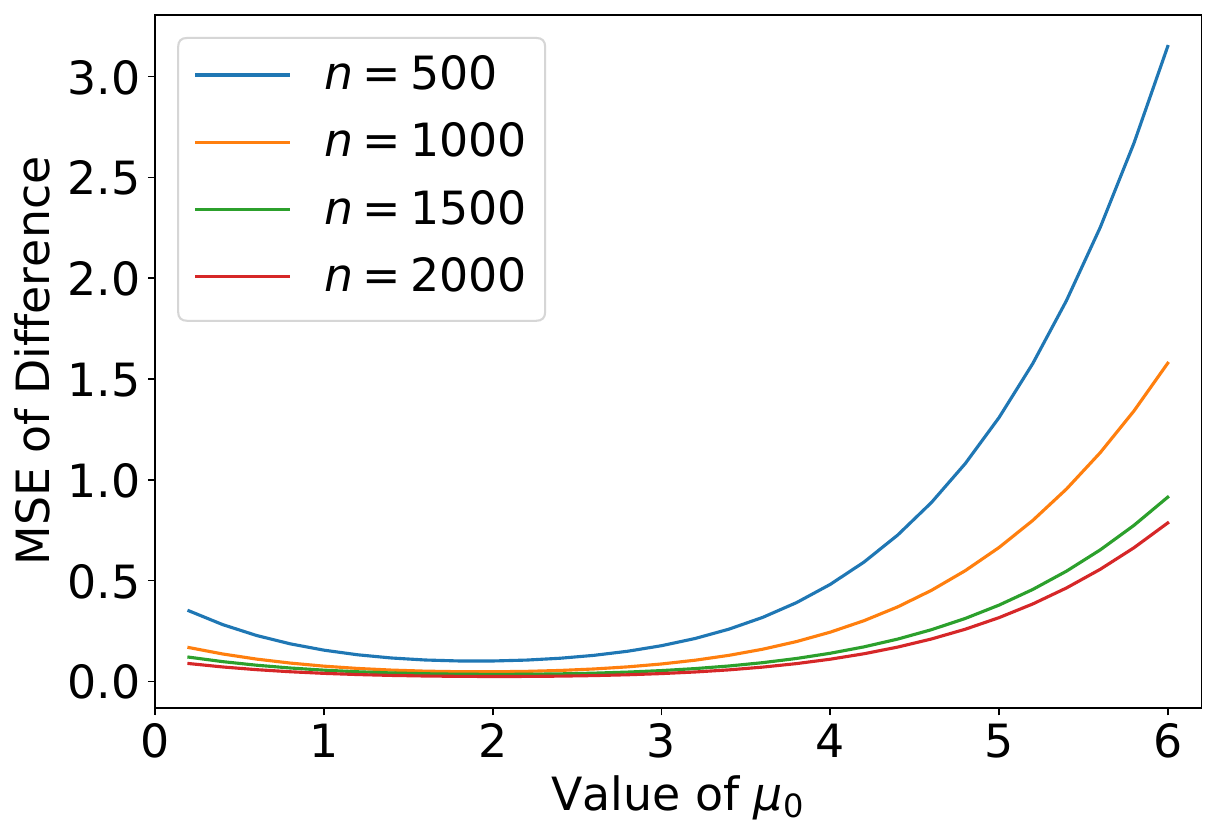}\label{fig_simu:frontier_asymp_msebigmu4}}
\caption{Simulation results with increasing $n$. $x$-axis in all figures shows different values of $\mu_0$. In Figures~\ref{fig_simu:frontier_asymp_mseratiosmallmu3}-\ref{fig_simu:frontier_asymp_mseratiobigmu4}, $y$-axis shows $\sum_{b=1}^{ 1000}(\hat{\sigma}_b^2/\sigma_{0,b}^2-1)^2/{ 1000}$. In Figures~\ref{fig_simu:frontier_asymp_mseSRsmallmu3}-\ref{fig_simu:frontier_asymp_mseSRbigmu4}, $y$-axis shows $\sum_{b=1}^{ 1000}(\frac{\mu_0}{\hat{\sigma}_b}-\frac{\mu_0}{\sigma_{0,b}})^2/{ 1000}$. And in Figures~\ref{fig_simu:frontier_asymp_msesmallmu3}-\ref{fig_simu:frontier_asymp_msebigmu4}, $y$-axis shows $\sum_{b=1}^{ 1000}(\hat{\sigma}_b^2-\sigma_{0,b}^2)^2/{ 1000}$. }
\label{fig_simu:frontier_asymp}
\end{figure}

Similar to the previous Subsection~\ref{subsubsec:simu_asymp_SR}, we investigate the effect of varying sample size $n$ and data dimension $p$ on  volatility predictions on the efficient frontier. Again we set $c=1/2$ or $c=3/2$, with $n = 500, 1000, 1500, 2000$, $\bSigma = \bSigma_3$, and $\bmu=\bmu_3$ or $\bmu_4$. 
Following the same procedure as described above in Section~\ref{sec:simu_frontier}, for each $(n,p)$ and  given $\mu_0$, we obtain  $1000$ values of $\sigma_0^2$ and $\hat{\sigma}^2$, denoted as $\sigma_{0,b}^2$ and $\hat{\sigma}^2_b$ for $b\in [1000]$. The outcomes are subsequently illustrated in Figure~\ref{fig_simu:frontier_asymp}. The MSE's of the ratio $\hat{\sigma}^2/\sigma_{0}^2-1$, the SR difference $\mu_0/\hat\sigma-\mu_0/\sigma_0$, and the variance difference $\hat{\sigma}^2-\sigma_{0}^2$ are illustrated in Figure~\ref{fig_simu:frontier_asymp} as we change the sample size $n$.
In Figure~\ref{fig_simu:frontier_asymp}, we see that the discrepancy between $ \sigma_0 $ and $ \hat{\sigma} $ decreases as the sample size $ n $ increases as expected. 
The ratio consistently remains around 1, as shown in Figures~\ref{fig_simu:frontier_asymp_mseratiosmallmu3}-\ref{fig_simu:frontier_asymp_mseratiobigmu4}, supporting the ratio consistency in Theorem~\ref{thm:frontier}.
When examining Sharpe ratio differences, Figures~\ref{fig_simu:frontier_asymp_mseSRsmallmu3} and \ref{fig_simu:frontier_asymp_mseSRbigmu3} show that the differences increase as $ \mu_0 $ grows, whereas Figures~\ref{fig_simu:frontier_asymp_mseSRsmallmu4} and \ref{fig_simu:frontier_asymp_mseSRbigmu4} show small differences even for larger $ \mu_0 $. This aligns with our theoretical result in Item 1 of Theorem~\ref{thm:frontier} as $ \cA_{rr} $ is unbounded for $ \bmu = \bmu_3 $ and bounded for $ \bmu = \bmu_4 $.
Regarding differences between $ \hat{\sigma}^2 $ and  $ \sigma_0^2$, Figures~\ref{fig_simu:frontier_asymp_msesmallmu3} and \ref{fig_simu:frontier_asymp_msebigmu3} show small differences for larger $ \mu_0 $ when $ \cA_{rr} $ is unbounded so that $ \mu_0 $ can take a large value and still ensures the convergence of $ \hat{\sigma}^2 $ to $ \sigma_0^2 $. In contrast, Figures~\ref{fig_simu:frontier_asymp_msesmallmu4} and \ref{fig_simu:frontier_asymp_msebigmu4} show larger differences for larger $\mu_0$ when $ \cA_{rr} $ is bounded. This matches our result in Item 2 of Theorem~\ref{thm:frontier}.

\section{Real Data Analysis}
\label{sec:realdata}
In this section, we conduct real data experiments to demonstrate the effectiveness and applicability of our methodology across various scenarios. 
The MV portfolio in Section~\ref{subsec:meanvariance} and the efficient frontier in Section~\ref{subsec:realdata_frontier} are both constructed using daily stock returns from the constituents of S\&P500 assuming we know the out-of-sample $\bmu$. In addition, we consider the global minimum variance portfolio in Appendix \ref{subsec:globalmin}, which has the advantage
of not requiring the knowledge of $\bmu$, and calibrated models in Appendix \ref{sec:appendixB2} to show the effect of estimating $\bmu$ in the MV portfolio.  Last but not least, we consider the MV portfolio with historical sample mean $\hbmu$ in Appendix~\ref{subsec:meanvarianceunknown} following Remark \ref{remark:unknown_mu}.

After deleting stocks with missing values, we end up $p=365$ stocks in total. 
Portfolios are built using historical data spanning different time length, including one year, two years, and four years, and rebalanced monthly to reflect changing market dynamics. 
To illustrate the portfolio construction, we give an example of building up the portfolio for Jan 2017 using one-year data. We utilize historical data from Jan 2016 to Dec 2016 to determine the optimized portfolio positions (i.e. the allocation vector $\wb^*$). These positions are then held for the entire Jan 2017. We refer to Jan 2017 as the current {\it testing month} and calculate the excess returns of the portfolio $\wb^*$ for each trading day in Jan 2017. We repeat the procedure in a rolling fashion for all testing months spanning from Jan 2013 to Jun 2023 and record the daily returns for each trading day.  The candidate set for the regularization matrix $\cQ_1$ is defined as $\cQ_1 = \{q \cdot \hbSigma_{pre}, q \in [1:30]/10\}$, where $\hbSigma_{pre}$ represents the sample covariance matrix of returns calculated from Jan 2004 to Dec 2008, which is a pre-training period not overlapping with data for portfolio construction and evaluation. 
We also consider a second candidate set $\cQ_2 = \{q \cdot \Ib_p, q \in [1:30]/10\}$, where $\Ib_p$ is the identity matrix of dimension $p$.

\subsection{Mean variance portfolio}
\label{subsec:meanvariance}
In this section, we construct the MV portfolios in order to achieve high Sharpe ratio. Recall that the MV portfolio allocation vector is  $\wb\propto (\hbSigma+\Qb)^{-1}\bmu$. We set its scaling by assuming $\sum_{i}|\wb_i|=1$, meaning that the book of the portfolio is one unit. 
The weight vector $\wb$ is 
\begin{align}
 \wb= {(\hbSigma+\Qb)^{-1}\bmu}/{\|(\hbSigma+\Qb)^{-1}\bmu\|_1}.
 \label{eq:meanvar_port}
\end{align}
Based on Theorem~\ref{thm:main_theorem}, we can utilize $\hat{SR}$ to optimize  $\Qb\in\cQ_1$ or $\Qb\in\cQ_2$ to yield a high $SR$.
We present the steps for our MV portfolio analysis.
\begin{enumerate}[leftmargin=*]
    \item For each testing month, we observe the historical daily returns $\Rb\in\RR^{n\times p}$, where $n$ is the total number of trading days with one-, two-, or four-year data and $p=365$ is the number of selected stocks. Calculate $\hbSigma$ and the portfolio $\wb$ by \eqref{eq:meanvar_port}.
    \item For each testing month, we run experiments for all candidate  $\Qb\in\cQ$ and also consider no regularization, i.e.  $\Qb=\0$, where we have $\wb\propto \hbSigma^{+}\bmu$ and $\hbSigma^{+}$ is the pseudo inverse, and the optimized  $\Qb^*\in\cQ$ using the estimation in Theorem~\ref{thm:main_theorem}. Here, 
    \begin{align*}
         \Qb^*=\argmax_{\Qb\in\cQ}\bigg(1-\frac{c}{p}\tr\hbSigma(\hbSigma+\Qb)^{-1}\bigg)\cdot\frac{\bmu^\top(\hbSigma+\Qb)^{-1}\bmu}{\sqrt{\bmu^\top(\hbSigma+\Qb)^{-1}\hbSigma(\hbSigma+\Qb)^{-1}\bmu}}.
    \end{align*}
    \item  We roll the procedure above for all testing months. Note that the value  $\Qb^*$ changes from month to month. 
    With the weight vector $\wb$ using all  $\Qb\in\cQ$, $\Qb=\0$ or $\Qb^*$, we can then compute the portfolio returns for each trading day in the testing month.
    \item We report the realized Sharpe ratio of daily portfolio returns over the future three years. 
\end{enumerate}
\begin{figure}[t!]
    \centering
    \subfloat[One year, $\cQ = \cQ_1, c>1$]{\includegraphics[width=0.33\textwidth]{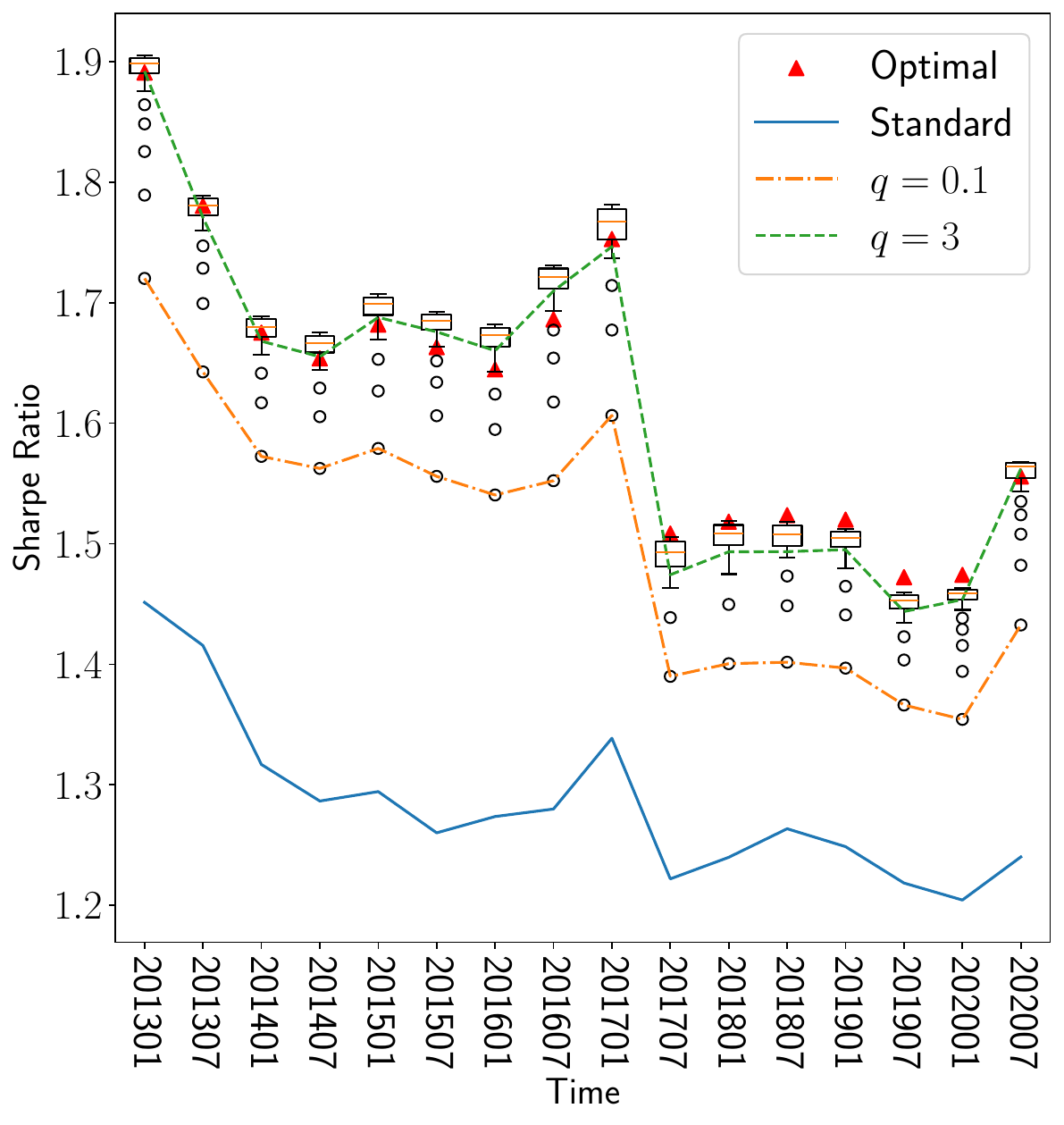}\label{fig_realmean:oneyear}}
    \subfloat[Two years, $\cQ = \cQ_1, c<1$]{\includegraphics[width=0.33\textwidth]{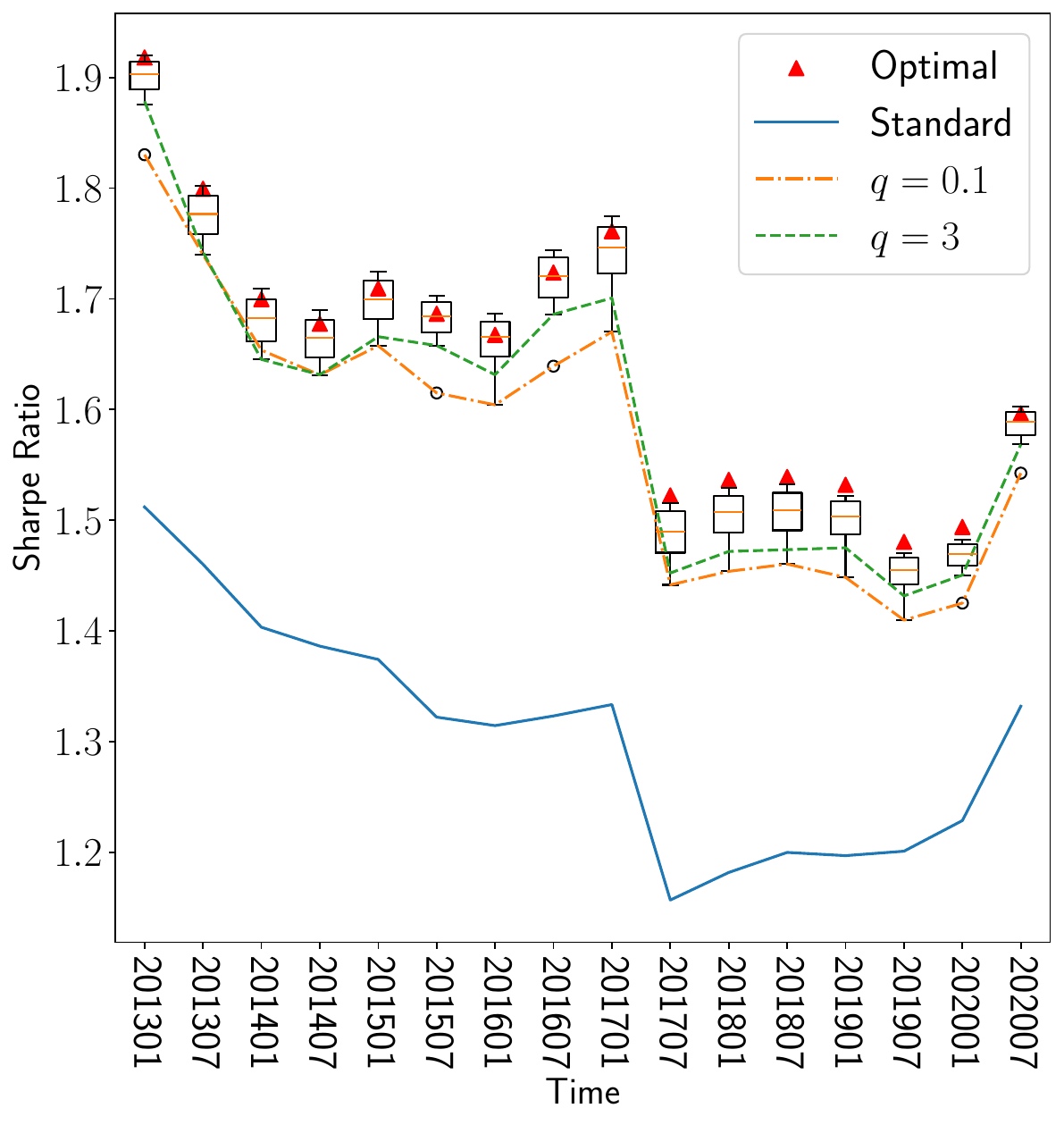}\label{fig_realmean:twoyear}}
    \subfloat[Four years, $\cQ = \cQ_1, c<1$]{\includegraphics[width=0.33\textwidth]{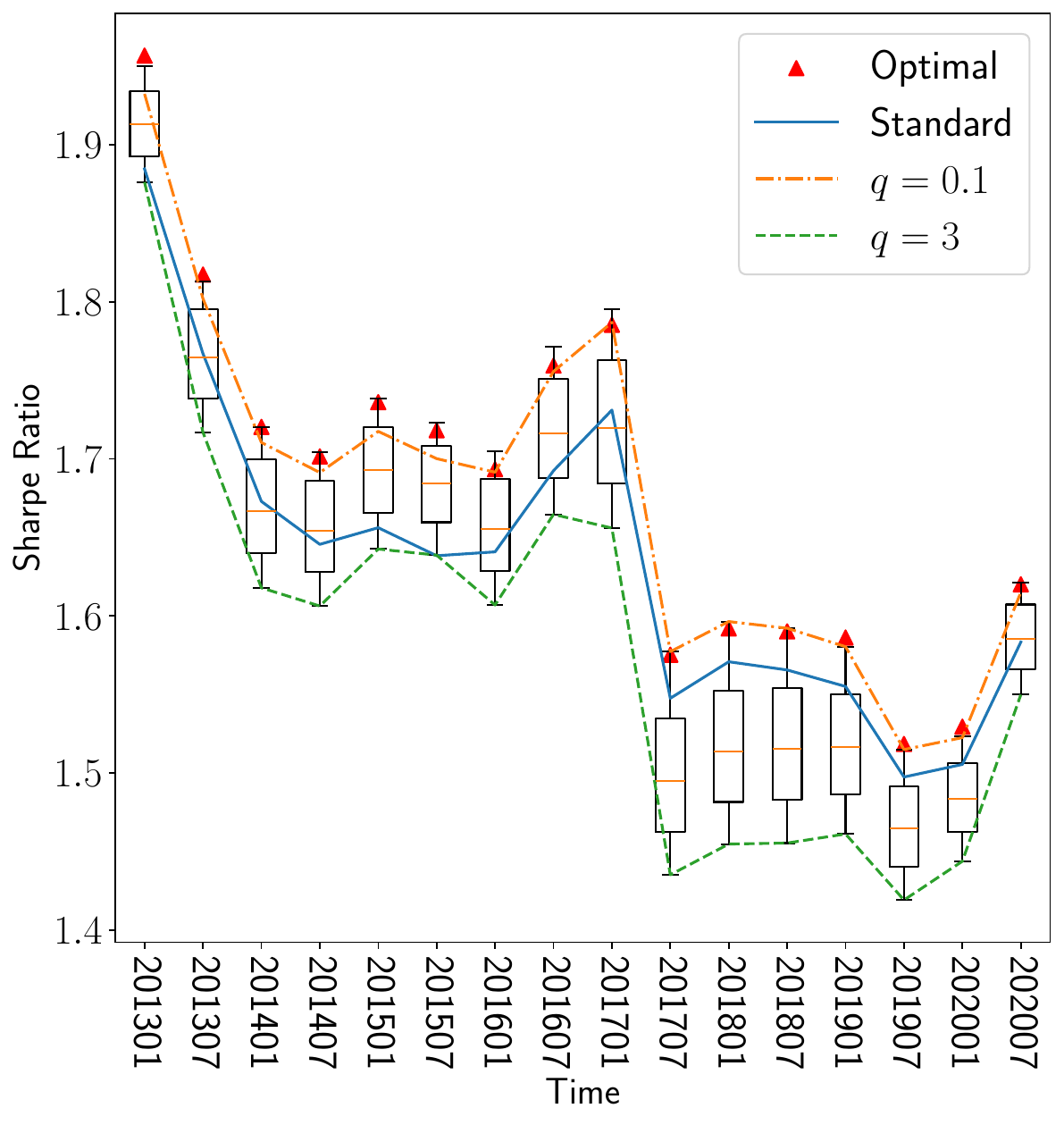}\label{fig_realmean:fouryear}}

    \subfloat[One year, $\cQ = \cQ_2, c>1$]{\includegraphics[width=0.33\textwidth]{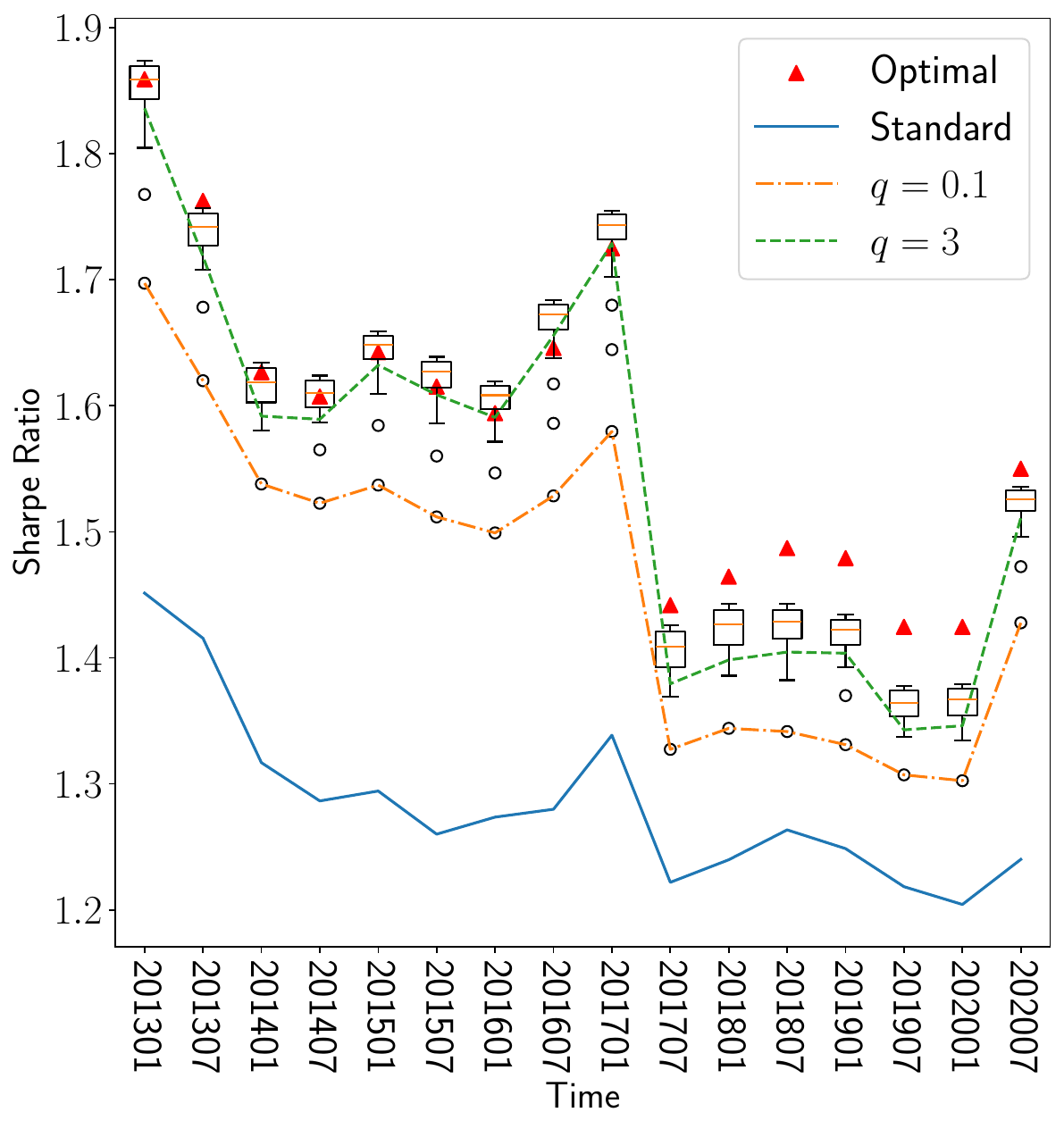}\label{fig_realmean:oneyear1}}
    \subfloat[Two years, $\cQ = \cQ_2, c<1$]{\includegraphics[width=0.33\textwidth]{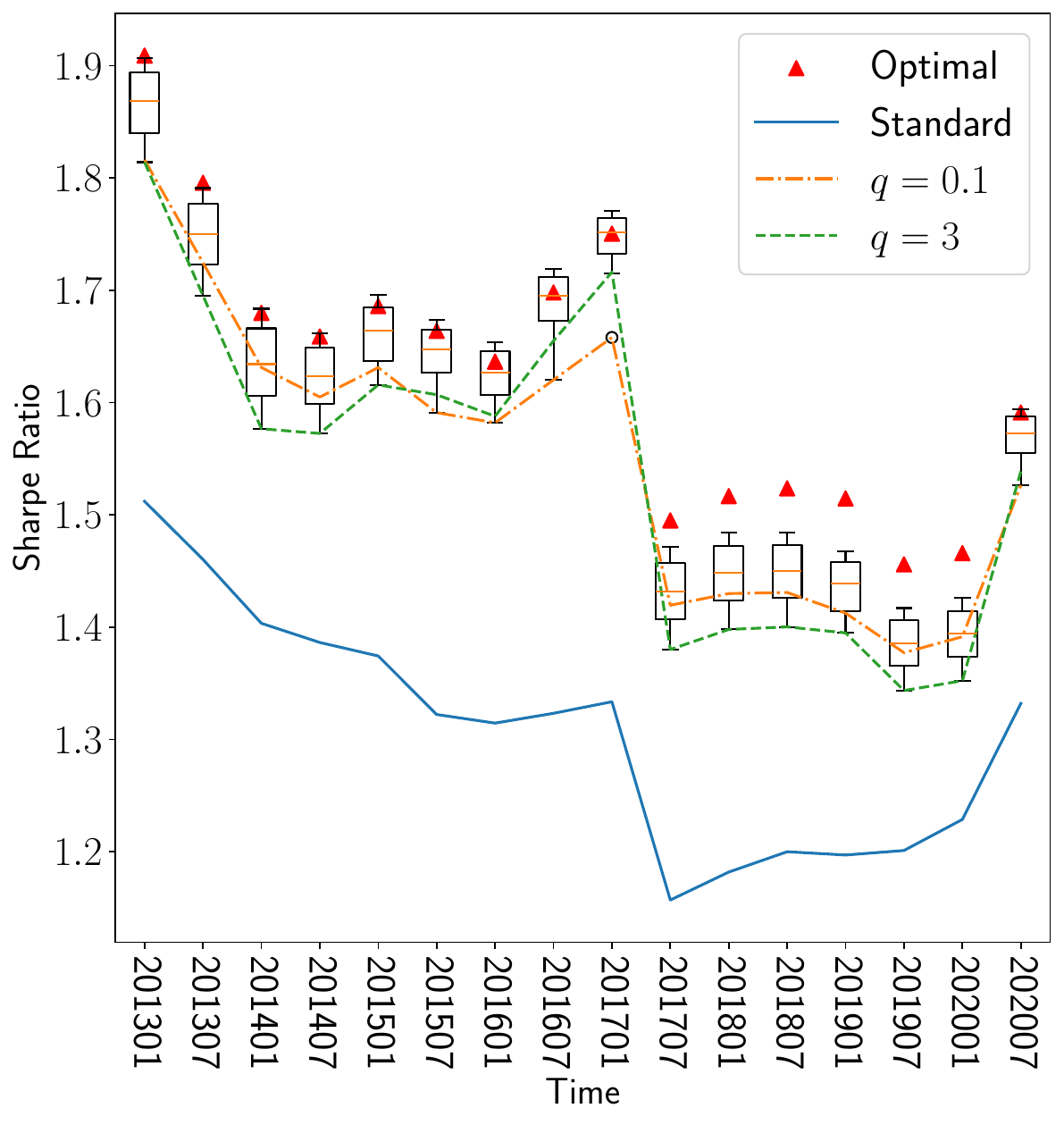}\label{fig_realmean:twoyear1}}
    \subfloat[Four years, $\cQ = \cQ_2, c<1$]{\includegraphics[width=0.33\textwidth]{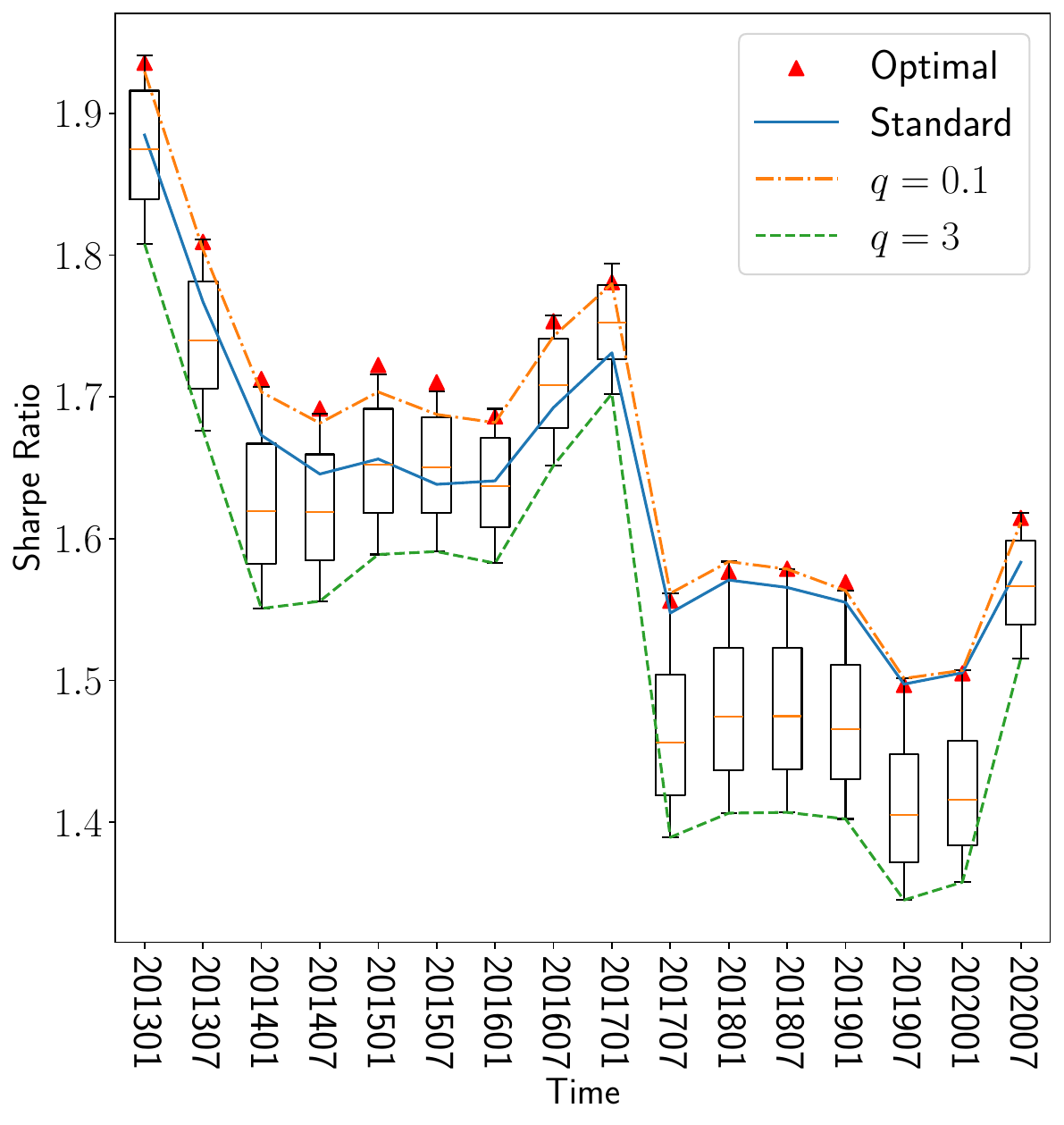}\label{fig_realmean:fouryear1}}
\caption{SR of mean-variance portfolios. The x-axis labels the rolling period, while the y-axis represents the out-of-sample SR of the portfolio returns every three years in the future. The blue solid line, orange dash-dot line and green dash line correspond to the SR with  $q=0$, $q=0.1$ (minimum $q$) and $q=3$ (maximum $q$) respectively. The boxplot displays SR for all  $\Qb\in\cQ$, and the red triangle indicates the SR under our optimized  $\Qb^*\in\cQ$.}
    \label{fig_realmean}
\end{figure}
In this experiment, we use the mean from each out-of-sample testing month as $\bmu$ to construct MV portfolio.   Readers may refer to Appendix~\ref{subsec:meanvarianceunknown} for real data experiments with unknown $\bmu$.
Figures \ref{fig_realmean:oneyear}-\ref{fig_realmean:fouryear} present the results obtained  with  $\Qb\in\cQ_1$, while Figures \ref{fig_realmean:oneyear1}-\ref{fig_realmean:fouryear1} present the results obtained with  $\Qb\in\cQ_2$. Obviously, when sample size is small, if we apply  $\Qb=\0$ (the blue curve) or   $q=0.1$ (the orange curve,  minimum $q$), we fail to attain the best possible SR.  If we use   $q=3$ (the green curve,  maximum $q$), it is also a bad choice when sample size is large. However, the optimized  $\Qb^*\in\cQ$ (the red triangles) consistently demonstrates superior performance throughout the experiments, regardless of the candidate sets and the varying lengths of historical data. From the results, we see clear advantage of actively optimized  $\Qb^*\in\cQ$ for each testing month over some ad-hoc constant  $\Qb\in\cQ$ for all testing months.   Moreover, the optimal regularization $\Qb^* \in \cQ_1$ performs better than  $\Qb^* \in \cQ_2$. For instance, comparing Figures \ref{fig_realmean:fouryear} and \ref{fig_realmean:fouryear1}, using the identical blue lines as the benchmark, the optimal SR values (red triangles) obtained from $\cQ_1$ are visibly higher than those obtained from $\cQ_2$. 

\subsection{Efficient frontier}
\label{subsec:realdata_frontier}
In this section, we test the real data application of Theorem~\ref{thm:frontier} to efficient frontier evaluation. For each specific $\mu_0$, we need to estimate the standard deviation of the regularized MV portfolio, presented in \eqref{eq:frontier2}. Here are the steps of our real data analysis.
\begin{enumerate}[leftmargin=*]
\item For each testing month, centralize the observed historicla daily returns $\Rb\in\RR^{n\times p}$, where $n$ corresponds to one-, two-, or four-year data and $p=365$. Then, calculate the sample covariance $\hbSigma$
from the observed $\Rb$. Let $\rb$ be the average return vector in the testing month. The optimal portfolio is given by $\wb^*=\gb+\mu_0\hb$, where $\gb$ and $\hb$ is given by \eqref{eq:def_gh}.
\item We run experiments for all  $\Qb$ in the candidate sets, including  $\Qb=\0$ and the optimized  $\Qb^*$. By Theorem~\ref{thm:frontier}, the optimized  $\Qb^*$  is obtained by minimizing   $\frac{(\gb+\mu_0\hb)^\top\hbSigma(\gb+\mu_0\hb)}{(1-c/p\cdot\tr \hbSigma(\hbSigma+\Qb)^{-1})^2}$ over all  $\Qb\in\cQ$. The case  $\Qb=\0$ replaces  $(\hbSigma+\Qb)^{-1}$ with $\hbSigma^+$ in $\wb^*$. 
\item  We repeat Step 2 for each testing month from Jan 2013 to Jun 2023. This allows us to collect daily portfolio returns for each  $\Qb$.  
\item We calculate the standard deviation of the daily returns for each  $\Qb$, including  $\Qb=\0$ and  $\Qb^*$, over the ten-year period. This allows us to generate a boxplot of the standard deviation for each candidate  $\Qb$ in $\cQ_1$ or $\cQ_2$ given $\mu_0$. 
\end{enumerate}
\begin{figure}[t!]
    \centering
    \subfloat[One year, $\cQ=\cQ_1$]{\includegraphics[width=0.33\textwidth]{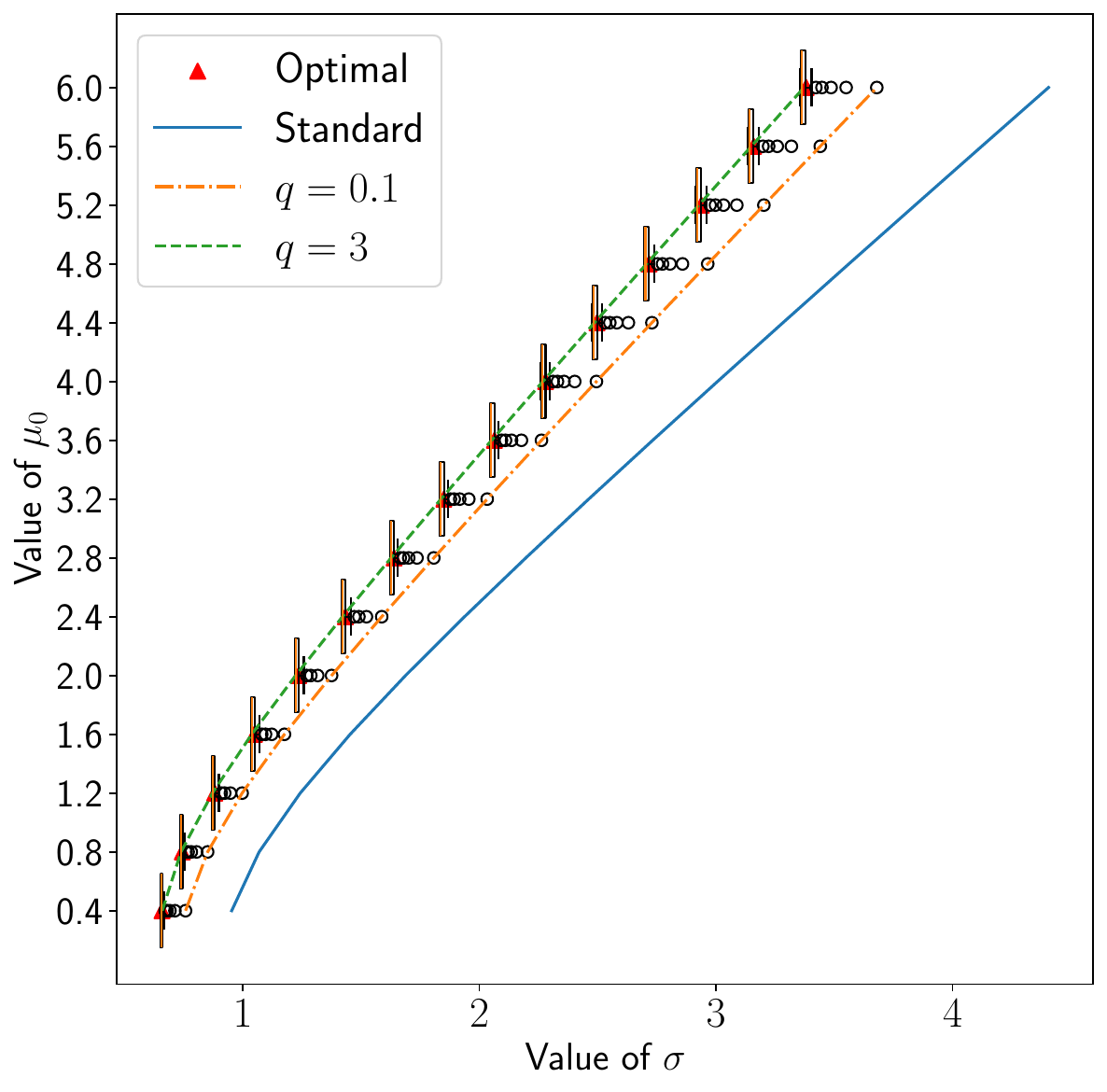}\label{fig_frontier:1yearQ1}}
    \subfloat[Two years, $\cQ=\cQ_1$]{\includegraphics[width=0.33\textwidth]{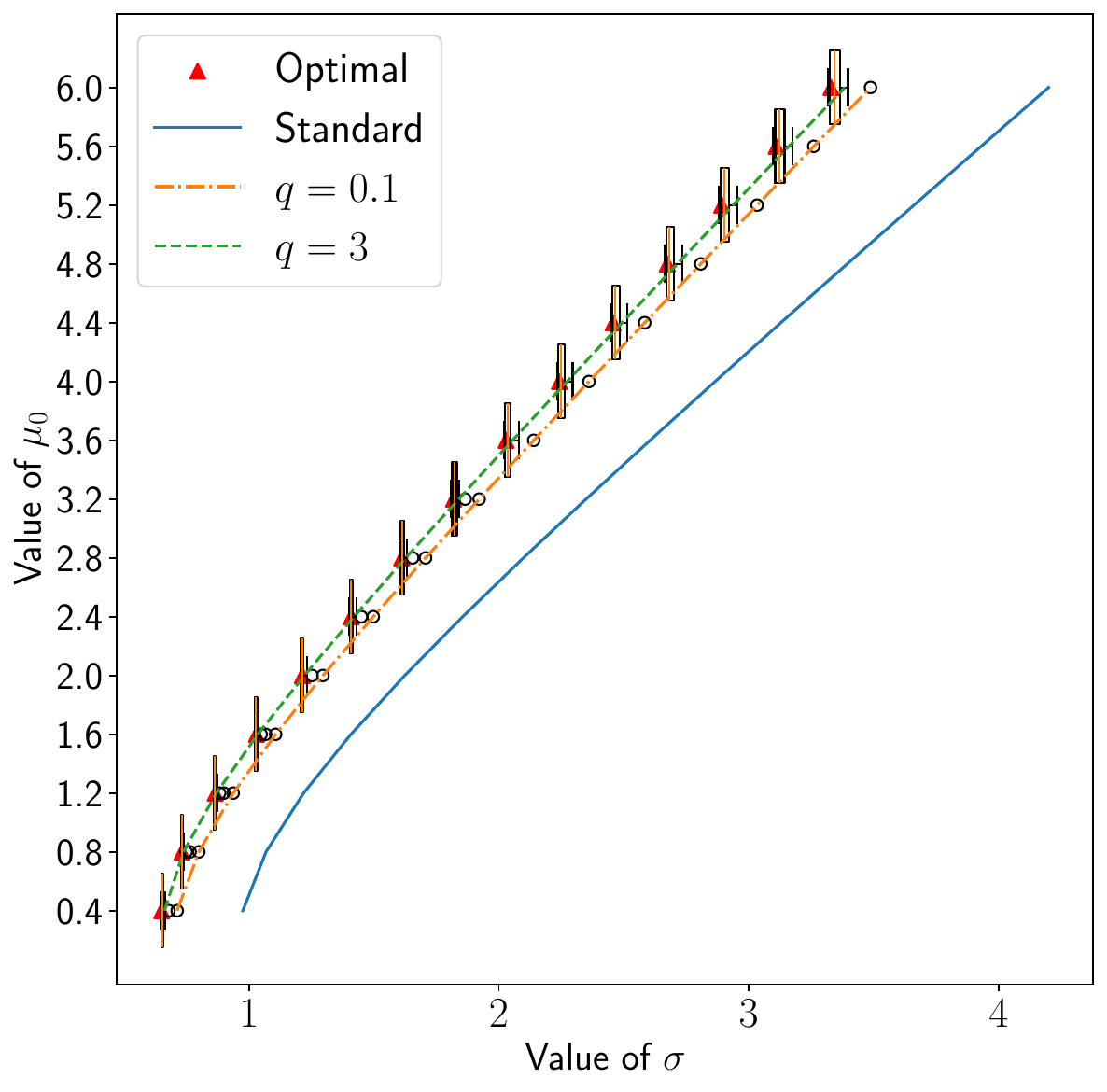}\label{fig_frontier:2yearQ1}}
    \subfloat[Four years, $\cQ=\cQ_1$]{\includegraphics[width=0.33\textwidth]{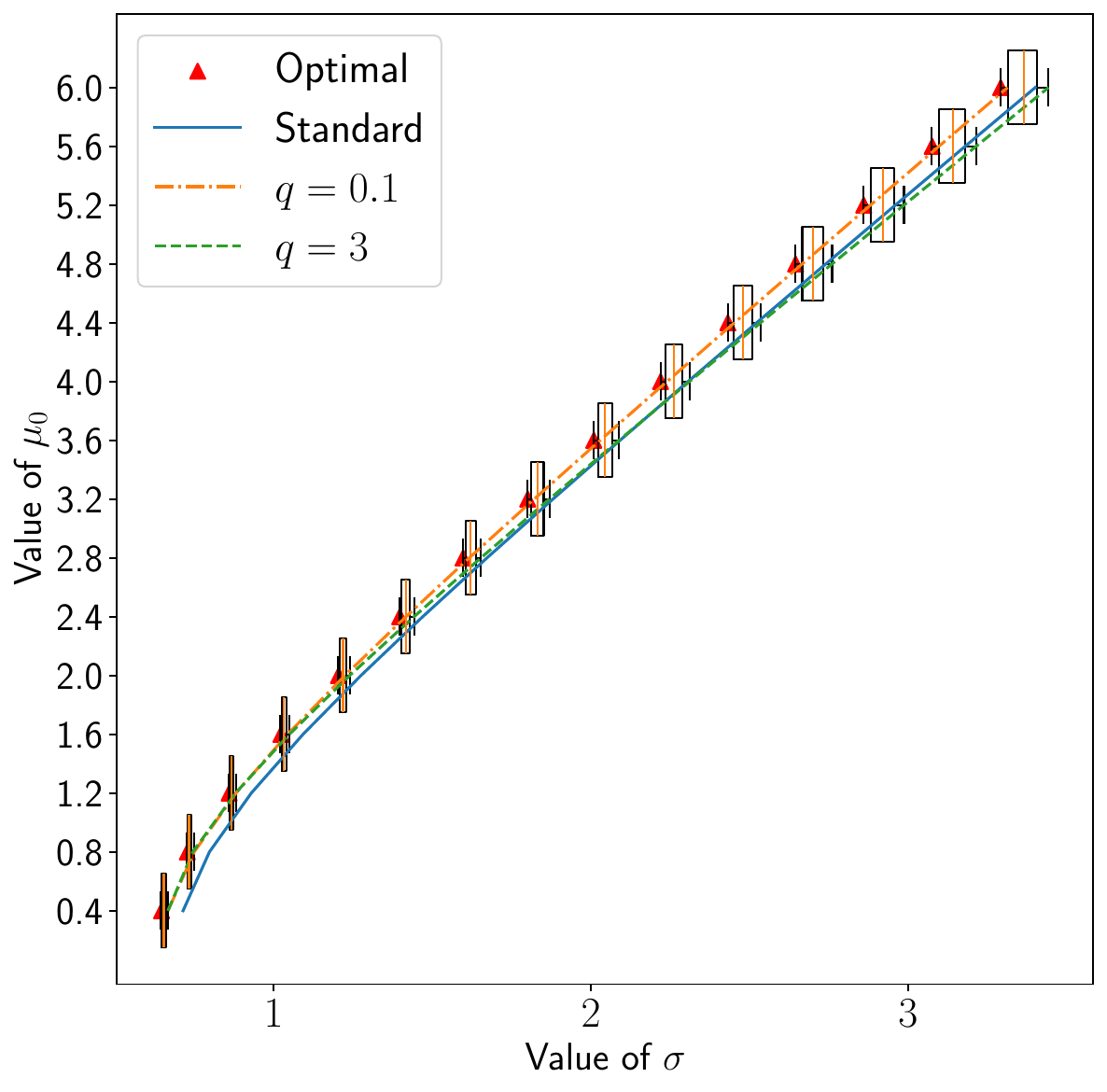}\label{fig_frontier:4yearQ1}}

    \subfloat[One year, $\cQ=\cQ_2$]{\includegraphics[width=0.33\textwidth]{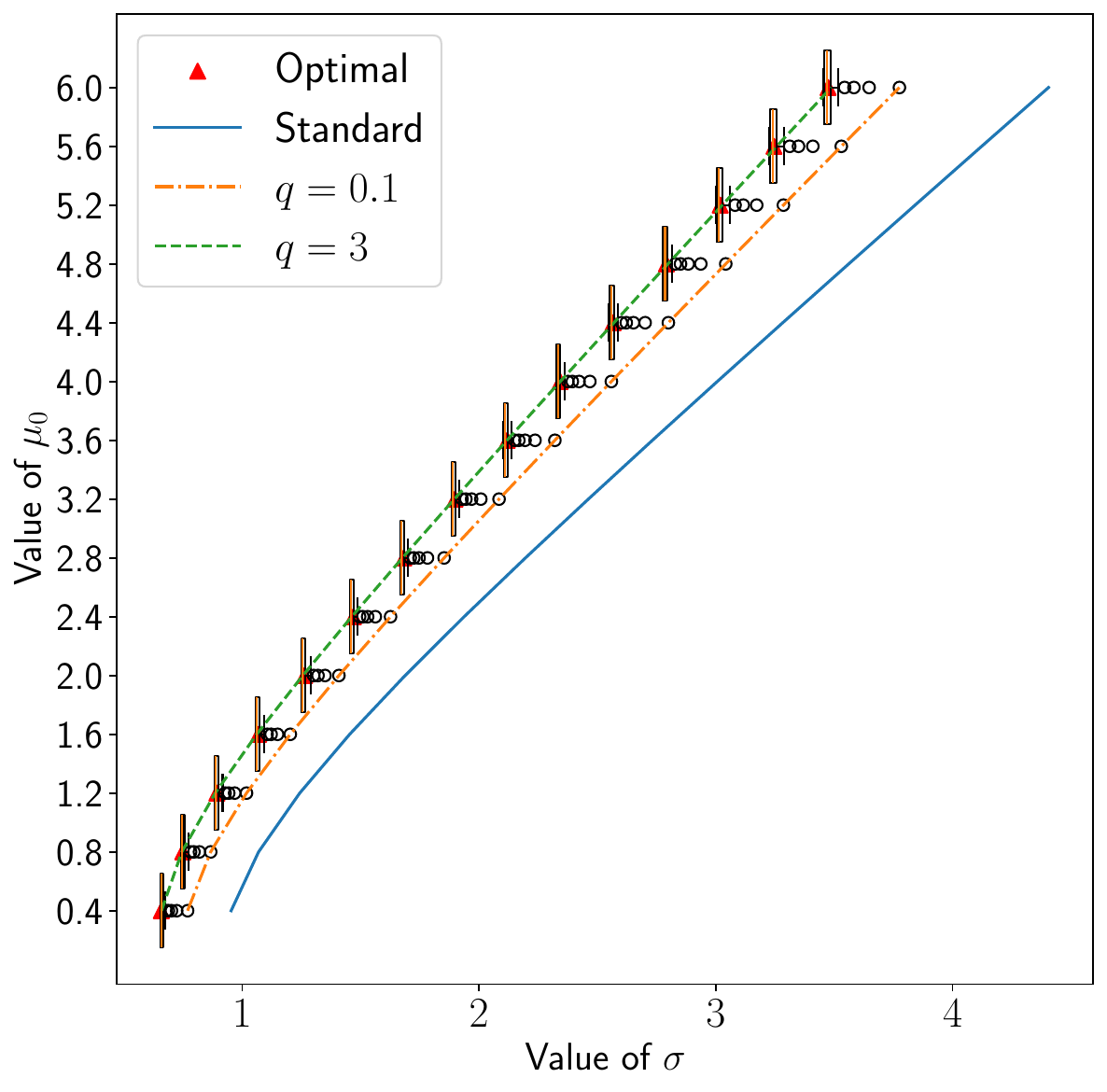}\label{fig_frontier:1yearQ2}}
    \subfloat[Two years, $\cQ=\cQ_2$]{\includegraphics[width=0.33\textwidth]{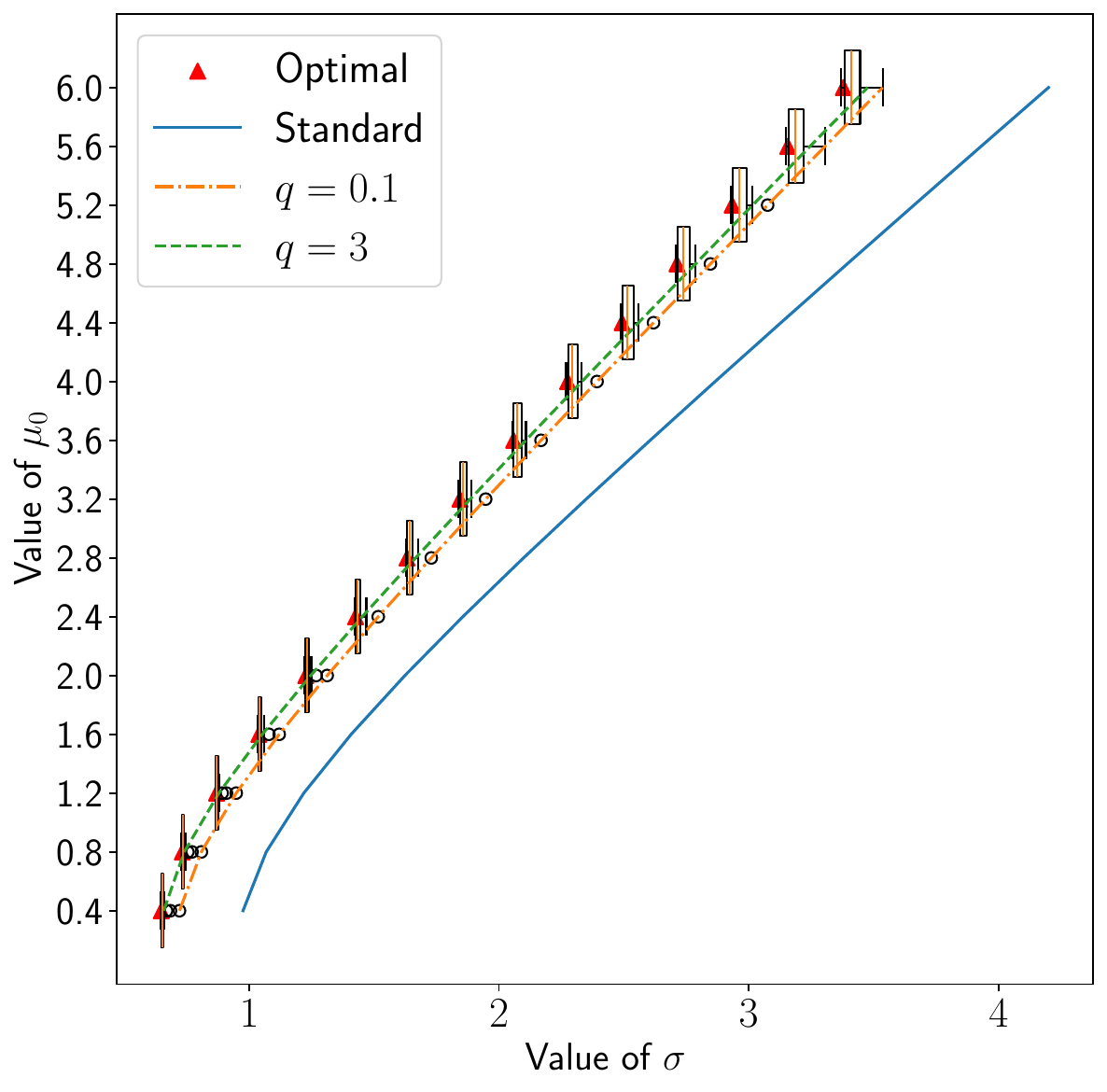}\label{fig_frontier:2yearQ2}}
    \subfloat[Four years, $\cQ=\cQ_2$]{\includegraphics[width=0.33\textwidth]{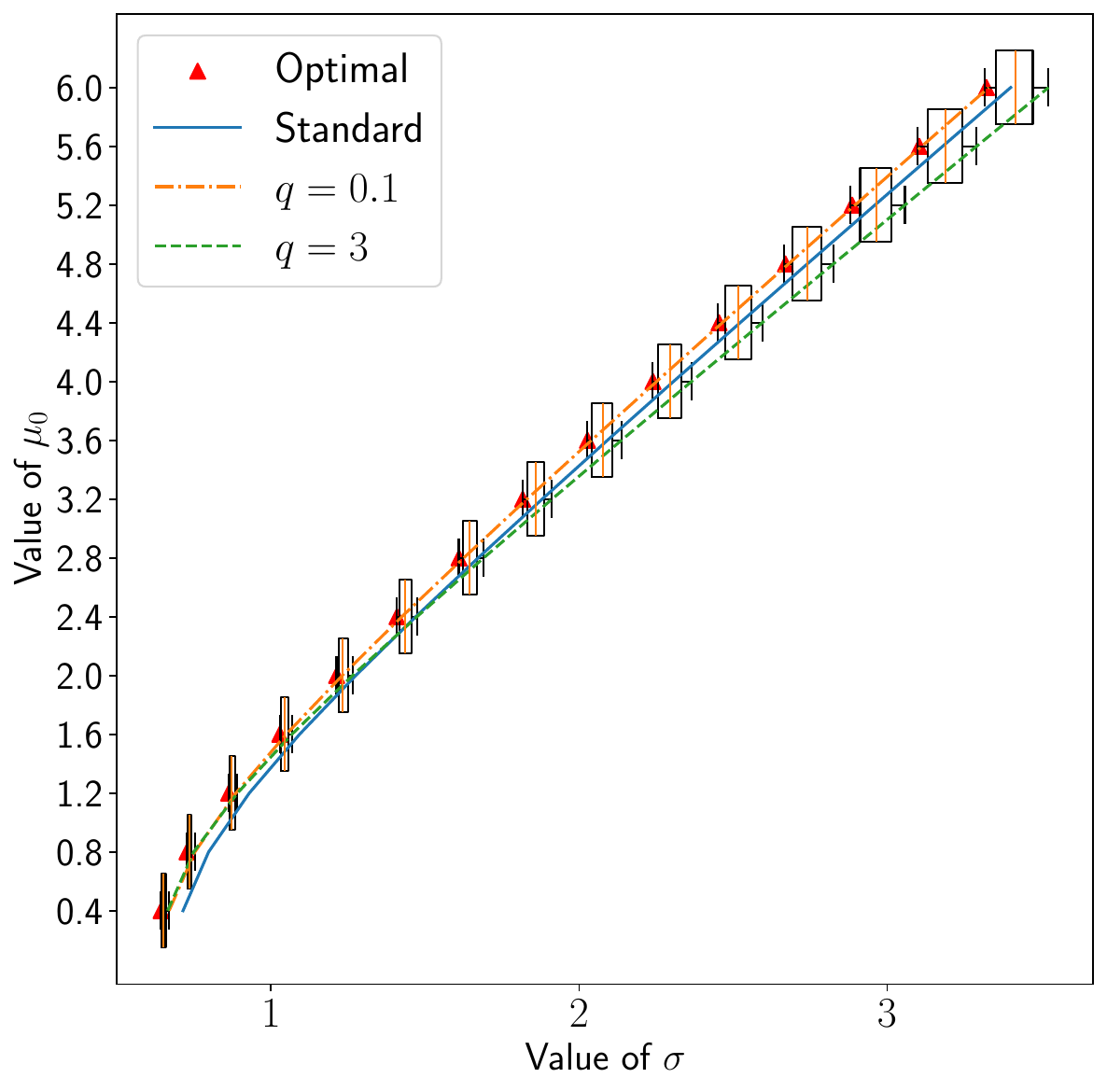}\label{fig_frontier:4yearQ2}}
\caption{The corrected efficient frontier. The x-axis represents the values of $\sigma$, while the y-axis represents  values of $\mu_0$. The blue solid line, orange dash-dot line and green dash line in the figure correspond to   $q=0$, $q=0.1$ (minimum $q$) and $q=3$ (maximum $q$) respectively. The boxplot displays the values of $\sigma$ for all  $\Qb\in\cQ$, and the red triangle indicates the volatility of the portfolio returns under our optimized  $\Qb^*\in\cQ$.}
    \label{fig_frontier:Q}
\end{figure}
As shown in Figure~\ref{fig_frontier:Q}, the standard approach with no regularization ($\Qb=\0$) substitutes $(\hbSigma+\Qb)^{-1}$ with $\hbSigma^+$ and exhibits a  larger volatility in comparison to the regularized portfolios for any given $\mu_0$. Furthermore, as illustrated in Figures~\ref{fig_frontier:1yearQ2} and \ref{fig_frontier:4yearQ2}, the minimum standard deviation occurs when $q$ is approximately  $3$ and $0.1$ respectively. However, our corrected estimator consistently identify the optimal  $\Qb\in\cQ_2$ and result in a smaller standard deviation. This again demonstrates the practical usefulness of Theorem~\ref{thm:frontier} in real efficient frontier recovery of the regularized portfolios.

\section{Conclusion}
\label{sec:conclusion}
In this work,  we propose a novel in-sample method to estimate the out-of-sample Sharpe ratio in high-dimensional portfolio optimization with or without a risk-free asset, under an unknown population covariance matrix $\bSigma$. 
By random matrix theory, we show that the proposed estimator shares the same asymptotic properties  with the out-of-sample   Sharpe ratio, even when  $\|\bSigma\|_{\op}$ is unbounded.  Specifically, we only need the  bounded $\|\bSigma/p\|_{\tr}$, which is more applicable to the real  financial data. We have also validated our theory through simulations and demonstrated its practical application through real data experiments. 


In future work, it would be  valuable to explore the  asymptotic distribution of the proposed estimator and  develop statistical inference, including studying the estimator’s correlation with the out-of-sample Sharpe ratio. A robust inference methodology would further enhance understanding and provide robust tools for evaluating the Sharpe ratio.


\bibliographystyle{agsm}

\bibliography{ref}

\newpage

\appendix
\bigskip
\begin{center}
{\large\bf Supplementary Material for ``Estimation of Out-of-Sample Sharpe Ratio for High Dimensional Portfolio Optimization''}
\end{center}

\begin{abstract}

This supplementary material contains additional experimental  results and all the proofs for theorems in the main text.

\end{abstract}

\section{Additional Numerical Experiments}
\label{sec:appendixA}
\subsection{Alternative settings}
\label{subsec:simu_alternative}
In this section, we change the values of the population covariance matrix $\bSigma$, the regularization matrix $\Qb$ and the mean vector $\bmu$ separately to assess the accuracy of Theorem~\ref{thm:main_theorem}, while keeping the values of the other quantities the same as those in Section~\ref{subsec:basic_simulation}. This provides comparison as we evaluate the effects of $\bSigma$, $\Qb$ and $\bmu$ separately. 
We investigate these changes for $(n,p)=(1500,750)$ and $(n,p)=(1500,2250)$ as in the basic settings. We still generate $\Rb\in\RR^{n\times p}$ for  1000 times and follow the procedure in Section~\ref{subsec:basic_simulation} to get $SR(\Qb)$ and $\hat{SR}(\Qb)$.

\medskip
\noindent\textbf{Results with different $\bSigma$'s.} Here we fix $\bmu=\bmu_0$ and $\Qb=q\cdot\Qb_0$. To validate our high dimensional statistical correction, we select two other covariance structures, denoted as $\bSigma_1$ and $\bSigma_2$. 
We construct $\bSigma_1$ as the diagonal matrix $\bSigma_1=\diag(\lambda_1,\dots,\lambda_p)$, where $\{\lambda_i\}_{i=1}^p$ is the sequence described in Section~\ref{subsec:basic_simulation}. This choice of $\bSigma_1$ represents the situation where we have no factor structure but well-conditioned spectrum.
We construct $\bSigma_2$ as $\bSigma_2 = \bSigma_0+\bxi_1\bxi_1^\top+\bxi_2\bxi_2^\top$, where $\bSigma_0$ is defined in Section~\ref{subsec:basic_simulation}, and $\bxi_1$ and $\bxi_2$ are two random generated vectors that are orthogonal to each other and also to the vector $\one$. The norm of $\bxi_1$ and $\bxi_2$ is chosen such that $\|\bxi_1\|_2^2=\|\bxi_2\|_2^2=p$. This choice of $\bSigma_2$ allows more factors besides the market factor, e.g. the Size and Value factors in Fama-French 3 factor model \citep{fama1993common}.
From a mathematical perspective, the assumptions on $\bSigma_1$ and $\bSigma_2$ cover a broader range of cases to verify our theory. The simulation procedure is similar to Section~\ref{subsec:basic_simulation}, except that we now replace the population covariance  $\bSigma_0$ with $\bSigma=\bSigma_1$ or $\bSigma=\bSigma_2$.
\begin{figure}[t]
    \centering
    \subfloat[$\bSigma=\bSigma_1,c=1/2$]{\includegraphics[width=0.24\textwidth]{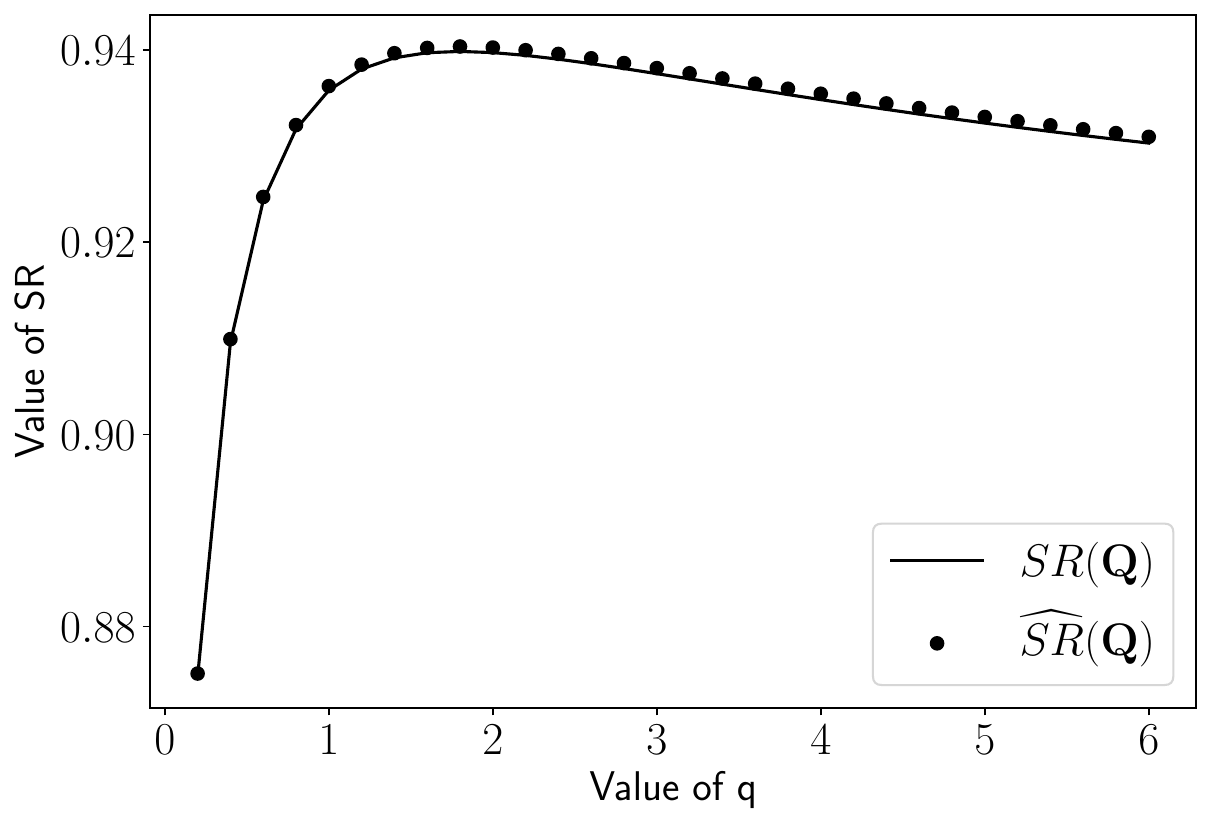}\label{fig_simu:sigma1n1500c<1}}
    \subfloat[$\bSigma=\bSigma_2,c=1/2$]{\includegraphics[width=0.24\textwidth]{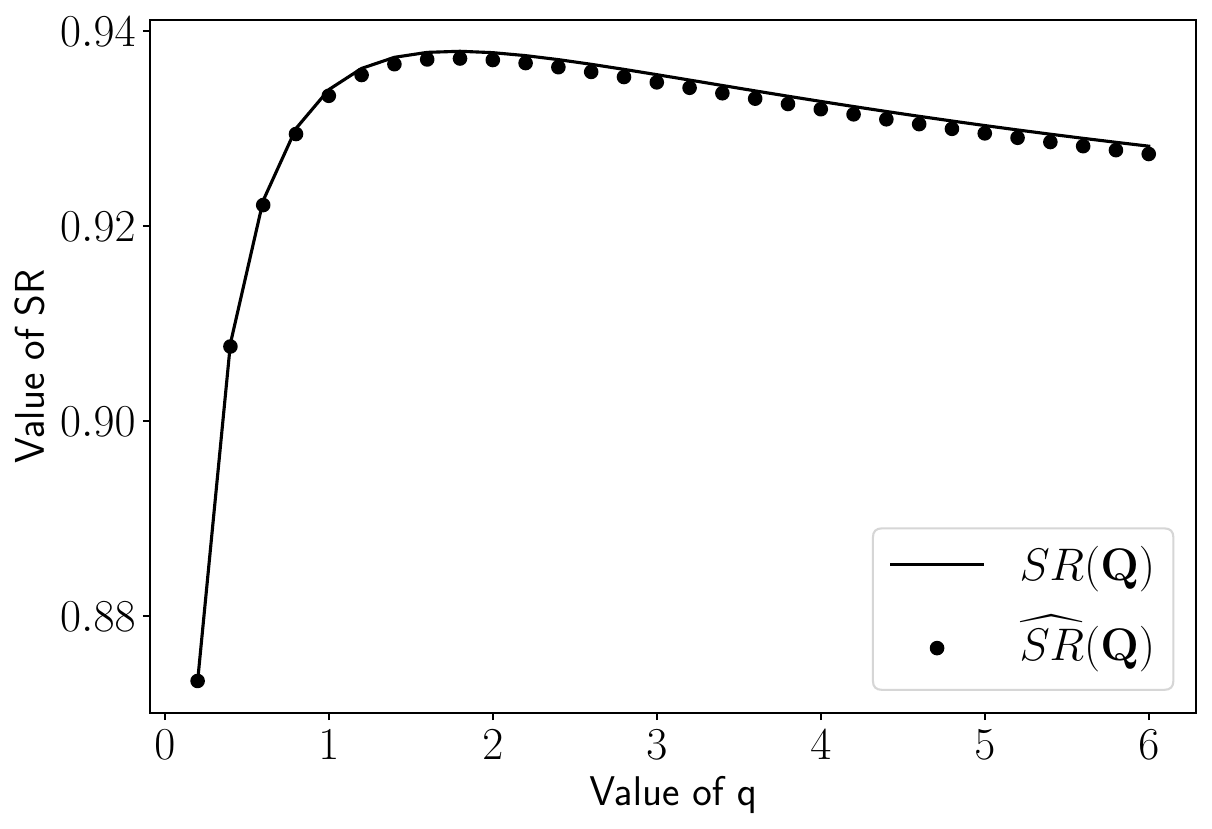}\label{fig_simu:sigma2n1500c<1}}
    \subfloat[$\bSigma=\bSigma_1,c=3/2$]{\includegraphics[width=0.24\textwidth]{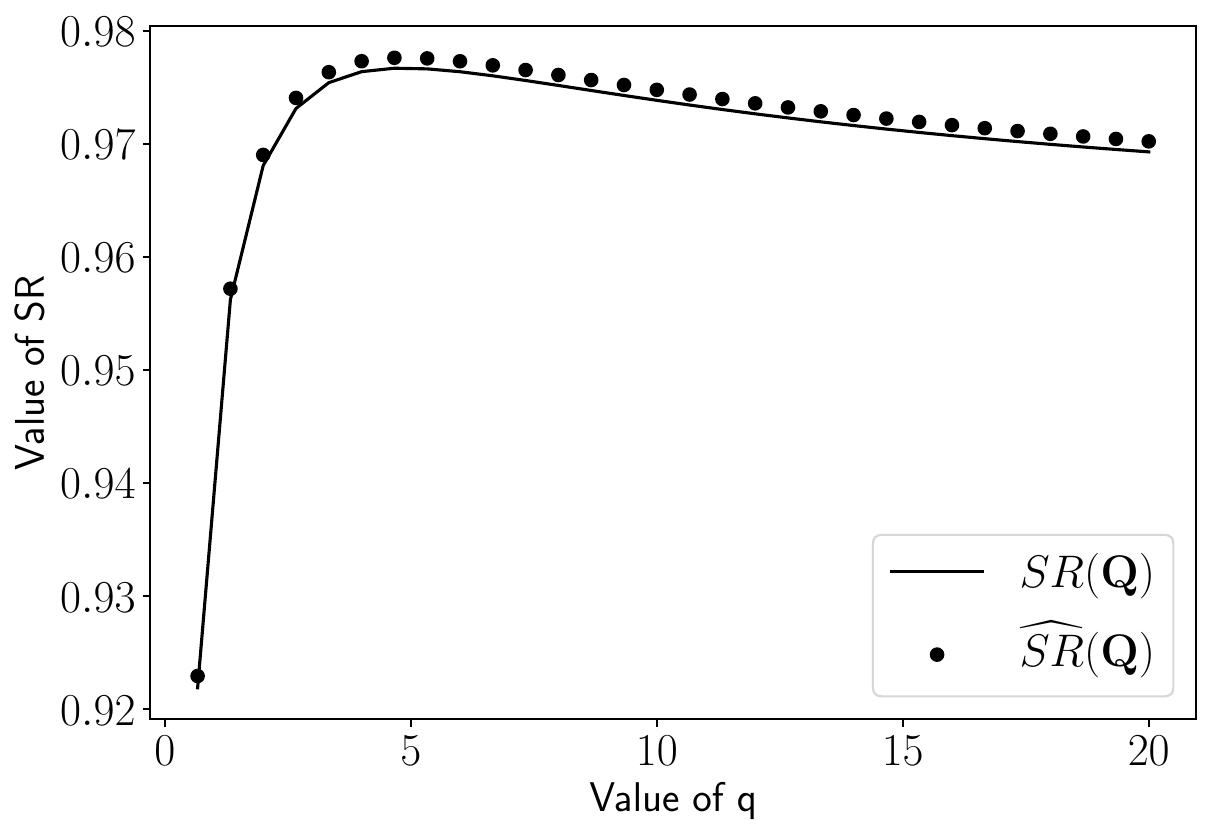}\label{fig_simu:sigma1n1500c>1}}
    \subfloat[$\bSigma=\bSigma_2,c=3/2$]{\includegraphics[width=0.24\textwidth]{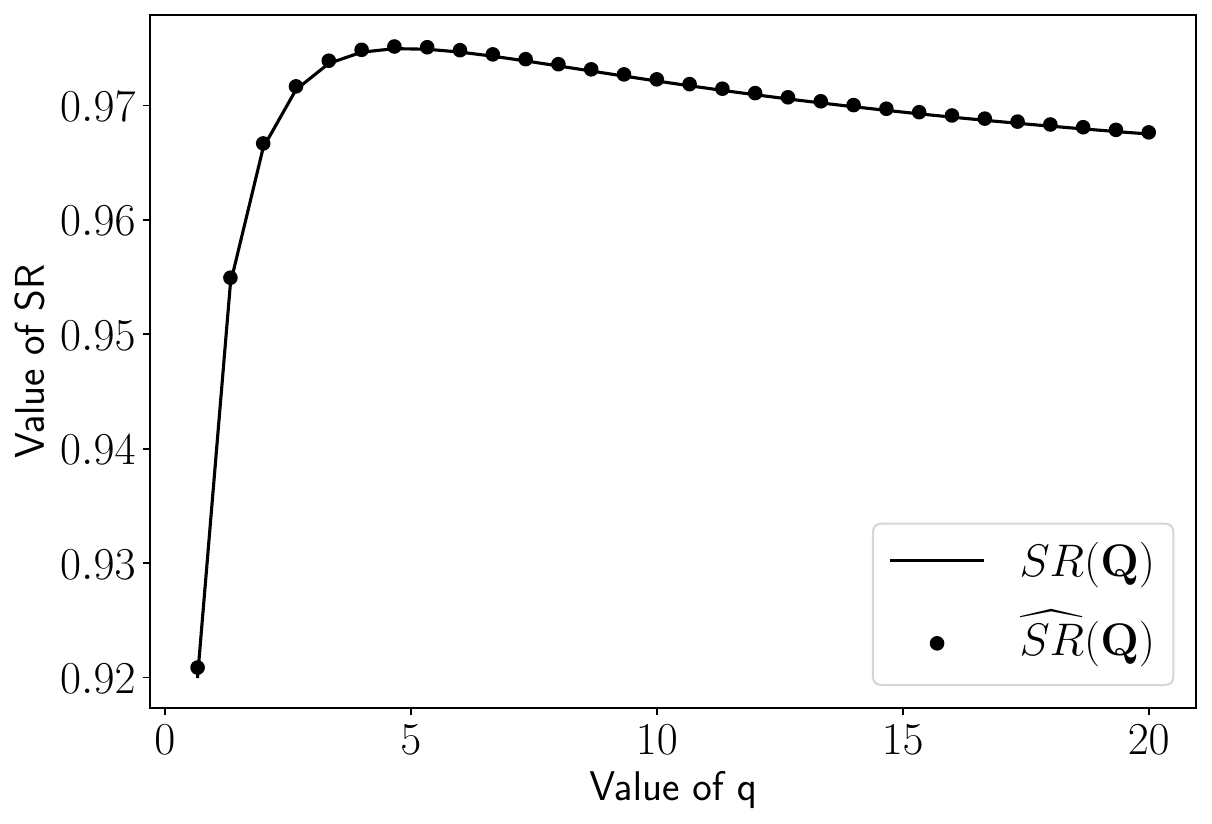}\label{fig_simu:sigma2n1500c>1}}
    \caption{Simulation results with different $\bSigma$'s. The x-axis and the y-axis follow exactly Figure~\ref{fig_simu:Basen1500}.
    }
    \label{fig_simu:bSigma}
\end{figure}

The results presented in Figure~\ref{fig_simu:bSigma} show that our estimated $\hat{SR}(\Qb)$, obtained without access to the population covariance matrix, again closely aligns with the true $SR(\Qb)$ for both $\bSigma_1$ and $\bSigma_2$ in both the cases of $c=1/2$ and $c=3/2$. This observation demonstrates that our high dimensional in-sample correction effectively captures the out-of-sample Sharpe ratio in all 3 cases in Assumption~\ref{ass:assump3}, no matter if we have bounded eigenvalues or a few diverging spikes. 

\medskip
\noindent\textbf{Results with different $\Qb$'s.}  Here we fix $\bmu=\bmu_0$ and $\bSigma=\bSigma_0$. 
In the basic settings, we define $\Qb=q\cdot\Qb_0=q\cdot \diag(3,...,3,1,...,1)$, where the number of $3$'s and $1$'s are both $p/2$. To provide a comparison, we introduce three additional matrices,  $\Qb_1=0.1\Qb_0+q\cdot \diag(\lambda_1,...,\lambda_p)$ where $\lambda_i$ is defined in Section~\ref{subsec:basic_simulation};  $\Qb_2=0.5\Ib_p+q\Qb_0$  and $\Qb_3=q\bSigma_0$.   The simulation process is similar to that in Section~\ref{subsec:basic_simulation}, except that we now replace the regularization term $q \cdot \Qb_0$ with $\Qb = \Qb_1$, $\Qb = \Qb_2$  or $\Qb=\Qb_3$. Additionally, we plot the line $SR_{\max} = \sqrt{\bmu_0^\top \bSigma_0^{-1} \bmu_0}$ to compare the values of $SR(\Qb)$ and $SR_{\max}$, which provides further validation of the findings in Appendix~\ref{sec:discussion_largest_SR} later and offers insight into how $\Qb$ can be shaped to approximate the highest achievable Sharpe ratio.

\begin{figure}[t]
    \centering
    \subfloat[$\Qb=\Qb_1,c=1/2$]{\includegraphics[width=0.32\textwidth]{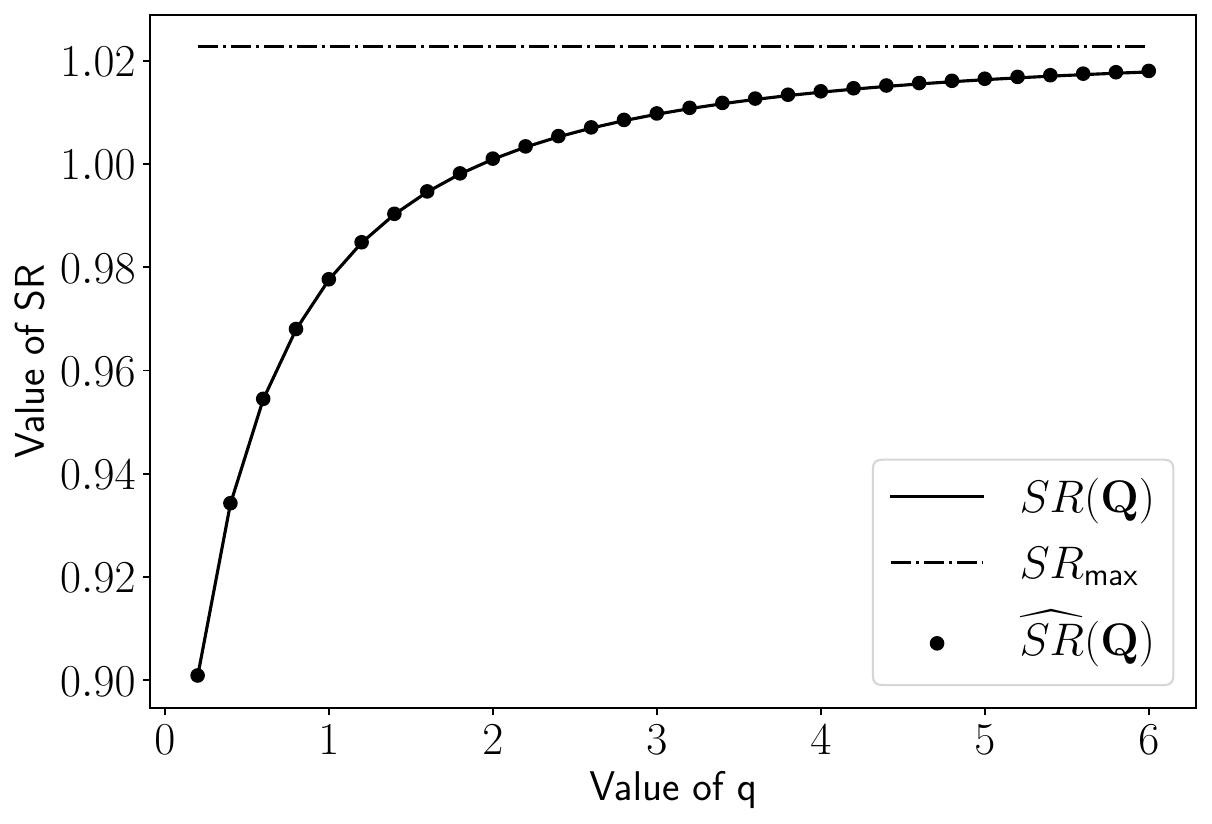}\label{fig_simu:Q1n1500c<1}}
    \subfloat[$\Qb=\Qb_2,c=1/2$]{\includegraphics[width=0.32\textwidth]{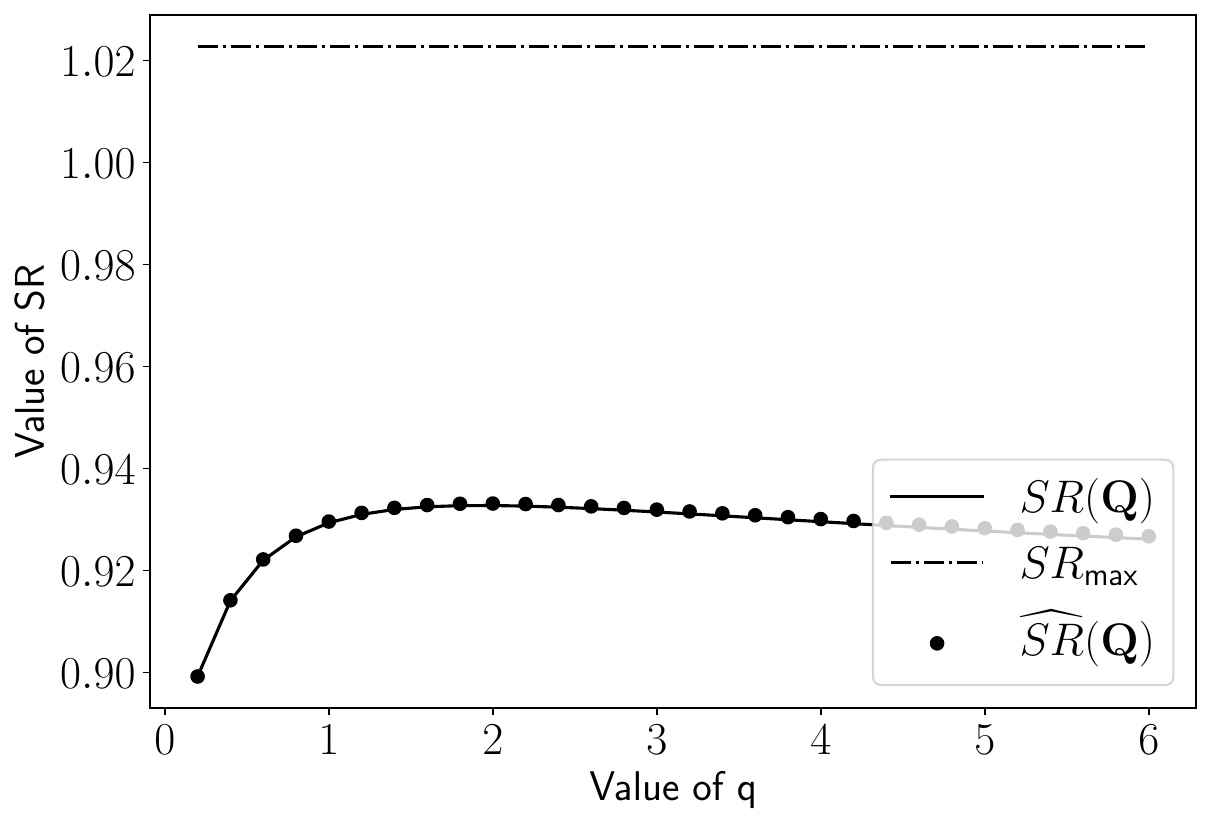}\label{fig_simu:Q2n1500c<1}}
    \subfloat[$\Qb=\Qb_3,c=1/2$]{\includegraphics[width=0.32\textwidth]{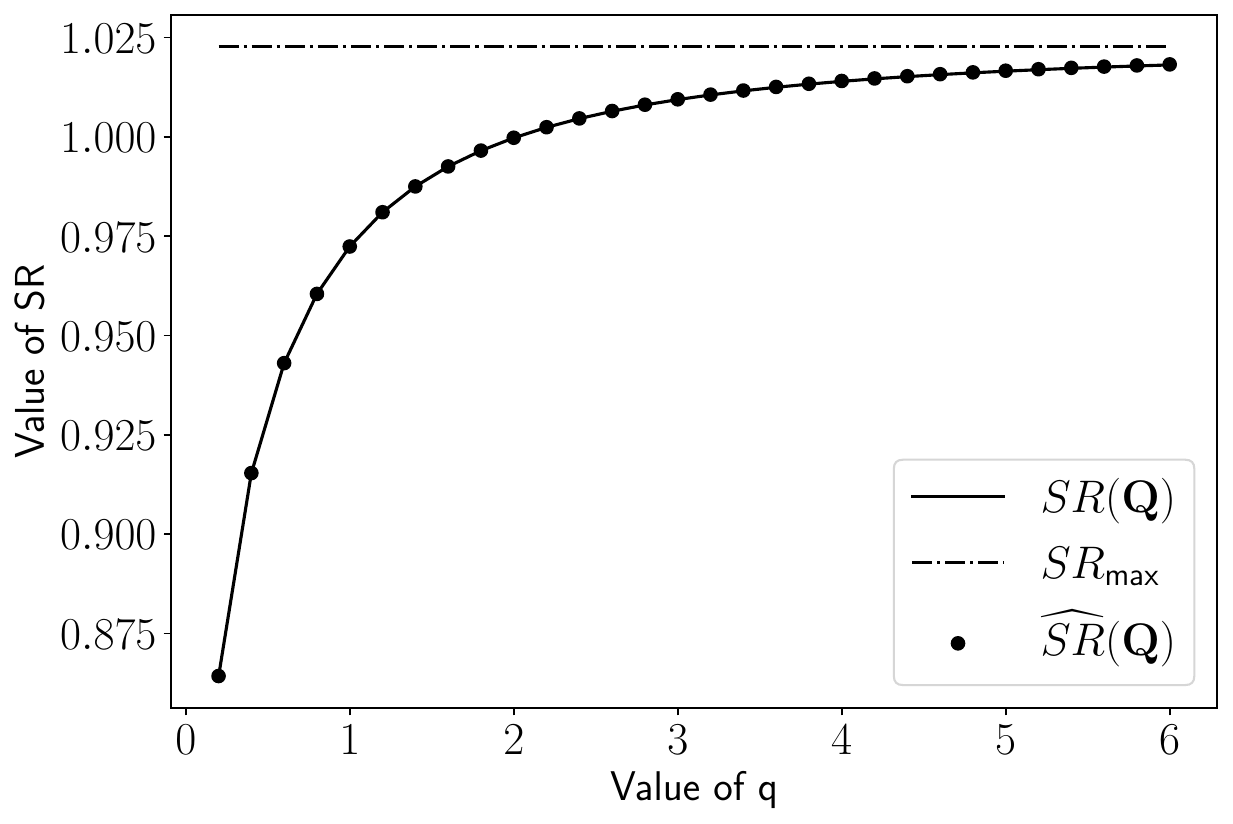}\label{fig_simu:Q3n1500c<1}}
    
    \subfloat[$\Qb=\Qb_1,c=3/2$]{\includegraphics[width=0.32\textwidth]{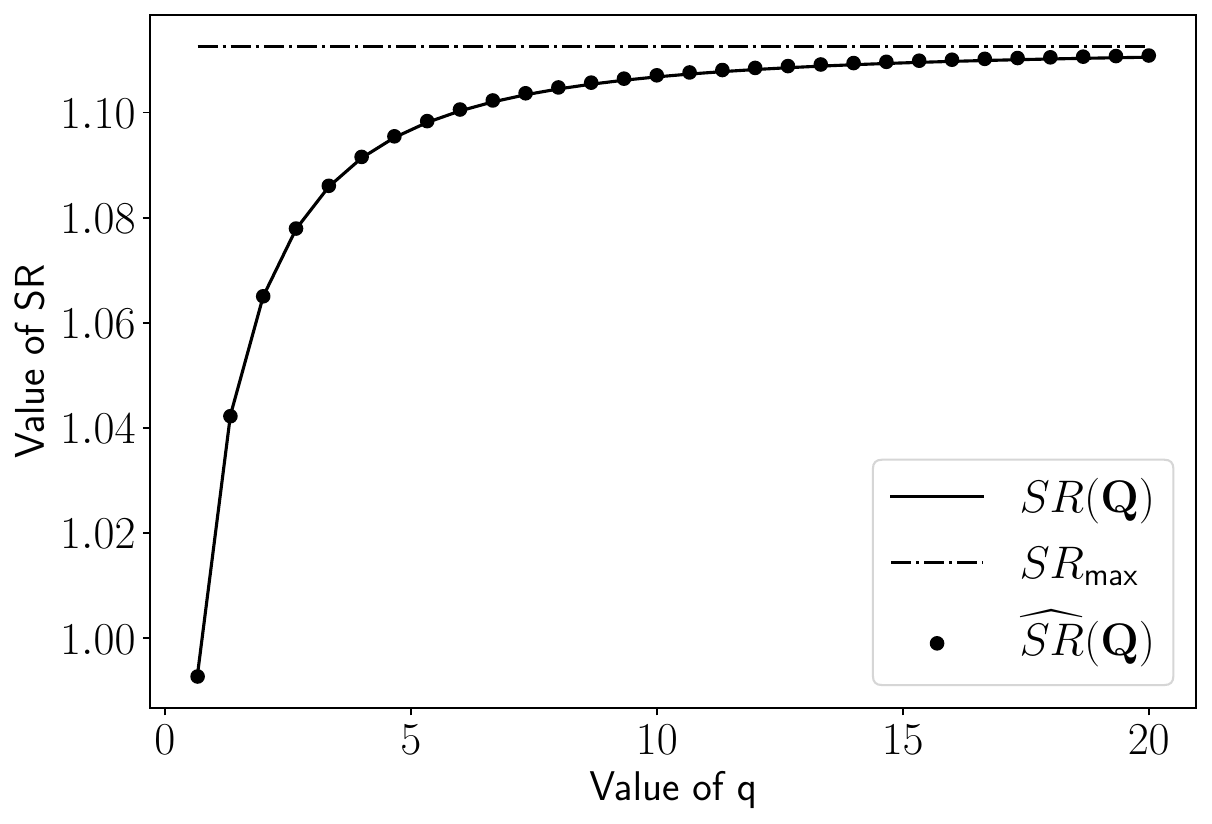}\label{fig_simu:Q1n1500c>1}}
    \subfloat[$\Qb=\Qb_2,c=3/2$]{\includegraphics[width=0.32\textwidth]{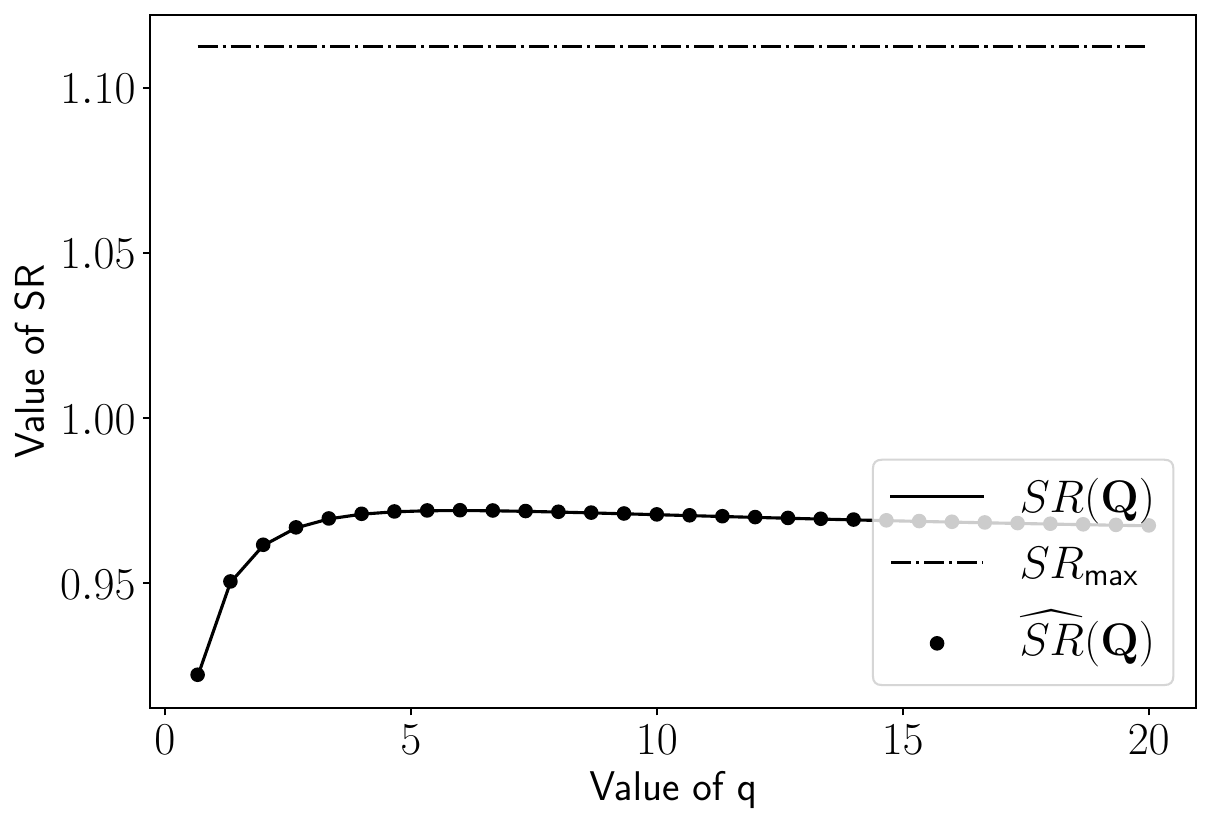}\label{fig_simu:Q2n1500c>1}}
    \subfloat[$\Qb=\Qb_3,c=3/2$]{\includegraphics[width=0.32\textwidth]{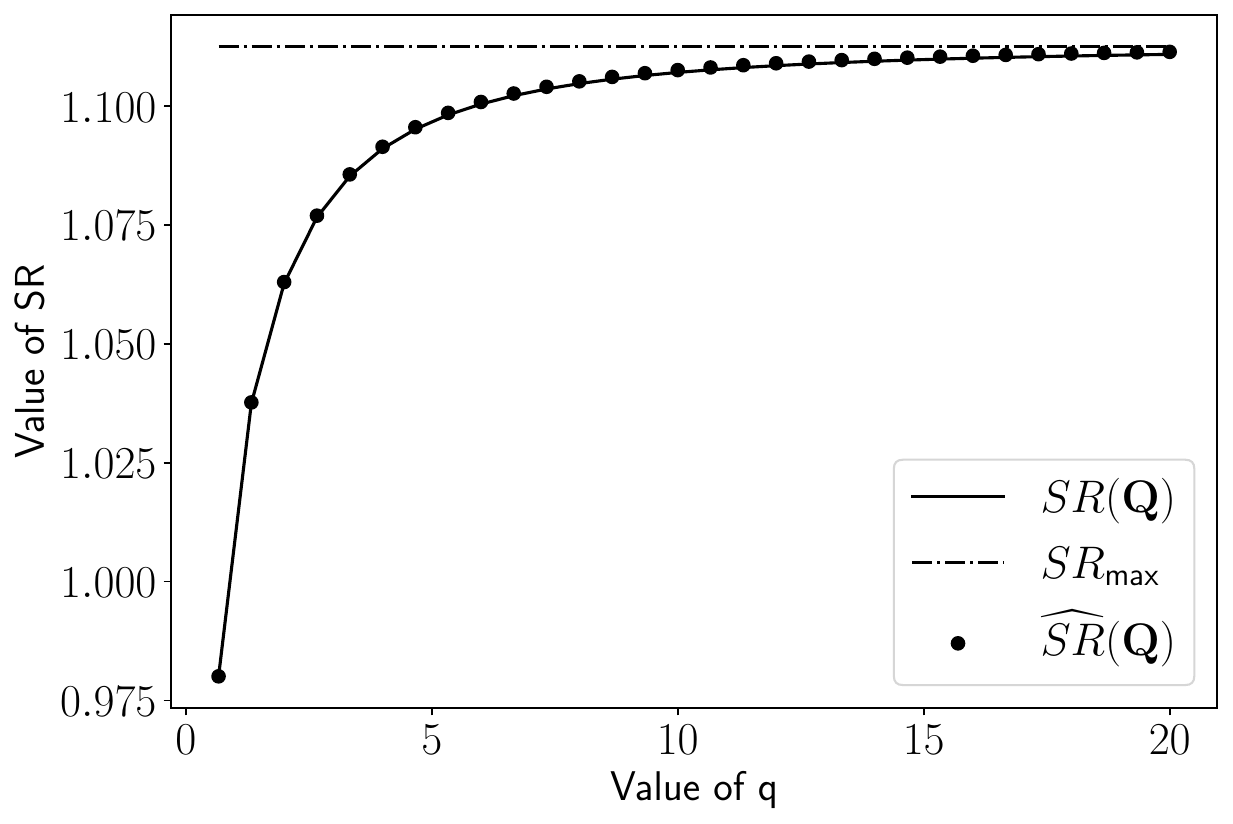}\label{fig_simu:Q3n1500c>1}}
    \caption{Simulation results with different $\Qb$'s. The x-axis and the y-axis follow exactly Figure~\ref{fig_simu:Basen1500}.
    }
    \label{fig_simu:Qb}
\end{figure}

 In Figure~\ref{fig_simu:Qb}, it is clear that $\hat{SR}(\Qb)$ aligns closely with $SR(\Qb)$ for both $\Qb = \Qb_1$ and $\Qb = \Qb_2$ when $c = 1/2$ and $c = 3/2$. However, the structure of $\Qb$ significantly affects the maximally achievable Sharpe ratio. For instance, in the settings of $\Qb = q\cdot \Qb_0$ or $\Qb = \Qb_2$, the gap between the highest achievable Sharpe ratio $SR(\Qb)$ and the theoretically maximal Sharpe ratio $SR_{\max}$ is clearly noticeable. 
Interestingly, in the case of $\Qb = \Qb_1$  or $\Qb=\Qb_3$, as we adjust $\Qb$ to approximate the structure of the population covariance matrix $\bSigma$ via its residual covariance $\diag(\lambda_1,...,\lambda_p)$,  or make $\Qb$ exactly the same as $\bSigma$, the maximal Sharpe ratio that we can achieve approaches $SR_{\max}$ as the tuning parameter $q$ increases. This observation not only supports our findings in Appendix~\ref{sec:discussion_largest_SR} but also suggests a practical approach for designing $\Qb$. Specifically, by constructing $\Qb$ to be structurally similar or proportionally aligned with $\bSigma$, we can enhance the Sharpe ratio and achieve values closer to $SR_{\max}$. Therefore, for optimal performance, it is advisable to design $\Qb$ within the candidate set $\cQ$ to reflect the true structure of $\bSigma$ or its corresponding residual covariance. If we do not have prior information about the covariance structure, we may either use some independent data to estimate the residual covariance or find good proxies such as the trading volume to capture the rough level of the residual variances.

\medskip
\noindent\textbf{Results with different $\boldsymbol{\mu}$'s.} Here we fix $\bSigma=\bSigma_0$ and $\Qb=q\cdot\Qb_0$ and introduce different settings for the mean vectors, $\bmu_1$ and $\bmu_2$. 
For $\bmu_1$, we assume that each element follows an independent uniform distribution, $\Unif(-\sqrt{2/p},\sqrt{2/p})$. This assumption ensures that each individual asset has a nonzero expected return. 
For $\bmu_2$, we add the vector $2\cdot\1_p$ to $\bmu_1$, that is $\bmu_2=\bmu_1+2\cdot\1_p$. Recall that the vector $\one$ represents the market factor, so $\bmu_2$ induces a portfolio that not only bets on ``alpha'' signals but also holds a long position in the market factor. 
The simulation procedure is similar to Section~\ref{subsec:basic_simulation}, except that we now replace the mean vector  $\bmu_0$ with $\bmu=\bmu_1$ or $\bmu=\bmu_2$.
\begin{figure}[t]
    \centering
    \subfloat[$\bmu=\bmu_1,c=1/2$]{\includegraphics[width=0.24\textwidth]{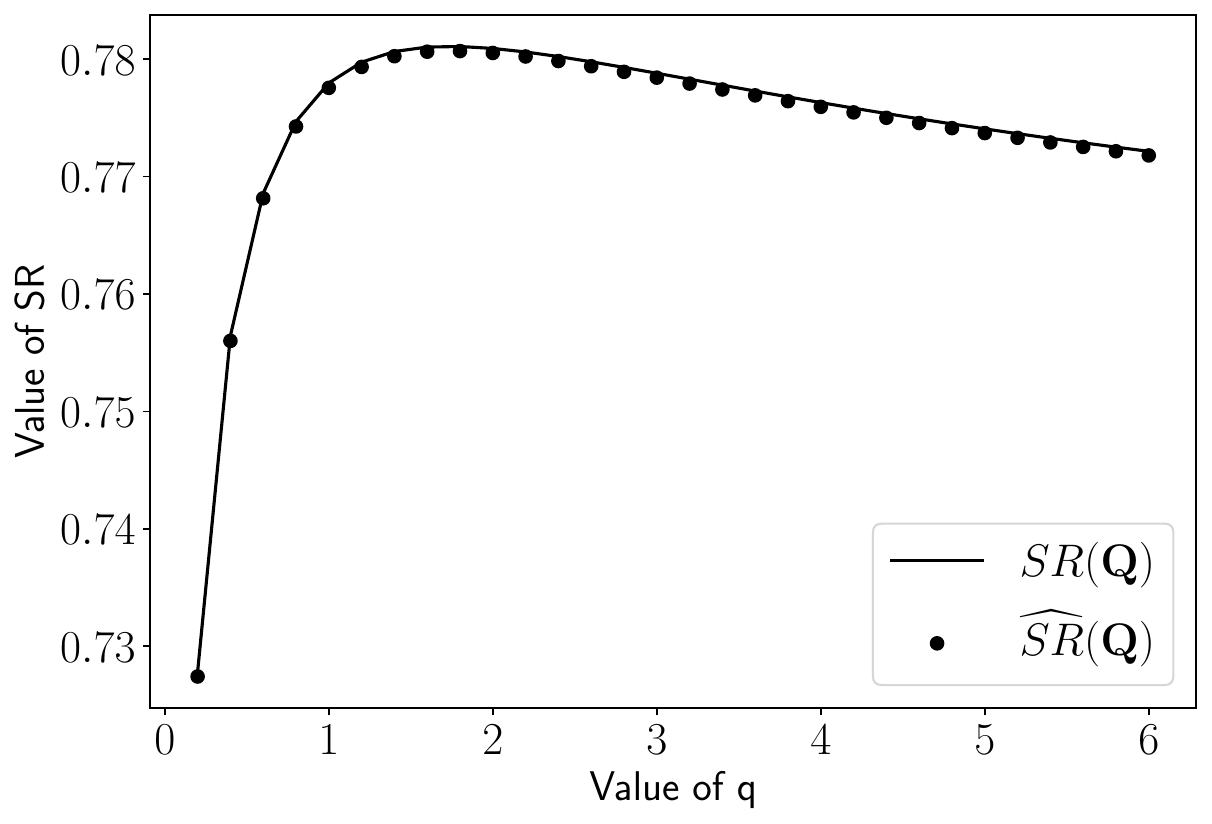}\label{fig_simu:mu1n1500c<1}}
    \subfloat[$\bmu=\bmu_2,c=1/2$]{\includegraphics[width=0.24\textwidth]{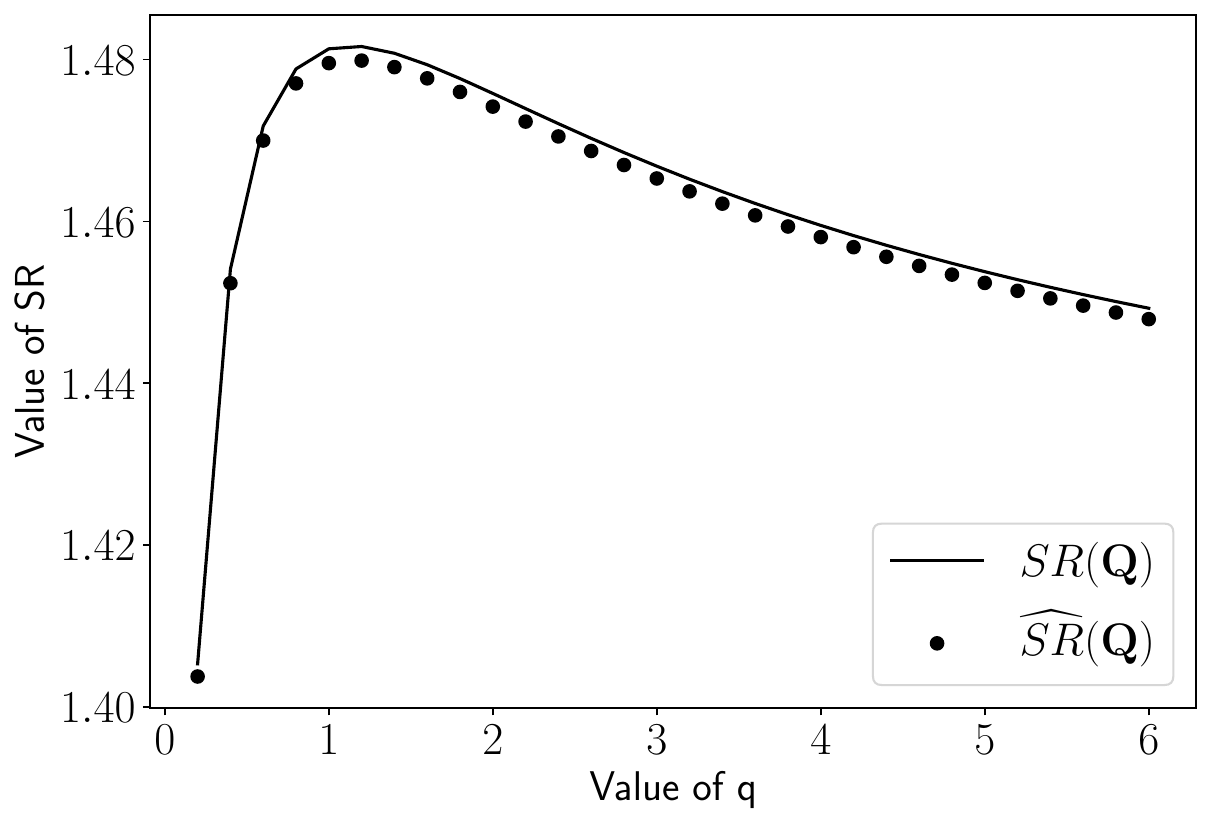}\label{fig_simu:mu2n1500c<1}}
    \subfloat[$\bmu=\bmu_1,c=3/2$]{\includegraphics[width=0.24\textwidth]{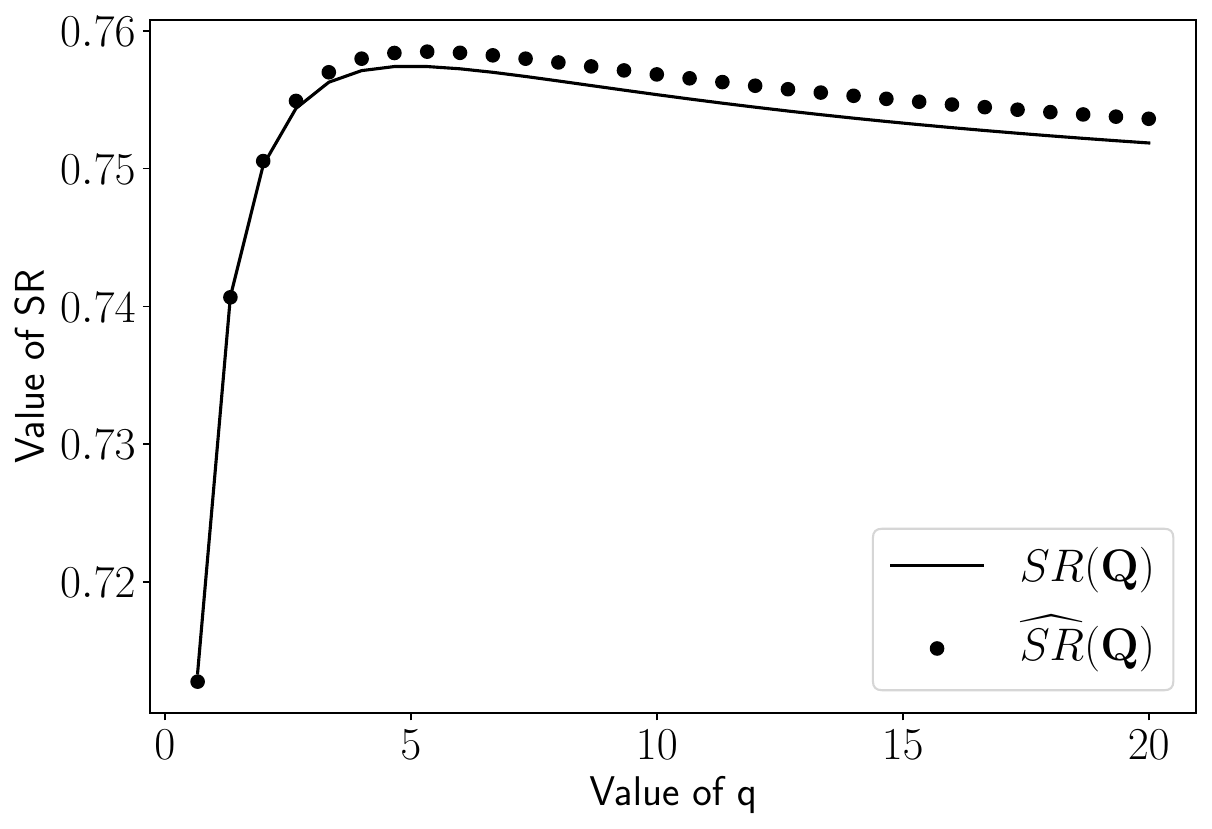}\label{fig_simu:mu1n1500c>1}}
    \subfloat[$\bmu=\bmu_2,c=3/2$]{\includegraphics[width=0.24\textwidth]{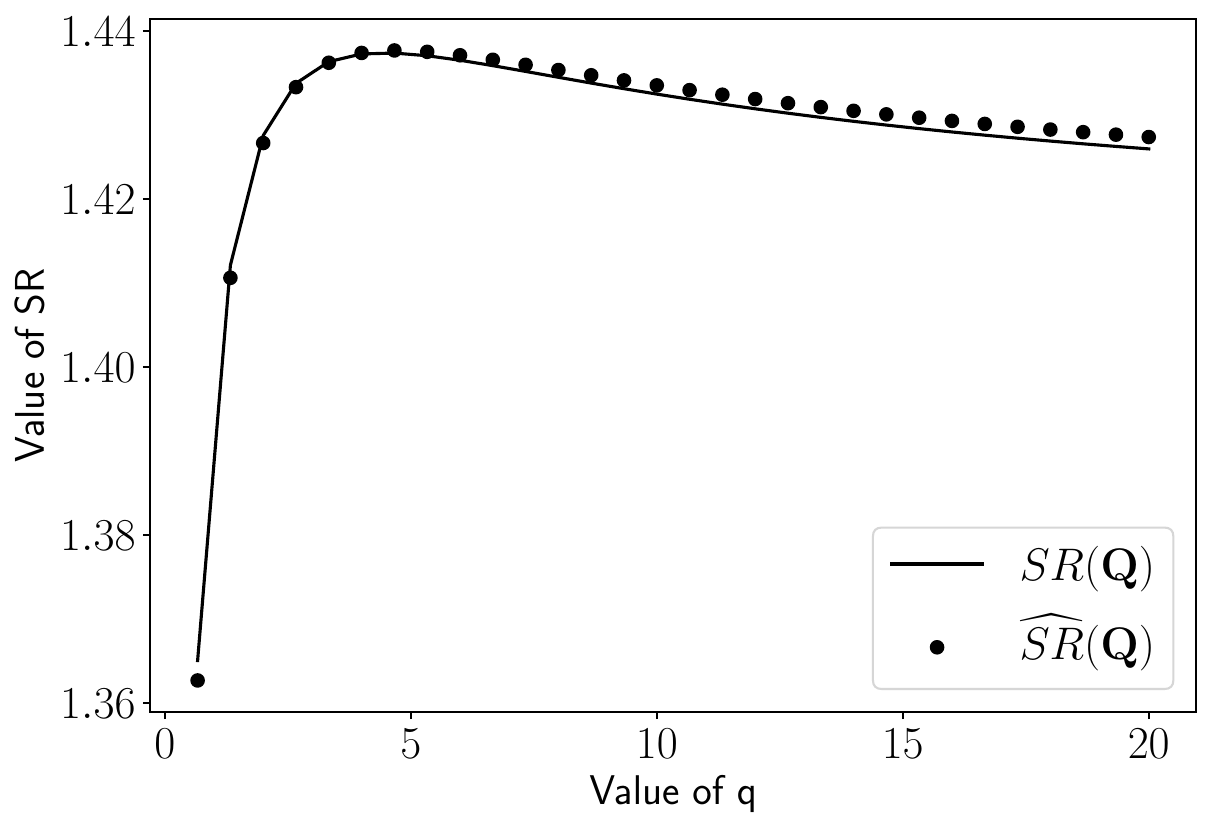}\label{fig_simu:mu2n1500c>1}}
    \caption{Simulation results with different $\bmu$'s. The x-axis and the y-axis follow exactly Figure~\ref{fig_simu:Basen1500}.
    }
    \label{fig_simu:bmu}
\end{figure}

The results presented in Figure~\ref{fig_simu:bmu} demonstrate that our estimated $\hat{SR}(\Qb)$ closely aligns with the true $SR(\Qb)$ for both $\bmu_1$ and $\bmu_2$. This comparison across different values of $\bmu$ further highlights the robustness and accuracy of our approach in estimating $SR(\Qb)$.
Furthermore, when comparing the cases of $\bmu=\bmu_1$ and $\bmu=\bmu_2$, we see a significant difference in $SR(\Qb)$. This observation indicates that for different mean returns, Sharpe ratio can vary dramatically, but our method can accurately evaluate the out-of-sample Sharpe ratio regardless of the structure of the mean returns.


\subsection{Experiments with unknown mean vector}
\label{subsec:simu_unknown}
In this section, we present additional experiments on the situation where the mean vector $\bmu$ is unknown. Recall $\hbSigma$ and $\hbmu$ are sample covariance and sample mean, respectively. The out-of-sample Sharpe ratio is given by $SR(\Qb)=\frac{\hbmu^\top(\hbSigma+\Qb)^{-1}\bmu}{\sqrt{\hbmu^\top(\hbSigma+\Qb)^{-1}\bSigma(\hbSigma+\Qb)^{-1}\hbmu}}$, and the estimation $\hat{SR}(\Qb)$ is given by
\begin{align*}
    \hat{SR}(\Qb)=\frac{\hbmu^\top(\hbSigma+\Qb)^{-1}\hbmu-\frac{\tr(\hbSigma+\Qb)^{-1}\hbSigma}{n-\tr(\hbSigma+\Qb)^{-1}\hbSigma}}{\sqrt{\hbmu^\top(\hat{\bSigma}+\Qb)^{-1}\hbSigma(\hat{\bSigma}+\Qb)^{-1}\hbmu}}\cdot\big(1-\frac{c}{p} \tr\hbSigma(\hbSigma+\Qb)^{-1}\big).
\end{align*}
See Appendix \ref{sec:unknown_mu} for more details on the theoretical justification.
We continue to fix the population matrix at $ \bSigma = \bSigma_0 $ and conduct experiments across three scenarios with mean vectors $ \bmu = \bmu_0 $, $ \bmu_1 $, and $ \bmu_2 $. For the regularization matrix $ \Qb $, we define it as $ \Qb_1 = q \cdot \diag(\lambda_1, \ldots, \lambda_p) $ and $\Qb_2=q\cdot\bSigma_0$, where $\bSigma_0$ and each $ \lambda_i $ are specified in Section~\ref{subsec:basic_simulation}. Additionally, we plot the line corresponding to the maximal Sharpe ratio, given by $ SR_{\max} = \sqrt{\bmu^\top \bSigma^{-1} \bmu} $, to assess whether designing $ \Qb $ as approximately proportional to $ \bSigma $ can still lead to an approximation of $ SR_{\max} $ when the mean vector $ \bmu $ is unknown. Based on the discussion in Section~\ref{sec:unknown_SRmax}, we also plot the value $SR_{\text{L}}=\frac{SR_{\max}^2}{\sqrt{SR_{\max}^2+c}}$ where $SR(q\cdot\bSigma_0)$ approximate with large $q$. The constant $c$ here also reflects the high-dimensional nature of the problem: if $p$ is much smaller than the sample size $n$, it may be possible to approximate the maximal Sharpe ratio even when $\bmu$ is unknown. 
\begin{figure}[t]
    \centering
    \subfloat[$\bmu=\bmu_0,c=1/2$]{\includegraphics[width=0.32\textwidth]{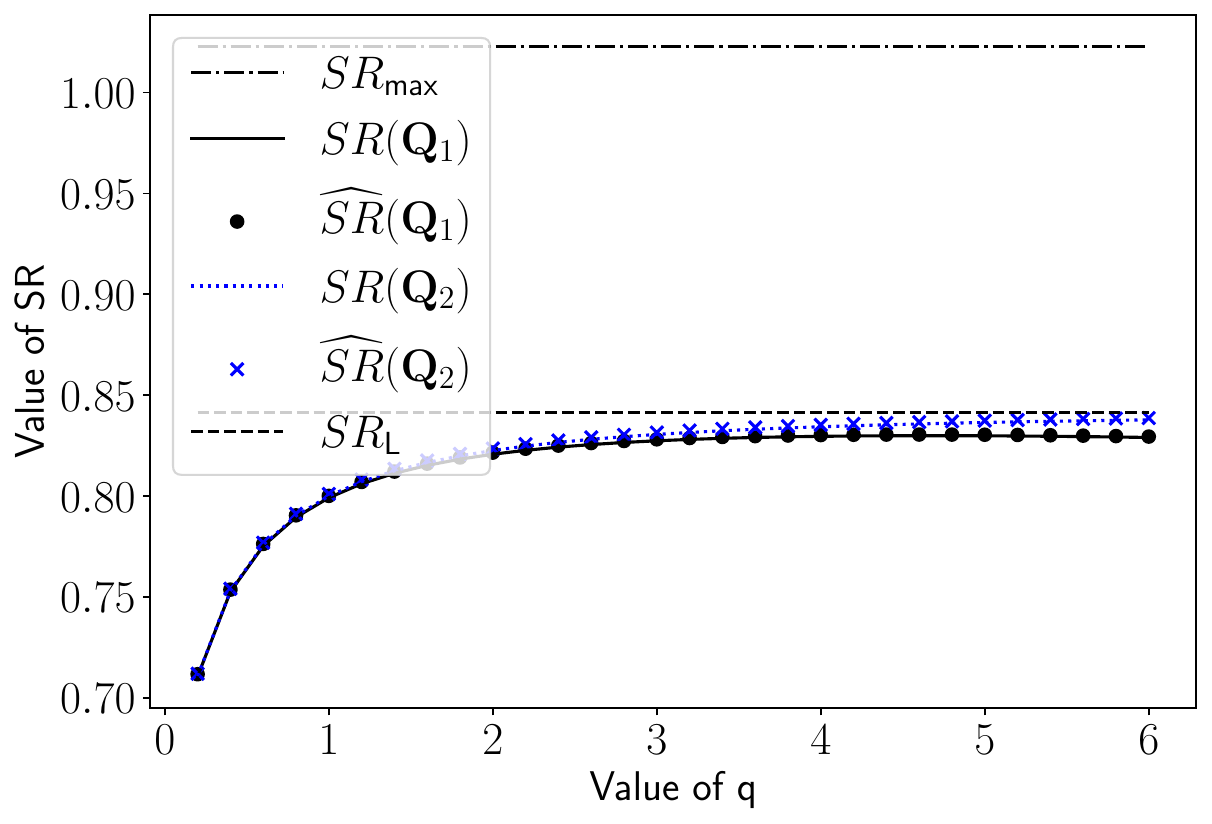}\label{fig_simu:unknown_mubase_c<1}}
     \subfloat[$\bmu=\bmu_1,c=1/2$]{\includegraphics[width=0.32\textwidth]{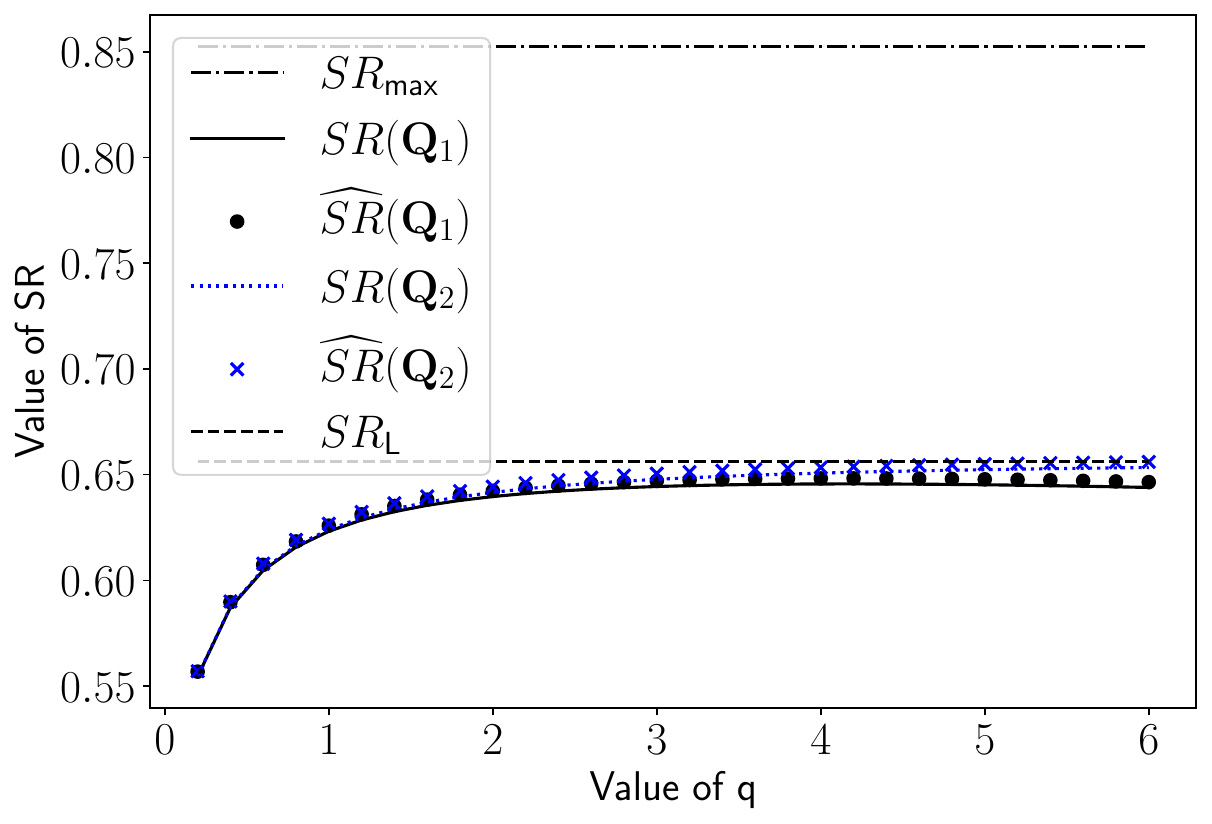}\label{fig_simu:unknown_mu1_c<1}}
      \subfloat[$\bmu=\bmu_2,c=1/2$]{\includegraphics[width=0.32\textwidth]{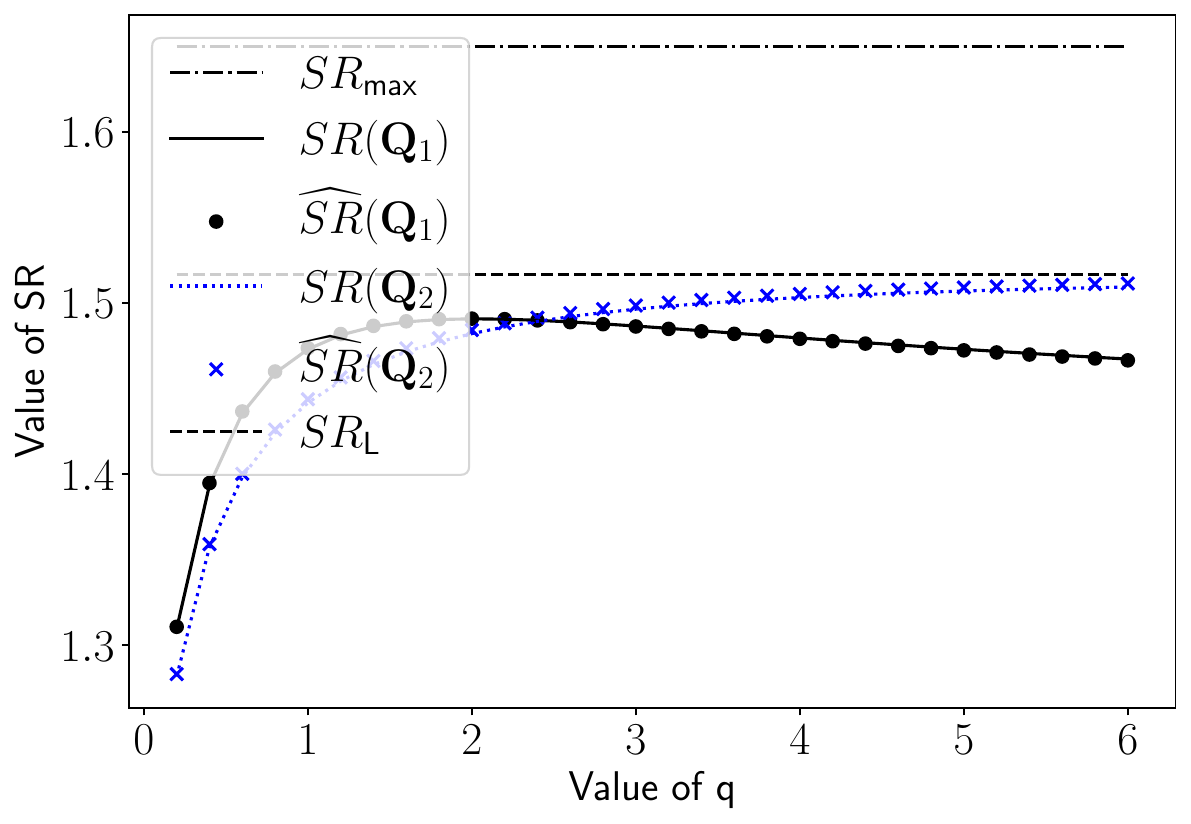}\label{fig_simu:unknown_mu2_c<1}}

       \subfloat[$\bmu=\bmu_0,c=3/2$]{\includegraphics[width=0.32\textwidth]{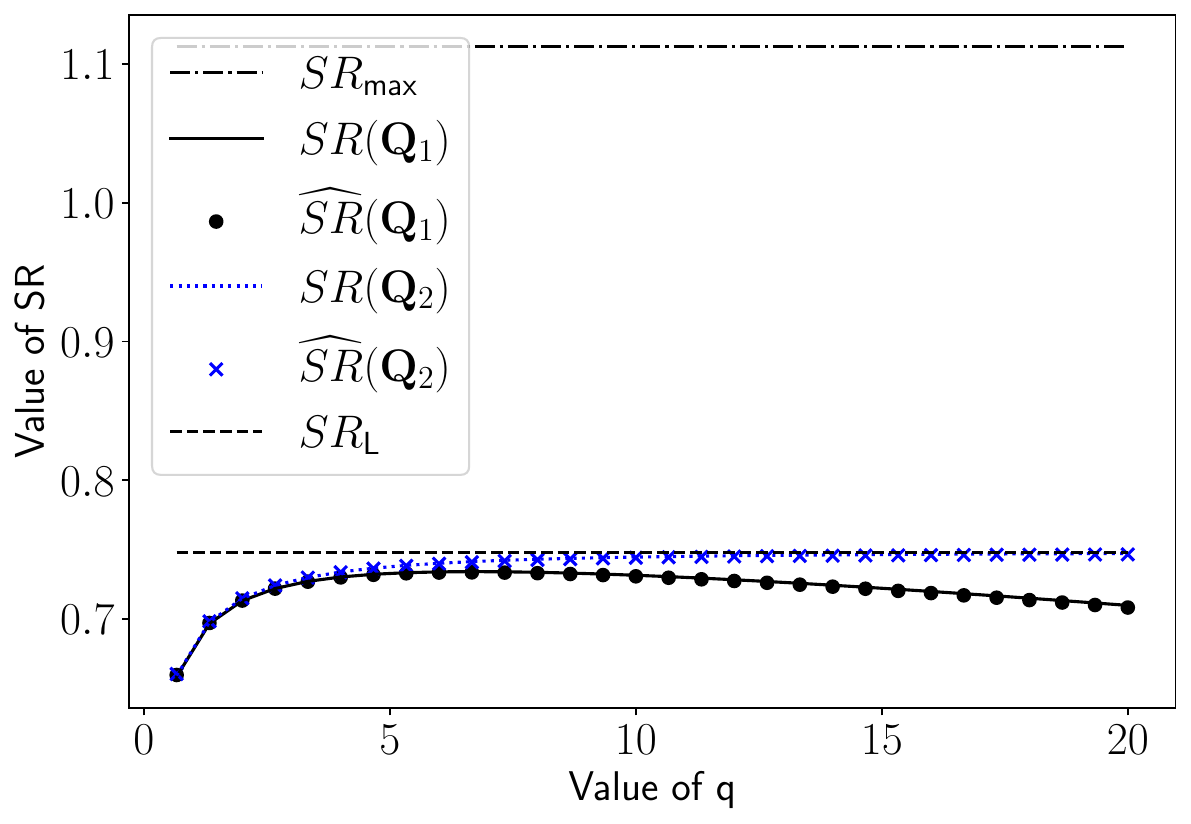}\label{fig_simu:unknown_mubase_c>1}}
        \subfloat[$\bmu=\bmu_1,c=3/2$]{\includegraphics[width=0.32\textwidth]{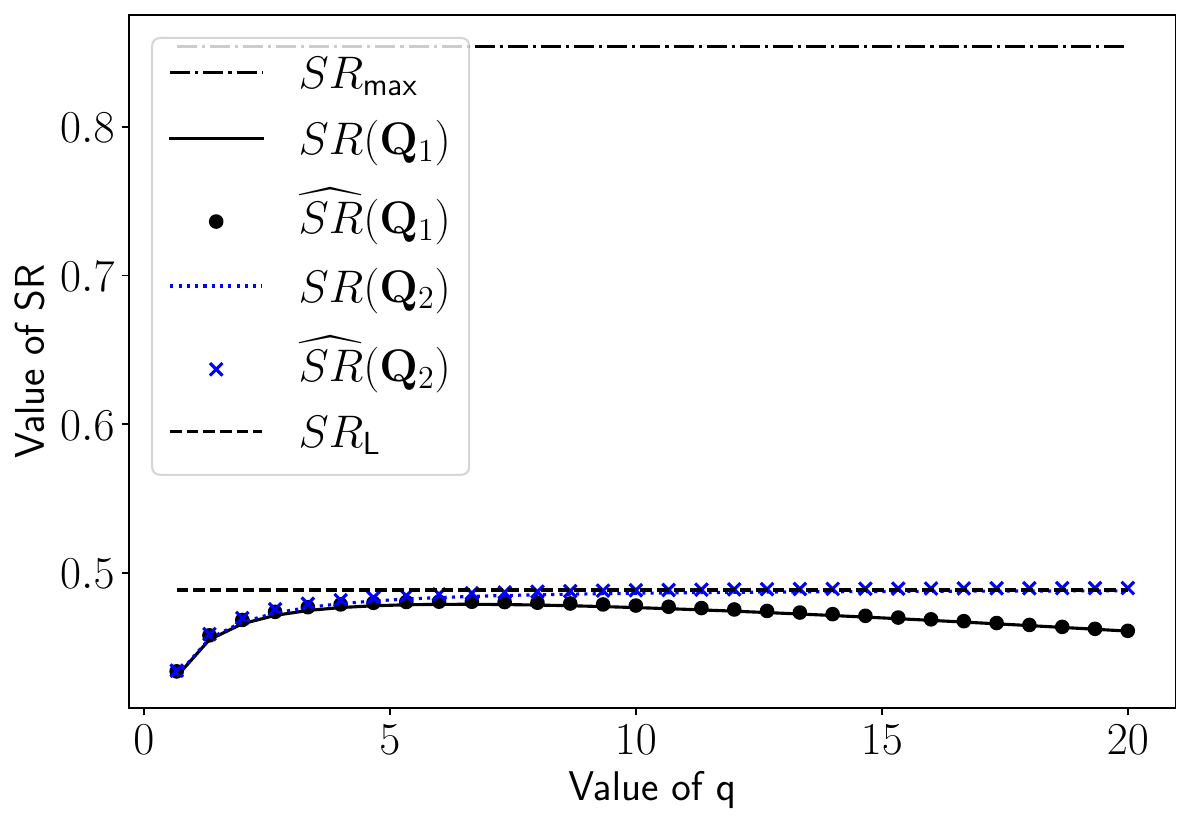}\label{fig_simu:unknown_mu1_c>1}}
         \subfloat[$\bmu=\bmu_2,c=3/2$]{\includegraphics[width=0.32\textwidth]{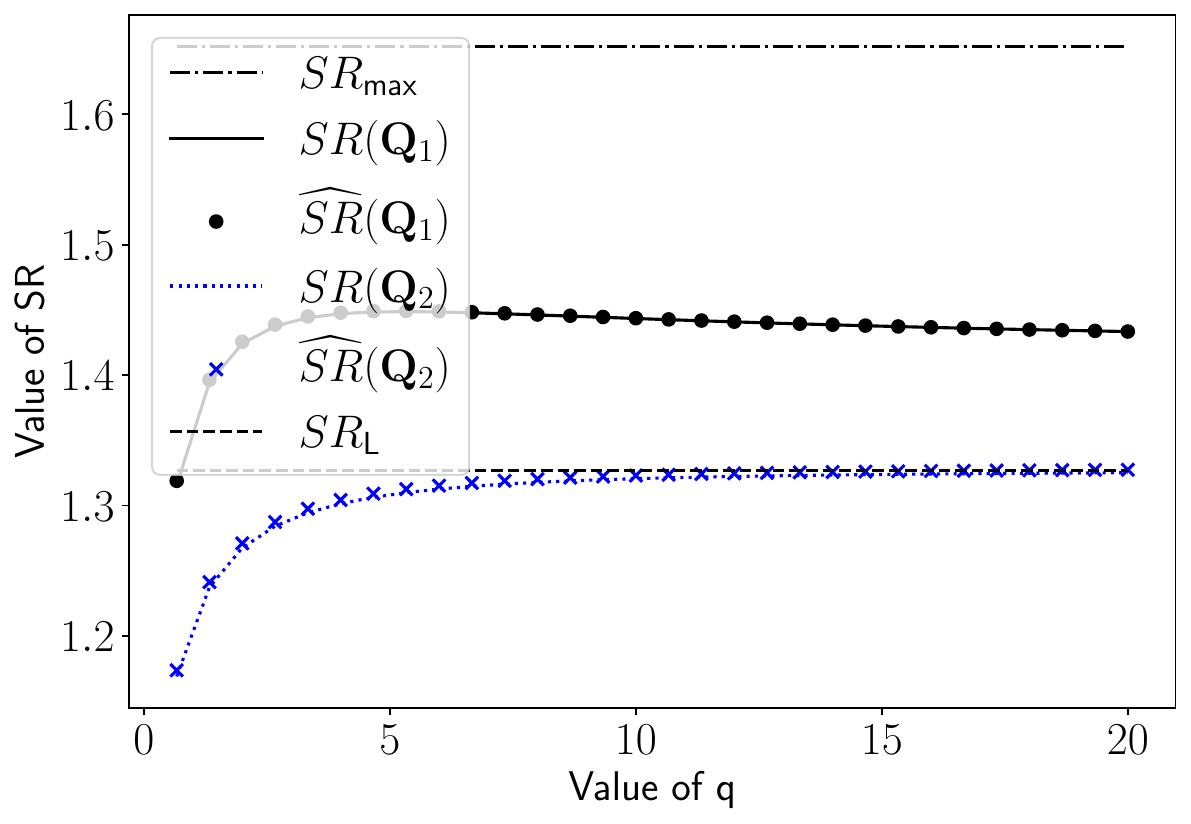}\label{fig_simu:unknown_mu2_c>1}}
    \caption{ Simulation results with $\bmu$ unknown. The x-axis and the y-axis follow exactly Figure~\ref{fig_simu:Basen1500}.}
    \label{fig_simu:unknownmu}
\end{figure}

The results presented in Figure~\ref{fig_simu:unknownmu} illustrate that when the mean vector $ \bmu $ is unknown, our proposed estimator $ \hat{SR}(\Qb) $  closely aligns with the true Sharpe ratio $ SR(\Qb) $. However, a particularly intriguing observation arises: in scenarios where $ \bmu $ is known, a larger $ q $ aligned with $ \bSigma $ seems beneficial for optimal portfolio construction; unlike the scenario where $ \bmu $ is known, when we set $ \Qb = q \bSigma$ in the case unknown $\bmu$, an increase in the tuning parameter $ q $ maintains a gap between $ SR_{\max} $ and the optimally achievable $ SR(\Qb) $.  Notably, the optimal value of $ q $ may even appear in the intermediate range when $\Qb$ has a slightly different form (change $\Qb_2$ to $\Qb_1$), suggesting that simply increasing $ q $ does not always yield better results.
This gap also indicates a challenge in approximating $ SR_{\max} $ when using both sample mean and sample covariance together. In the context of an unknown $ \bmu $, it appears that the best portfolio construction strategy may no longer involve using the true $ \bSigma$. As direct evidence, Figure~\ref{fig_simu:unknown_mu2_c>1} shows that with $\Qb = q \Qb_1$, where $\Qb_1$ slightly differs from the true covariance matrix $\bSigma$, the performance of $SR(\Qb_1)$ can even outperform the performance when $\Qb=q\cdot\bSigma$. This raises fascinating questions about the effectiveness of traditional approaches and suggests that when faced with uncertainty regarding $\bmu$, relying on the true covariance matrix might not lead to the best outcomes. Instead, an alternative approach may be necessary to navigate this uncertainty effectively. We leave this problem for future investigation.

\subsection{Additional experiments on regularization matrix}
\label{sec:comparsion_optimal}
In this section, we conduct additional experiments by optimizing the entire matrix $\hat{\Qb}=\arg\max_{\Qb}\hat{SR}(\Qb)$ (i.e. without restricting \(\Qb\in\cQ\)). The experiment is designed to support our claim below Theorem~\ref{thm:main_theorem} that optimizing \(\hat{SR}(\Qb)\) over the entire space can lead to overfitting to the in-sample data, which will potentially inflate \(\hat{SR}(\hat{\Qb})\) and create a discrepancy between the estimated \(\hat{SR}(\hat{\Qb})\) and the true \(SR(\hat{\Qb})\).

Recall that
\[
\hat{SR}(\Qb)=\left(1-\frac{c}{p}\operatorname{tr}\left[\hbSigma(\hbSigma+\Qb)^{-1}\right]\right)
\cdot \frac{\bmu^\top(\hbSigma+\Qb)^{-1}\bmu}{\sqrt{\bmu^\top(\hbSigma+\Qb)^{-1}\hbSigma(\hbSigma+\Qb)^{-1}\bmu}}.
\]
To provide an example of overfitting where optimizing \(\hat{SR}(\Qb)\) without constraints, we specifically pick up two  pairs of $(n,p)$: $(500,250) $ and $ (1000,500)$ (other pairs are omitted due to heavy computational burden). We adopt the same settings as in Section~\ref{subsec:basic_simulation} (\(\bmu=\bmu_0\) and \(\bSigma=\bSigma_0\)),  generate $\Rb\in\RR^{n\times p}$ for  20 times and optimize $\Qb$ in  $\hat{SR}(\Qb)$ to obtain $\hat{\Qb}$ each time. To ensure that \(\Qb\) remains positive definite during optimization, we decompose $\Qb$ as \(\Qb=\Lb\Lb^\top\) and then optimize $\Lb$, where \(\Lb\) is a lower-triangular matrix. Additionally, we consider an alternative case by constraining \(\Qb\) to be diagonal.

\begin{table}[t]
\caption{ Comparison for $SR(\hat{\Qb})$ and  $\hat{SR}(\hat{\Qb})$. The mean gives the average value and the range gives the minimum and maximum values over the 20 independent trials.}
\label{table:comparison}
\scalebox{0.9}{ 
\begin{tabular}{cccccc}
\hline
\multirow{2}{*}{$(n,p)$} & \multirow{2}{*}{$SR_{\max}$} & \multicolumn{4}{c}{Optimization over full $\Qb$}\\
        &             & mean of $SR(\hat{\Qb})$ & Range of $SR(\hat{\Qb})$ & mean of $\hat{SR}(\hat{\Qb})$ & Range of $\hat{SR}(\hat{\Qb})$ \\ \cline{2-6} 
$(500,250)$              & 0.923       & 0.643  & [0.604,0.694]        & 1.299        & [1.178,1.454]              \\
$(1000,500)$             & 1.123       & 0.791  & [0.738,0.824]        & 1.513        & [1.406,1.578]              \\ \hline
\multirow{2}{*}{$(n,p)$} & \multirow{2}{*}{$SR_{\max}$} & \multicolumn{4}{c}{Optimization over diagonal $\Qb$}             \\
        &             & mean of $SR(\hat{\Qb})$ & Range of $SR(\hat{\Qb})$ & mean of $\hat{SR}(\hat{\Qb})$ & Range of $\hat{SR}(\hat{\Qb})$ \\ \cline{2-6} 
$(500,250)$              & 0.923       & 0.770  & [0.715,0.818]        & 0.967        & [0.909,1.056]              \\
$(1000,500)$             & 1.123       & 0.944  & [0.912,0.979]        & 1.146        & [1.081,1.224]              \\ \hline
\end{tabular}}
\end{table}
As shown in Table~\ref{table:comparison}, when we optimize the full matrix \( \Qb \), the value of $\hat{SR}(\hat{\Qb})$ is even much higher than $SR_{\max}$. Despite the optimal matrix \( \hat{\Qb} \) having a high value in \( \hat{SR}(\hat{\Qb}) \), the performance of the true objective \( SR(\hat{\Qb}) \) is poor,  remaining significantly lower than $SR_{\max}$. 
This discrepancy indicates that the matrix \( \hat{\Qb} \), which is optimized based on \( \hat{SR}(\Qb) \), is not suitable to obtain a high  out-of-sample Sharpe ratio.  As discussed below Theorem~\ref{thm:main_theorem}, the primary reason for discrepancy is that $\hat{\Qb}$ is overfitted to  in-sample data. We further conduct experiments by restricting \(\Qb\) to be diagonal. Although this restriction reduces the discrepancy relative to the full matrix optimization, it does not eliminate it completely. With less flexibility in the optimization of \(\Qb\), the value of \( \hat{SR}({\Qb}) \) decreases, but the performance of \( SR(\hat{\Qb}) \) improves. In Appendix~\ref{sec:optimalQ}, we prove that if \(\Qb\) is drawn from a fixed, finite‑dimensional family (fixed degree of freedom), then  \(\hat{SR}(\hat{\Qb})\) and \(SR(\hat{\Qb})\) are asymptotically consistent under mild  conditions. It will be interesting to further investigate the relation between discrepancy and the flexibility of optimization  in future work.

\section{Additional Real Data Experiments}
\label{sec:appendixB}

\subsection{Global minimum variance portfolio}
\label{subsec:globalmin}
In this section, we consider the global minimum variance portfolio, which has the advantage of no requirement about the knowledge of $\bmu$. Our construction of the global minimum variance portfolio is defined as  $\wb\propto (\hbSigma+\Qb)^{-1}\one$, where we simply apply the regularization $\Qb \in\cQ_1$ or $\Qb\in\cQ_2$, where $\cQ_1$ and $\cQ_2$ are defined at beginning of Section~\ref{sec:realdata}.  Since the total weights sum up to 1, the exact weight vector $\wb$ takes the form:
\begin{align}
     \wb= \frac{(\hbSigma+\Qb)^{-1}\one}{\one^\top(\hbSigma+\Qb)^{-1}\one}.
    \label{eq:minvar_port}
\end{align}
Note that the scale is determined slightly differently from the mean variance portfolio, for the reason that   $\one^\top(\hbSigma+\Qb)^{-1}\one$ is always positive in the global minimum variance portfolio, but   $\one^\top(\hbSigma+\Qb)^{-1}\bmu$ may become negative and unstable in the MV portfolio. Consequently, the variance of the portfolio can be expressed as:
\begin{align*}
     \Var_{\tilde\Rb}(\wb^\top\tilde\Rb)=\frac{\one^\top(\hbSigma+\Qb)^{-1}\bSigma(\hbSigma+\Qb)^{-1}\one}{(\one^\top(\hbSigma+\Qb)^{-1}\one)^2},
\end{align*}
where $\tilde\Rb$ is an out-of-sample return vector with  mean $\bmu$ and covariance  $\bSigma$.
Theorem~\ref{thm:main_theorem} provides an in-sample method to estimate
\begin{align*}
     1/\sqrt{\Var_{\tilde\Rb}(\wb^\top\xb)}=\frac{\one^\top(\hbSigma+\Qb)^{-1}\one}{\sqrt{\one^\top(\hbSigma+\Qb)^{-1}\bSigma(\hbSigma+\Qb)^{-1}\one}},
\end{align*}
which exactly has the form of \eqref{eq:def_SRQ}, replacing $\bmu$ with $\one$.
To minimize the out-of-sample $\Var_{\tilde\Rb}(\wb^\top\xb)$, Theorem~\ref{thm:main_theorem} can help us determine the best $\Qb^*$ from a finite candidate set. Below we give the steps of our real data analysis for the global minimum portfolio.

\begin{figure}[t!]
    \centering
    \subfloat[One year, $\cQ = \cQ_1, c>1$]{\includegraphics[width=0.33\textwidth]{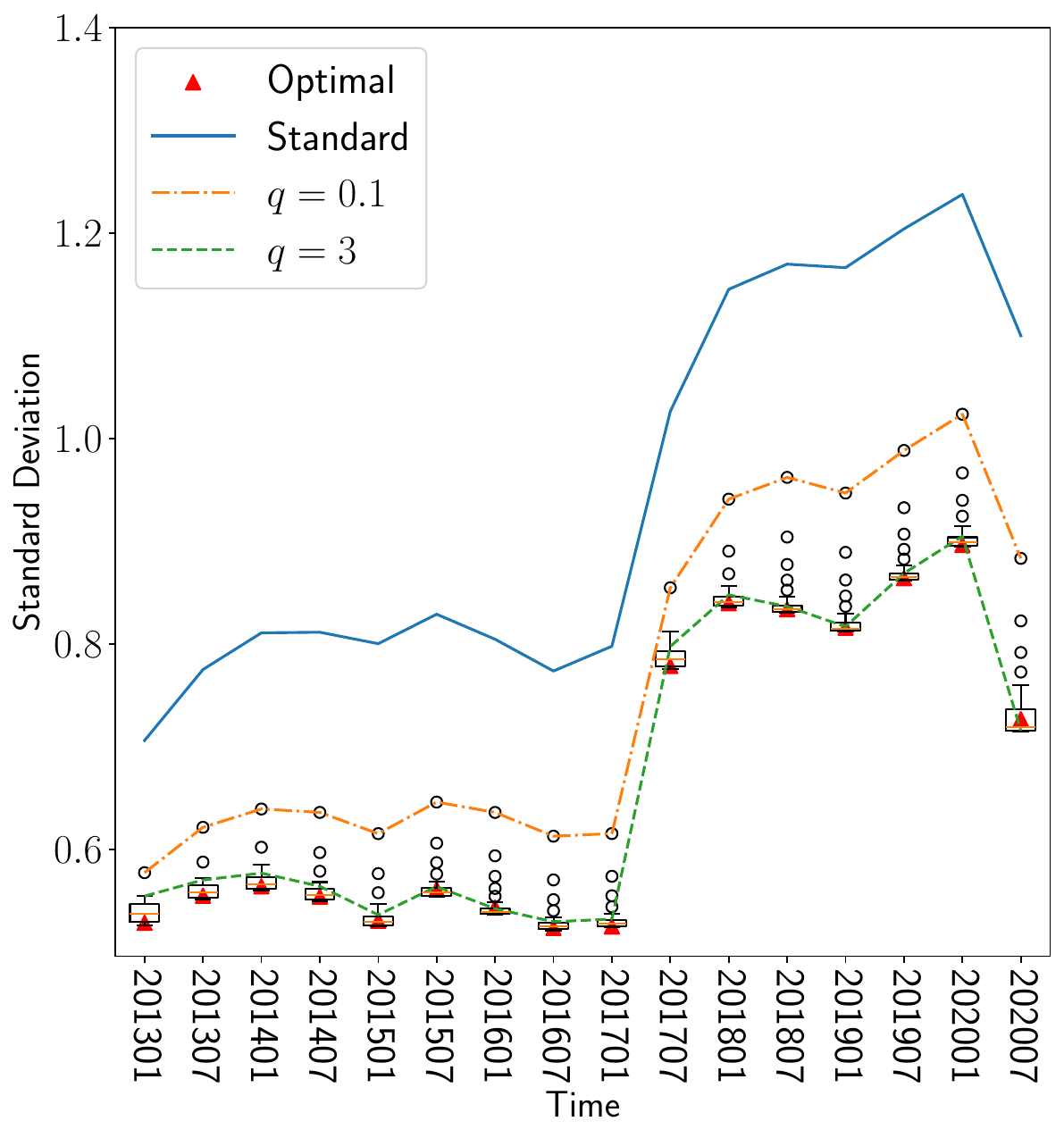}\label{fig_realsd:oneyear}}
    \subfloat[Two years, $\cQ = \cQ_1, c<1$]{\includegraphics[width=0.33\textwidth]{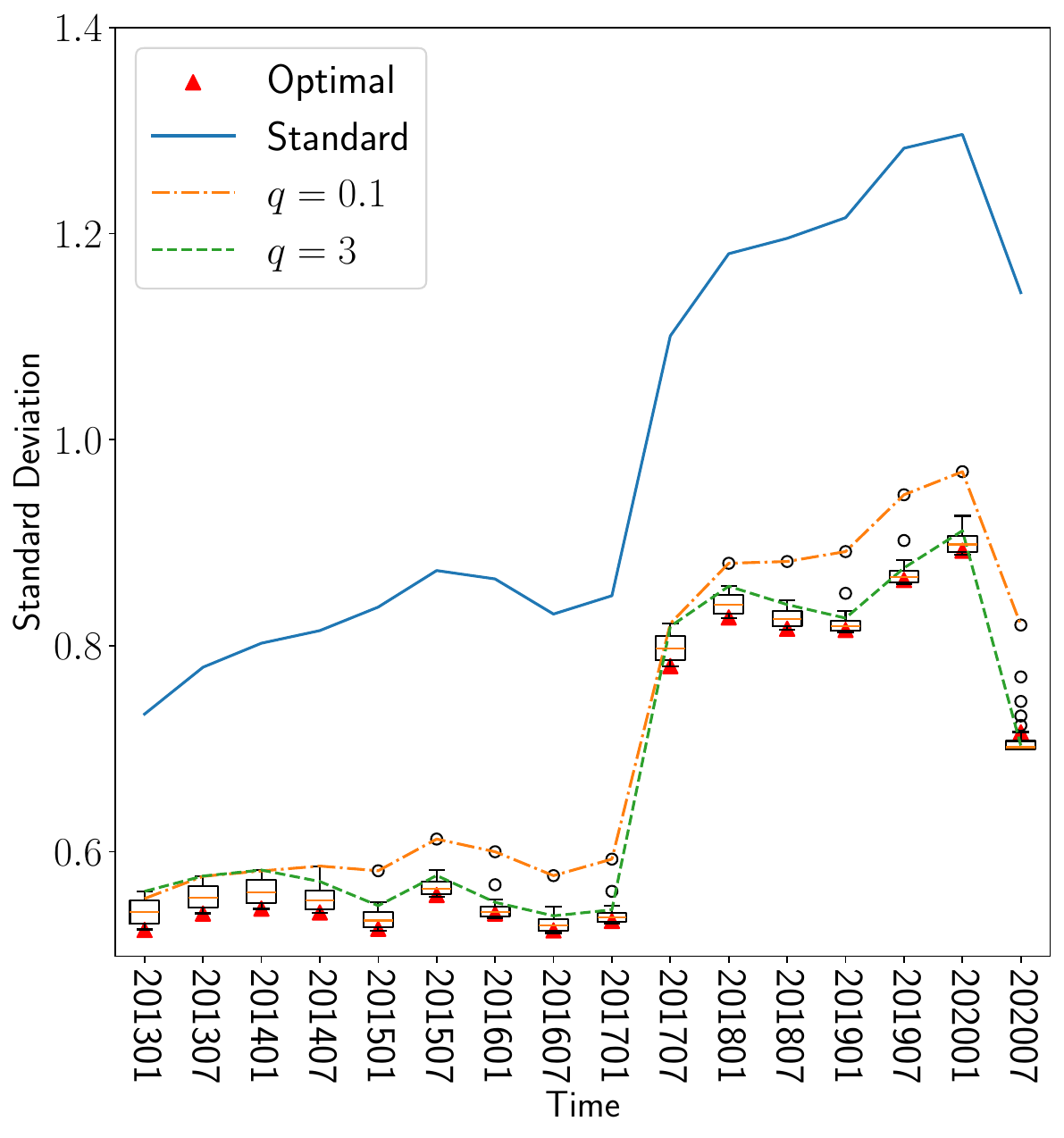}\label{fig_realsd:twoyear}}
    \subfloat[Four years, $\cQ = \cQ_1, c<1$]{\includegraphics[width=0.33\textwidth]{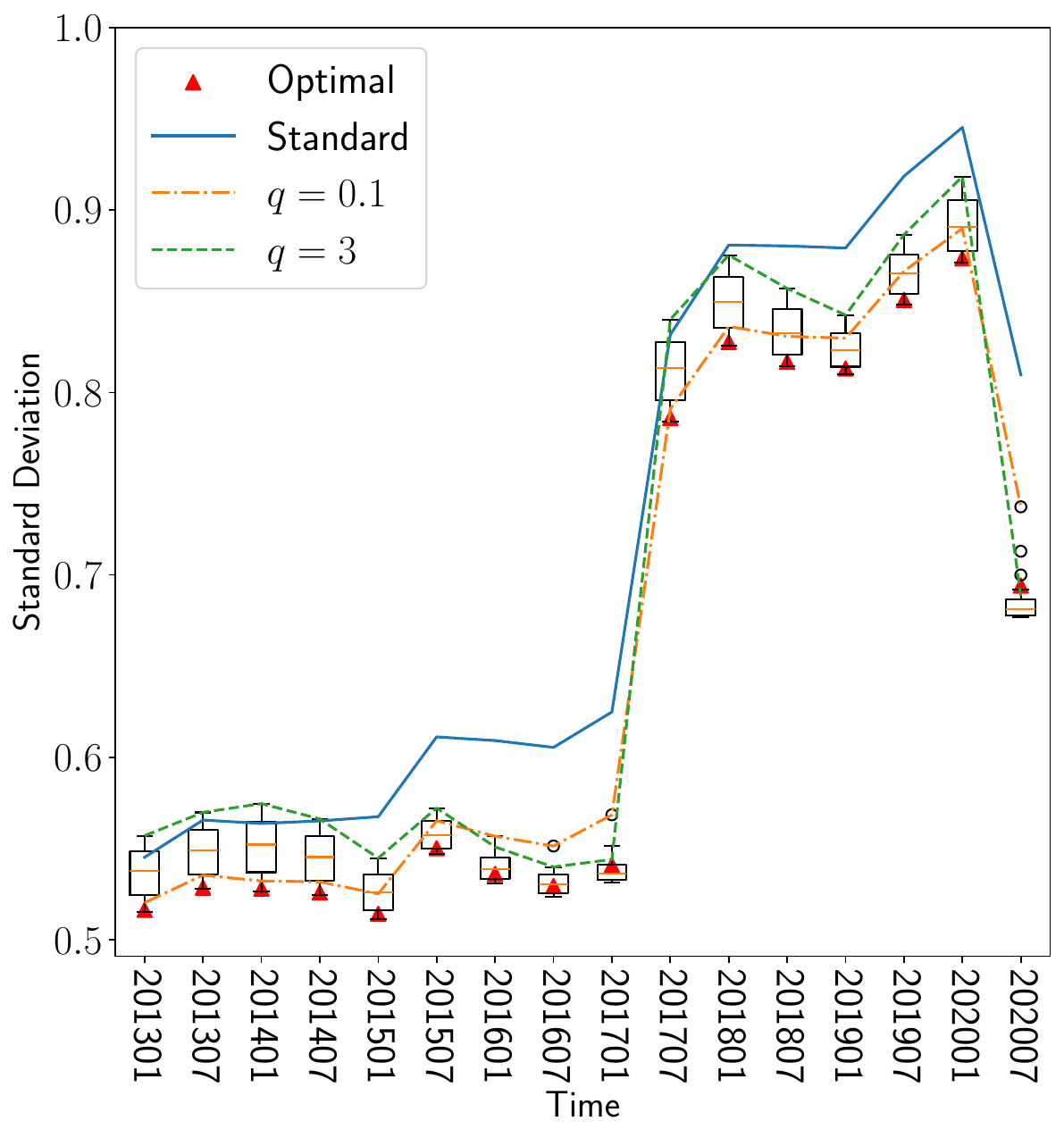}\label{fig_realsd:fouryear}}

    \subfloat[One year, $\cQ = \cQ_2, c>1$]{\includegraphics[width=0.33\textwidth]{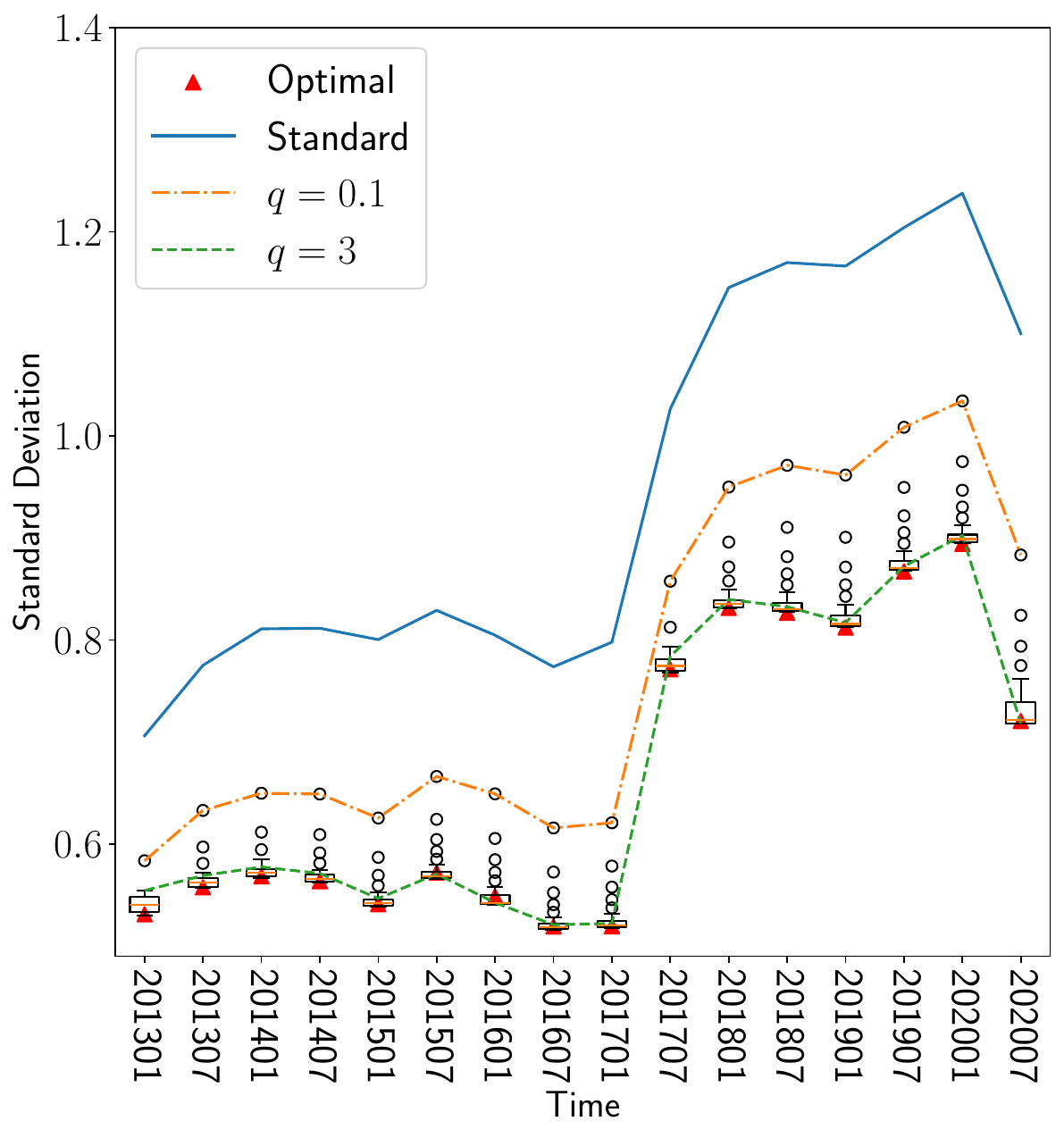}\label{fig_realsd:oneyear1}}
    \subfloat[Two years, $\cQ = \cQ_2, c<1$]{\includegraphics[width=0.33\textwidth]{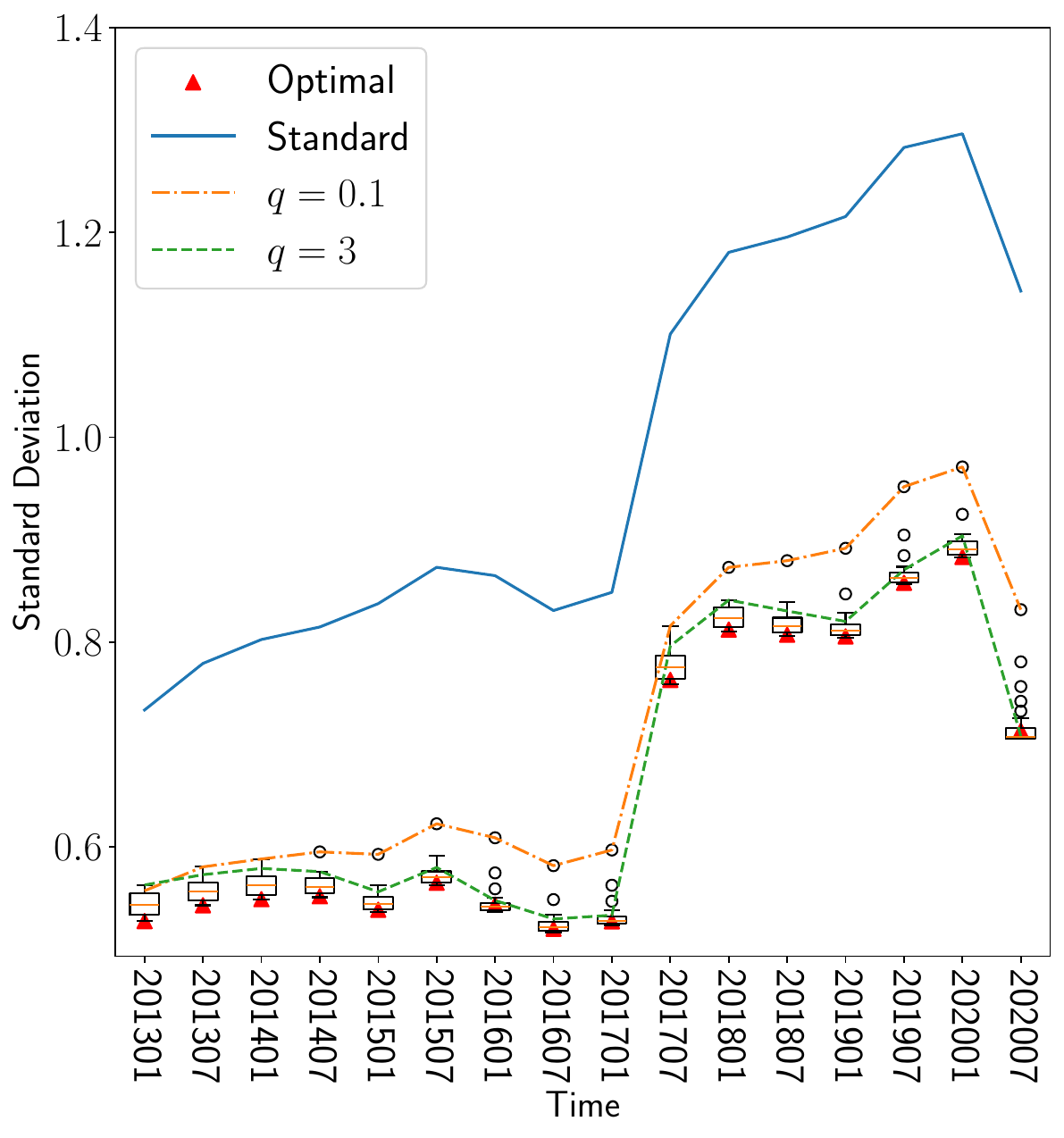}\label{fig_realsd:twoyear1}}
    \subfloat[Four years, $\cQ = \cQ_2, c<1$]{\includegraphics[width=0.33\textwidth]{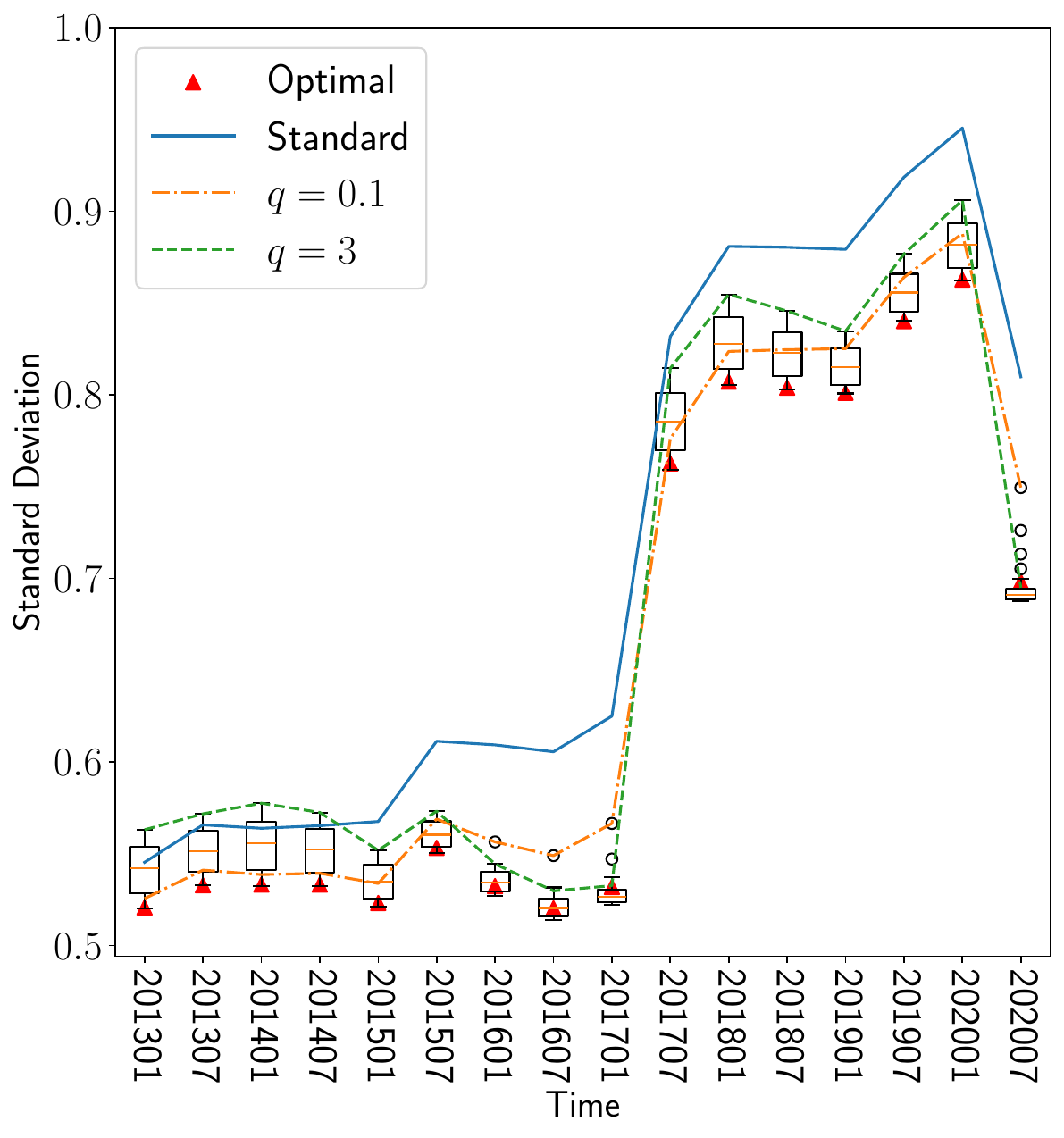}\label{fig_realsd:fouryear1}}
    
    \caption{Standard deviation (volatility) of global minimum variance portfolios. The x-axis represents the rolling time, while the y-axis represents the standard deviation of the portfolio returns every three years. The blue solid line, orange dash-dot line and green dash line in the figure correspond to the  standard deviation with  $q=0$, $q=0.1$ (minimum $q$) and $q=3$ (maximum $q$) respectively. The boxplot displays the standard deviation for all  $\Qb\in\cQ$, and the red triangle indicates the standard deviation  under optimized  $\Qb^*\in\cQ$.}
    \label{fig_realsd}
\end{figure}

\begin{enumerate}[leftmargin=*]
    \item For each testing month, we observe the daily returns $\Rb\in\RR^{n\times p}$, where $n$ is the total number of trading days with one-, two-, and four-year historical data length and $p=365$ is the number of selected stocks. Calculate the sample covariance matrix $\hbSigma$ and the portfolio weight $\wb$ is given by \eqref{eq:minvar_port}. 
    \item  For each testing month, we run experiments for all candidate  $\Qb$ and also consider no regularization, i.e.  $\Qb=\zero$, where we have $\wb\propto \hbSigma^{+}\one$ and $\hbSigma^{+}$ is the pseudo inverse, and the optimized  $\Qb^*\in\cQ$ using the estimation in Theorem~\ref{thm:main_theorem}. Here, 
    \begin{align*}
         \Qb^*=\argmax_{\Qb\in\cQ}\bigg(1-\frac{c}{p}\tr\hbSigma(\hbSigma+\Qb)^{-1}\bigg)\cdot\frac{\one^\top(\hbSigma+\Qb)^{-1}\one}{\sqrt{\one^\top(\hbSigma+\Qb)^{-1}\hbSigma(\hbSigma+\Qb)^{-1}\one}}.
    \end{align*}
    \item  We roll the procedure above for all testing months. Note that the value of  $\Qb^*$ changes from month to month. 
    With the weight vector $\wb$ using all  $\Qb\in\cQ$, $\Qb=\zero$ or $\Qb^*$, we can then compute the portfolio returns for each trading day in the testing month.
    \item We report the standard deviation of the portfolio returns every three years, as shown in Figure~\ref{fig_realsd}. For instance, the results labeled on the x-axis as 202001 in Figure~\ref{fig_realsd} represent the standard deviation of out-of-sample portfolio returns from Jan 2020 to Dec 2022 with 36 testing months. When the x-axis label updates to 202007, the results reflect the standard deviation of returns from Jul 2020 to Jun 2023. 
\end{enumerate}

Figure~\ref{fig_realsd:oneyear}-\ref{fig_realsd:fouryear} present the results obtained using historical data of varying lengths with  $\Qb\in\cQ_1$, and Figure~\ref{fig_realsd:oneyear1}-\ref{fig_realsd:fouryear1} use  $\Qb\in\cQ_2$. From Figure~\ref{fig_realsd}, we can make a few conclusions. Firstly, with no regularization ( $\Qb=\zero$, the blue curve), the out-of-sample volatility is very high. So adding proper regularization to protect the condition of the sample covariance is absolutely necessary.
Secondly, as shown in Figure~\ref{fig_realsd}, for each of the 6 cases, if we fix one choice $\Qb\in\cQ$ for all testing months, for most months  $q=0.1$ (the orange curve,  minimum tuning parameter) achieves the smallest out-of-sample volatility, but occasionally  $q=3$ (the green curve,  minimum tuning parameter) outperforms  $q=0.1$. If the volatility happens to be monotone in the range of $q$ values, the boxplot will be in the middle between the orange and green curves, but if the best $q$ is neither the minimum $q$ nor the maximum $q$, we also see cases that the boxplot is beneath both the orange and green curves. 
This indicates that applying a constant level of regularization for all months such as  $q = 0.1$ (minimum tuning parameter) to the portfolio may not be suitable for all time periods. 
So if we want to update  $\Qb$ from month to month, it is important to have a in-sample method to estimate the out-of-sample volatility. 
Thirdly, from all 6 cases in Figure~\ref{fig_realsd}, we see that the portfolios with optimized  $\Qb^*$ (the red triangles), automatically updated every month, consistently exhibit a smaller standard deviation regardless of the candidate sets and the varying lengths of historical data.
So we can conclude that our method can consistently and flexibly select an appropriate  $\Qb\in\cQ$ based on the data itself, leading to smaller out-of-sample variances of the global minimum variance portfolio and demonstrating significant advantages over naively choosing a fixed  $\Qb\in\cQ$ for all time points. The last side observation is that the volatility level jumps to a higher level starting from the label 201707, which actually represents Jul 2017 to Jun 2020, the first time period covering the COVID outbreak. 

\subsection{Calibrated models with real data}
\label{sec:appendixB2}
In the previous experiments for mean-variance portfolios, we made the assumption that the mean vector $\bmu$ in S\&P500 is known.   In this section, we hope to look into how our methods are affected when $\bmu$ is unobserved. To this end, we use real data to calibrate a model for $\bmu$ and then use the calibrated model to simulate synthetic data. 
Here are the steps.
\begin{enumerate}[leftmargin=*]
    \item Using solely the returns in the testing month $j$, we have the out-of-sample average return vector $\bmu_j$ and the out-of-sample estimated covariance matrix $\tilde\bSigma_j=\tilde{\Xb}_j^\top\tilde{\Xb}_j/\tilde{n}_j$. Here, $\tilde{\Xb}_j$ is the return matrix for the testing month $j$, and $\tilde{n}_j$ is the number of trading days in that month. Additionally, we define the in-sample empirical covariance matrix $\hbSigma_j=\Xb_j^\top\Xb_j/n_j$, where $\Xb_j$ is the return matrix consisting of five years prior to the testing month $j$, and $n_j$ is the total number of trading days in those five years.
    \item Using the well-defined $\bmu_j$, $\tilde\bSigma_j$, and $\hbSigma_j$, we generate returns of $21$ trading days for each trading month. The returns for the trading days in the month $j$ are drawn from the distribution $\cN(\check{\bmu}_j,\check{\bSigma}_j)$, where:
\begin{equation}
\label{eq:realdatacalibrate}
    \begin{split}
    \check{\bmu}_j=\rho\cdot\bmu_j+\sqrt{1-\rho^2}\cdot\cN(0,\sigma^2_{\bmu_j}\Ib),\quad\check{\bSigma}_j=0.95\cdot\hbSigma_j+\sqrt{1-0.95^2}\cdot \tilde\bSigma_j.
    \end{split}
\end{equation}
Here, $\rho>0$ is a tuning parameter that controls the level of perturbation to $\bmu_j$, and $\sigma_{\bmu_j}$ is the standard deviation of the elements in $\bmu_j$.

\begin{figure}[t!]
    \centering
    \subfloat[$\cQ_1$, two years]{\includegraphics[width=0.24\textwidth]{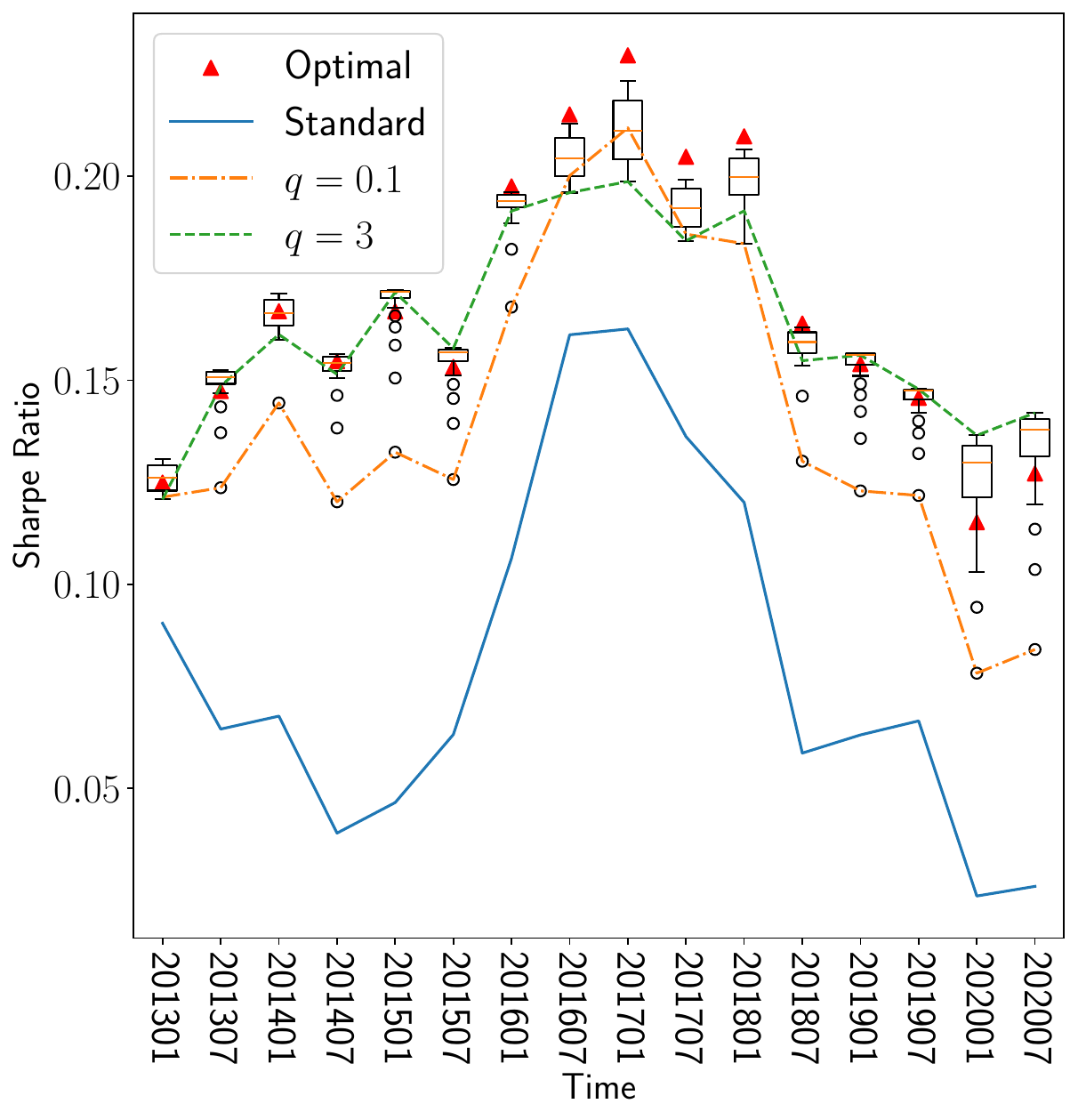}\label{fig_calibrate:largesimuoneyear2-0.05}}
    \subfloat[$\cQ_2$, two years]{\includegraphics[width=0.24\textwidth]{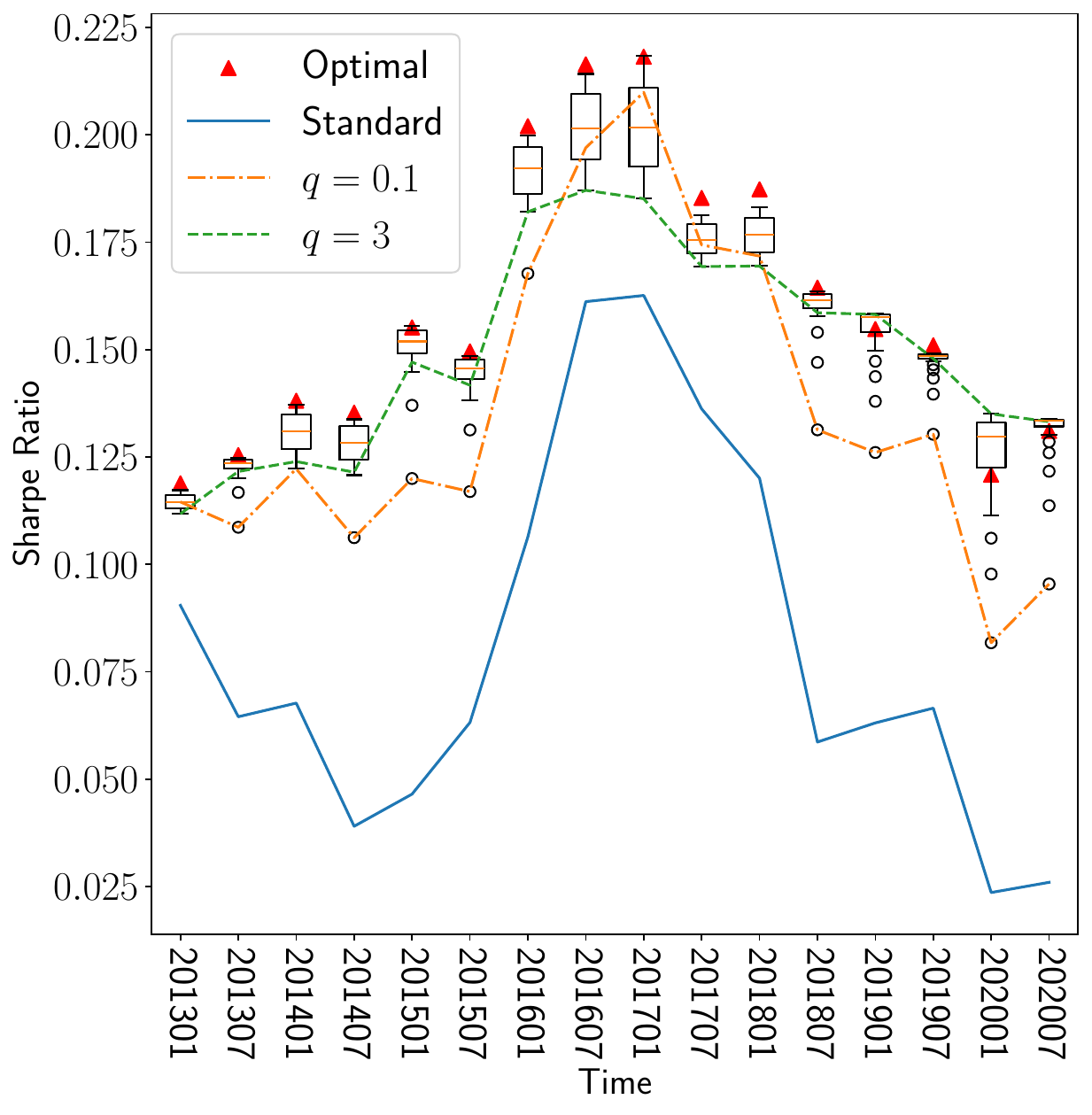}\label{fig_calibrate:smallsimuoneyear2-0.05}}
    \subfloat[$\cQ_1$, four years]{\includegraphics[width=0.24\textwidth]{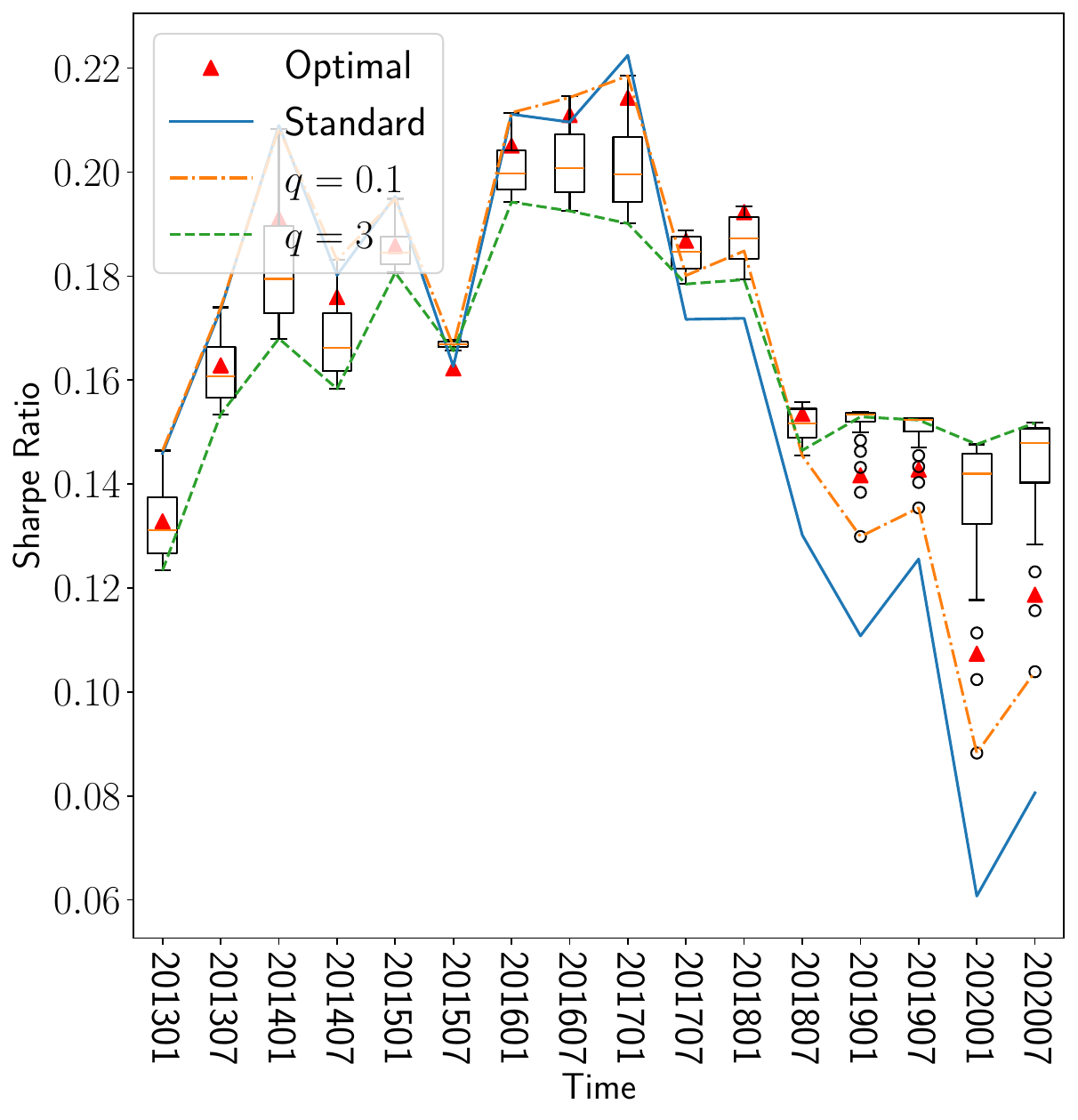}\label{fig_calibrate:largesimuoneyear4-0.05}}
    \subfloat[$\cQ_2$, four years]{\includegraphics[width=0.24\textwidth]{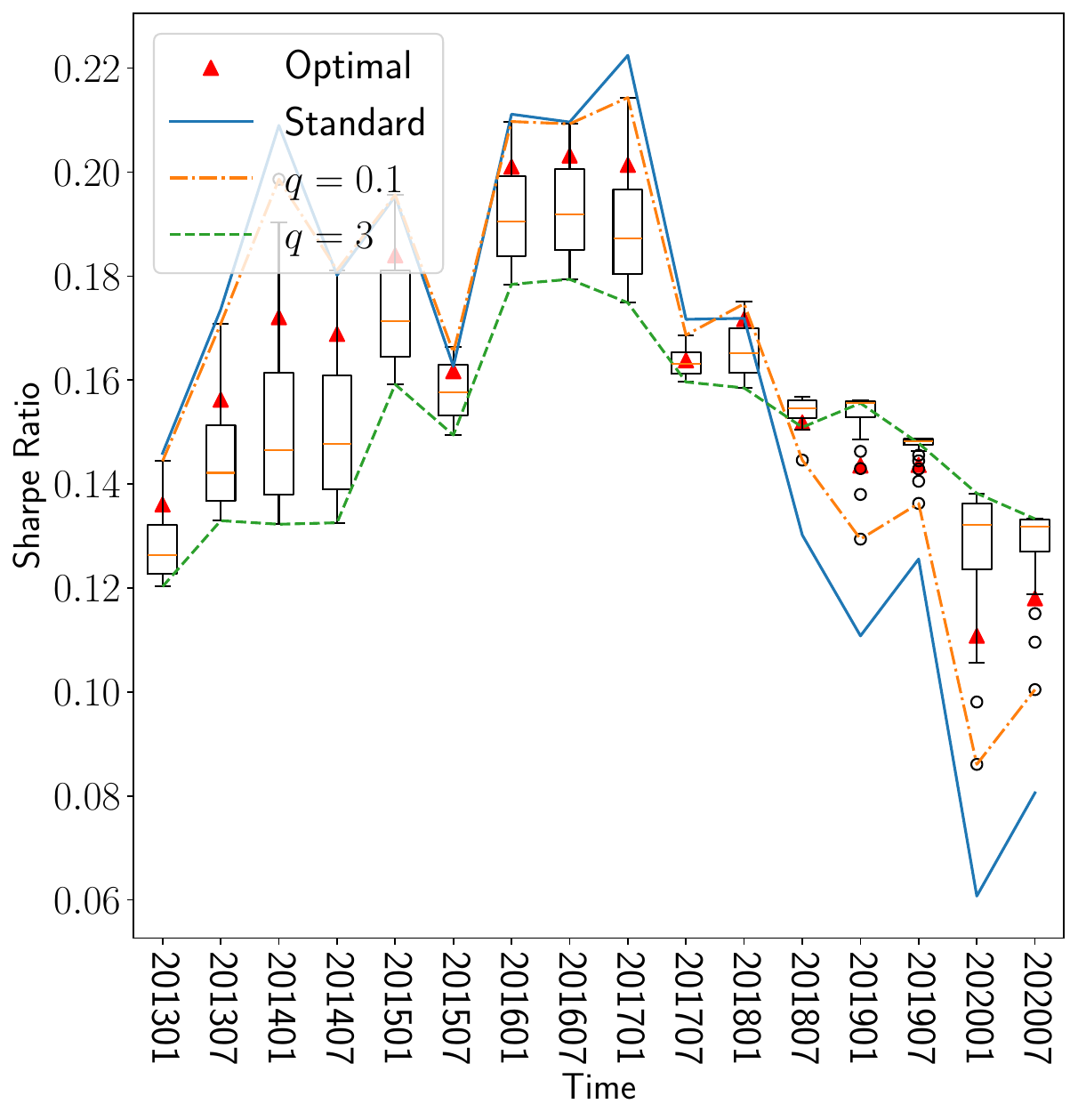}\label{fig_calibrate:smallsimuoneyear4-0.05}}

    \subfloat[$\cQ_1$, two years]{\includegraphics[width=0.24\textwidth]{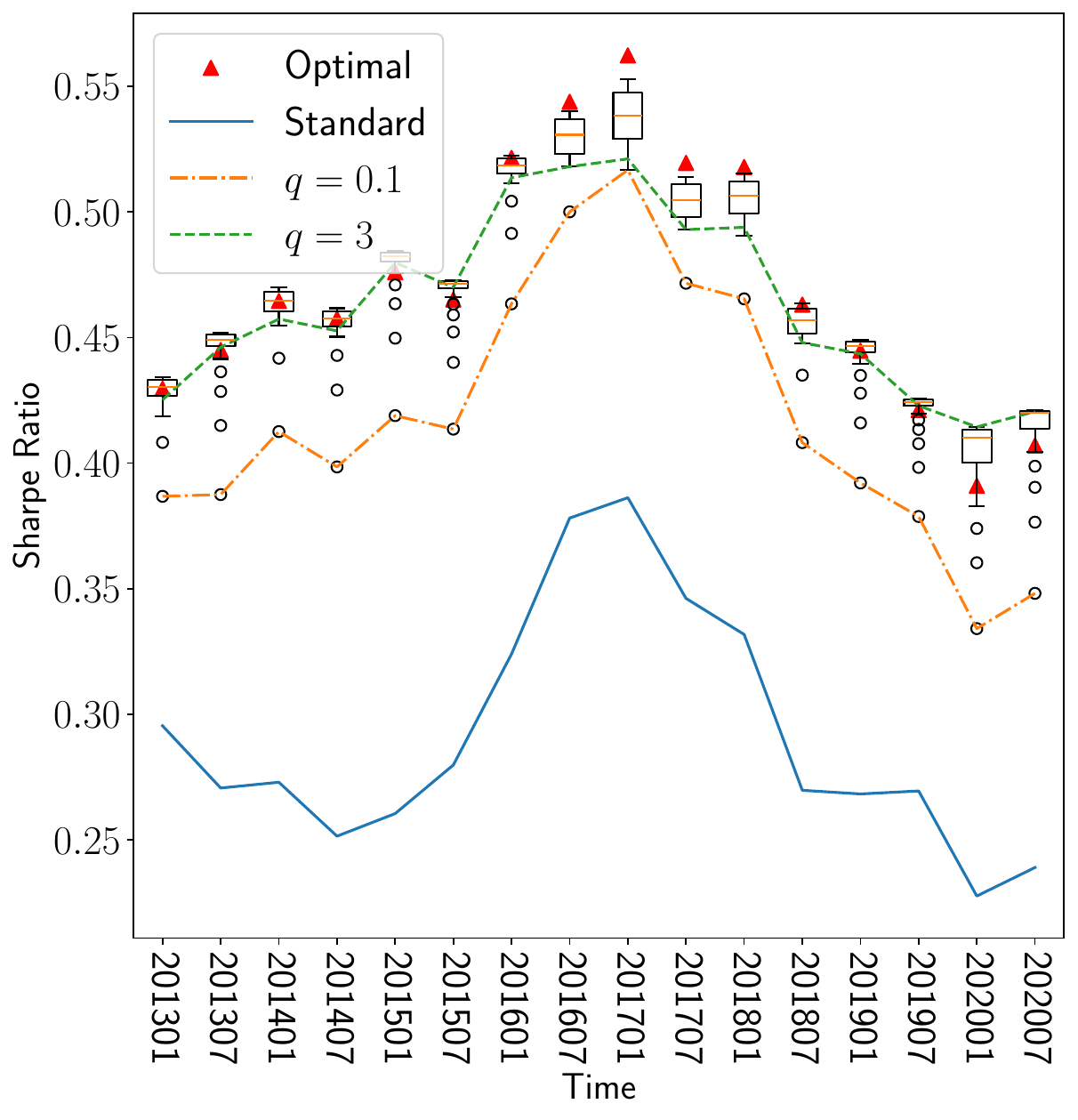}\label{fig_calibrate:largesimuoneyear2-0.2}}
    \subfloat[$\cQ_2$, two years]{\includegraphics[width=0.24\textwidth]{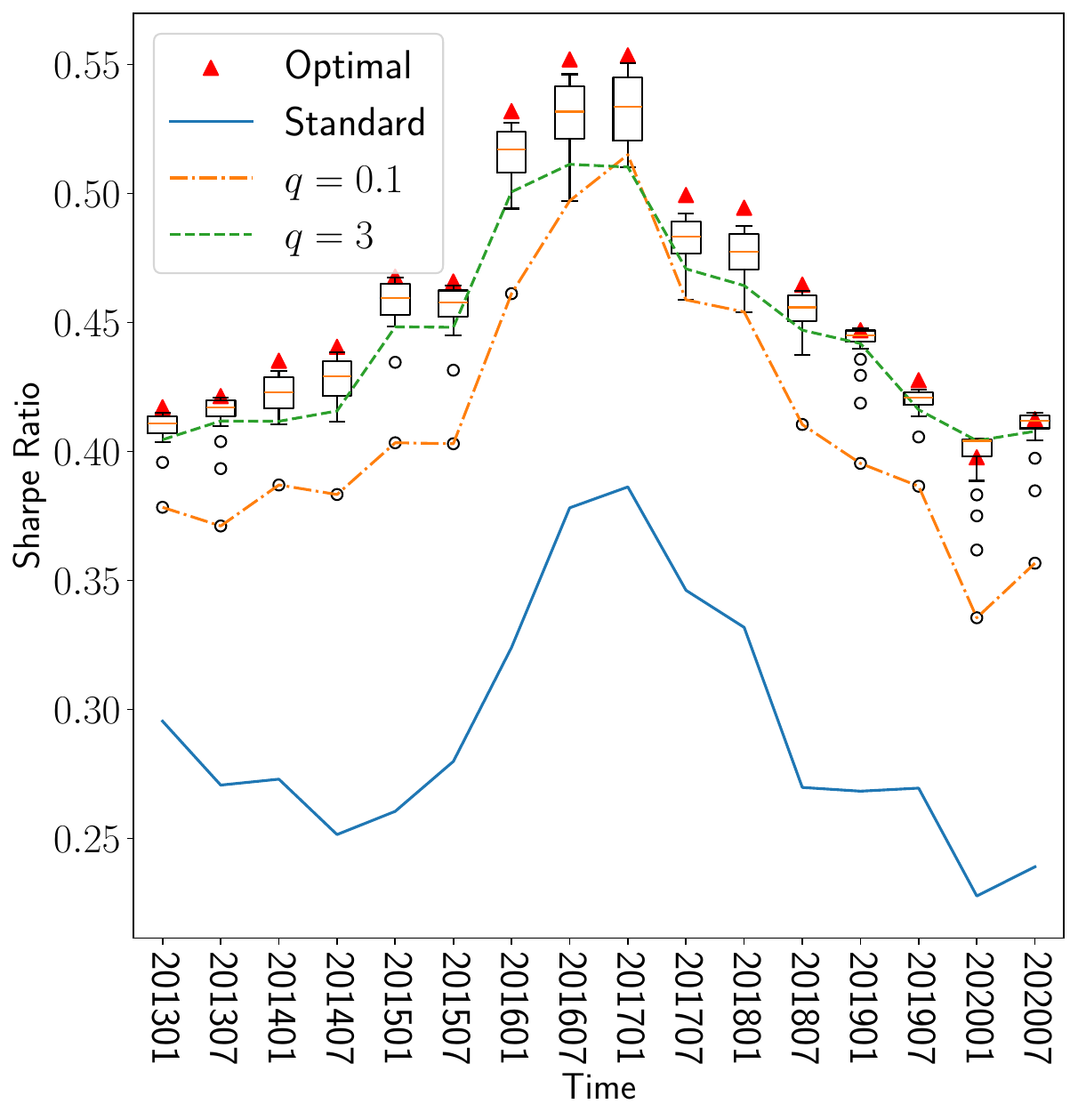}\label{fig_calibrate:smallsimuoneyear2-0.2}}
    \subfloat[$\cQ_1$, four years]{\includegraphics[width=0.24\textwidth]{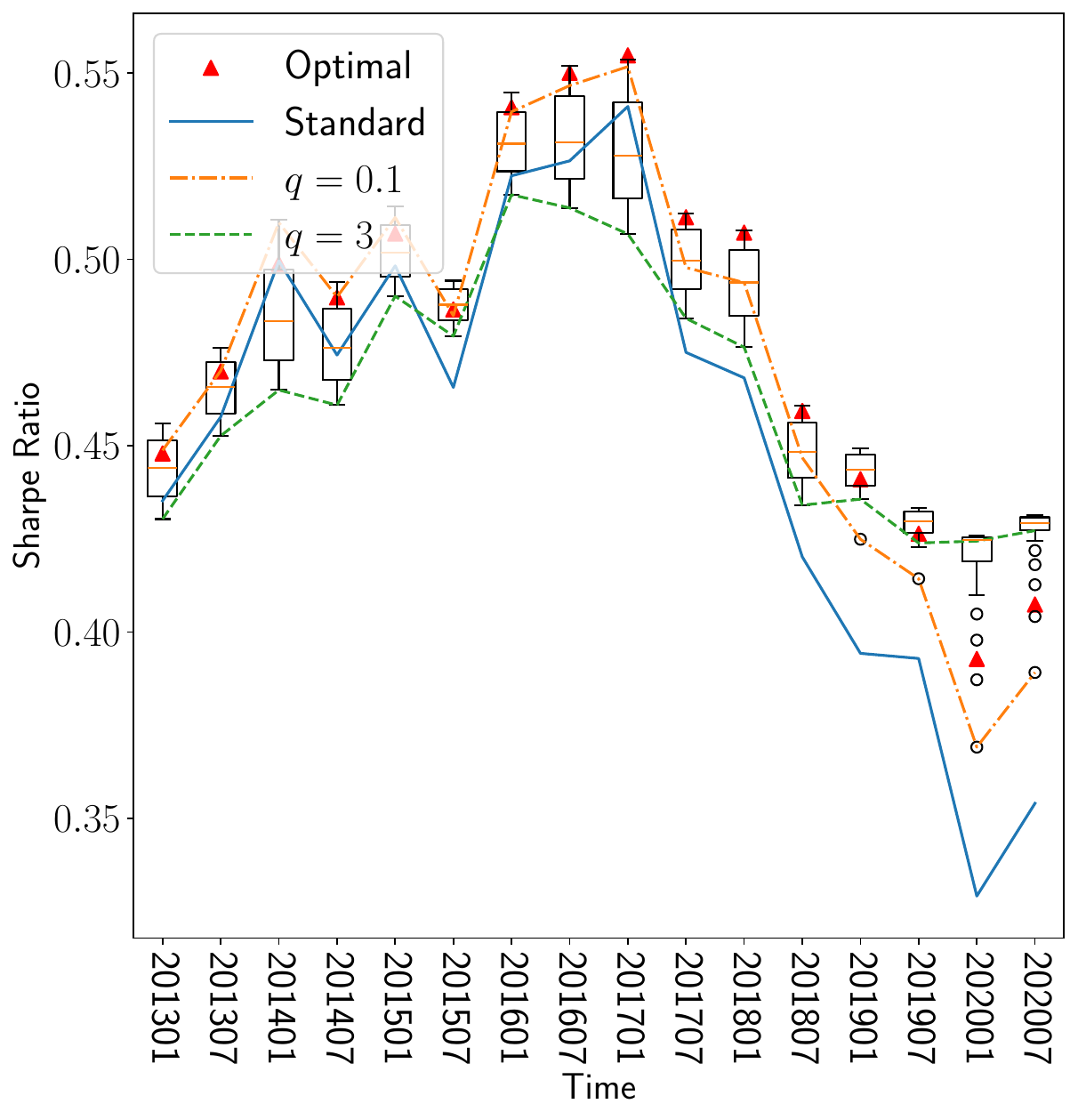}\label{fig_calibrate:largesimuoneyear4-0.2}}
    \subfloat[$\cQ_2$, four years]{\includegraphics[width=0.24\textwidth]{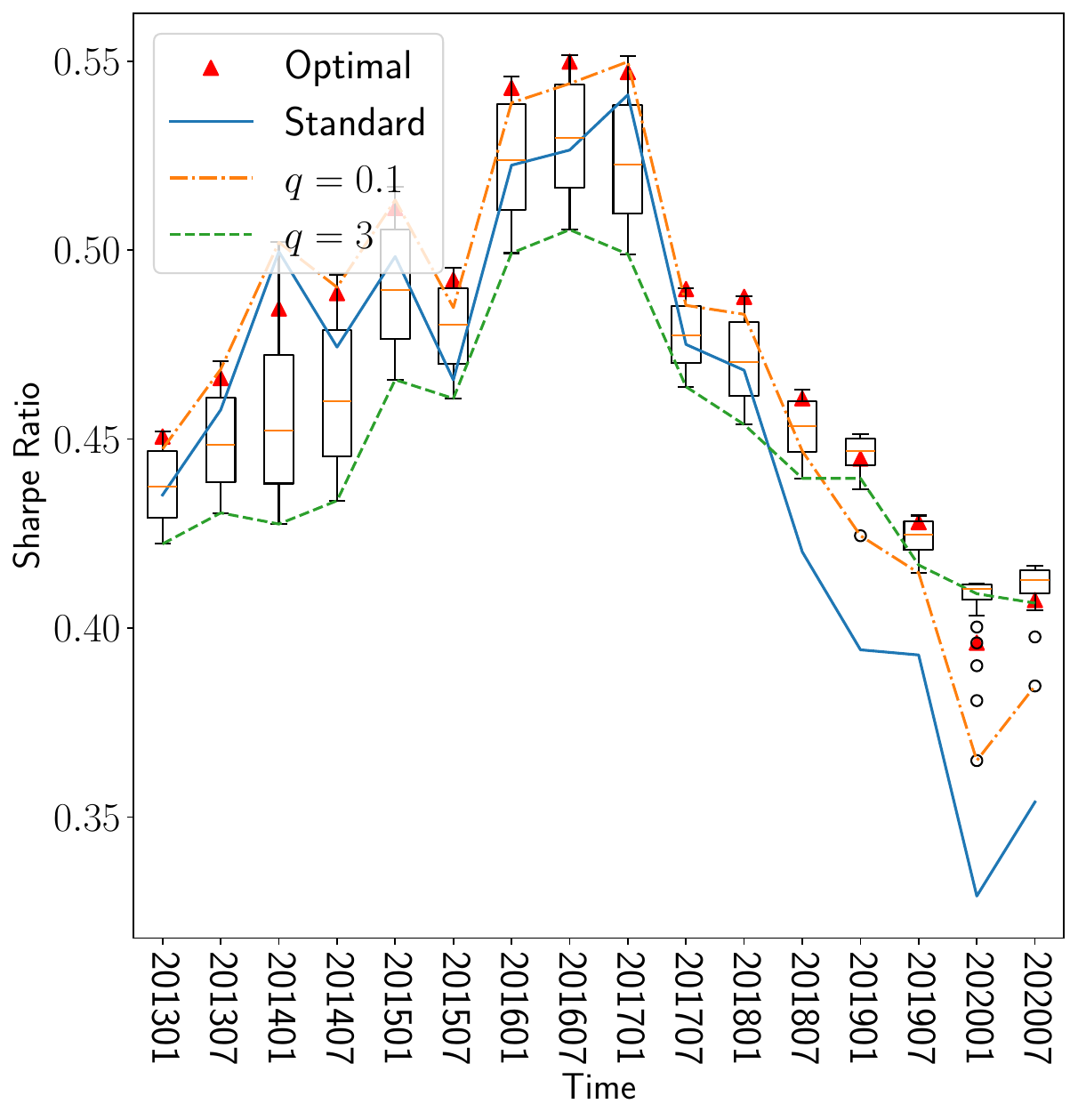}\label{fig_calibrate:smallsimuoneyear4-0.2}}

    \subfloat[$\cQ_1$, two years]{\includegraphics[width=0.24\textwidth]{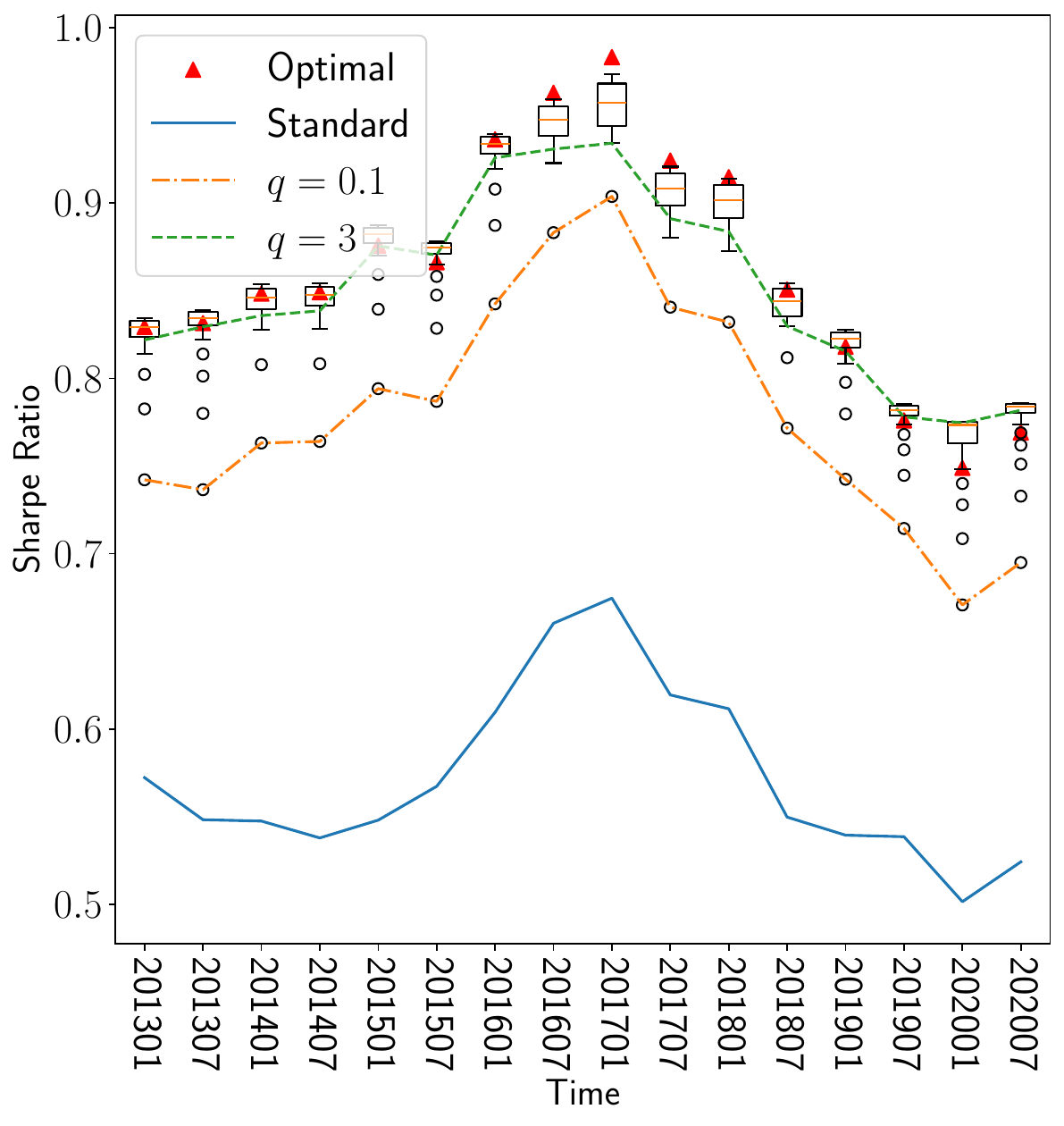}\label{fig_calibrate:largesimuoneyear2-0.4}}
    \subfloat[$\cQ_2$, two years]{\includegraphics[width=0.24\textwidth]{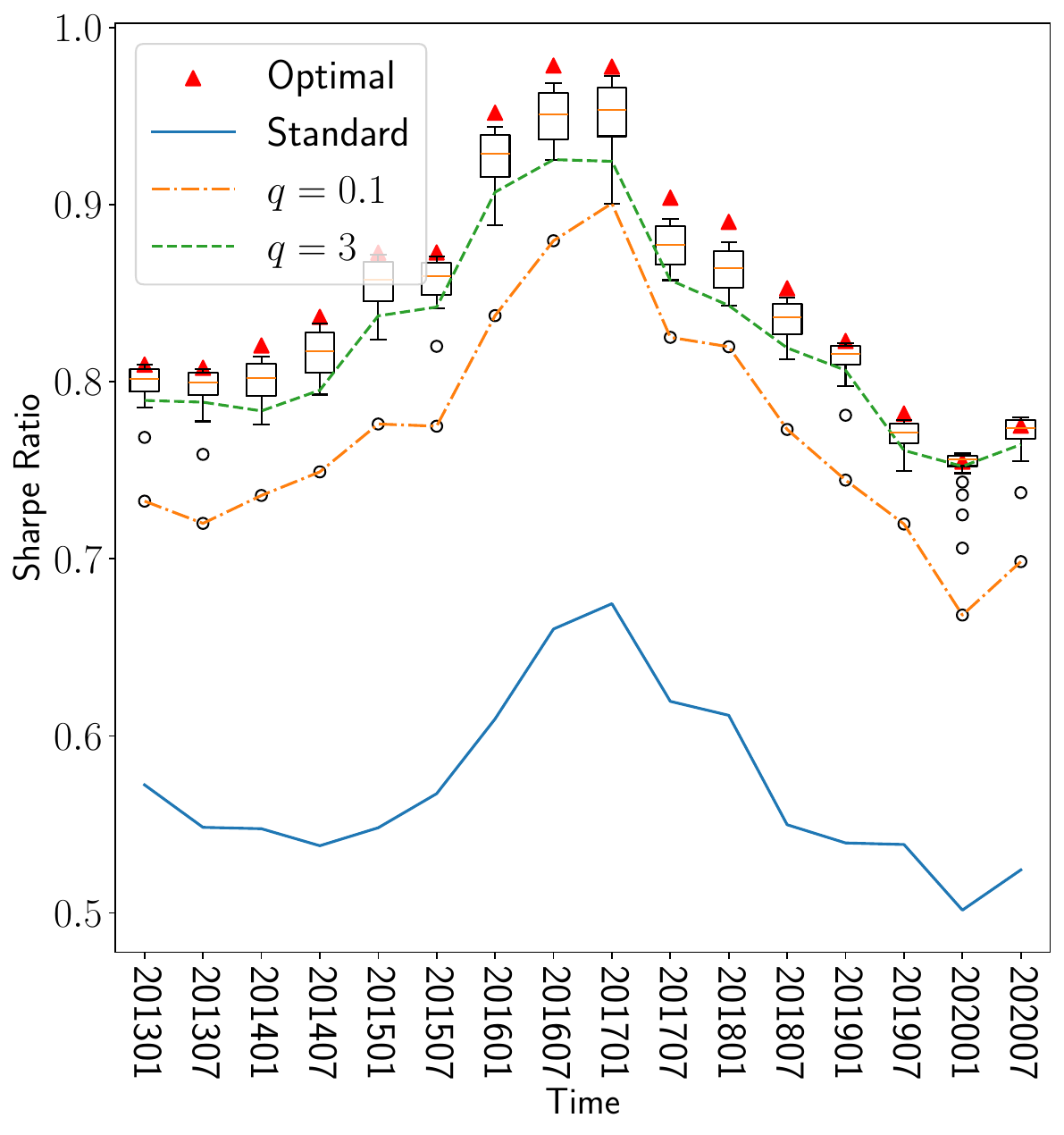}\label{fig_calibrate:smallsimuoneyear2-0.4}}
    \subfloat[$\cQ_1$, four years]{\includegraphics[width=0.24\textwidth]{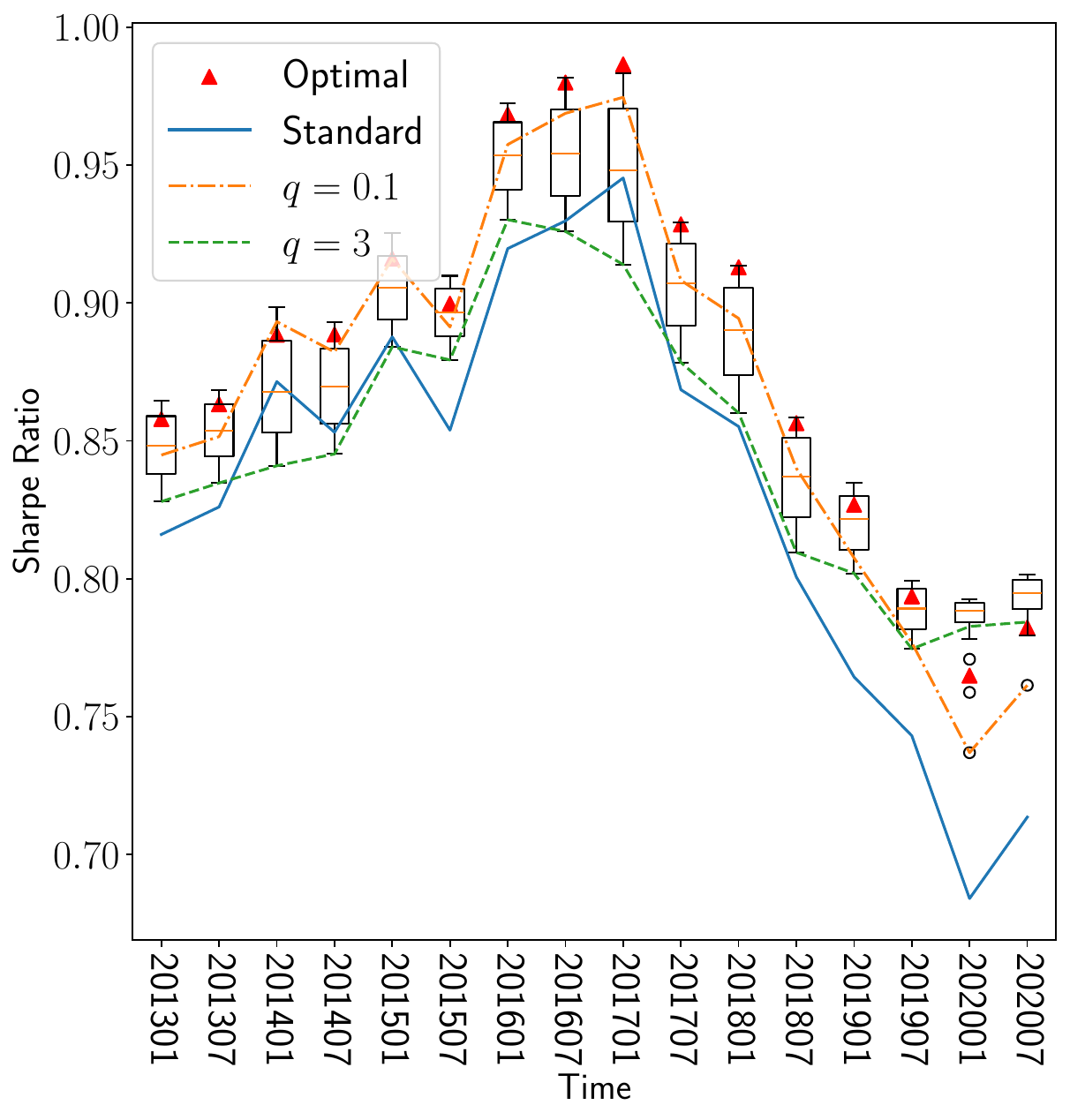}\label{fig_calibrate:largesimuoneyear4-0.4}}
    \subfloat[$\cQ_2$, four years]{\includegraphics[width=0.24\textwidth]{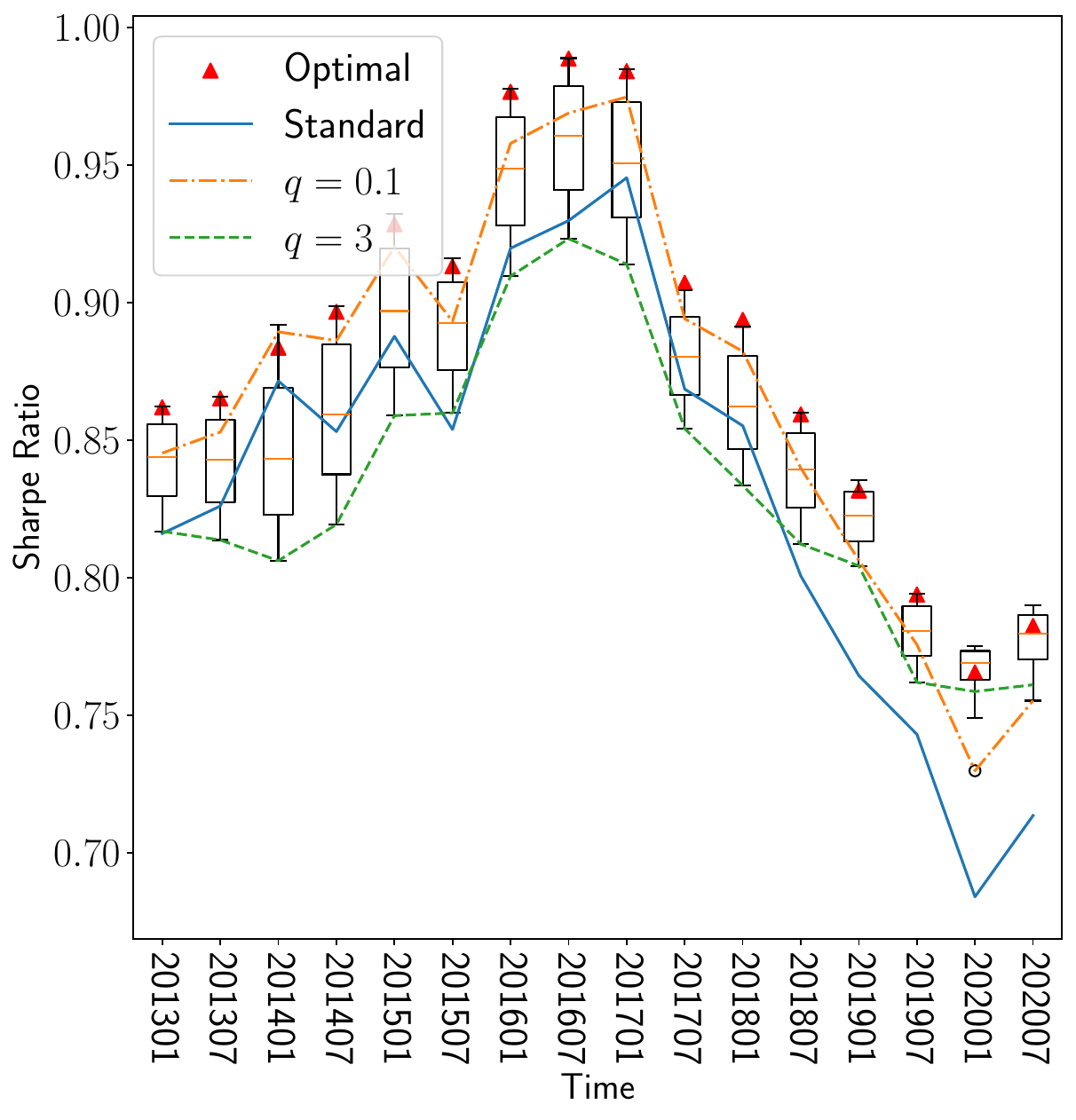}\label{fig_calibrate:smallsimuoneyear4-0.4}}
    
    \caption{Sharpe ratio of mean-variance portfolio with calibrated models. The three rows of plots correspond to $\rho=0.05, 0.2, 0.4$ respectively. The x-axis represents the rolling time, while the y-axis represents the Sharpe ratio of portfolio returns every three years.}
    \label{fig_calibratesimurho}
\end{figure}

\item To use the calibrated model, we assume that the mean vector $\bmu_j$ is known and used in portfolio construction, but the actual true mean of the returns for the month $j$ is $\check{\bmu}_j$ which is unknown. We use different values of $\rho$ to generate different calibrated models, where a larger $\rho$ reflects a higher signal level of predicting $\check{\bmu}_j$ using $\bmu_j$. 
The allocation vector for the testing month $j$ is then constructed as  $\wb\propto (\hbSigma+\Qb)^{-1}\bmu_j$, and the remaining steps follow the same procedure in Section~\ref{subsec:meanvariance}.
\end{enumerate}

We present experimental results for $\rho=0.05$, $0.2$, and $0.4$ in three rows of plots respectively in Figure~\ref{fig_calibratesimurho}. As indicated in equation \eqref{eq:realdatacalibrate}, smaller values of $\rho$ result in larger perturbations on the mean vector $\bmu_j$. 
$\rho=0.05$, $0.2$, and $0.4$, the mean cosine values of the angle between $\bmu_j$ and $\check{\bmu}_j$ are $0.0428$, $0.1368$, and $0.2605$, respectively, with corresponding standard deviations of $0.2676$, $0.2649$, and $0.2562$.
We observed that as $\rho$ increases, the Sharpe ratio for the portfolio gradually increases. This makes intuitive sense as higher prediction power leads to better portfolio out-of-sample performance. 
Specifically, for $\rho=0.2$ and $0.4$, the optimized  $\Qb^*$ based on our in-sample SR estimator always displays the best performance among the candidate  $\Qb$ in $\cQ=\cQ_1$ or $\cQ=\cQ_2$. For $\rho=0.05$, we may not always achieve the best out-of-sample SR due to the low signal level ($\check{\bmu}_j$ becomes nearly orthogonal to $\bmu_j$). However, we still see comparable performances. 

\subsection{Mean variance portfolio with unknown mean vector}
\label{subsec:meanvarianceunknown}
In this section, we investigate the MV portfolio construction where $\bmu$ is unknown. When $\bmu$ is unknown, we define the MV portfolio allocation vector as $\wb \propto (\hbSigma + \Qb)^{-1} \hbmu$. The portfolio is scaled under the assumption $\sum_{i} |\wb_i| = 1$, which ensures that the total exposure of the portfolio sums to one unit. Thus, the weight vector $\wb$ can be expressed as:
\begin{align}
 \wb= (\hbSigma+\Qb)^{-1}\hbmu/\|(\hbSigma+\Qb)^{-1}\hbmu\|_1.
 \label{eq:meanunknownvar_port}
\end{align}
In this experiment, we change the set $\cQ_1 = \{q \cdot \diag(\hbSigma_{pre}-\lambda_1\bxi_1\bxi_1^\top), q \in [1:30]/10\}$, where $\hbSigma_{pre}$ represents the sample covariance matrix of returns calculated from Jan 2004 to Dec 2008, and $\lambda_1$ and $\bxi_1$ are the largest eigenvalues and corresponding eigenvectors of $\hbSigma_{pre}$, and we still keep the set $\cQ_2$. According to Theorem~\ref{thm:unknown_mu}, for each $\Qb \in \cQ$, we can leverage $\hat{SR}$ to identify the optimal regularization matrix $\Qb$ within $\cQ_1$ or $\cQ_2$ that maximizes the Sharpe ratio ($SR$). We present the steps for our MV portfolio analysis.
\begin{enumerate}[leftmargin=*]
    \item For each testing month, we observe the excess daily returns $\Rb\in\RR^{n\times p}$, where $n$ is the total number of trading days with one-, two-, and four-year historical data length and $p=365$ is the number of selected stocks. Calculate the sample mean $\hbmu$, the sample covariance matrix $\hbSigma$ and the portfolio weight $\wb$ given by \eqref{eq:meanunknownvar_port}.
    \item For each testing month, we run experiments for all candidate  values and also consider no regularization, i.e. $\Qb=\zero$, where we have $\wb\propto \hbSigma^{+}\hbmu$ and $\hbSigma^{+}$ is the pseudo inverse, and the optimized $\Qb^*\in\cQ$ using the estimation in Theorem~\ref{thm:unknown_mu}. Here, 
    \begin{align*}
        \Qb^*=\argmax_{\Qb^*\in\cQ}\bigg(1-\frac{c}{p}\tr\hbSigma(\hbSigma+\Qb)^{-1}\bigg)\cdot\frac{\hbmu^\top(\hbSigma+\Qb)^{-1}\hbmu-\frac{\tr\hbSigma(\hbSigma+\Qb)^{-1}}{n-\tr\hbSigma(\hbSigma+\Qb)^{-1}}}{\sqrt{\hbmu^\top(\hbSigma+\Qb)^{-1}\hbSigma(\hbSigma+\Qb)^{-1}\hbmu}}.
    \end{align*}
    \item  We roll the procedure above for all testing months. Note that the optimal $\Qb^*$ changes from month to month. 
    With the weight vector $\wb$ using all $\Qb\in\cQ$, $\Qb=\zero$ or $\Qb^*$, we can then compute the portfolio returns for each trading day in the testing month.
    \item We report the realized Sharpe ratio of daily portfolio returns over the future three years. 
\end{enumerate}
\begin{figure}[t!]
    \centering
    \subfloat[One year, $\cQ = \cQ_1, c>1$]{\includegraphics[width=0.33\textwidth]{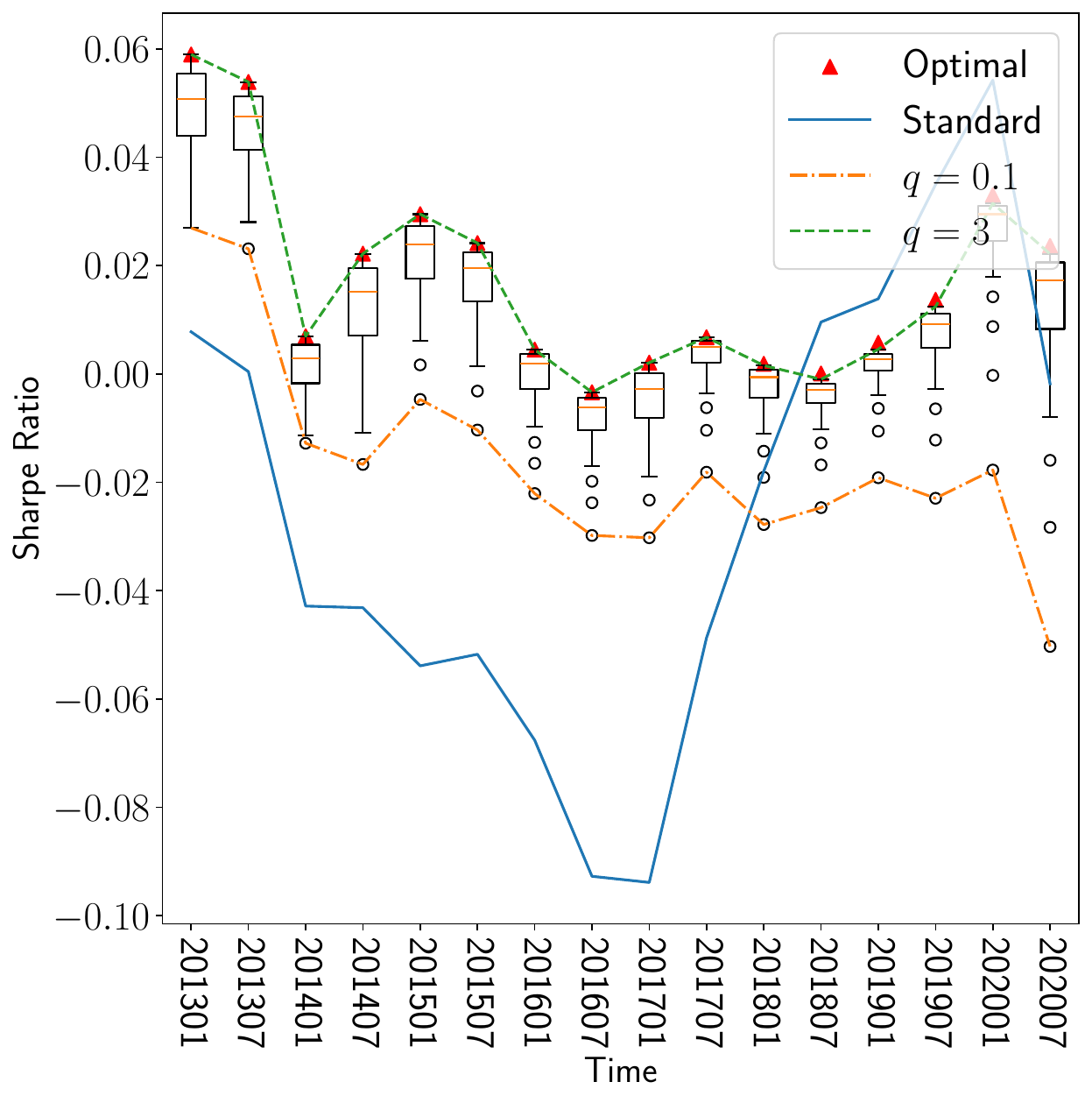}\label{fig_realmeanunknown:oneyear}}
    \subfloat[Two years, $\cQ = \cQ_1, c<1$]{\includegraphics[width=0.33\textwidth]{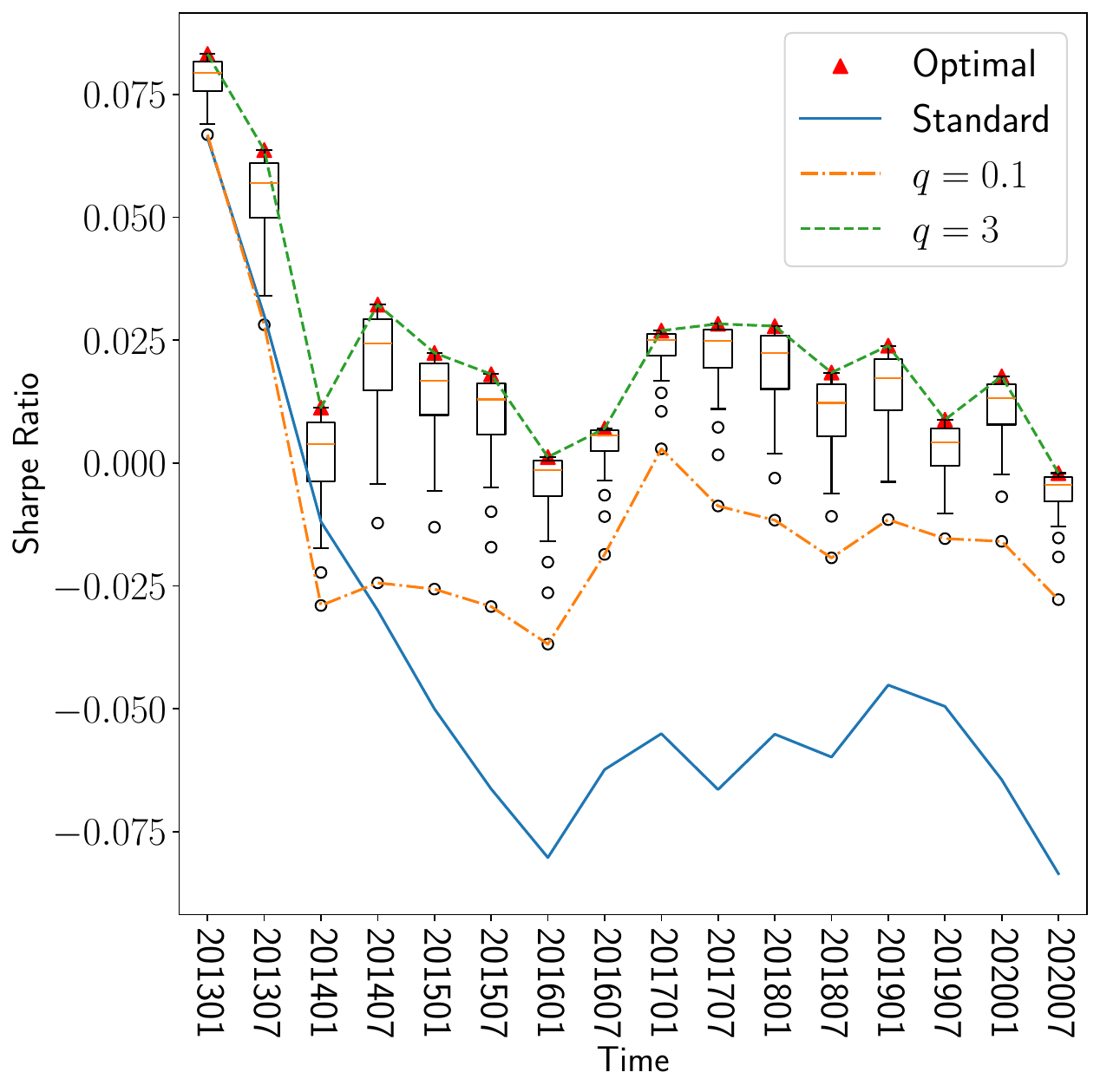}\label{fig_realmeanunknown:twoyear}}
    \subfloat[Four years, $\cQ = \cQ_1, c<1$]{\includegraphics[width=0.33\textwidth]{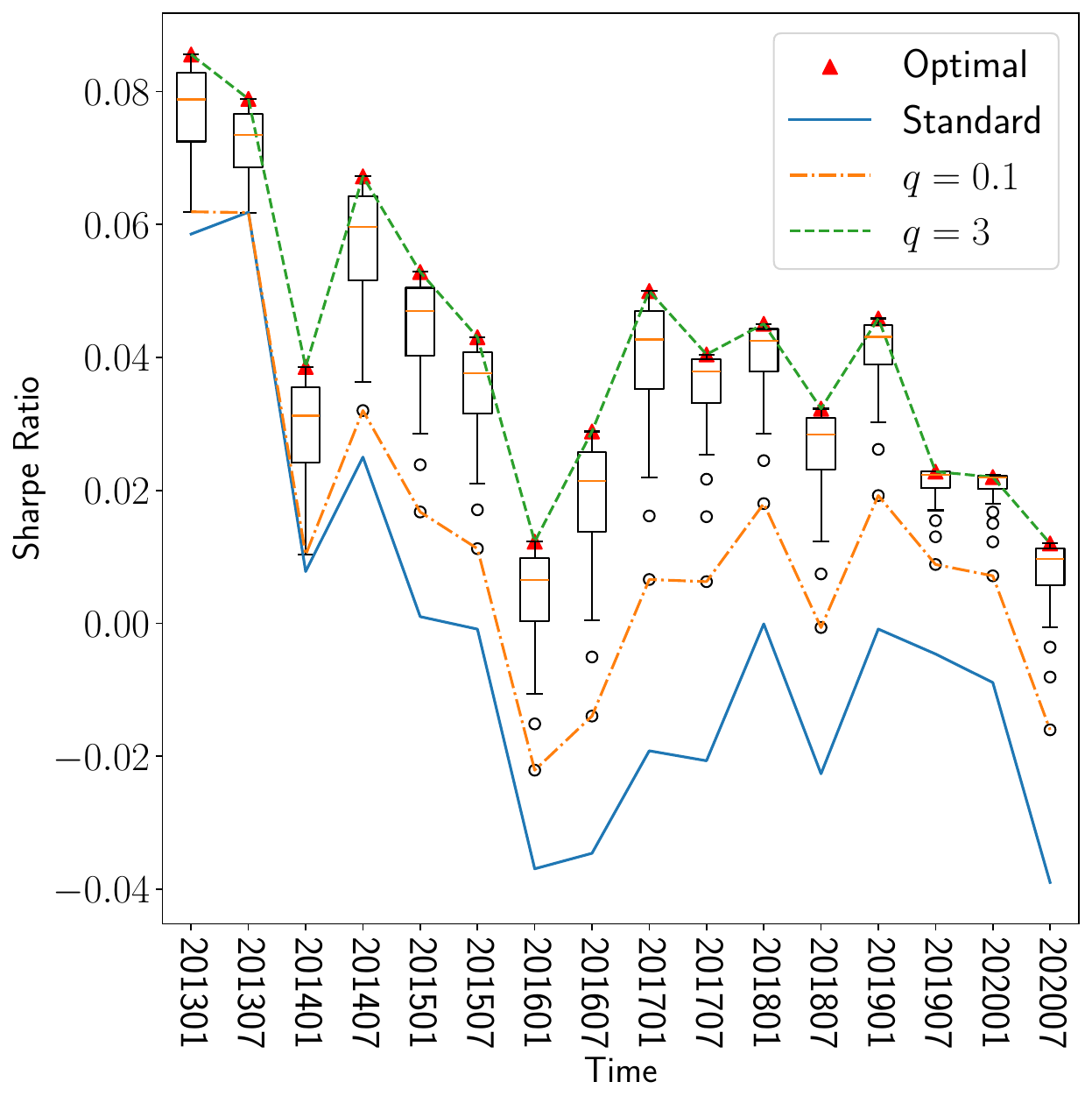}\label{fig_realmeanunknown:fouryear}}

    \subfloat[One year, $\cQ = \cQ_2, c>1$]{\includegraphics[width=0.33\textwidth]{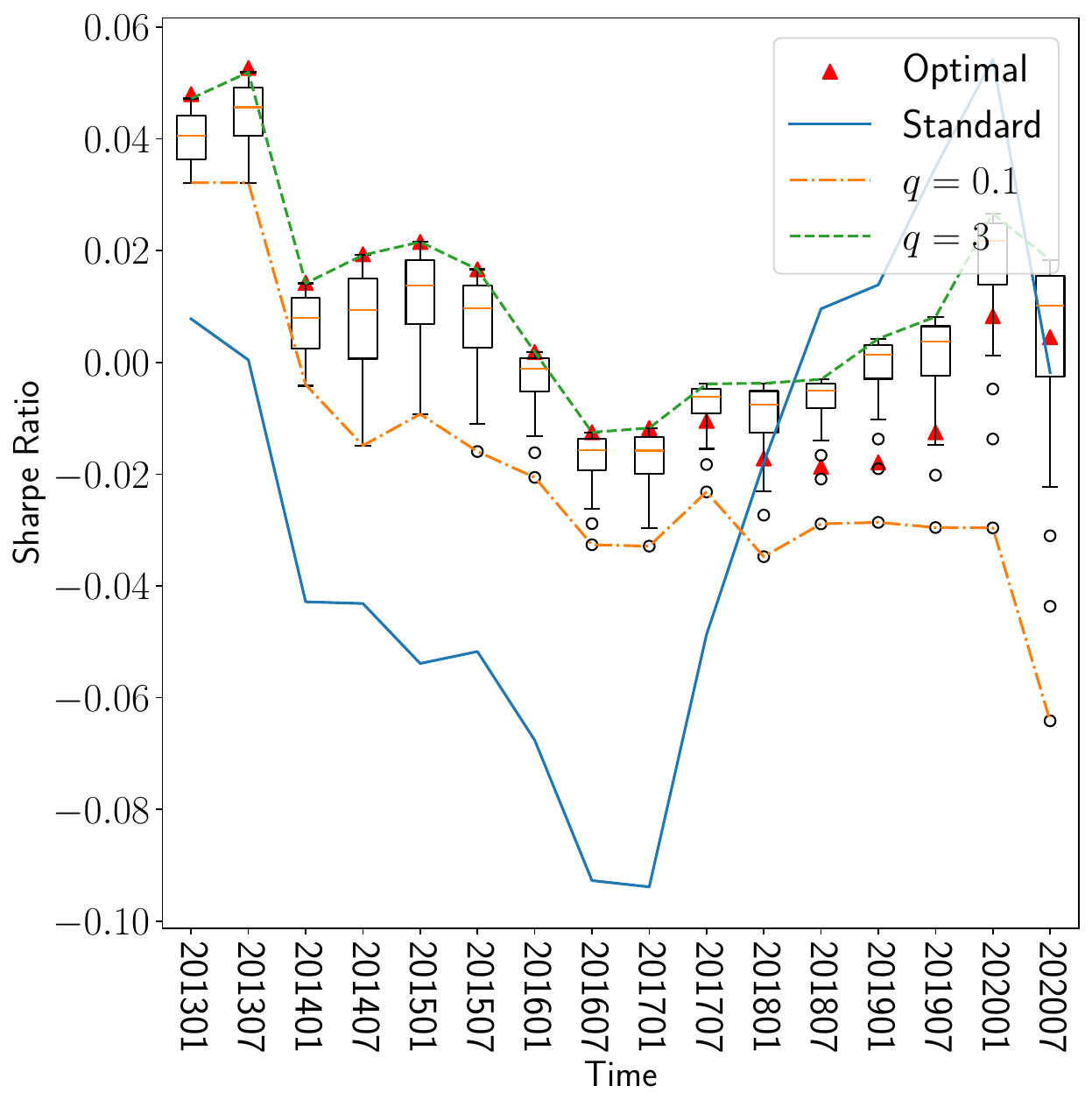}\label{fig_realmeanunknown:oneyear1}}
    \subfloat[Two years, $\cQ = \cQ_2, c<1$]{\includegraphics[width=0.33\textwidth]{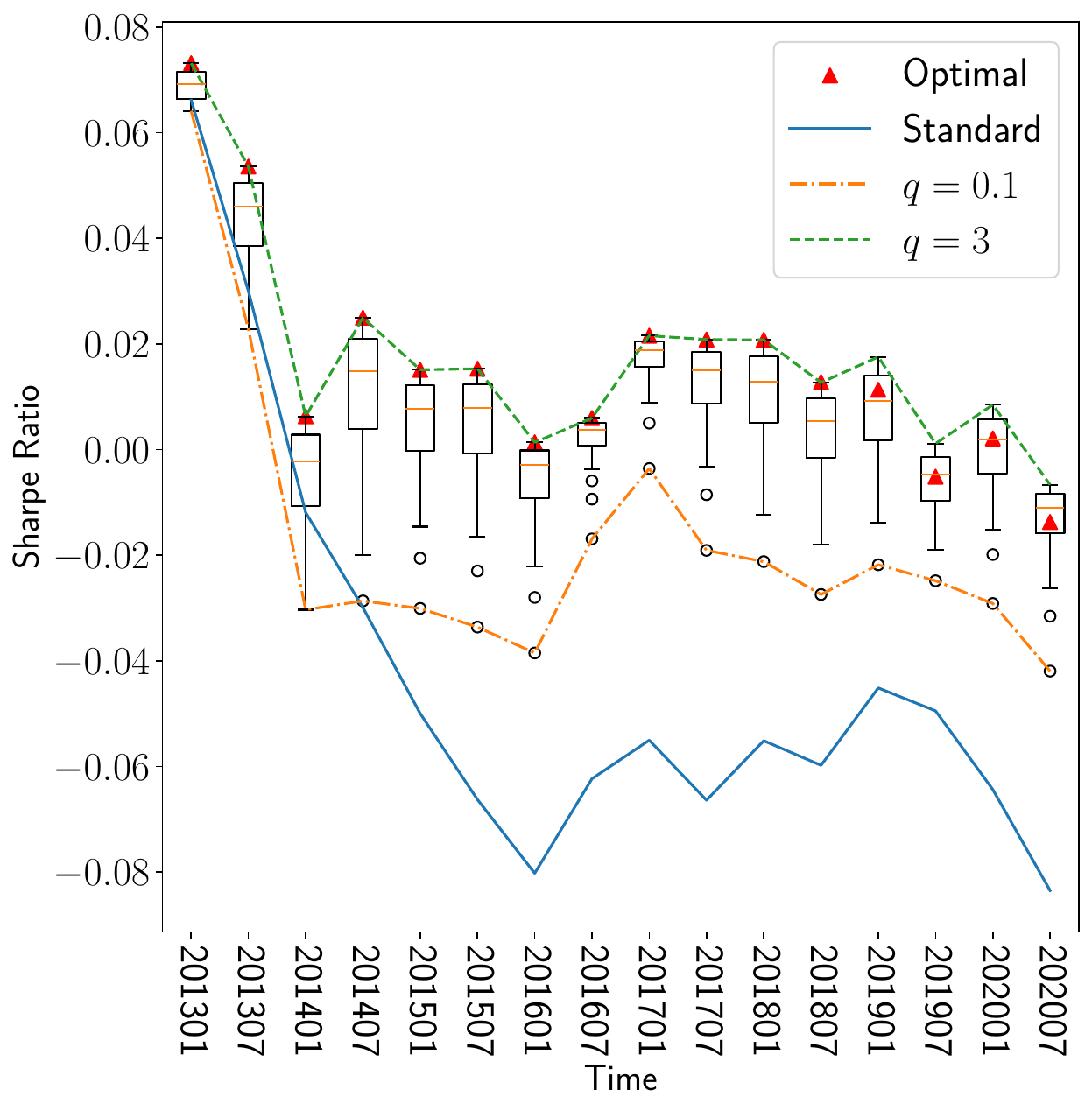}\label{fig_realmeanunknown:twoyear1}}
    \subfloat[Four years, $\cQ = \cQ_2, c<1$]{\includegraphics[width=0.33\textwidth]{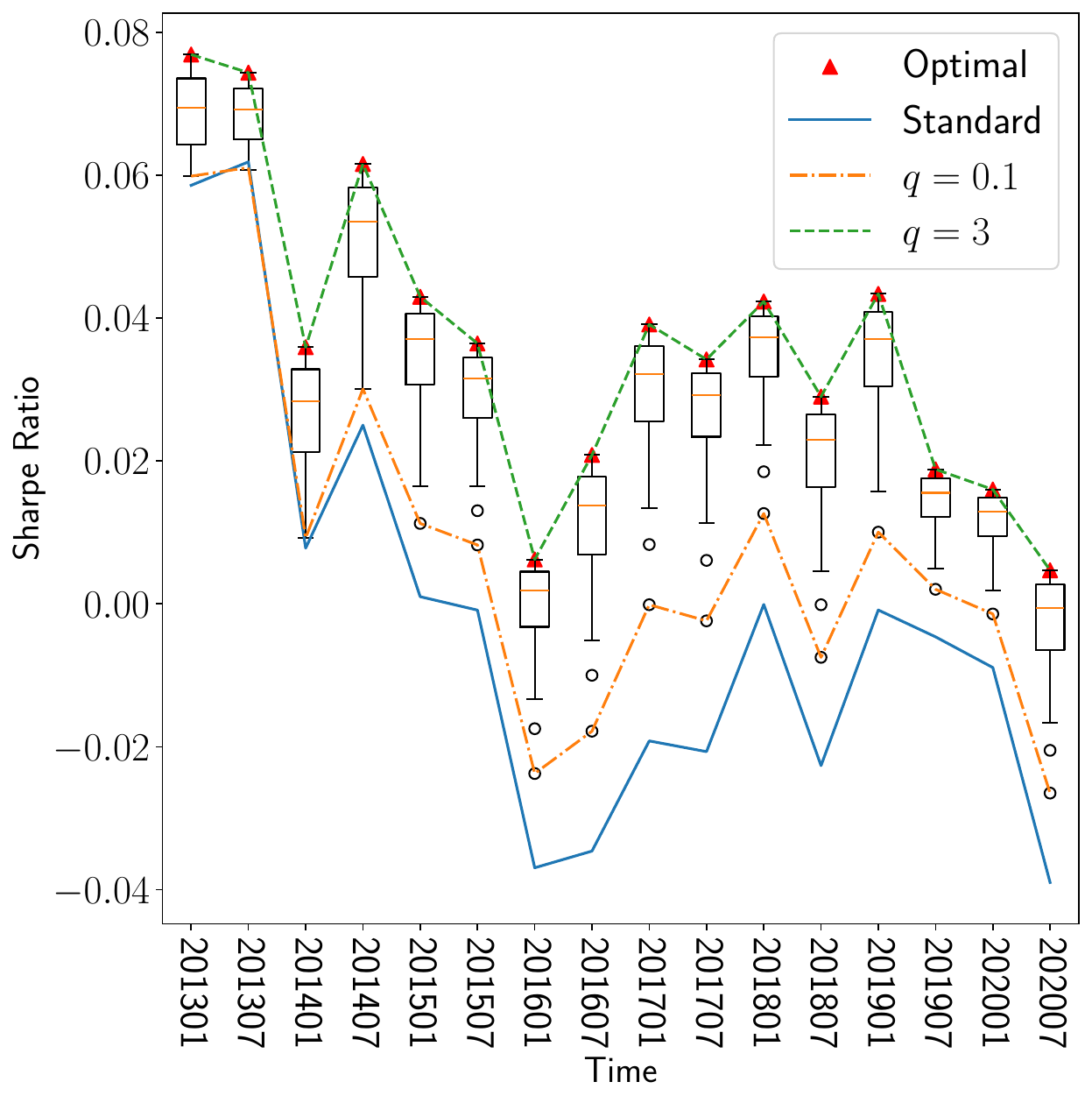}\label{fig_realmeanunknown:fouryear1}}
    
    \caption{ SR of mean-variance portfolios. The x-axis labels the rolling period, while the y-axis represents the out-of-sample SR of the portfolio returns every three years in the future. The blue solid line, orange dash-dot line and green dash line correspond to the SR with $q=0$, $q=0.1$ (minimum $q$) and $q=3$ (maximum $q$) respectively. The boxplot displays the Sharpe ratio for all $\Qb\in\cQ$, and the red triangle indicates the SR under our optimized $\Qb^*\in\cQ$.}
    \label{fig_realmeanunknown}
\end{figure}
Here we apply the sample mean and covariance matrix to construct the portfolio weight $\wb$. Different from Section~\ref{subsec:meanvariance}, we do not access to the mean vector $\bmu$ in the testing month. 
Figures \ref{fig_realmeanunknown:oneyear} to \ref{fig_realmeanunknown:fouryear} present the results for $\Qb \in \cQ_1$, while Figures \ref{fig_realmeanunknown:oneyear1} to \ref{fig_realmeanunknown:fouryear1} display the outcomes for $\Qb \in \cQ_2$. The blue curve represents the scenario where $\Qb = \zero$. The results clearly highlight the benefits of an actively optimized $\Qb^* \in \cQ$ for each testing month compared to a static, fixed $\Qb \in \cQ$ used uniformly across all months.
Moreover, using the sample mean $\hbmu$ leads to a significant drop in the Sharpe ratio compared to the scenario where $\bmu$ is known. This observation underscores the importance of accurate $\bmu$ estimation in practical settings. It suggests that relying the sample mean to estimate $\bmu$ may  be insufficient.

\section{Supporting Lemmas}
\label{sec:preliminary}
In this section, we will present several lemmas that play a crucial role in determining the asymptotic behavior of $T_{n,1}(\Qb)$ and $T_{n,2}(\Qb)$. These lemmas are instrumental in our analysis and provide valuable insights into the problem at hand.  By Assumption~\ref{ass:assump2}, readers may keep in mind that $\lambda_{\min}(\Qb)$ is bounded away from $0$, and we let $c_{\Qb}$ be a constant that is smaller than $\lambda_{\min}(\Qb)$. The special case $\Qb=\zero$ when $c<1$ is excluded in this section, the proof of Theorem~\ref{thm:main_theorem} when $\Qb=\zero$ is simple and we directly show it in Section~\ref{sec:appendix_proof_thm_bound}.

\begin{lemma}[Theorem 1 in  \cite{rubio2011spectral}]
\label{lemma:keylemma}
Given $\Xb\in\RR^{n\times p}$ with $\Xb=\Zb\bSigma^{\frac12}$ where the elements of $\Zb\in\RR^{n\times p}$ are i.i.d  with  zero mean, variance 1 and finite $8+\varepsilon$ order moment for some $\varepsilon>0$,  $\bSigma^{\frac12}\in\RR^{p\times p}\geq0$ is a  non negative Hermitian matrix with $\|\bSigma\|_{\op}$ bounded. For any deterministic $\Ab\in\RR^{p\times p}$ with bounded trace  norm ($\tr(\sqrt{\Ab\Ab^\top})\leq C$ with respect to $n$) and Hermitian non-negative matrix $\Sbb\in\RR^{p\times p}$, define 
\begin{align*}
    m_n(z)=\tr \bigg(\frac{\Xb^\T\Xb}{n}+\Sbb+z\Ib_p\bigg)^{-1}\Ab,
\end{align*}
and assume that $p,n$ tends to infinity with proportion $p/n\to c>0$. Then for fixed $z\in\bbC- (-\infty,0]$, 
\begin{align}
\label{eq:convergenceinlemma}
    m_n(z)-\tr\bigg( \frac{\bSigma}{1+s_0(z)}+\Sbb+z\Ib_p \bigg)^{-1}\Ab~\cvas~ 0.
\end{align}
Here, $s_0(z)$ is the Stieltjes transform of a positive finite measure and  uniquely solves the following equation:
\begin{align*}
    s_0(z)=\frac{c}{p}\tr\bSigma\bigg(\frac{\bSigma}{1+s_0(z)}+\Sbb+z\Ib_p\bigg)^{-1}.
\end{align*}
\end{lemma}
Note that Corollary 1.11 in \cite{widder1938stieltjes} gives that the convergence of such Stieltjes transform  in \eqref{eq:convergenceinlemma} is indeed internally closed uniform convergence, combined with the fact that the functions above are all analytic, we have the following corollary:
\begin{corollary}
\label{corol:keycorol}
In the condition of Lemma~\ref{lemma:keylemma}, it still holds the convergence  that
\begin{align}
\label{eq:convergenceincorol}
    \frac{\rmd m_n(z)}{\rmd z}-\frac{\rmd ~\tr\Big(\frac{\bSigma}{1+s_0(z)}+\Sbb+z\Ib_p \Big)^{-1}\Ab}{\rmd z}~\cvas~ 0.
\end{align}
\end{corollary} Lemma~\ref{lemma:keylemma} and Corollary~\ref{corol:keycorol} are limited to the constraint $\|\bSigma\|_{\op}$ is bounded. We first extend the results above to the case where $\|\bSigma\|_{\op}$ is unbounded. We give the following definition first. 
\begin{definition}
\label{def:s0-f(z)}
Define $s_0(z)$ by the equation
\begin{align*}
    s_0(z)=\frac{c}{p}\tr~\bSigma\bigg( \frac{\bSigma}{1+s_0(z)}+\Qb +z\bSigma\bigg)^{-1}.
\end{align*}
Then let \begin{align*}
    &f_n(z)=\tr\bigg(\frac{\Zb^\top\Zb}{n}+\bSigma^{-\frac{1}{2}}\Qb \bSigma^{-\frac{1}{2}}+z\Ib\bigg)^{-1}\bSigma^{-\frac{1}{2}}\Ab \bSigma^{-\frac{1}{2}},\\
    &\overline{f}_n(z)=\tr\bigg(\frac{\Ib}{1+s_0}+ \bSigma^{-\frac{1}{2}}\Qb \bSigma^{-\frac{1}{2}}+z\Ib \bigg)^{-1}\bSigma^{-\frac{1}{2}}\Ab \bSigma^{-\frac{1}{2}}.
\end{align*}
\end{definition}
The existence and uniqueness of $s_0(z)$ can be found in \citet{rubio2011spectral}. By Definition~\ref{def:s0-f(z)} above, we have the following lemma.
\begin{lemma}
\label{lemma:unbounded_op_znot=0}
Suppose that Assumption~\ref{ass:assump1}, \ref{ass:assump2} and \ref{ass:assump3}  hold, let $f_n(z)$ and $\overline{f}_n(z)$ be defined in Definition~\ref{def:s0-f(z)} with $\|\Ab\|_{\tr}\leq C<+\infty$, then for any fixed $z\in\bbC\backslash(-\infty,0]$, it holds that
\begin{align*}
    f_n(z)-\overline{f}_n(z)~\cvas~ 0, \quad f'_n(z)-\overline{f}'_n(z)~\cvas~ 0 .
\end{align*}
\end{lemma}
\begin{proof}[Proof of Lemma~\ref{lemma:unbounded_op_znot=0}]
Lemma~\ref{lemma:keylemma} and Corollary~\ref{corol:keycorol} directly supports Lemma~\ref{lemma:unbounded_op_znot=0}.
\end{proof}
From Lemma~\ref{lemma:unbounded_op_znot=0} above, we directly conclude that for fixed $z\in\bbC\backslash(-\infty,0]$, 
\begin{align*}
 \tr\bigg(\frac{\Xb^\top\Xb}{n}+\Qb +z\bSigma\bigg)^{-1}\Ab  -  \tr\bigg(\frac{\bSigma}{1+s_0(z)}+\Qb +z\bSigma\bigg)^{-1}\Ab ~\cvas~0,
\end{align*}
and the convergence also holds for the derivatives of $z$. However, this conclusion is not sufficient for us to prove Theorem~\ref{thm:main_theorem}, as it is crucial for us to consider the case when $z=0$. For the Case 1 in Assumption~\ref{ass:assump3}, by $\lambda_{\min}(\Qb)\geq c_{\Qb}$, we can easily let the matrix $\Sbb$ in Lemma~\ref{lemma:keylemma} by $\Sbb=\Qb-c_{\Qb}/2\Ib$, and $z=c_{\Qb}/2\Ib$, then we can extend the convergence to $z=0$ when the Case 1 in Assumption~\ref{ass:assump3} holds. 
In the following lemmas (Lemma~\ref{lemma:constant_level}-\ref{lemma:z=0converge_c>1}), we aim to extend the convergence to $z=0$ in the case 2 and 3 in Assumption~\ref{ass:assump3}. To prove so, we limit the range of $z$ into $[0,+\infty)$. We have the following lemma.
\begin{lemma}
\label{lemma:constant_level}
Suppose that Assumption~\ref{ass:assump1}, \ref{ass:assump2} and \ref{ass:assump3} hold, when $z\in[0,+\infty)$, $s_0(z)$ defined in Definition~\ref{def:s0-f(z)} exists and is unique. It holds that
$
    0<s_0(z)\leq {c\|\bSigma/p\|_{\tr}}/{\lambda_{\min}(\Qb)}.
$ 
For $z\in[0,+\infty)$, we further define $s_{1,\bSigma}(z)$ and $s_{1,\Qb}(z)$ by
\begin{align*}
    &s_{1,\bSigma}(z)=\bigg(\frac{s_{1,\bSigma}(z)}{(1+s_0(z))^2}-1\bigg)\frac{c}{p}\tr\bSigma\bigg( \frac{ \bSigma}{1+s_0(z)}+ \Qb +z\bSigma  \bigg)^{-1}
     \bSigma 
    \bigg( \frac{ \bSigma}{1+s_0(z)}+ \Qb +z\bSigma  \bigg)^{-1},\\
    &s_{1,\Qb}(z)=\frac{c}{p}\tr~\bSigma\bigg( \frac{\bSigma}{1+s_0(z)}+\Qb+z\bSigma \bigg)^{-1}\bigg(\frac{s_{1,\Qb}(z)}{(1+s_0(z))^2}\bSigma-\Qb\bigg)\bigg( \frac{\bSigma}{1+s_0(z)}+\Qb+z\bSigma \bigg)^{-1},
\end{align*}
then $s_{1,\bSigma}(z)$ and $s_{1,\Qb}(z)$ exist, and it holds that 
\begin{align*}
    -s_0(z)(1+s_0(z))^2\leq s_{1,\bSigma}(z)<0,\quad -s_0(z)\leq s_{1,\Qb}(z)<0.
\end{align*}
\end{lemma}
Lemma~\ref{lemma:constant_level} shows that the values of $s_0(z)$, $s_{1,\bSigma}(z)$ and $s_{1,\Qb}(z)$ range in a constant level. The proof can be found in Section~\ref{sec:prooflemmaconstantlevel}. With Lemma~\ref{lemma:constant_level}, we give the following lemmas which play a key role in the analysis under Case 2 and 3.

\begin{lemma}
 \label{lemma:z=0converge}
Suppose that Assumption~\ref{ass:assump1}, \ref{ass:assump2} and \ref{ass:assump3} Case 2 hold ($p/n\to c<1$), let $f_n(z)$ and $\overline{f}_n(z)$ be defined in Definition~\ref{def:s0-f(z)}, where $\Ab$ is a deterministic matrix with $\|\Ab\|_{\tr}\leq C<+\infty$. For $z\in[0,+\infty)$, it holds that 
\begin{align*}
    \lim\limits_{n\to +\infty}\lim\limits_{z\to 0^{+}} f_n(z)=\lim\limits_{z\to 0^{+}} \lim\limits_{n\to +\infty} f_n(z), \quad 
    \lim\limits_{n\to +\infty}\lim\limits_{z\to 0^{+}} \overline{f}_n(z)=\lim\limits_{z\to 0^{+}} \lim\limits_{n\to +\infty} \overline{f}_n(z).
\end{align*}
Therefore, for any fixed $z\in\bbC\backslash(-\infty,0)$, it holds that
\begin{align*}
    f_n(z)-\overline{f}_n(z)~\cvas~ 0, \quad f'_n(z)-\overline{f}'_n(z)~\cvas~ 0 .
\end{align*}
\end{lemma}
Lemma~\ref{lemma:z=0converge} established the key asymptotic behavior in Case 2 in Assumption~\ref{ass:assump3}. 
However, a notable distinction in Case 3 in Assumption~\ref{ass:assump3} is the presence of $p/n\to c\geq1$. Consequently, the proof technique utilized in Lemma~\ref{lemma:z=0converge} is no longer applicable due to the inability to uniformly bound $|f'_n(z)|$ by a trivial constant. Nevertheless, we can still address the case $z=0$ if we make the additional assumption that the number of eigenvalues $K$ tending to infinity remains fixed and $\lambda_1/\lambda_K^2\leq C$:

\begin{lemma}
 \label{lemma:z=0converge_c>1}
Suppose that Assumption~\ref{ass:assump1}, \ref{ass:assump2} and \ref{ass:assump3} Case 3 hold ($p/n\to c\geq1$), let $f_n(z)$ and $\overline{f}_n(z)$ be defined in Definition~\ref{def:s0-f(z)}, where $\Ab$ is a deterministic matrix with $\|\Ab\|_{\tr}\leq C<+\infty$. For $z\in[0,+\infty)$, it holds that 
\begin{align*}
    \lim\limits_{n\to +\infty}\lim\limits_{z\to 0^{+}} f_n(z)=\lim\limits_{z\to 0^{+}} \lim\limits_{n\to +\infty} f_n(z), \quad 
    \lim\limits_{n\to +\infty}\lim\limits_{z\to 0^{+}} \overline{f}_n(z)=\lim\limits_{z\to 0^{+}} \lim\limits_{n\to +\infty} \overline{f}_n(z).
\end{align*}
Therefore, for any fixed $z\in\bbC\backslash(-\infty,0)$, it holds that
\begin{align*}
    f_n(z)-\overline{f}_n(z)~\cvas~ 0, \quad f'_n(z)-\overline{f}'_n(z)~\cvas~ 0 .
\end{align*}
\end{lemma}
Combined with Corollary~\ref{corol:keycorol}, Lemma~\ref{lemma:z=0converge} and \ref{lemma:z=0converge_c>1}, we have the following Proposition.
\begin{proposition}
\label{prop:all_case_converge}
Suppose that Assumption~\ref{ass:assump1}, \ref{ass:assump2} and \ref{ass:assump3} hold, let $f_n(z)$ and $\overline{f}_n(z)$ be defined in Definition~\ref{def:s0-f(z)}, where $\Ab$ is a deterministic matrix with $\|\Ab\|_{\tr}\leq C<+\infty$. For any fixed $z\in\bbC\backslash(-\infty,0)$, it holds that 
\begin{align*}
    f_n(z)-\overline{f}_n(z)~\cvas~ 0, \quad f'_n(z)-\overline{f}'_n(z)~\cvas~ 0 .
\end{align*}
\end{proposition}
\begin{proof}[Proof of Proposition~\ref{prop:all_case_converge}]
For Case 1 in Assumption~\ref{ass:assump3}, we have $\lambda_{\min}(\Qb)\geq c_{\Qb}>0$ in Lemma~\ref{lemma:keylemma} and Corollary~\ref{corol:keycorol} above.  We may separate $\Qb+z\Ib$ into $\Qb-c_{\Qb}/2\Ib+(z+c_{\Qb}/2)\Ib $, therefore the equation still holds when $z=0$.
For Case 2 and Case 3 in Assumption~\ref{ass:assump3}, the results directly hold from Lemma~\ref{lemma:z=0converge} and \ref{lemma:z=0converge_c>1}.
\end{proof}
The difference between Proposition~\ref{prop:all_case_converge} and Lemma~\ref{lemma:unbounded_op_znot=0} is the range of $z$. Proposition~\ref{prop:all_case_converge} covers the case of the convergence when $z=0$. Note that when diverging spikes exist in the matrix $\bSigma$, $|f_n(z)|$ and $|\overline{f}_n(z)|$ may tends to $0$ as $n$ tends to infinity,  we also give the ratio convergence on $f_n(z)$ and $\overline{f}_n(z)$. We have the following lemma. 
\begin{lemma}
\label{lemma:stronger_type_convergence}
Suppose that Assumption~\ref{ass:assump1}, \ref{ass:assump2} and \ref{ass:assump3} hold. Recall the definitions of $f_n(z)$ and $\overline{f}_n(z)$ from Definition~\ref{def:s0-f(z)} where we choose $\Ab$ to be a semi-positive definite matrix. Then it holds that 
\begin{align*}
    f_n(z)/\overline{f}_n(z)~\cvas~1,\quad f'_n(z)/\overline{f}'_n(z) ~\cvas~1.
\end{align*}
\end{lemma}
The proof of Lemma~\ref{lemma:stronger_type_convergence} is given in Section~\ref{sec:prooflemmastronger_type_convergence}. 
Under these lemmas, we are able to prove Theorem~\ref{thm:main_theorem}. We give the proof of Theorem~\ref{thm:main_theorem} in Section~\ref{sec:appendix_proof_thm_bound}.

\subsection{Proof of Lemma~\ref{lemma:constant_level}}
\label{sec:prooflemmaconstantlevel}
We write $s_0(z)$, $s_{1,\bSigma}(z)$ and $s_{1,\Qb}(z)$ as $s_0$, $s_{1,\bSigma}$ and $s_{1,\Qb}$ for simplicity.
We first prove the existence and uniqueness of $s_0$ under the certain condition. Readers may refer to the section 4.1 and 4.2 in \cite{rubio2011spectral} for the existence and uniqueness of $s_0$ for $z\in\bbC_+$. With similar proof, the existence and uniqueness of $s_0$ when $z\in\bbC_-$ still hold. We prove that when $z\in[0,+\infty)$, $s_0>0$ has the unique solution. To prove so, denote $\ub=\{u_i>0\}$ by the eigenvalue of $\bSigma^{-\frac12}\Qb\bSigma^{-\frac12} +z\Ib$, and $g(s;\ub)$ by $\frac{c}{p}\sum_{i=1}^p \frac{1}{\frac{1}{1+s}+u_i}$, hence $s_0$ satisfies
\begin{align*}
    s_0=g(s_0;\ub)&=\frac{c}{p}\sum_{i=1}^p \frac{1}{\frac{1}{1+s_0}+u_i}=\frac{c}{p}\sum_{i=1}^p\frac{1}{u_i}\cdot \frac{(1+s_0)u_i}{1+u_i(1+s_0)}\\
    &=\frac{c}{p}\sum_{i=1}^p \frac{1}{u_i}\cdot\bigg(1-\frac{1}{1+u_i(1+s_0)} \bigg)>0,
\end{align*}
when $s=0$ we have $g(s;\ub)-s>0$, when $s\to+\infty$, we have $g(s;\ub)-s<0$, the existence of $s_0$ then comes from Bolzano-Cauchy theorem. As for the uniqueness, assume $g(s;\ub)-s$ exists two roots $s_1$ , $s_2$, and $s_1<s_2$. Note that $g(s;\ub)-s$ is a concave function, we conclude that $g(s;\ub)-s$ increases then decreases or just decreases. Note that $g(0;\ub)>0$ and $g(s;\ub)-s$ exists roots, we conclude that $g(s;\ub)-s$ decreases  in the neighbourhood of $s_1$ and decreases in $[s_1,+\infty)$. Hence, in the interval $[s_1,s_2]$, $g(s;\ub)-s\equiv 0$, which violates to the condition $g(s;\ub)-s$ is concave. Therefore we prove the uniqueness of $s_0$. 
For the bound of $s_0$, it holds that
\begin{align*}
    s_0= \frac{c}{p}\tr \bSigma \bigg(\frac{\bSigma}{1+s_0}+\Qb+z\bSigma \bigg)^{-1}\leq c\|\bSigma/p \|_{\tr}/\lambda_{\min}(\Qb).
\end{align*}
Here, the  inequality comes from $\tr(\Ab\Bb)\leq \|\Ab\|_{\tr}\|\Bb\|_{\op}$ and $\lambda_{\min}(\Ab+\Bb)\geq \lambda_{\min}(\Ab)+\lambda_{\min}(\Bb)$ for non negative matrix $\Ab$ and $\Bb$.

We now give the analysis of $s_{1,\bSigma}$. Denote $u_i>0$ by the eigenvalues of $\bSigma^{-\frac{1}{2}}\Qb \bSigma^{-\frac{1}{2}}+z\Ib$, then we have
\begin{align*}
    s_0&=\frac{c}{p}\sum_{i=1}^p\frac{(1+s_0)}{(1+s_0)\cdot u_i+1},\\
    s_{1,\bSigma}&=\left(\frac{s_{1,\bSigma}}{\left(1+s_0\right)^2}-1\right) \cdot \frac{c}{p} \cdot \sum_{i=1}^p \frac{\left(1+s_0\right)^2}{\left(1+\left(1+s_0\right)\cdot u_i\right)^2}.
\end{align*}
Define $m=\frac{c}{p}\tr\big(\frac{\Ib}{1+s_0}+ \bSigma^{-\frac{1}{2}}\Qb \bSigma^{-\frac{1}{2}}+z\Ib\big)^{-2}/(1+s_0)^2>0$, then we have 
\begin{align*}
    m&=\frac{c}{p} \cdot \sum_{i=1}^p \frac{1}{(1+\left(1+s_0)\cdot u_i\right)^2}\leq \frac{c}{p(1+s_0)} \cdot \sum_{i=1}^p \frac{\left(1+s_0\right)}{1+\left(1+s_0\right)\cdot u_i} =  \frac{s_0}{1+s_0},\\
    s_{1,\bSigma}&=\big(s_{1,\bSigma}-(1+s_0)^2\big)\cdot m=\frac{m(1+s_0)^2}{m-1},
\end{align*}
where the first equation is by the definition of $m$ and $s_0$, and the second equation is by \eqref{eq:defofs1z_c>1}. The existence of $s_{1,\bSigma}$ comes from the second equality. Moreover from $0<m\leq s_0/(1+s_0)$, we conclude that $s_{1,\bSigma}$ exists, $s_{1,\bSigma}<0$ and $-s_{1,\bSigma}\leq s_0\cdot (1+s_0)^2$. This completes bound for $s_{1,\bSigma}$.

We next bound for $s_{1,\Qb}$. Denote by $u_i$ ($u_i>0$) the eigenvalue of the matrix $\bSigma^{-\frac{1}{2}}\Qb\bSigma^{-\frac{1}{2}}+z\Ib$, by the definition of $s_0$, we can easily see that
\begin{align*}
    s_0=\frac{c}{p}\sum_{i=1}^p\frac{(1+s_0)}{(1+s_0)\cdot u_i+1}.
\end{align*}
By the definition of $s_{1,\Qb}$, we can further express $s_{1,\Qb}$ as 
\begin{align}
s_{1,\Qb}&=\frac{c}{p}\tr~\bSigma\bigg( \frac{\bSigma}{1+s_0}+\Qb \bigg)^{-1}\bigg(\frac{s_{1,\Qb}}{(1+s_0)^2}\bSigma-\Qb\bigg)\bigg( \frac{\bSigma}{1+s_0}+\Qb \bigg)^{-1}\nonumber\\
&=s_{1,\Qb}\cdot \frac{c}{p}\sum_{i=1}^p \frac{1}{(1+(1+s_0)u_i)^2}-\frac{c}{p}\tr~\bSigma\bigg( \frac{\bSigma}{1+s_0}+\Qb \bigg)^{-1}\Qb\bigg( \frac{\bSigma}{1+s_0}+\Qb \bigg)^{-1}\nonumber\\
&=s_{1,\Qb}\cdot \frac{c}{p}\sum_{i=1}^p \frac{1}{(1+(1+s_0)u_i)^2}-s_0+\frac{c}{p}\tr~\bSigma\bigg( \frac{\bSigma}{1+s_0}+\Qb \bigg)^{-1}\frac{\bSigma}{1+s_0}\bigg( \frac{\bSigma}{1+s_0}+\Qb \bigg)^{-1}\nonumber\\
&=s_{1,\Qb}\cdot \frac{c}{p}\sum_{i=1}^p \frac{1}{(1+(1+s_0)u_i)^2}-s_0+\frac{c}{p}\sum_{i=1}^p \frac{1+s_0}{(1+(1+s_0)u_i)^2}.\label{eq:s1Q_s0_express}
\end{align}
Here, the second equality comes from the decomposition of the matrix $\Ib/(1+s_0)+\bSigma^{-\frac{1}{2}}\Qb\bSigma^{-\frac{1}{2}}$, the third equality comes from $\Qb=\Qb+\bSigma/(1+s_0)-\bSigma/(1+s_0)$ and the definition of $s_0$, and the last equality still comes from the decomposition of the matrix $\Ib/(1+s_0)+\bSigma^{-\frac{1}{2}}\Qb\bSigma^{-\frac{1}{2}}$. 

Define  $m=\frac{c}{p}\tr\big(\frac{\Ib}{1+s_0}+ \bSigma^{-\frac{1}{2}}\Qb \bSigma^{-\frac{1}{2}}\big)^{-2}/(1+s_0)^2$, then we have
\begin{align*}
    m&=\frac{c}{p} \cdot \sum_{i=1}^p \frac{1}{(1+(1+s_0)\cdot u_i)^2}\leq \frac{c}{p(1+s_0)} \cdot \sum_{i=1}^p \frac{1+s_0}{1+(1+s_0)\cdot u_i} = \frac{s_0}{1+s_0}.
\end{align*}
By the inequality above, the coefficient of $s_{1,\Qb}$ in the right side in \eqref{eq:s1Q_s0_express} is defined as $m$ and is smaller than $1$, we conclude the existence of $s_{1,\Qb}$. We first prove that $s_{1,\Qb}<0$. By \eqref{eq:s1Q_s0_express}, if $s_{1,\Qb}>0$, then we have
\begin{align*}
    s_{1,\Qb}&=m\cdot s_{1,\Qb}-s_0+\frac{c}{p}\sum_{i=1}^p \frac{1+s_0}{(1+(1+s_0)u_i)^2}\\
    &\leq s_{1,\Qb}\cdot \frac{s_0}{1+s_0}-s_0+\frac{c}{p}\sum_{i=1}^p \frac{1+s_0}{1+(1+s_0)u_i}\\
    &=s_{1,\Qb}\cdot\frac{s_0}{1+s_0}<s_{1,\Qb},
\end{align*}
where the first inequality is by $m\leq s_0/(1+s_0)$ and $1+(1+s_0)u_i\geq 1$, the second equality is by the definition of $s_0$ and the last inequality is by $s_0/(1+s_0)<1$. We have a contradictory, therefore we conclude that $s_{1,\Qb}<0$. To prove $s_{1,\Qb}\geq -s_0$, from \eqref{eq:s1Q_s0_express} we have
\begin{align*}
    1+s_0+s_{1,\Qb}=1+\frac{c}{p}\sum_{i=1}^p \frac{1+s_0+s_{1,\Qb}}{(1+(1+s_0)u_i)^2}=1+m\cdot(1+s_0+s_{1,\Qb}).
\end{align*}
From $0<m<1$ we conclude that $1+s_0+s_{1,\Qb}>0$. Therefore we have that
\begin{align*}
    1+s_0+s_{1,\Qb}=1+m\cdot(1+s_0+s_{1,\Qb})\geq 1, 
\end{align*}
from which we conclude that $s_{1,\Qb}\geq -s_0$.

\subsection{Proof of Lemma~\ref{lemma:z=0converge}}
\label{sec:prooflemmaz=0converge}
We need to prove that  the limit $n\to+\infty$ and $z\to0^+$ can be exchanged in $f_n(z)$. We note that
\begin{align*}
    |f_n(z)|&=\bigg|\tr\bigg(\frac{\bSigma^{\frac{1}{2}}\Zb^\top\Zb\bSigma^{\frac{1}{2}}}{n}+\Qb +z\bSigma\bigg)^{-1}\Ab\bigg| \leq \frac{\big\| \Ab\big\|_{\tr}}{\lambda_{\min}(\Qb)+z\lambda_{\min}(\bSigma)}\leq \frac{\big\| \Ab\big\|_{\tr}}{\lambda_{\min}(\Qb)},
\end{align*}
where  $\lambda_{\min }(\cdot)$ denote the  smallest  eigenvalues.  We have $|f_n(z)|$ is uniformly bounded by $\|\Ab\|_{\tr}/\lambda_{\min}(\Qb)$ when $z\in [0,+\infty)$. Hence we easily get that  $|f_n(z)-f_n(0)|\leq 2\|\Ab\|_{\tr}/\lambda_{\min}(\Qb)$. Moreover, when $p/n\to c<1$, we have
\begin{align*}
    |f_n'(z)|= 
    \tr\bigg[\bigg(\frac{\Zb^\top\Zb}{n}+\bSigma^{-\frac{1}{2}}\Qb \bSigma^{-\frac{1}{2}}+z\Ib \bigg)^{-2}  \bigg] \bSigma^{-\frac{1}{2}}\Ab\bSigma^{-\frac{1}{2}} \leq \frac{\|\bSigma^{-\frac{1}{2}}\Ab\bSigma^{-\frac{1}{2}}  \|_{\tr}}{\lambda_{\min}^2(\Zb^\top\Zb/n)}\leq \frac{4(\| \bSigma^{-1}\|_{\op}\|\Ab\|_{\tr})}{(1-\sqrt{c})^4}. 
\end{align*}
The second inequality holds almost surely for large enough $n$, by the Bai-Yin theorem \citep{bai1993limit} that $\lambda_{\min}(\Zb^\top\Zb/n)\geq  (1-\sqrt{c})^2/2$ almost surely for all $n$ large enough. As its derivatives are bounded, the sequence $f_n, n=1,2,3, \ldots$ is equicontinuous, and by the Arzela-Ascoli theorem, we deduce that $f_n$ converges uniformly to its limit. By the Moore-Osgood theorem, we can exchange limits (as $n, p \rightarrow \infty$ and $z \rightarrow 0^{+}$).

We next prove that the limit $n\to+\infty$ and $z\to0^+$ can be exchanged in $\overline{f}_n(z)$. Note that 
\begin{align*}
    \overline{f}_n(z)&=\tr\bigg(\frac{\Ib}{1+s_0}+ \bSigma^{-\frac{1}{2}}\Qb \bSigma^{-\frac{1}{2}}+z\Ib \bigg)^{-1}\bSigma^{-\frac{1}{2}}\Ab \bSigma^{-\frac{1}{2}}\\
    &= \tr \bigg(\frac{\bSigma}{1+s_0}+ \Qb +z\bSigma \bigg)^{-1}\Ab \leq \|\Ab\|_{\tr}/\lambda_{\min}(\Qb), 
\end{align*}
we have $\overline{f}_n(z)$ is uniformly bounded by $\|\Ab\|_{\tr}/\lambda_{\min}(\Qb)$ when $z\in [0,+\infty)$.  Recall that 
\begin{align*}
    s_0(z)=\frac{c}{p}\tr \bigg(\frac{\Ib}{1+s_0(z)}+ \bSigma^{-\frac{1}{2}}\Qb \bSigma^{-\frac{1}{2}}+z\Ib \bigg)^{-1},
\end{align*}
when $z\geq 0$, Lemma~\ref{lemma:constant_level} gives that 
\begin{align*}
    0<s_0(z)\leq  c\|\bSigma/p \|_{\tr}/\lambda_{\min}(\Qb).
\end{align*}
Combined the truth $s_0(z)$ is analytic from \citet{rubio2011spectral}, we have 
\begin{align*}
    \overline{f}_n'(z)&=\tr \bigg( \frac{\bSigma}{1+s_0(z)}+\Qb+z\bSigma \bigg)^{-1}\bigg(\frac{s_1(z)}{(1+s_0(z))^2}\bSigma-\bSigma\bigg)\bigg( \frac{\bSigma}{1+s_0(z)}+\Qb+z\bSigma \bigg)^{-1}\Ab\\
    &=\bigg(\frac{s_1(z)}{(1+s_0(z))^2}-1\bigg)\tr \bigg(\frac{\Ib}{1+s_0(z)}+ \bSigma^{-\frac{1}{2}}\Qb \bSigma^{-\frac{1}{2}}+z\Ib\bigg)^{-2}\bSigma^{-\frac{1}{2}}\Ab \bSigma^{-\frac{1}{2}},
\end{align*}
where $s_1(z)=s_0'(z)$ satisfies
\begin{align}
    s_1(z)&=\frac{c}{p}\tr\bSigma\bigg( \frac{ \bSigma}{1+s_0(z)}+ \Qb +z\bSigma  \bigg)^{-1}
    \bigg(\frac{s_1(z)}{(1+s_0(z))^2} \bSigma- \bSigma \bigg)
    \bigg( \frac{ \bSigma}{1+s_0(z)}+ \Qb +z\bSigma  \bigg)^{-1}\label{eq:defofs1z}.
\end{align}
Lemma~\ref{lemma:constant_level} further gives that $s_1(z)$ is uniformly bounded, we can see 
\begin{align*}
    |\overline{f}'_n(z)|\leq \bigg|\frac{s_1(z)}{(1+s_0(z))^2}-1 \bigg|\cdot c\cdot (1+s_0(z))^2
\end{align*}
is uniformly bounded. Apply the similar argument above which proved the exchange in $f_n(z)$, we obtain that the limit $n\to+\infty$ and $z\to0^+$ can be exchanged in $\overline{f}_n(z)$. We complete the proof that $n\to+\infty$ and $z\to0^+$ can be exchanged in $f_n(z)$ and $\overline{f}_n(z)$. 

Therefore, we have
\begin{align*}
0&=\lim_{z\to0^+}\lim_{n\to+\infty}\overline{f}_n(z)-f_n(z)=\lim_{n\to+\infty}\lim_{z\to0^+}\overline{f}_n(z)-f_n(z)\\
&=\lim_{n\to+\infty}\overline{f}_n(0)-f_n(0).
\end{align*}
Here, the first equality comes from Lemma~\ref{lemma:unbounded_op_znot=0}.  The internally closed uniform convergence of $\overline{f}_n(z)-f_n(z)$ is hence extended to the area $\bbC\backslash(-\infty,0)$, which also indicates that for fixed $z\in\bbC\backslash(-\infty,0)$, $\overline{f}'_n(z)-f'_n(z)~\cvas~0$.

\subsection{Proof of Lemma~\ref{lemma:z=0converge_c>1}}
When $z\in[0,+\infty)$, from Lemma~\ref{lemma:constant_level} we can obtain the existence and uniqueness of $s_0>0$. We also have $|f_n(z)|$ and $|\overline{f}_n(z)|$ are uniformly bounded by $\|\Ab\|_{\tr}/\lambda_{\min}(\Qb)$, and $|s_0(z)|$ is bounded by $c\|\bSigma/p \|_{\tr}/\lambda_{\min}(\Qb)$. We still define $s_1(z)=s_0'(z)$ which satisfies
\begin{align}
    s_1(z)&=\frac{c}{p}\tr\bSigma\bigg( \frac{ \bSigma}{1+s_0(z)}+ \Qb +z\bSigma  \bigg)^{-1}
    \bigg(\frac{s_1(z)}{(1+s_0(z))^2} \bSigma- \bSigma \bigg)
    \bigg( \frac{ \bSigma}{1+s_0(z)}+ \Qb +z\bSigma  \bigg)^{-1}\nonumber\\
    &=\bigg(\frac{s_1(z)}{(1+s_0(z))^2}-1\bigg)\cdot \frac{c}{p}\tr\bigg(\frac{\Ib}{1+s_0(z)}+ \bSigma^{-\frac{1}{2}}\Qb \bSigma^{-\frac{1}{2}}+z\Ib\bigg)^{-2},\label{eq:defofs1z_c>1}
\end{align}
then 
\begin{align*}
    \overline{f}_n'(z)&=\tr \bigg( \frac{\bSigma}{1+s_0(z)}+\Qb+z\bSigma \bigg)^{-1}\bigg(\frac{s_1(z)}{(1+s_0(z))^2}\bSigma-\bSigma\bigg)\bigg( \frac{\bSigma}{1+s_0(z)}+\Qb+z\bSigma \bigg)^{-1}\Ab\\
    &=\bigg(\frac{s_1(z)}{(1+s_0(z))^2}-1\bigg)\tr \bigg(\frac{\Ib}{1+s_0(z)}+ \bSigma^{-\frac{1}{2}}\Qb \bSigma^{-\frac{1}{2}}+z\Ib\bigg)^{-2}\bSigma^{-\frac{1}{2}}\Ab \bSigma^{-\frac{1}{2}}.
\end{align*}
To complete the proof, we only need to show $z\to 0^+$ and $n\to+\infty$ can be exchanged in $f_n(z)$ and $\overline{f}_n(z)$. Hence it remains to prove that $|f_n(z)|$ and $|f'_n(z)|$ are uniformly bounded. 

We first prove that $|\overline{f}'_n(z)|$ is uniformly bounded. To prove so, we have 
\begin{align*}
    |\overline{f}'_n(z)|\leq |s_1(z)-(1+s_0(z))^2|\cdot \| \bSigma^{-\frac{1}{2}}\Ab \bSigma^{-\frac{1}{2}}\|_{\tr}.
\end{align*}
It is equal to prove that $|s_1(z)|$ is uniformly bounded.   By Lemma~\ref{lemma:constant_level} we complete the proof that $|\overline{f}'_n(z)|$ is uniformly bounded.

We next prove that $|f'_n(z)|$ is uniformly bounded almost surely for large enough $n$.  Recall that
\begin{align*}
    f_n'(z)&= 
    \tr\bigg[\bigg(\frac{\bSigma^{\frac{1}{2}}\Zb^\top\Zb\bSigma^{\frac{1}{2}}}{n}+\Qb +z\bSigma \bigg)^{-1} \bSigma \bigg(\frac{\bSigma^{\frac{1}{2}}\Zb^\top\Zb\bSigma^{\frac{1}{2}}}{n}+\Qb +z\bSigma \bigg)^{-1} \bigg] \Ab.
\end{align*}
 By Assumption~\ref{ass:assump2}, when the constant $c_2\not=0$, it is easy to bound $f'_n(z)$ by $\frac{\|\bSigma^{-1/2}\Ab\bSigma^{-1/2}\|_{\tr}}{c'^2c_2^2}$ with some constant $c',c_2>0$. We then consider the case $c_2=0$ and $c'\Ib\leq\Qb\leq C'\Ib$.
If we can prove that for any fixed vector $\ub$ with $\|\ub\|=1$, 
\begin{align}
\label{eq:uniformbound_vector}
    \ub^\top\bigg(\frac{\bSigma^{\frac{1}{2}}\Zb^\top\Zb\bSigma^{\frac{1}{2}}}{n}+\Qb +z\bSigma \bigg)^{-1} \bSigma \bigg(\frac{\bSigma^{\frac{1}{2}}\Zb^\top\Zb\bSigma^{\frac{1}{2}}}{n}+\Qb +z\bSigma \bigg)^{-1} \ub 
\end{align}
is uniformly bounded almost surely by constant $C$ for large enough $n$,  we will get that 
\begin{align*}
|f_n'(z)|&=\sum_{i=1}^{p}\lambda_{\Ab,i}\bxi_{\Ab,i}^\top\bigg(\frac{\bSigma^{\frac{1}{2}}\Zb^\top\Zb\bSigma^{\frac{1}{2}}}{n}+\Qb +z\bSigma \bigg)^{-1} \bSigma \bigg(\frac{\bSigma^{\frac{1}{2}}\Zb^\top\Zb\bSigma^{\frac{1}{2}}}{n}+\Qb +z\bSigma \bigg)^{-1} \bxi_{\Ab,i}\\
&\leq C\sum_{i=1}^p \lambda_{\Ab,i}=C\|\Ab\|_{\tr}<+\infty
\end{align*}
almost surely for large enough $n$. Here, $\lambda_{\Ab,i}$ is the eigenvalue of matrix $\Ab$ and  $\bxi_{\Ab,i}$ is the corresponding eigenvector. Hence, the main target now is to deduce the uniform bound of \eqref{eq:uniformbound_vector}. We separate the proof of the uniform bound of \eqref{eq:uniformbound_vector} into the case  $z\leq \sqrt{\frac{\lambda_{\min}(\Qb)}{18c(\sum_{j=1}^K\sqrt{\lambda_j})^2}}$  and  $z\geq \sqrt{\frac{\lambda_{\min}(\Qb)}{18c(\sum_{j=1}^K\sqrt{\lambda_j})^2}}$. 

When  $z\geq \sqrt{\frac{\lambda_{\min}(\Qb)}{18c(\sum_{j=1}^K\sqrt{\lambda_j})^2}}$,  it can be seen that \eqref{eq:uniformbound_vector} is smaller than the value when we let $z= \sqrt{\frac{\lambda_{\min}(\Qb)}{18c(\sum_{j=1}^K\sqrt{\lambda_j})^2}}$, 
hence it remained for us to complete the proof when $z\leq \sqrt{\frac{\lambda_{\min}(\Qb)}{18c(\sum_{j=1}^K\sqrt{\lambda_j})^2}}$.

When  $z\leq \sqrt{\frac{\lambda_{\min}(\Qb)}{18c(\sum_{j=1}^K\sqrt{\lambda_j})^2}}$, by Assumption~\ref{ass:assump3}, the number of eigenvalues tending to infinity remains fixed,  we can separate $\bSigma$ by  
\begin{align*}
    \bSigma=\bGamma_s \bLambda_s \bGamma_s^\top + \bGamma_b \bLambda_b \bGamma_b^\top,
\end{align*}
here $\bGamma_s\in \RR^{p\times K}$ is the eigenvectors corresponding to the diagnal matrix $\bLambda_s\in\RR^{K\times K}$ with $\lambda_{\min}(\bLambda_s)\to +\infty$ as $n\to+\infty$,  $\bGamma_b\in \RR^{p\times (p-K)}$ is the eigenvectors corresponding to the diagnal matrix $\bLambda_b\in\RR^{(p-K)\times (p-K)}$ with $\|\bLambda_b\|_{\op}$ uniformly bounded. To give the uniform bound of \eqref{eq:uniformbound_vector}, we can easily see that it is equal to give the uniform bound of 
\begin{align}
   \label{eq:uniformbound_vector_s} 
   \ub^\top\bigg(\frac{\bSigma^{\frac{1}{2}}\Zb^\top\Zb\bSigma^{\frac{1}{2}}}{n}+\Qb +z\bSigma \bigg)^{-1} \bGamma_s \bLambda_s \bGamma_s^\top\bigg(\frac{\bSigma^{\frac{1}{2}}\Zb^\top\Zb\bSigma^{\frac{1}{2}}}{n}+\Qb +z\bSigma \bigg)^{-1} \ub 
\end{align}
due to 
\begin{align*}
    \ub^\top\bigg(\frac{\bSigma^{\frac{1}{2}}\Zb^\top\Zb\bSigma^{\frac{1}{2}}}{n}+\Qb +z\bSigma \bigg)^{-1} \bGamma_b \bLambda_b \bGamma_b^\top\bigg(\frac{\bSigma^{\frac{1}{2}}\Zb^\top\Zb\bSigma^{\frac{1}{2}}}{n}+\Qb +z\bSigma \bigg)^{-1} \ub\leq \|\bLambda_b\|_{\op}/\lambda_{\min}^2(\Qb)
\end{align*}
is uniformly bounded. Now, for the reason that $K$ is fixed, from \eqref{eq:uniformbound_vector_s}, we  only need to prove that for any deterministic $\ub$ with $\| \ub\|=1$,
\begin{align}
\label{eq:uniformbound_vector_single}
\sqrt{\lambda_j}\bigg|\bxi_j^\top\bigg(\frac{\bSigma^{\frac{1}{2}}\Zb^\top \Zb \bSigma^{\frac{1}{2}}}{n}+\Qb+z\bSigma\bigg)^{-1} \ub\bigg|
\end{align}
is uniformly bounded for $n$ large enough. Here, $\lambda_j\in \bLambda_s$, and $\bxi_j$ is the eigenvector corresponding to $\lambda_j$.  
Similar to the decomposition of $\bSigma$, we decompose $\bSigma^{\frac{1}{2}}\Zb^\top \Zb \bSigma^{\frac{1}{2}}/n$ into:
\begin{align*}
    \frac{\bSigma^{\frac{1}{2}}\Zb^\top \Zb \bSigma^{\frac{1}{2}}}{n}=\hbGamma_s \hbLambda_s\hbGamma_s^\top+\hbGamma_b\hbLambda_b\hbGamma_b^\top,
\end{align*}
and rewrite $\tQb=\Qb+z\bSigma$, then \eqref{eq:uniformbound_vector_single} can be written as 
\begin{align*}
    &\sqrt{\lambda_j}\big| \bxi_j^\top \big(\hbGamma_s \hbLambda_s\hbGamma_s^\top+\hbGamma_b\hbLambda_b\hbGamma_b^\top+\tQb \big)^{-1}\ub\big|\\
    &\quad=\sqrt{\lambda_j}\left| \bxi_j^\top\left( \big[\hbGamma_s,\hbGamma_b\big]\begin{bmatrix}
        \hbLambda_s+\hbGamma_s^\top\tQb \hbGamma_s & \hbGamma_s^\top\tQb \hbGamma_b\\
        \hbGamma_b^\top\tQb \hbGamma_s&\hbLambda_b+\hbGamma_b^\top\tQb \hbGamma_b
    \end{bmatrix}\begin{bmatrix}
        \hbGamma_s^\top\\
        \hbGamma_b^\top
    \end{bmatrix} \right)^{-1}\ub\right|\\
    &\quad =\sqrt{\lambda_j}\left| \bxi_j^\top \big[\hbGamma_s,\hbGamma_b\big]\begin{bmatrix}
        \hbLambda_s+\hbGamma_s^\top\tQb \hbGamma_s & \hbGamma_s^\top\tQb \hbGamma_b\\
        \hbGamma_b^\top\tQb \hbGamma_s&\hbLambda_b+\hbGamma_b^\top\tQb \hbGamma_b
    \end{bmatrix}^{-1}\begin{bmatrix}
        \hbGamma_s^\top\\
        \hbGamma_b^\top
    \end{bmatrix} \ub\right|.
\end{align*}
We further define $\Ab=\hbLambda_s+\hbGamma_s^\top\tQb \hbGamma_s$, $\Bb=\hbGamma_s^\top\tQb \hbGamma_b$ and $\Db=\hbLambda_b+\hbGamma_b^\top\tQb \hbGamma_b$, by the inverse of block matrices, we have
\begin{align*}
    \begin{bmatrix}
        \Ab&\Bb\\
        \Bb^\top &\Db
    \end{bmatrix}^{-1}=\begin{bmatrix}
        (\Ab-\Bb\Db^{-1}\Bb^\top)^{-1}&-(\Ab-\Bb\Db^{-1}\Bb^\top)^{-1}\Bb\Db^{-1}\\
        -\Db^{-1}\Bb^\top (\Ab-\Bb\Db^{-1}\Bb^\top)^{-1}&\Db^{-1}+\Db^{-1}\Bb^{\top}(\Ab-\Bb\Db^{-1}\Bb^\top)^{-1}\Bb\Db^{-1}
    \end{bmatrix},
\end{align*}
hence \eqref{eq:uniformbound_vector_single} can be written as 
\begin{align}
\label{eq:uniformbound_vector_single1}
\eqref{eq:uniformbound_vector_single}=\sqrt{\lambda_j}\big|\bxi_j^\top \hbGamma_s \bDelta_{ss}\hbGamma_s^\top\ub+\bxi_j^\top \hbGamma_b \bDelta_{bs}\hbGamma_s^\top\ub+\bxi_j^\top \hbGamma_s \bDelta_{sb}\hbGamma_b^\top\ub+\bxi_j^\top \hbGamma_b \bDelta_{bb}\hbGamma_b^\top\ub \big|,
\end{align}
where
\begin{align*}
    &\bDelta_{ss}=(\Ab-\Bb\Db^{-1}\Bb^\top)^{-1},\\
    &\bDelta_{bs}=\bDelta_{sb}^\top=-(\Ab-\Bb\Db^{-1}\Bb^\top)^{-1}\Bb\Db^{-1},\\
&\bDelta_{bb}=\Db^{-1}+\Db^{-1}\Bb^{\top}(\Ab-\Bb\Db^{-1}\Bb^\top)^{-1}\Bb\Db^{-1}.
\end{align*}
We will bound the operator norm of the three terms above. We first  note that $ \lambda_{\min}(\Db)\geq \lambda_{\min}(\hbGamma_b^\top\tQb \hbGamma_b)\geq \lambda_{\min}(\tilde{\Qb})\geq \lambda_{\min}(\Qb)$, we have
$\|\Db^{-1}\|_{\op}\leq 1/\lambda_{\min}(\Qb)$ is uniformly bounded. We next bound $\| \Bb\|_{\op}$. Indeed, we have
\begin{align*}
\|\Bb\|_{\op}&=\big\| \hbGamma_s^\top (\Qb+z\bSigma)\hbGamma_b\big\|_{\op}\leq \|\Qb\|_{\op}+z\big\| \hbGamma_s^\top(\bGamma_s\bLambda_s\bGamma_s^\top+\bGamma_b\bLambda_b\bGamma_b^\top)\hbGamma_b\big\|_{\op}\\
&\leq \|\Qb\|_{\op}+z\|\bLambda_b\|_{\op}+z\sum_{j=1}^K \lambda_j \big\|\bxi_j^\top \hbGamma_b\big\|,
\end{align*}
where the first inequality is by $\|\Ab+\Bb\|_{\op}\leq \|\Ab\|_{\op}+\|\Bb\|_{\op}$, and the second inequality is by $\|\Ab\Bb\|_{\op}\leq \|\Ab\|_{\op}\|\Bb\|_{\op}$ and $\|\hbGamma_s \|_{\op},\|\bGamma_s\|_{\op},\|\hbGamma_b \|_{\op},\|\bGamma_b\|_{\op}\leq 1$. We remind readers here $\bxi_j\in\bGamma_s$. As demonstrated in \citet{bao2022statistical}'s Theorem~2.5, for $n$ large enough, almost surely we have 
\begin{align}
\label{eq:boundofbxiiGammab}
    \lambda_j\|\bxi_j\hbGamma_b\|^2\leq 2\lambda_j\bigg(1-\frac{\lambda_j^2-c}{\lambda_j(\lambda_j+c)} \bigg)=\frac{2c(\lambda_j+1)}{\lambda_j+c}\leq 2c. 
\end{align}
Therefore, when  $z\leq \sqrt{\frac{\lambda_{\min}(\Qb)}{18c(\sum_{j=1}^K\sqrt{\lambda_j})^2}}$, we have almost surely 
\begin{align*}
\|\Bb\|_{\op}\leq \|\Qb\|_{\op}+ (z\cdot \sum_{j=1}^{K}\sqrt{3c\lambda_j})=O(1).
\end{align*}
for $n$ large enough. Hence when  $z\leq \sqrt{\frac{\lambda_{\min}(\Qb)}{18c(\sum_{j=1}^K\sqrt{\lambda_j})^2}}$,
\begin{align*}
    \frac{\lambda_{\min}(\Ab-\Bb\Db^{-1}\Bb^\top)}{\lambda_{K}}\geq \frac{2\lambda_K/3-\|\Bb\|_{\op}^2/\lambda_{\min}(\Qb)}{\lambda_K}\geq \frac{1}{2}
\end{align*}
almost surely for all $n$ large enough. Here, the first inequality utilizes $\hlambda_K\geq 2\lambda_K/3$ almost surely for large enough $n$.  We prove that when  $z\leq \sqrt{\frac{\lambda_{\min}(\Qb)}{18c(\sum_{j=1}^K\sqrt{\lambda_j})^2}}$ for any $j\in[1:K]$, 

\begin{align}
\label{eq:deltass}
    \|\bDelta_{ss} \|_{\op}\cdot \sqrt{\lambda_j}=\sqrt{\lambda_j}/\lambda_{\min}(\Ab-\Bb\Db^{-1}\Bb^\top)\leq 2\sqrt{\lambda_j}/\lambda_K\leq 2C.
\end{align}
almost surely for all $n$ large enough. Note that    $\|\Bb\|_{\op}=O(1)$  as long as  $z\leq \sqrt{\frac{\lambda_{\min}(\Qb)}{18c(\sum_{j=1}^K\sqrt{\lambda_j})^2}}$ for $n$ large enough, we conclude that 
 
\begin{align}
\label{eq:deltasb}
    \|\bDelta_{bs}\|_{\op} \cdot \sqrt{\lambda_j}\leq \sqrt{\lambda_j}\|\bDelta_{ss}\|_{\op}\|\Bb\|_{\op}\leq 3C
\end{align}
almost surely for all $n$ large enough when  $z\leq \sqrt{\frac{\lambda_{\min}(\Qb)}{18c(\sum_{j=1}^K\sqrt{\lambda_j})^2}}$  for some large enough absolute $C>0$. As for $\|\bDelta_{bb}\|_{\op}$, we have
\begin{align}
\label{eq:deltabb}
    \| \bDelta_{bb}\|_{\op}\leq \frac{1}{\lambda_{\min}(\Qb)}+\frac{1}{\lambda_{\min}^2(\Qb)}\cdot \|\Bb\|_{\op}^2/\lambda_K\cdot \|\bDelta_{ss}\|_{\op}\cdot\lambda_K\leq\frac{4}{3\lambda_{\min}(\Qb)}
\end{align}
almost surely for all $n$ large enough when  $z\leq \sqrt{\frac{\lambda_{\min}(\Qb)}{18c(\sum_{j=1}^K\sqrt{\lambda_j})^2}}$ . Combined \eqref{eq:uniformbound_vector_single1} with \eqref{eq:boundofbxiiGammab}-\eqref{eq:deltabb},  we conclude that when  $z\leq \sqrt{\frac{\lambda_{\min}(\Qb)}{18c(\sum_{j=1}^K\sqrt{\lambda_j})^2}}$,

\begin{align*}
\eqref{eq:uniformbound_vector_single}&\leq  \sqrt{\lambda_j} \cdot (\|\bDelta_{ss}\|_{\op}+2\|\bDelta_{bs}\|_{\op}+\| \bxi_j \hbGamma_b\|\cdot \|\bDelta_{bb}\|_{\op})\\
&\leq 6C
\end{align*}
almost surely for all $n$ large enough and for some large enough absolute constant $C>0$. We hence conclude that when  $z\leq \sqrt{\frac{\lambda_{\min}(\Qb)}{18c(\sum_{j=1}^K\sqrt{\lambda_j})^2}}$, $|f'_n(z)|$ is uniformly bounded. Hence we have that $|f'_n(z)|$ is uniformly bounded when $z\geq0$. After similar argument in the proof of Lemma~\ref{lemma:z=0converge} completes the proof of Lemma~\ref{lemma:z=0converge_c>1}. 

\subsection{Proof of Lemma~\ref{lemma:stronger_type_convergence}}
\label{sec:prooflemmastronger_type_convergence}
Before we give the proof, we prove the following lemma.
\begin{lemma}
\label{lemma:stronger_lemma}
For each $n\in N_+$ and $i\in[n]$, define sequence $a_{i,n},b_{i,n}\geq 0$ satisfy
\begin{align*}
\lim_{n\to+\infty}a_{i,n}/b_{i,n}\to 1
\end{align*}
for any $i\leq n$. Then it holds that
\begin{align*}
    \lim_{n\to+\infty}\sum_{i=1}^n a_{i,n}/\sum_{i=1}^n b_{i,n} =1.
\end{align*}
\end{lemma}
\begin{proof}[\textbf{Proof of Lemma~\ref{lemma:stronger_lemma}}]
For any $\varepsilon>0$, there exists $N$, when $n>N$, $|a_{i,n}-b_{i,n}|\leq \varepsilon b_{i,n}$ for any $i\leq n$, therefore for $n>N$ we have
\begin{align*}
    \sum_{i=1}^n a_{i,n}/\sum_{i=1}^n b_{i,n}\leq \sum_{i=1}^n (1+\varepsilon)b_{i,n}/\sum_{i=1}^n b_{i,n}= 1+\varepsilon.
\end{align*}
Here, the first inequality is by $a_{i,n}\leq (1+\varepsilon)b_{i,n}$. Similarly we have $\sum_{i=1}^n a_{i,n}/\sum_{i=1}^n b_{i,n}\geq 1-\varepsilon$. This completes the proof of Lemma~\ref{lemma:stronger_lemma}.
\end{proof}

\begin{proof}[\textbf{Proof of Lemma~\ref{lemma:stronger_type_convergence}}]
Without loss of generality, we can assume $\|\Ab\|_{\tr}$ is bounded. Otherwise, we can over $\|\Ab\|_{\tr}$ in both numerator and denominator.

By Proposition~\ref{prop:all_case_converge}, we have that for any given $z\in \bbC\backslash(-\infty,0)\bigcup \{0\}$, 
\begin{align*}
    |f_n(z)-\overline{f}_n(z)|~\cvas~0,\quad |f'_n(z)-\overline{f}'_n(z)|~\cvas~0.
\end{align*}
It is also clear that $|f_n(z)|$, $|\overline{f}_n(z)|$, $|f'_n(z)|$ and $|\overline{f}'_n(z)|$ are all bounded. 

When Case 1 in Assumption~\ref{ass:assump3} holds ($\|\bSigma\|_{\op}$ is bounded), we decompose $\Ab$ by $\Ab=\sum_{i=1}^p\lambda_{i,\Ab}\bxi_{i,\Ab}\bxi_{i,\Ab}^\top$, where $\lambda_{i,\Ab}$ is the eigenvalue and $\bxi_{i,\Ab}$ is the corresponding eigenvector. Then $f_n(z)$ and $\overline{f}_n(z)$ can be expressed as 
\begin{align*}
    f_n(z)&=\sum_{i=1}^p \lambda_{i,\Ab}\bxi_{i,\Ab}^\top\bigg(\frac{\bSigma^{\frac{1}{2}}\Zb^\top\Zb\bSigma^{\frac{1}{2}}}{n}+ \Qb +z\bSigma \bigg)^{-1} \bxi_{i,\Ab},\\
    \overline{f}_n(z)&=\sum_{i=1}^p \lambda_{i,\Ab}\bxi_{i,\Ab}^\top\bigg(\frac{\bSigma}{1+s_0}+ \Qb +z\bSigma \bigg)^{-1} \bxi_{i,\Ab}.
\end{align*}
Let $a_{i,n}=\lambda_{i,\Ab}\bxi_{i,\Ab}^\top\big(\frac{\bSigma^{\frac{1}{2}}\Zb^\top\Zb\bSigma^{\frac{1}{2}}}{n}+ \Qb +z\bSigma \big)^{-1} \bxi_{i,\Ab}$ and $b_{i,n}=\lambda_{i,\Ab}\bxi_{i,\Ab}^\top\big(\frac{\bSigma}{1+s_0}+ \Qb +z\bSigma \big)^{-1} \bxi_{i,\Ab}$, we have that $a_{i,n}/\lambda_{i,\Ab}-b_{i,n}/\lambda_{i,\Ab}~\cvas~0$ and  $b_{i,n}/\lambda_{i,\Ab}$ has the constant lower bound,  hence $a_{i,n}/b_{i,n}~\cvas~1$ for any $i\leq p$. From Lemma~\ref{lemma:stronger_lemma}, we have 
$f_n(z)/\overline{f}_n(z)~\cvas~1$.

When Case 2 and 3 in Assumption~\ref{ass:assump3} hold, we have
\begin{align*}
{f}_n(z)&=\tr\bigg(\frac{\Zb^\top\Zb}{n}+ \bSigma^{-\frac{1}{2}}\Qb \bSigma^{-\frac{1}{2}}+z\Ib \bigg)^{-1}\bSigma^{-\frac{1}{2}}\Ab \bSigma^{-\frac{1}{2}},\\
    \overline{f}_n(z)&=\tr\bigg(\frac{\Ib}{1+s_0}+ \bSigma^{-\frac{1}{2}}\Qb \bSigma^{-\frac{1}{2}}+z\Ib \bigg)^{-1}\bSigma^{-\frac{1}{2}}\Ab \bSigma^{-\frac{1}{2}}.
\end{align*}
Decompose $\bSigma^{-\frac{1}{2}}\Ab\bSigma^{-\frac{1}{2}}$ by $\bSigma^{-\frac{1}{2}}\Ab\bSigma^{-\frac{1}{2}}=\sum_{i=1}^p\lambda_{i}\bxi_{i}\bxi_i^\top$, where $\lambda_i$ is the eigenvalue and $\bxi_i$ is the corresponding eigenvector. Hence $f_n(z)$ and $\overline{f}_n(z)$ can be expressed as 
\begin{align*}
    f_n(z)&=\sum_{i=1}^p \lambda_i \bxi_i^\top\bigg(\frac{\Zb^\top\Zb}{n}+ \bSigma^{-\frac{1}{2}}\Qb \bSigma^{-\frac{1}{2}}+z\Ib \bigg)^{-1} \bxi_i,\\
    \overline{f}_n(z)&=\sum_{i=1}^p\lambda_i \bxi_i^\top\bigg(\frac{\Ib}{1+s_0}+ \bSigma^{-\frac{1}{2}}\Qb \bSigma^{-\frac{1}{2}}+z\Ib \bigg)^{-1} \bxi_i.
\end{align*}
Similarly, let $a_{i,n}=\lambda_i \bxi_i^\top\big(\frac{\Zb^\top\Zb}{n}+ \bSigma^{-\frac{1}{2}}\Qb \bSigma^{-\frac{1}{2}}+z\Ib \big)^{-1} \bxi_i$ and $b_{i,n}=\lambda_i\bxi_i^\top \big(\frac{\Ib}{1+s_0}+ \bSigma^{-\frac{1}{2}}\Qb \bSigma^{-\frac{1}{2}}+z\Ib \big)^{-1} \bxi_i$, we have that $a_{i,n}/\lambda_i-b_{i,n}/\lambda_i~\cvas~0$ and  $b_{i,n}/\lambda_i$ has the constant lower bound, hence $a_{i,n}/b_{i,n}~\cvas~1$ for any $i\leq p$. From Lemma~\ref{lemma:stronger_lemma}, we have 
$f_n(z)/\overline{f}_n(z)~\cvas~1$. The proof for $f'_n(z)/\overline{f}'_n(z)~\cvas~1$ is similar to $f_n(z)/\overline{f}_n(z)~\cvas~1$, and we thus omit it. 
\end{proof}


\section{Proof of Theorem~\ref{thm:main_theorem}}
\label{sec:appendix_proof_thm_bound}
Let $s_0(z)$ be defined in Definition~\ref{def:s0-f(z)}, and set $s_0=s_0(0)$. We rewrite Theorem~\ref{thm:main_theorem} and extend it as follows:
\begin{theorem}
    \label{thm:main_theoremadd}
    Suppose Assumptions~\ref{def:X}, \ref{ass:assump1}, \ref{ass:assump2} and \ref{ass:assump3} hold. For any $\Qb\in\cQ$, a good estimator $\hat{SR}(\Qb)$ for $SR(\Qb)$ which is defined in \eqref{eq:definition_SR}, is given as follows. 
    \begin{align*}
        \hat{SR}(\Qb)=\frac{T_{n,1}(\Qb)}{\sqrt{\big|\hat{T}_{n,2}(\Qb)\big|}},\quad \text{where}\quad \hat{T}_{n,2}(\Qb)=\frac{\tr(\hat{\bSigma}+\Qb)^{-1}\hbSigma(\hat{\bSigma}+\Qb)^{-1}\Ab}{\big(1-\frac{c}{p} \tr\hbSigma(\hbSigma+\Qb)^{-1}\big)^2}.
    \end{align*}
The following properties hold:
\begin{enumerate}
    \item If $\|\Ab\|_{\tr} $ is bounded, it holds that $\hat{T}_{n,2}(\Qb)-T_{n,2}(\Qb)~\cvas~0$.
    \item If $\Ab$ is semi-positive definite, it holds that $\hat{T}_{n,2}(\Qb)/T_{n,2}(\Qb)~\cvas~1 $ and $ \hat{SR}(\Qb)/SR(\Qb)~\cvas~1$.
    \item If $\Ab$ is semi-positive definite and $\tr\Big(\frac{\bSigma}{1+s_0}+\Qb \Big)^{-1}\Ab$ is bounded, it holds that $ \hat{SR}(\Qb)-SR(\Qb)~\cvas~0 $.
\end{enumerate}
\end{theorem}

In order to prove Theorem~\ref{thm:main_theoremadd}, we give the precise asymptotic of $T_{n,1}(\Qb)$ and $T_{n,2}(\Qb)$, and then construct a statistics only using in-sample data which share the same asymptotic as $T_{n,1}(\Qb)$ and $T_{n,2}(\Qb)$, finally we use these statistics to estimate $SR(\Qb)$. Define
\begin{align*}
    g(\bSigma)=\tr\bigg( \frac{\bSigma}{1+s_0}+\Qb \bigg)^{-1}\bSigma\bigg( \frac{\bSigma}{1+s_0}+\Qb \bigg)^{-1}\Ab.
\end{align*}
Recall the definitions of $T_{n,1}(\Qb)$ and $T_{n,2}(\Qb)$ in \eqref{eq:Tn1} and \eqref{eq:Tn2} that
\begin{align*}
    T_{n,1}(\Qb)=\tr (\hbSigma+\Qb)^{-1}\Ab,\quad T_{n,2}(\Qb)=\tr (\hbSigma+\Qb)^{-1}\bSigma(\hbSigma+\Qb)^{-1}\Ab.
\end{align*}
The proposition below gives us the precise asymptotic of $T_{n,1}(\Qb)$ and $T_{n,2}(\Qb)$ and is proved in Appendix~\ref{appsec:proofthm1}.
\begin{proposition}
\label{prop:converge_Tn1Tn2}
Suppose that Assumptions in Theorem~\ref{thm:main_theoremadd} hold. Let $T_{n,1}(\Qb)$ and $T_{n,2}(\Qb)$ defined in \eqref{eq:Tn1} and \eqref{eq:Tn2} respectively, then if $\|\Ab\|_{\tr}$ is bounded, it holds that
\begin{align*}
    &T_{n,1}(\Qb)-\tr\bigg( \frac{\bSigma}{1+s_0}+\Qb \bigg)^{-1}\Ab~\cvas~0,\\
    &T_{n,2}(\Qb)-\frac{(1+s_0)^2-s_{1,\bSigma}}{(1+s_0)^2}g(\bSigma)~\cvas~0.
\end{align*}
Here, $s_0$ and $s_{1,\bSigma}$ uniquely solve the equations
\begin{align*}
    s_0=\frac{c}{p}\tr~\bSigma\bigg( \frac{\bSigma}{1+s_0}+\Qb \bigg)^{-1},\quad s_{1,\bSigma}=\bigg(\frac{s_{1,\bSigma}}{(1+s_0)^2}-1\bigg)\cdot\frac{c}{p}\tr~\bSigma\bigg( \frac{\bSigma}{1+s_0}+\Qb\bigg)^{-1}\bSigma\bigg( \frac{\bSigma}{1+s_0}+\Qb\bigg)^{-1}.
\end{align*}
Moreover, if $\Ab$ is  semi-positive definite, then 
\begin{align*}
    &T_{n,1}(\Qb)/\tr\bigg( \frac{\bSigma}{1+s_0}+\Qb \bigg)^{-1}\Ab~\cvas~1,\\
    &T_{n,2}(\Qb)/\bigg(\frac{(1+s_0)^2-s_{1,\bSigma}}{(1+s_0)^2}g(\bSigma)\bigg)~\cvas~1.
\end{align*}
\end{proposition}
The next proposition further indicates how to construct an estimator for $g(\bSigma)$ only using $\hat\bSigma$ instead of $\bSigma$ to achieve the same asymptotic of $T_{n,1}(\Qb)$ and $T_{n,2}(\Qb)$.
\begin{proposition}
\label{prop:converge_gSigma}
    Let $s_0$  be defined in Proposition~\ref{prop:converge_Tn1Tn2}. For any matrix $\Ab$ with $\|\Ab\|_{\tr}$ bounded, define 
\begin{align*}
    \hat{g}(\hbSigma)=\hat{T}_1-\hat{T}_2=\tr(\hbSigma+\Qb)^{-1}\hbSigma(\hbSigma+\Qb)^{-1}\Ab,
\end{align*}
where $\hat{T}_1=\tr(\hat{\bSigma}+\Qb)^{-1}\Ab$, $\hat{T}_2=\tr(\hat{\bSigma}+\Qb)^{-1}\Qb(\hat{\bSigma}+\Qb)^{-1}\Ab$ and $\hbSigma=\Xb^\top\Xb/n$. Then under the same conditions as Proposition~\ref{prop:converge_Tn1Tn2}, 
\begin{align*}
    \hat{g}(\hbSigma)-\frac{(1+s_0+ s_{1,\Qb})g(\bSigma)}{(1+s_0)^2}~\cvas~0,
\end{align*}
where $s_{1,\Qb}$ solves the following equation
\begin{align*}
    s_{1,\Qb}=\frac{c}{p}\tr~\bSigma\bigg( \frac{\bSigma}{1+s_0}+\Qb \bigg)^{-1}\bigg(\frac{s_{1,\Qb}}{(1+s_0)^2}\bSigma-\Qb\bigg)\bigg( \frac{\bSigma}{1+s_0}+\Qb \bigg)^{-1}.
\end{align*}
Moreover, if $\Ab$ is semi-positive definite, then
\begin{align*}
    \frac{(1+s_0)^2\hat{g}(\hbSigma)/g(\bSigma)}{1+s_0+ s_{1,\Qb}}~\cvas~1.
\end{align*}
\end{proposition}
Proposition~\ref{prop:converge_gSigma} constructs the estimation for the key part in the asymptotic of $T_{n,2}(\Qb)$ given in Proposition~\ref{prop:converge_Tn1Tn2}, and is proved in Appendix~\ref{appsec:proofprop2}.   Now, the only thing remains for us is to construct estimators for $s_0$, $s_{1,\bSigma}$ and $s_{1,\Qb}$ respectively. The next proposition established such precise estimations.
\begin{proposition}
\label{prop:converge_s}
Let $s_0$, $s_{1,\bSigma}$ and $s_{1,\Qb}$ be defined in Proposition~\ref{prop:converge_Tn1Tn2} and Proposition~\ref{prop:converge_gSigma} respectively. Define 
\begin{align*}
    f_1(\hat{\bSigma})=\frac{1}{p}\tr\Qb(\hat{\bSigma}+\Qb)^{-1},\quad f_2(\hat{\bSigma})=\frac{1}{p}\tr\Qb(\hat{\bSigma}+\Qb)^{-1}\Qb(\hat{\bSigma}+\Qb)^{-1},
\end{align*}
then under the same condition as Proposition~\ref{prop:converge_Tn1Tn2},
 \begin{align*}
     &\frac{1}{s_0}\cdot\frac{c\big(1-f_1(\hat{\bSigma})\big)}{1+c\big(f_1(\hat{\bSigma})-1\big)}~\cvas~1,\quad \frac{1}{s_{1,\Qb}}\cdot c\cdot\frac{f_2(\hat{\bSigma})-f_1(\hat{\bSigma})}{\big(1+c(f_1(\hat{\bSigma})-1)\big)^2}~\cvas~1,\\
     &\frac{1}{s_{1,\bSigma}}\cdot c\cdot\frac{\big(-1+2f_1(\hat{\bSigma})-f_2(\hat{\bSigma})+c(f_1(\hat{\bSigma})-1)^2\big) }{\big(1+c(f_1(\hat{\bSigma})-1)\big)^4}~\cvas~1.
 \end{align*}
 Moreover, it holds from the conclusion above that
 \begin{align*}
    \frac{1+s_0+s_{1,\Qb}}{(1+s_0)^2-s_{1,\bSigma}} /\bigg(1-\frac{c}{p}\tr\hbSigma\big(\hbSigma+\Qb\big)^{-1} \bigg)^2\cvas1.
\end{align*}
\end{proposition}
Proposition~\ref{prop:converge_s} establishes the statistics only related to $\Xb$ to estimate the constants, and is proved in Appendix~\ref{appsec:proofprop3}. With Proposition~\ref{prop:converge_Tn1Tn2}-\ref{prop:converge_s} above, we give the proof of Theorem~\ref{thm:main_theoremadd}.

\begin{proof}[\textbf{Proof of Theorem~\ref{thm:main_theoremadd}}]

Below we first prove the three conclusions when $\Qb\neq \zero$.

\noindent\textbf{\underline{The first conclusion.}} If $\|\Ab\|_{\tr} $ is bounded, from Proposition~\ref{prop:converge_Tn1Tn2} and \ref{prop:converge_gSigma}, we conclude
    \begin{align*}
        T_{n,2}(\Qb)-\frac{(1+s_0)^2-s_{1,\bSigma}}{(1+s_0)^2}g(\bSigma)~\cvas~0,\quad \frac{(1+s_0)^2\hat{g}(\hbSigma)}{1+s_0+s_{1,\Qb}}-g(\bSigma)~\cvas~0.
    \end{align*}
    Here, $\hat{g}(\hbSigma)=\tr(\hbSigma+\Qb)^{-1}\hbSigma(\hbSigma+\Qb)^{-1}\Ab $. Therefore we have
    \begin{align}
    \label{eq:Tn2_construct}
        T_{n,2}(\Qb)-\frac{(1+s_0)^2-s_{1,\bSigma}}{1+s_0+s_{1,\Qb}}\hat{g}(\hbSigma)~\cvas~0. 
    \end{align}   
Proposition~\ref{prop:converge_s} further gives us the estimation of $s_0$, $s_{1,\bSigma}$ and $s_{1,\Qb}$. Note that $s_0=s_0(0)$, $s_{1,\bSigma}=s_{1,\bSigma}(0)$, where $s_0(z)$ and $s_{1,\bSigma}(z)$ are defined in Lemma~\ref{lemma:constant_level}. We have $s_0>0$, $-s_0(1+s_0)^2<s_{1,\bSigma}<0$ and $1\leq 1+s_0+s_{1,\Qb}\leq 1+s_0$ are all bounded. Easy to conclude that $\frac{(1+s_0)^2-s_{1,\bSigma}}{(1+s_0)^2} $, $\frac{1+s_0+s_{1,\Qb}}{(1+s_0)^2} $, $ \frac{(1+s_0)^2-s_{1,\bSigma}}{1+s_0+s_{1,\Qb}}=\Theta(1)$. From Proposition~\ref{prop:converge_s} we also have
    \begin{align}
    \label{eq:Tn2_coef_construct}
       \frac{1}{(1+c(f_1(\hat\Sigma)-1))^2}/ \frac{(1+s_0)^2-s_{1,\bSigma}}{1+s_0+s_{1,\Qb}}~\cvas~1.
    \end{align}
    Thus combined \eqref{eq:Tn2_construct} with \eqref{eq:Tn2_coef_construct}  we have
    \begin{align*}
    &T_{n,2}(\Qb)-\frac{\hat{g}(\hbSigma)}{(1+c(f_1(\hat\Sigma)-1))^2}~\cvas~0.
\end{align*}
This directly indicates that
\begin{align*}
    \hat{T}_{n,2}(\Qb)-T_{n,2}(\Qb)~\cvas~0.
\end{align*}

\noindent\textbf{\underline{The second conclusion.}} When $\Ab$ is semi-positive definite, from Proposition~\ref{prop:converge_Tn1Tn2} and \ref{prop:converge_gSigma}, we conclude that
    \begin{align*}
        \frac{(1+s_0)^2-s_{1,\bSigma}}{(1+s_0)^2}g(\bSigma)/T_{n,2}(\Qb)~\cvas~1,\quad \frac{(1+s_0)^2\hat{g}(\hbSigma)}{1+s_0+s_{1,\Qb}}/g(\bSigma)~\cvas~1.
    \end{align*}
Hence we have 
\begin{align}
\label{eq:ness_over1}
    \frac{(1+s_0)^2-s_{1,\bSigma}}{1+s_0+s_{1,\Qb}}\cdot \frac{\hat{g}(\hbSigma)}{T_{n,2}(\Qb)}~\cvas~1.
\end{align}
It only requires to prove that 
\begin{align}
\label{eq:constant_over1}
    \frac{1+s_0+s_{1,\Qb}}{(1+s_0)^2-s_{1,\bSigma}}\cdot \frac{1}{(1+c(f_1(\hat\Sigma)-1))^2}~\cvas~1.
\end{align}
Proposition~\ref{prop:converge_s} shows that \eqref{eq:constant_over1} holds.
Combined \eqref{eq:constant_over1} with \eqref{eq:ness_over1}, we have
\begin{align*}
   \frac{1}{(1+c(f_1(\hat\Sigma)-1))^2}\cdot \frac{\hat{g}(\hbSigma)}{T_{n,2}(\Qb)}~\cvas~1,
\end{align*}
which completes the proof of the second conclusion in Theorem~\ref{thm:main_theoremadd}.

\noindent\textbf{\underline{The third conclusion.}} As for the last conclusion in Theorem~\ref{thm:main_theoremadd}, we both have $\Ab$ is semi-positive definite and $\Ab$ is bounded. Therefore by $\hat{SR}(\Qb)/SR(\Qb)~\cvas~1$, we can conclude that
\begin{align*}
    \big(\hat{SR}(\Qb)-SR(\Qb)\big)/SR(\Qb)~\cvas~0.
\end{align*}
As long as  we prove that $SR(\Qb)$ is almost surely upper bounded, we can easily conclude that $\hat{SR}(\Qb)-SR(\Qb)~\cvas~0$. It remains for us to prove that $SR(\Qb)$ is almost surely upper bounded. By Proposition~\ref{prop:converge_Tn1Tn2}, we have
\begin{align*}
    \frac{\tr\Big( \frac{\bSigma}{1+s_0}+\Qb \Big)^{-1}\Ab}{SR(\Qb)\cdot\sqrt{\frac{(1+s_0)^2-s_{1,\bSigma}}{(1+s_0)^2}\cdot g(\bSigma)}}~\cvas~1.
\end{align*}
Also note that by Lemma~\ref{lemma:constant_level}, the adjust factor  $\frac{(1+s_0)^2-s_{1,\bSigma}}{(1+s_0)^2}$ has the lower bound $1$ and the upper bound $1+s_0$,  it only requires us to prove that $\tr\big( \frac{\bSigma}{1+s_0}+\Qb \big)^{-1}\Ab/g(\bSigma)$ is upper bounded.  Recall that $g(\bSigma)=\tr\big( \frac{\bSigma}{1+s_0}+\Qb \big)^{-1}\bSigma\big( \frac{\bSigma}{1+s_0}+\Qb \big)^{-1}\Ab$, we have that
\begin{align}
    g(\bSigma)&=\tr\bigg( \frac{\bSigma}{1+s_0}+\Qb \bigg)^{-1}\bSigma\bigg( \frac{\bSigma}{1+s_0}+\Qb \bigg)^{-1}\Ab\nonumber\\
    &=(1+s_0)\tr\bigg( \frac{\bSigma}{1+s_0}+\Qb \bigg)^{-1}\frac{\bSigma}{1+s_0}\bigg( \frac{\bSigma}{1+s_0}+\Qb \bigg)^{-1}\Ab\nonumber\\
    &=(1+s_0)\bigg[\tr\bigg( \frac{\bSigma}{1+s_0}+\Qb \bigg)^{-1}\Ab -\tr\bigg( \frac{\bSigma}{1+s_0}+\Qb \bigg)^{-1}\Qb\bigg( \frac{\bSigma}{1+s_0}+\Qb \bigg)^{-1}\Ab\bigg].\label{eq:thirdconclude1}
\end{align}
Here, the last equality is by $\bSigma/(1+s_0)=\bSigma/(1+s_0)+\Qb-\Qb$. By the semi-positive definite of $\bSigma$, $\Qb$ and $\Ab$,  we can easily conclude that
\begin{align}
    \tr\bigg( \frac{\bSigma}{1+s_0}+\Qb \bigg)^{-1}&\Qb\bigg( \frac{\bSigma}{1+s_0}+\Qb \bigg)^{-1}\Ab\nonumber\\
    & \leq \lambda_{\max} \bigg(\Qb\bigg( \frac{\bSigma}{1+s_0}+\Qb \bigg)^{-1} \bigg) \cdot \tr\bigg( \frac{\bSigma}{1+s_0}+\Qb \bigg)^{-1}\Ab\nonumber\\
    &\leq\frac{\lambda_{\max}(\bSigma^{-\frac{1}{2}}\Qb\bSigma^{-\frac{1}{2}})}{\lambda_{\max}(\bSigma^{-\frac{1}{2}}\Qb\bSigma^{-\frac{1}{2}}+\Ib/(1+s_0))} \cdot \tr\bigg( \frac{\bSigma}{1+s_0}+\Qb \bigg)^{-1}\Ab\nonumber\\
    &\leq \frac{C(1+s_0)}{C(1+s_0)+1}\cdot \tr\bigg( \frac{\bSigma}{1+s_0}+\Qb \bigg)^{-1}\Ab.\label{eq:thirdconclude2}
\end{align}
Here, $\lambda_{\max}(\cdot) $ and $\lambda_{\min}(\cdot)$ denote by the largest and smallest eigenvalue of the matrix. $C$ is the maximal eigenvalue of $\bSigma^{-1/2}\Qb\bSigma^{-1/2}$.   The first inequality is by $\tr(\Ab\Bb)\leq \lambda_{\max}(\Ab)\tr\Bb$ for $\Ab,\Bb\geq0$, the second inequality is by $\lambda_{\max} \Big(\Qb\Big( \frac{\bSigma}{1+s_0}+\Qb \Big)^{-1} \Big)=\lambda_{\max} \Big(\Big( \frac{\Ib}{1+s_0}+\bSigma^{-\frac{1}{2}}\Qb\bSigma^{-\frac{1}{2}} \Big)^{-1}\cdot \bSigma^{-\frac{1}{2}}\Qb\bSigma^{-\frac{1}{2}} \Big)$. Therefore, plugging \eqref{eq:thirdconclude2} into \eqref{eq:thirdconclude1}, we have
\begin{align*}
    g(\bSigma)&\geq \frac{1+s_0}{C(1+s_0)+1}\cdot \tr\bigg( \frac{\bSigma}{1+s_0}+\Qb \bigg)^{-1}\Ab 
\end{align*}
which indicates that 
\begin{align*}
    \frac{\tr\Big( \frac{\bSigma}{1+s_0}+\Qb \Big)^{-1}\Ab}{\sqrt{g(\bSigma)}}&\leq \sqrt{\frac{C(1+s_0)+1}{1+s_0}}\cdot\sqrt{\tr\bigg( \frac{\bSigma}{1+s_0}+\Qb \bigg)^{-1}\Ab }\\
    &<+\infty.
\end{align*}
Therefore we conclude that $SR(\Qb)$ is almost surely bounded, which completes the proof of the third conclusion in Theorem~\ref{thm:main_theoremadd}.

As for  the case $\Qb=\zero$ with $c<1$. 
The proofs of Lemma~\ref{lemma:constant_level} and \ref{lemma:z=0converge} when $\Qb=\zero$ are slightly different. The solution of $s_0$ in Lemma~\ref{lemma:constant_level} can be obtained directly by the cancellation of $\bSigma$, and  $|f_n(z)|$ in Lemma~\ref{lemma:z=0converge} can also be bounded by $2\|\bSigma^{-\frac{1}{2}}\Ab\bSigma^{-\frac{1}{2}}\|_{\tr}/(1-\sqrt{c})^2$. Then when $\Qb=\zero$ under $c<1$, it is clear that $1+s_0~\cvas~1/(1-c)$ and $\frac{(1+s_0)^2-s_{1,\bSigma}}{1+s_0+s_{1,\Qb}}~\cvas~1/(1-c)^2$ after simple algebra calculation.
Following very similar steps, we can easily prove the same three conclusions. So we omit some details here.  
Readers may note that 
\begin{align*}
    1-\frac{c}{p}\tr\hbSigma(\hbSigma+\Qb)^{-1}~\cvas~1-c
\end{align*}
as $n$ tends to infinity when $\Qb=\zero$. The general form of $\hat{T}_{n,2}(\Qb)$
\begin{align*}
    \hat{T}_{n,2}(\Qb)=\frac{\tr(\hat{\bSigma}+\Qb)^{-1}\hbSigma(\hat{\bSigma}+\Qb)^{-1}\Ab}{\big(1-\frac{c}{p} \tr\hbSigma(\hbSigma+\Qb)^{-1}\big)^2}
\end{align*}
covers the case $\Qb=\zero$. 
This completes the proof of Theorem~\ref{thm:main_theorem}.
\end{proof}

\subsection{Proof of Proposition~\ref{prop:converge_Tn1Tn2}}
\label{appsec:proofthm1}

Based on the previous result in Section~\ref{sec:preliminary}, We give the following proposition which covers Proposition~\ref{prop:converge_Tn1Tn2} by taking $\lambda=0$ and $\Bb=\bSigma$.
\begin{proposition}
\label{prop:prop1proof}
Define the positive matrix $\Bb$ taking value on $\bSigma$ or $\Qb$, and let
\begin{align*}
    T_{n,1}(\lambda,\Bb)&=\tr\bigg(\frac{\Xb^\T\Xb}{n}+\Qb+\lambda\Bb\bigg)^{-1}\Ab,\\
     T_{n,2}(\lambda,\Bb)&=\tr\bigg(\frac{\Xb^\T\Xb}{n}+\Qb+\lambda\Bb\bigg)^{-1}\Bb\bigg(\frac{\Xb^\T\Xb}{n}+\Qb+\lambda\Bb\bigg)^{-1}\Ab.
\end{align*}
Under the same condition as Proposition~\ref{prop:converge_Tn1Tn2}, if $\lambda\geq0$, it holds that
\begin{align*}
        &T_{n,1}(\lambda,\Bb)-\tr\bigg( \frac{\bSigma}{1+s_{0,\Bb}(\lambda)}+\Qb+\lambda\Bb \bigg)^{-1}\Ab~\cvas~0,\\
        &T_{n,2}(\lambda,\Bb)-\tr\bigg( \frac{\bSigma}{1+s_{0,\Bb}(\lambda)}+\Qb+\lambda\Bb \bigg)^{-1}\bigg(\Bb-\frac{s_{1,\Bb}(\lambda)}{(1+s_{0,\Bb}(\lambda))^2}\bSigma\bigg)\bigg( \frac{\bSigma}{1+s_{0,\Bb}(\lambda)}+\Qb+\lambda\Bb \bigg)^{-1}\Ab~\cvas~0.
\end{align*}
Here, $s_{0,\Bb}(\lambda)$ and $s_{1,\Bb}(\lambda)$ uniquely solve the following equations
\begin{align*}
    &s_{0,\Bb}(\lambda)=\frac{c}{p}\tr\bSigma\bigg( \frac{ \bSigma}{1+s_{0,\Bb}(\lambda)}+ \Qb +\lambda \Bb  \bigg)^{-1},\\
    &s_{1,\Bb}(\lambda)=\frac{c}{p}\tr\bSigma\bigg( \frac{ \bSigma}{1+s_{0,\Bb}(\lambda)}+ \Qb +\lambda \Bb  \bigg)^{-1}
    \bigg(\frac{s_{1,\Bb}(\lambda)}{(1+s_{0,\Bb}(\lambda))^2} \bSigma- \Bb \bigg)
    \bigg( \frac{ \bSigma}{1+s_{0,\Bb}(\lambda)}+ \Qb +\lambda \Bb  \bigg)^{-1}.
\end{align*}
Moreover, if $\Ab$ is semi-positive definite and $\Bb=\bSigma$, it holds that
\begin{align*}
        &T_{n,1}(\lambda,\Bb)/\tr\bigg( \frac{\bSigma}{1+s_{0,\Bb}(\lambda)}+\Qb+\lambda\Bb \bigg)^{-1}\Ab~\cvas~1,\\
        &T_{n,2}(\lambda,\Bb)/\tr\bigg( \frac{\bSigma}{1+s_{0,\Bb}(\lambda)}+\Qb+\lambda\Bb \bigg)^{-1}\bigg(\Bb-\frac{s_{1,\Bb}(\lambda)}{(1+s_{0,\Bb}(\lambda))^2}\bSigma\bigg)\bigg( \frac{\bSigma}{1+s_{0,\Bb}(\lambda)}+\Qb+\lambda\Bb \bigg)^{-1}\Ab~\cvas~1.
\end{align*}
\end{proposition}

\begin{proof}[\textbf{Proof of Proposition~\ref{prop:prop1proof}}]
It is easy to check that $|T_{n,1}(\lambda,\Bb)|$ is bounded. 
For the convergence of $T_{n,1}(\lambda,\Bb)$, we can easy to see that 
\begin{align*}
    &T_{n,1}(\lambda,\Bb=\bSigma)=\tr\bSigma^{-\frac12} \bigg(\frac{\Zb^\T\Zb}{n}+\bSigma^{-\frac12}\Qb\bSigma^{-\frac12}+\lambda\Ib_p\bigg)^{-1}\bSigma^{-\frac12}\Ab,\\
    &T_{n,1}(\lambda,\Bb=\Qb)=\tr \bigg(\frac{\Xb^\T\Xb}{n}+\Qb+\lambda\Qb\bigg)^{-1}\Ab.
\end{align*}
Here $\Xb=\bSigma^{\frac{1}{2}}\Zb$. If $\Bb=\bSigma$, from Proposition~\ref{prop:all_case_converge} we can see 
\begin{align*}
    &T_{n,1}(\lambda,\Bb)-\tr\bigg( \frac{\bSigma}{1+s_{0,\Bb}(\lambda)}+\Qb+\lambda\Bb \bigg)^{-1}\Ab~\cvas~0,\\
    &\frac{\partial T_{n,1}(\lambda,\Bb)}{\partial \lambda}-\frac{\partial \tr \Big(\frac{\bSigma}{1+s_{0,\Bb}(\lambda)}+\Qb+\lambda\Bb\Big)^{-1}\Ab}{\partial\lambda}~\cvas~0
\end{align*}
holds for $\lambda\geq 0$. Here, we have 
\begin{align*}
    \frac{\rmd s_{0,\Bb}(\lambda)}{\rmd \lambda}=\frac{\rmd \frac{c}{p}\tr\bSigma\Big( \frac{ \bSigma}{1+s_{0,\Bb}(\lambda)}+ \Qb +\lambda \Bb  \Big)^{-1}}{\rmd \lambda}.
\end{align*}

If $\Bb=\Qb$, we can still have that for $z\in\bbC\backslash(-\infty,0]$, 
\begin{align*}
    T_{n,1}(z,\Bb)-\tr\bigg( \frac{\bSigma}{1+s_{0,\Bb}(z)}+\Qb+z\Bb \bigg)^{-1}\Ab~\cvas~0,
\end{align*}
and the convergence can be extended to $z=0$. The analytic function is internally uniform convergence, therefore the derivatives also converge. We conclude that when $\lambda\geq0$,
\begin{align*}
    &T_{n,1}(\lambda,\Bb)-\tr\bigg( \frac{\bSigma}{1+s_{0,\Bb}(\lambda)}+\Qb+\lambda\Bb \bigg)^{-1}\Ab~\cvas~0,\\
    &\frac{\partial T_{n,1}(\lambda,\Bb)}{\partial \lambda}-\frac{\partial \tr \Big(\frac{\bSigma}{1+s_{0,\Bb}(\lambda)}+\Qb+\lambda\Bb\Big)^{-1}\Ab}{\partial\lambda}~\cvas~0.
\end{align*}
Direct calculations give us 
\begin{align*}
    &\frac{\partial T_{n,1}(\lambda,\Bb)}{\partial \lambda}=-T_{n,2}(\lambda,\Bb),\\
    &\frac{\partial \tr \Big(\frac{\bSigma}{1+s_{0,\Bb}(\lambda)}+\Qb+\lambda\Bb\Big)^{-1}\Ab}{\partial\lambda}\\
    &\qquad=-\tr\bigg( \frac{\bSigma}{1+s_{0,\Bb}(\lambda)}+\Qb+\lambda\Bb \bigg)^{-1}\bigg(\Bb-\frac{s'_{0,\Bb}(\lambda)}{(1+s_{0,\Bb}(\lambda))^2}\bSigma\bigg)\bigg( \frac{\bSigma}{1+s_{0,\Bb}(\lambda)}+\Qb+\lambda\Bb \bigg)^{-1}\Ab,\\
    & s'_{0,\Bb}(\lambda)=\frac{c}{p}\tr\bSigma\bigg( \frac{ \bSigma}{1+s_{0,\Bb}(\lambda)}+ \Qb +\lambda \Bb  \bigg)^{-1}
    \bigg(\frac{s'_{0,\Bb}(\lambda)}{(1+s_{0,\Bb}(\lambda))^2} \bSigma- \Bb \bigg)
    \bigg( \frac{ \bSigma}{1+s_{0,\Bb}(\lambda)}+ \Qb +\lambda \Bb  \bigg)^{-1}.
\end{align*}
By the definition of $s_{1,\Bb}(\lambda)$ in Proposition~\ref{prop:prop1proof} we have $s'_{0,\Bb}(\lambda)=s_{1,\Bb}(\lambda)$, therefore the convergence of $T_{n,2}(\lambda,\Bb)$ is proved by 
\begin{align*}
    &\bigg|T_{n,2}(\lambda)-\tr\bigg( \frac{\bSigma}{1+s_{0,\Bb}(\lambda)}+\Qb+\lambda\Bb \bigg)^{-1}\bigg(\Bb-\frac{s'_0(\lambda)}{(1+s_{0,\Bb}(\lambda))^2}\bSigma\bigg)\bigg( \frac{\bSigma}{1+s_{0,\Bb}(\lambda)}+\Qb+\lambda\Bb \bigg)^{-1}\Ab\bigg|\\
    &\quad = \bigg|\frac{\partial T_{n,1}(\lambda,\Bb)}{\partial \lambda}-\frac{\partial \tr \Big(\frac{\bSigma}{1+s_{0,\Bb}(\lambda)}+\Qb+\lambda\Bb\Big)^{-1}\Ab}{\partial\lambda}\bigg| ~\cvas~0.
\end{align*}
This completes the first conclusion of Proposition~\ref{prop:prop1proof}. When $\Ab$ is semi-positive definite and $\Bb=\bSigma$, we have $T_{n,1}(\lambda,\Bb)=f_n(\lambda)$ and $T_{n,2}(\lambda,\Bb)=f'_n(\lambda)$. The conclusion directly comes from Lemma~\ref{lemma:stronger_type_convergence}. 
This completes the proof of Proposition~\ref{prop:prop1proof}.
\end{proof}
By taking $\Bb=\bSigma$, Proposition~\ref{prop:converge_Tn1Tn2} holds directly from Proposition~\ref{prop:prop1proof}.

\subsection{Proof of Proposition~\ref{prop:converge_gSigma}}
\label{appsec:proofprop2}
It is easy to get that the constants $s_0$ defined in Proposition~\ref{prop:converge_Tn1Tn2} and $s_{1,\Bb}$ ($\Bb=\bSigma$ or $\Bb=\Qb$ defined in Proposition~\ref{prop:converge_Tn1Tn2} and \ref{prop:converge_gSigma}) are defined by $s_0=s_{0,\Bb}(0)$ and $s_{1,\Bb}=s_{1,\Bb}(0)$ in Proposition~\ref{prop:prop1proof}. 
We remind the readers that $s_{1,\Bb}$ solves the following equation:
\begin{align*}
    s_{1,\Bb}=\frac{c}{p}\tr\bSigma\bigg( \frac{ \bSigma}{1+s_{0,\Bb}(\lambda)}+ \Qb   \bigg)^{-1}
    \bigg(\frac{s_{1,\Bb}}{(1+s_0)^2} \bSigma- \Bb \bigg)
    \bigg( \frac{ \bSigma}{1+s_0}+ \Qb  \bigg)^{-1}.
\end{align*}
Moreover, from Proposition~\ref{prop:prop1proof}, if we let $\Bb=\Qb$ and $\lambda=0$, then 
\begin{align*}
    \tr \big(\hat{\bSigma}+\Qb\big)^{-1}\Qb\big(\hat{\bSigma}+\Qb\big)^{-1}\Ab-\tr\bigg( \frac{\bSigma}{1+s_0}+\Qb \bigg)^{-1}\bigg(\Qb-\frac{s_{1,\Qb}}{(1+s_0)^2}\bSigma\bigg)\bigg( \frac{\bSigma}{1+s_0}+\Qb \bigg)^{-1}\Ab~\cvas~0,
\end{align*}
by the definition $g(\bSigma)=\tr\Big( \frac{\bSigma}{1+s_0}+\Qb \Big)^{-1}\bSigma\Big( \frac{\bSigma}{1+s_0}+\Qb \Big)^{-1}\Ab$ and  $\hat{T}_2=\tr \big(\hat{\bSigma}+\Qb\big)^{-1}\Qb\big(\hat{\bSigma}+\Qb\big)^{-1}\Ab$, the equation above can be rewrited as
\begin{align}
\label{eq:gsigma1}
    \hat{T}_2-\tr\bigg( \frac{\bSigma}{1+s_0}+\Qb \bigg)^{-1}\Qb\bigg( \frac{\bSigma}{1+s_0}+\Qb \bigg)^{-1}\Ab+\frac{s_{1,\Qb}}{(1+s_0)^2}g(\bSigma)~\cvas~0.
\end{align}
\eqref{eq:gsigma1} gives an equation to estimate $g(\bSigma)$. However, another equation is required to eliminate $\tr\big( \bSigma/(1+s_0)+\Qb \big)^{-1}\Qb\big( \bSigma/(1+s_0)+\Qb \big)^{-1}\Ab$ in \eqref{eq:gsigma1} and get the desired estimation of $g(\bSigma)$. To get the other equation, from the definition of $g(\bSigma)$ we have
\begin{align}
    g(\bSigma)&=\tr\bigg( \frac{\bSigma}{1+s_0}+\Qb \bigg)^{-1}\bSigma\bigg( \frac{\bSigma}{1+s_0}+\Qb \bigg)^{-1}\Ab\nonumber\\
    &=(1+s_0)\tr\bigg( \frac{\bSigma}{1+s_0}+\Qb \bigg)^{-1}\frac{\bSigma}{1+s_0}\bigg( \frac{\bSigma}{1+s_0}+\Qb \bigg)^{-1}\Ab\nonumber\\
    &=(1+s_0)\bigg[\tr\bigg( \frac{\bSigma}{1+s_0}+\Qb \bigg)^{-1}\Ab-\tr\bigg( \frac{\bSigma}{1+s_0}+\Qb \bigg)^{-1}\Qb\bigg( \frac{\bSigma}{1+s_0}+\Qb \bigg)^{-1}\Ab \bigg].\label{eq:gSigma2}
\end{align}
Combing \eqref{eq:gSigma2} with $\tr\big(\hat{\bSigma}+\Qb\big)^{-1}\Ab-\tr\big( \bSigma/(1+s_0)+\Qb \big)^{-1}\Ab~\cvas~0$ in Proposition~\ref{prop:prop1proof}, we successfully obtain another estimation of $g(\bSigma)$:
\begin{align}
    &g(\bSigma)-(1+s_0)\bigg[\tr\bigg( \hat{\bSigma}+\Qb \bigg)^{-1}\Ab-\tr\bigg( \frac{\bSigma}{1+s_0}+\Qb \bigg)^{-1}\Qb\bigg( \frac{\bSigma}{1+s_0}+\Qb \bigg)^{-1}\Ab \bigg]\nonumber\\
    &\qquad\qquad=g(\bSigma)-(1+s_0)\bigg[\hat{T}_1-\tr\bigg( \frac{\bSigma}{1+s_0}+\Qb \bigg)^{-1}\Qb\bigg( \frac{\bSigma}{1+s_0}+\Qb \bigg)^{-1}\Ab  \bigg]~\cvas~0.\label{eq:gSigma3}
\end{align}
Here, we remind the readers that $\hat{T}_1=\tr\big(\hat{\bSigma}+\Qb\big)^{-1}\Ab$. From \eqref{eq:gsigma1} and \eqref{eq:gSigma3}, we could successfully eliminate $\tr\big( \bSigma/(1+s_0)+\Qb \big)^{-1}\Qb\big( \bSigma/(1+s_0)+\Qb \big)^{-1}\Ab$ and achieve that
\begin{align*}
&(1+s_0)\underbrace{\bigg\{\hat{T}_2-\tr\bigg( \frac{\bSigma}{1+s_0}+\Qb \bigg)^{-1}\Qb\bigg( \frac{\bSigma}{1+s_0}+\Qb \bigg)^{-1}\Ab +\frac{s_{1,\Qb}}{(1+s_0)^2}g(\bSigma)\bigg\}}_{\eqref{eq:gsigma1}}\\
&\qquad\qquad + \underbrace{g(\bSigma)-(1+s_0)\bigg[\hat{T}_1-\tr\bigg( \frac{\bSigma}{1+s_0}+\Qb \bigg)^{-1}\Qb\bigg( \frac{\bSigma}{1+s_0}+\Qb \bigg)^{-1}\Ab \bigg]}_{\eqref{eq:gSigma3}}~\cvas~0. 
\end{align*}
Therefore $\tr\big( \bSigma/(1+s_0)+\Qb \big)^{-1}\Qb\big( \bSigma/(1+s_0)+\Qb \big)^{-1}\Ab$ is eliminated and we conclude that $\frac{(1+s_0)^2(\hat{T}_1- \hat{T}_2)}{1+s_0+ s_{1,\Qb}}-g(\bSigma)~\cvas~0$. Similar to the proof of Lemma~\ref{lemma:stronger_type_convergence}, when $\Ab$ is semi-positive definite, we have $\frac{(1+s_0)^2(\hat{T}_1- \hat{T}_2)/g(\bSigma)}{1+s_0+ s_{1,\Qb}}~\cvas~1$.

\subsection{Proof of Proposition~\ref{prop:converge_s}}
\label{appsec:proofprop3}
We first recall that $ f_1(\hat{\bSigma})=\frac{1}{p}\tr\Qb(\hat{\bSigma}+\Qb)^{-1}, f_2(\hat{\bSigma})=\frac{1}{p}\tr(\hat{\bSigma}+\Qb)^{-1}\Qb(\hat{\bSigma}+\Qb)^{-1}\Qb$, and the definitions of $s_0$, $s_{1,\bSigma}$ and $s_{1,\Qb}$ are given by the solutions of the following equations: 
\begin{align*}
    &s_0=\frac{c}{p}\tr~\bSigma\bigg( \frac{\bSigma}{1+s_0}+\Qb \bigg)^{-1},\\
    &s_{1,\bSigma}=\bigg(\frac{s_{1,\bSigma}}{(1+s_0)^2}-1\bigg)\cdot\frac{c}{p}\tr~\bSigma\bigg( \frac{\bSigma}{1+s_0}+\Qb\bigg)^{-1}\bSigma\bigg( \frac{\bSigma}{1+s_0}+\Qb\bigg)^{-1},\\
    &s_{1,\Qb}=\frac{c}{p}\tr~\bSigma\bigg( \frac{\bSigma}{1+s_0}+\Qb \bigg)^{-1}\bigg(\frac{s_{1,\Qb}}{(1+s_0)^2}\bSigma-\Qb\bigg)\bigg( \frac{\bSigma}{1+s_0}+\Qb \bigg)^{-1}.
\end{align*}
Then similar to the proof of Proposition~\ref{prop:prop1proof}, for $\Ab$ semi-positive definite we  have
\begin{align*}
    \tr \bigg(\frac{\Xb^\T\Xb}{n}+\Qb\bigg)^{-1}\Ab/\tr \bigg( \frac{\bSigma}{1+s_0}+\Qb\bigg)^{-1}\Ab~\cvas~1.
\end{align*}
Let  $\Ab=\Qb/p$,  we have the following  equation:
\begin{align}
    f_1(\hat{\bSigma})/\bigg(\frac{1}{p}\tr~\Qb\bigg( \frac{\bSigma}{1+s_0}+\Qb \bigg)^{-1}\bigg)~\cvas~1,\label{eq:Cs1}
\end{align}
Moreover, if we let $\Ab=\Qb/p$ in the Proposition~\ref{prop:converge_gSigma}, we can also have
\begin{align}
    \frac{(1+s_0)^2\big[f_1(\hat{\bSigma})-f_2(\hat{\bSigma})\big]}{1+s_0+s_{1,\Qb}}/\bigg(\frac{1}{p}\tr~\bigg( \frac{\bSigma}{1+s_0}+\Qb\bigg)^{-1}\bSigma\bigg( \frac{\bSigma}{1+s_0}+\Qb\bigg)^{-1}\Qb\bigg)~\cvas~1.\label{eq:Cs2}
\end{align}
With equations \eqref{eq:Cs1} and \eqref{eq:Cs2}, we are now able to prove the conclusions in Proposition~\ref{prop:converge_s}. 

\noindent\textbf{\underline{The estimation of $s_0$.}} For $s_0=\frac{c}{p}\tr~\bSigma\Big( \frac{\bSigma}{1+s_0}+\Qb \Big)^{-1}$, we could further simplify the equations as
    \begin{align*}
        s_0&=\frac{c}{p}\tr~\bSigma\Big( \frac{\bSigma}{1+s_0}+\Qb \Big)^{-1}=c(1+s_0)\cdot\bigg[1-\frac{1}{p}\tr~\Qb\bigg( \frac{\bSigma}{1+s_0}+\Qb \bigg)^{-1}\bigg].
    \end{align*}
    Combined the equation above with \eqref{eq:Cs1}, we have 
    \begin{align*}
        \bigg(1-\frac{s_0}{c(1+s_0)}\bigg)/f_1(\hat{\bSigma})~\cvas~1,\quad\text{which indicates}\quad\frac{c\big(1-  f_1(\hat{\bSigma})\big)/s_0}{1-c\big(1-  f_1(\hat{\bSigma})\big)}~\cvas~1.
    \end{align*}

\noindent \textbf{\underline{The estimation of $s_{1,\Qb}$.}} For $s_{1,\Qb}=\frac{c}{p}\tr~\bSigma\Big( \frac{\bSigma}{1+s_0}+\Qb \Big)^{-1}\Big(\frac{s_{1,\Qb}}{(1+s_0)^2}\bSigma-\Qb\Big)\Big( \frac{\bSigma}{1+s_0}+\Qb \Big)^{-1}$, we could further simplify the equations as
    \begin{align*}
        s_{1,\Qb}&=\frac{cs_{1,\Qb}}{p}\tr~\frac{\bSigma}{1+s_0}\bigg( \frac{\bSigma}{1+s_0}+\Qb \bigg)^{-1}\frac{\bSigma}{1+s_0}\bigg( \frac{\bSigma}{1+s_0}+\Qb \bigg)^{-1}\\
        &\qquad-\frac{c}{p}\tr~\Qb\bigg( \frac{\bSigma}{1+s_0}+\Qb\bigg)^{-1}\bSigma\bigg( \frac{\bSigma}{1+s_0}+\Qb\bigg)^{-1}\\
        &=\frac{cs_{1,\Qb}}{p}\tr~\bigg(\frac{\bSigma}{1+s_0}+\Qb-\Qb\bigg)\bigg( \frac{\bSigma}{1+s_0}+\Qb \bigg)^{-1}\frac{\bSigma}{1+s_0}\bigg( \frac{\bSigma}{1+s_0}+\Qb \bigg)^{-1}\\
        &\qquad-\frac{c}{p}\tr~\Qb\bigg( \frac{\bSigma}{1+s_0}+\Qb\bigg)^{-1}\bSigma\bigg( \frac{\bSigma}{1+s_0}+\Qb\bigg)^{-1}\\
        &=\frac{cs_{1,\Qb}}{p}\tr~\frac{\bSigma}{1+s_0}\bigg( \frac{\bSigma}{1+s_0}+\Qb \bigg)^{-1}-\frac{cs_{1,\Qb}}{p}\tr~\Qb\bigg( \frac{\bSigma}{1+s_0}+\Qb \bigg)^{-1}\frac{\bSigma}{1+s_0}\bigg( \frac{\bSigma}{1+s_0}+\Qb \bigg)^{-1}\\
        &\qquad-\frac{c}{p}\tr~\Qb\bigg( \frac{\bSigma}{1+s_0}+\Qb\bigg)^{-1}\bSigma\bigg( \frac{\bSigma}{1+s_0}+\Qb\bigg)^{-1}\\
        &=\frac{cs_{1,\Qb}}{p}\tr~\frac{\bSigma}{1+s_0}\bigg( \frac{\bSigma}{1+s_0}+\Qb \bigg)^{-1}\\
        &\qquad-\bigg(\frac{cs_{1,\Qb}}{p(1+s_0)}+\frac{c}{p}\bigg)\cdot\tr~\Qb\bigg( \frac{\bSigma}{1+s_0}+\Qb\bigg)^{-1}\bSigma\bigg( \frac{\bSigma}{1+s_0}+\Qb\bigg)^{-1}\\
        &=cs_{1,\Qb}\bigg[1-\frac{1}{p}\tr~\Qb\bigg( \frac{\bSigma}{1+s_0}+\Qb \bigg)^{-1}\bigg]\\
        &\qquad-\bigg(\frac{cs_{1,\Qb}}{p(1+s_0)}+\frac{c}{p}\bigg)\cdot\tr~\Qb\bigg( \frac{\bSigma}{1+s_0}+\Qb\bigg)^{-1}\bSigma\bigg( \frac{\bSigma}{1+s_0}+\Qb\bigg)^{-1}\\
        &=\frac{s_0s_{1,\Qb}}{1+s_0}-\bigg(\frac{cs_{1,\Qb}}{p(1+s_0)}+\frac{c}{p}\bigg)\cdot\tr~\Qb\bigg( \frac{\bSigma}{1+s_0}+\Qb\bigg)^{-1}\bSigma\bigg( \frac{\bSigma}{1+s_0}+\Qb\bigg)^{-1}.
    \end{align*}
    Here, the last equality is by the definition of $s_0$. From the last equality above we have
    \begin{align}
        s_{1,\Qb}=-(1+s_0+s_{1,\Qb})\cdot \frac{c}{p}\tr~\Qb\bigg( \frac{\bSigma}{1+s_0}+\Qb\bigg)^{-1}\bSigma\bigg( \frac{\bSigma}{1+s_0}+\Qb\bigg)^{-1}.\label{eq:s1Qlastequality}
    \end{align}
    Combined \eqref{eq:s1Qlastequality} with  \eqref{eq:Cs2}, we conclude that 
    \begin{align}
    \label{eq:S1I}
        \frac{-c(1+s_0)^2(f_1(\hat{\bSigma})-f_2(\hat{\bSigma}))}{s_{1,\Qb}}~\cvas~1.
    \end{align}
    Note also that $\bigg(1+\frac{c(1-f_1(\hat{\bSigma}))}{1-c\big(1-f_1(\hat{\bSigma})\big)}\bigg)/(1+s_0)~\cvas~1$ holds, replacing $(1+s_0)$ with $\frac{1}{1-c\big(1-f_1(\hat{\bSigma})\big)}$ into \eqref{eq:S1I} gives
    \begin{align*}
        c\cdot\frac{(f_2(\hat{\bSigma})-f_1(\hat{\bSigma}))/s_{1,\Qb}}{\big(c\big(1-f_1(\hat{\bSigma})\big)-1\big)^2}~\cvas~1.
    \end{align*}
    
\noindent\textbf{\underline{The estimation of $s_{1,\bSigma}$.}} For $s_{1,\bSigma}=\Big(\frac{s_{1,\bSigma}}{(1+s_0)^2}-1\Big)\cdot\frac{c}{p}\tr~\bSigma\Big( \frac{\bSigma}{1+s_0}+\Qb\Big)^{-1}\bSigma\Big( \frac{\bSigma}{1+s_0}+\Qb\Big)^{-1}$, we could further simplify the equations as 
\begin{align*}
    s_{1,\bSigma}&=\bigg(\frac{s_{1,\bSigma}}{(1+s_0)^2}-1\bigg)\cdot\frac{c}{p}\tr~\bSigma\bigg( \frac{\bSigma}{1+s_0}+\Qb\bigg)^{-1}\bSigma\bigg( \frac{\bSigma}{1+s_0}+\Qb\bigg)^{-1}\\
    &=\frac{c(s_{1,\bSigma}-(1+s_0)^2)}{p}\tr~\frac{\bSigma}{1+s_0}\bigg( \frac{\bSigma}{1+s_0}+\Qb\bigg)^{-1}\frac{\bSigma}{1+s_0}\bigg( \frac{\bSigma}{1+s_0}+\Qb\bigg)^{-1}\\
    &=c(s_{1,\bSigma}-(1+s_0)^2)\bigg(1-\frac{1}{p}\tr~\Qb\bigg( \frac{\bSigma}{1+s_0}+\Qb\bigg)^{-1}\bigg)\\
    &\qquad-\frac{c(s_{1,\bSigma}-(1+s_0)^2)}{p(1+s_0)}\cdot\tr~\Qb\bigg( \frac{\bSigma}{1+s_0}+\Qb\bigg)^{-1}\bSigma\bigg( \frac{\bSigma}{1+s_0}+\Qb\bigg)^{-1}\\
    &=\frac{c(s_{1,\bSigma}-(1+s_0)^2)}{1+s_0}\bigg(\frac{s_0}{c}-\frac{1}{p}\tr~\Qb\bigg( \frac{\bSigma}{1+s_0}+\Qb\bigg)^{-1}\bSigma\bigg( \frac{\bSigma}{1+s_0}+\Qb\bigg)^{-1}\bigg).
\end{align*}
Here, the last equality uses the fact $1-\frac{1}{p}\tr~\Qb\Big( \frac{\bSigma}{1+s_0}+\Qb\Big)^{-1}=\frac{s_0}{c(1+s_0)}$. Combined the results above with \eqref{eq:Cs1} and \eqref{eq:Cs2}, we conclude that 
    \begin{align*}
       c[s_{1,\bSigma}-(1+s_0)^2]\cdot\bigg(\frac{s_0}{c(1+s_0)}-(1+s_0)\cdot\frac{f_1(\hat{\bSigma})-f_2(\hat{\bSigma})}{1+s_0+s_{1,\Qb}}\bigg)/s_{1,\bSigma}~\cvas~1,
    \end{align*}
Furthermore, by \eqref{eq:S1I} and replacing $(1+s_0)$ with $\frac{1}{1-c\big(1-f_1(\hat{\bSigma})\big)}$ we have
    \begin{align*}
       c\bigg[s_{1,\bSigma}-\frac{1}{\big(1-c\big(1-f_1(\hat{\bSigma})\big)\big)^2}\bigg]\cdot\Bigg(1-f_1(\hat{\bSigma})-\frac{f_1(\hat{\bSigma})-f_2(\hat{\bSigma})}{1+\frac{c(f_1(\hat{\bSigma})-f_2(\hat{\bSigma}))}{c(1-f_1(\hat{\bSigma}))-1}}\Bigg)/s_{1,\bSigma}~\cvas~1.
    \end{align*}
Solve this linear equation with $s_{1,\bSigma}$ and we conclude that
\begin{align*}
    c\cdot\frac{-1+2f_1(\hat{\bSigma})-f_2(\hat{\bSigma})+c(f_1(\hat{\bSigma})-1)^2 }{\big(1+c(f_1(\hat{\bSigma})-1)\big)^4\cdot s_{1,\bSigma}}~\cvas~1.
\end{align*}

Note that $s_0,|s_{1,\bSigma}|$ and $|s_{1,\Qb}|$ are bounded by Lemma~\ref{lemma:constant_level}, direct calculation gives 
 \begin{align*}
    \frac{1+s_0+s_{1,\Qb}}{(1+s_0)^2-s_{1,\bSigma}} /\bigg(1-\frac{c}{p}\tr\hbSigma\big(\hbSigma+\Qb\big)^{-1} \bigg)^2\cvas1.
\end{align*}
We complete the proof of Proposition~\ref{prop:converge_s}.

\section{Proof of Theorem~\ref{thm:frontier}}
\label{sec:proofthmfrontier}
In this section, we give the proof of Theorem~\ref{thm:frontier}. Before the proof, we give the following lemmas.
 The following lemma summarizes some results from Section~\ref{sec:preliminary} and Section~\ref{sec:appendix_proof_thm_bound}.

\begin{lemma}
\label{lemma:conclusions_previous}
Under the condition of Theorem~\ref{thm:frontier}, suppose that $\Ab$ is semi-positive, it holds that
\begin{align*}
    &\frac{(1+s_0)^2}{(1+s_0)^2-s_{1,\bSigma}}\cdot \frac{\tr(\hat{\bSigma}+\Qb)^{-1}\bSigma(\hat{\bSigma}+\Qb)^{-1}\Ab}{\tr\Big( \frac{\bSigma}{1+s_0}+\Qb \Big)^{-1}\bSigma\Big( \frac{\bSigma}{1+s_0}+\Qb \Big)^{-1}\Ab}~\cvas~1.\\
    &\frac{(1+s_0)^2}{1+s_0+s_{1,\Qb}}\cdot \frac{\tr(\hat{\bSigma}+\Qb)^{-1}\hbSigma(\hat{\bSigma}+\Qb)^{-1}\Ab}{\tr\Big( \frac{\bSigma}{1+s_0}+\Qb \Big)^{-1}\bSigma\Big( \frac{\bSigma}{1+s_0}+\Qb \Big)^{-1}\Ab}~\cvas~1.
\end{align*}
\end{lemma}
\begin{proof}[\textbf{Proof of Lemma~\ref{lemma:conclusions_previous}}]
The first result comes from Proposition~\ref{prop:converge_Tn1Tn2}, and the second result comes from Proposition~\ref{prop:converge_gSigma}.
\end{proof}
Lemma~\ref{lemma:constant_level} and \ref{lemma:conclusions_previous} jointly show that  $\tr(\hat{\bSigma}+\Qb)^{-1}\bSigma(\hat{\bSigma}+\Qb)^{-1}\Ab$, $\tr(\hat{\bSigma}+\Qb)^{-1}\hbSigma(\hat{\bSigma}+\Qb)^{-1}\Ab$ and $\tr\Big( \frac{\bSigma}{1+s_0}+\Qb \Big)^{-1}\bSigma\Big( \frac{\bSigma}{1+s_0}+\Qb \Big)^{-1}\Ab$ are in the same order.The following lemma indicates that  $\tr\Big( \frac{\bSigma}{1+s_0}+\Qb \Big)^{-1}\bSigma\Big( \frac{\bSigma}{1+s_0}+\Qb \Big)^{-1}\Ab=\Theta\Big( \tr\Big( \frac{\bSigma}{1+s_0}+\Qb \Big)^{-1}\Ab\Big)$. This equality will further simplify our analysis in the proof of Theorem~\ref{thm:frontier}.
\begin{lemma}
    \label{lemma:same_order_2_to_1}
Under the condition of Theorem~\ref{thm:frontier}, suppose that $\Ab$ is semi-positive, it holds that
\begin{align*}
\tr\bigg( \frac{\bSigma}{1+s_0}+\Qb \bigg)^{-1}\bSigma\bigg( \frac{\bSigma}{1+s_0}+\Qb \bigg)^{-1}\Ab=\Theta\bigg( \tr\bigg( \frac{\bSigma}{1+s_0}+\Qb \bigg)^{-1}\Ab\bigg).
\end{align*}
Moreover, if $\Ab=\bxi\bxi^\top$, then 
\begin{align*}
     \tr\bigg( \frac{\bSigma}{1+s_0}+\Qb \bigg)^{-1}\Ab=\Theta\big( \|\bSigma^{-\frac{1}{2}}\bxi\|_2^2\big).
\end{align*}
\end{lemma}
The proof of Lemma~\ref{lemma:same_order_2_to_1} can be found in Section~\ref{sec:prooflemma_same_order_2_to_1}. We give an additional lemma, which is applied to handle the dependence of $\Ab$ ($\Ab$ is rank-one and semi-positive in Theorem~\ref{thm:frontier}) and $\hbSigma$ in the proof of Theorem~\ref{thm:frontier}.
\begin{lemma}
\label{lemma:analysis_dependence}
Recall that
\begin{align*}
    \cA_{rr}=\rb^\top \bigg(\frac{\bSigma}{1+s_0}+\Qb\bigg)^{-1}\rb,\quad \cA_{r1}=\rb^\top \bigg(\frac{\bSigma}{1+s_0}+\Qb\bigg)^{-1}\one,\quad \cA_{11}=\one^\top \bigg(\frac{\bSigma}{1+s_0}+\Qb\bigg)^{-1}\one.
\end{align*}
Define 
\begin{align*}
             &\alpha_n=\rb^\top \big(\hbSigma+\Qb \big)^{-1}\rb-\mu_0\rb^\top \big(\hbSigma+\Qb \big)^{-1}\one,\quad \alpha_0=\cA_{rr}-\mu_0 \cA_{r1},\\
    &\beta_n=\mu_0\one^\top \big(\hbSigma+\Qb \big)^{-1}\one-\rb^\top \big(\hbSigma+\Qb \big)^{-1}\one,\quad \beta_0=\mu_0\cA_{11}-\cA_{r1}.
\end{align*}
Given some large enough constant $C$, then under the condition of Theorem~\ref{thm:frontier}, the following properties hold:
\begin{enumerate}
    \item If $\mu_0\leq C\sqrt{\cA_{rr}/\cA_{11}}$, it holds that $\frac{\alpha_n-\alpha_0}{\cA_{rr}}~\cvas~0,\quad
    \frac{\beta_n-\beta_0}{\sqrt{\cA_{rr}\cdot\cA_{11}}}~\cvas~0$.
    \item If  $\mu_0\geq C\sqrt{\cA_{rr}/\cA_{11}}$, it holds that $     \frac{\alpha_n-\alpha_0}{\mu_0\sqrt{\cA_{rr}\cA_{11}}}~\cvas~0,\quad
    \frac{\beta_n-\beta_0}{\mu_0\cA_{11}}~\cvas~0$.
\end{enumerate}
An additional conclusion is $D=\Theta(\cA_{rr}\cA_{11})$ where $D$ is defined in \eqref{eq:def_gh}.
\end{lemma}
The proof of Lemma~\ref{lemma:analysis_dependence} can be found in Section~\ref{sec:prooflemma_analysis_dependence}. With the lemmas above, we are now ready to prove Theorem~\ref{thm:frontier}.
\begin{proof}[\textbf{Proof of Theorem~\ref{thm:frontier}}]
We first give the proof of ratio consistency. Recall that
\begin{align}
\label{eq:def_gh_appendix}
\begin{array}{l@{~}l@{~}l@{~}l@{~}l}
&\gb=D^{-1}\big[B(\hbSigma+\Qb)^{-1}\one-A(\hbSigma+\Qb)^{-1}\rb\big],~\hb=D^{-1}\big[C(\hbSigma+\Qb)^{-1}\rb-A(\hbSigma+\Qb)^{-1}\one\big],\\
&A=\rb^\top(\hbSigma+\Qb)^{-1}\one,~B=\rb^\top(\hbSigma+\Qb)^{-1}\rb,~C=\one^\top(\hbSigma+\Qb)^{-1}\one,~D=BC-A^2,
\end{array}
\end{align}
and 
\begin{equation}
    \label{eq:frontierequation_appendix}
    \begin{split}
             &\sigma_0^2= \wb^{*\top}\bSigma\wb^*=(\gb+\mu_0\cdot\hb)^\top \bSigma (\gb+\mu_0\cdot\hb),\\ &\hat{\sigma}^2= (\gb+\mu_0\cdot\hb)^\top \hbSigma (\gb+\mu_0\cdot\hb)/(1-c/p\cdot \tr \hbSigma(\hbSigma+\Qb)^{-1})^2.
    \end{split}
\end{equation}
The first conclusion in Theorem~\ref{thm:frontier} states that $\hat{\sigma}^2/\sigma_0^2~\cvas~1$. By direct calculation of \eqref{eq:frontierequation_appendix}, we observe that
\begin{equation}
    \label{eq:frontierequation_appendix_fix}
    \begin{split}
             &\sigma_0^2= (\alpha_n \one+\beta_n\rb)^\top (\hbSigma+\Qb)^{-1}\bSigma(\hbSigma+\Qb)^{-1}(\alpha_n \one+\beta_n\rb)/D^2,\\ 
             &\hat{\sigma}^2=\frac{(\alpha_n \one+\beta_n\rb)^\top (\hbSigma+\Qb)^{-1}\hbSigma(\hbSigma+\Qb)^{-1}(\alpha_n \one+\beta_n\rb)}{D^2(1-c/p\cdot \tr \hbSigma(\hbSigma+\Qb)^{-1})^2}.
    \end{split}
\end{equation}
Here, we have
\begin{equation}
    \label{eq:alpha_beta_appendix}
    \begin{split}
             &\alpha_n=\rb^\top \big(\hbSigma+\Qb \big)^{-1}\rb-\mu_0\rb^\top \big(\hbSigma+\Qb \big)^{-1}\one,\\
    &\beta_n=\mu_0\one^\top \big(\hbSigma+\Qb \big)^{-1}\one-\rb^\top \big(\hbSigma+\Qb \big)^{-1}\one.
    \end{split}
\end{equation}
It is clear that $\alpha_n \one+\beta_n\rb$ is dependent with $\hbSigma$ and is not deterministic, hence we can not directly apply the results from Theorem~\ref{thm:main_theorem}. Fortunately, the dependence of $\alpha_n \one+\beta_n\rb$ and $\hbSigma$ is only due to the two numbers $\alpha_n$ and $\beta_n$. We further define 
\begin{align}
\label{eq:def_alphabeta0}
    \alpha_0=\cA_{rr}-\mu_0\cA_{r1},\quad \beta_0=\mu_0\cA_{11}-\cA_{r1}.
\end{align}
In the following steps of the proof of the ratio consistency, we mainly focus on the proof of the two key results below:
\begin{align}
\label{eq:target_frontier}
    \frac{(1+s_0)^2}{(1+s_0)^2-s_{1,\bSigma}}\cdot \frac{(\alpha_n \one+\beta_n\rb)^\top (\hbSigma+\Qb)^{-1}\bSigma(\hbSigma+\Qb)^{-1}(\alpha_n \one+\beta_n\rb)}{(\alpha_0 \one+\beta_0\rb)^\top (\bSigma/(1+s_0)+\Qb)^{-1}\bSigma(\bSigma/(1+s_0)+\Qb)^{-1}(\alpha_0 \one+\beta_0\rb)}~\cvas~1,
\end{align}
and 
\begin{align}
\label{eq:target_frontier1}
    \frac{(1+s_0)^2}{1+s_0+s_{1,\Qb}}\cdot \frac{(\alpha_n \one+\beta_n\rb)^\top (\hbSigma+\Qb)^{-1}\hbSigma(\hbSigma+\Qb)^{-1}(\alpha_n \one+\beta_n\rb)}{(\alpha_0 \one+\beta_0\rb)^\top (\bSigma/(1+s_0)+\Qb)^{-1}\bSigma(\bSigma/(1+s_0)+\Qb)^{-1}(\alpha_0 \one+\beta_0\rb)}~\cvas~1.
\end{align}

\medskip
\noindent\textbf{\underline{Proof of \eqref{eq:target_frontier}}}. To prove \eqref{eq:target_frontier}, 
we first give the order of the dominator $(\alpha_0 \one+\beta_0\rb)^\top (\bSigma/(1+s_0)+\Qb)^{-1}\bSigma(\bSigma/(1+s_0)+\Qb)^{-1}(\alpha_0 \one+\beta_0\rb) $. From Lemma~\ref{lemma:same_order_2_to_1}, it holds that
\begin{align}
    &(\alpha_0 \one+\beta_0\rb)^\top \bigg(\frac{\bSigma}{1+s_0}+\Qb\bigg)^{-1}\bSigma\bigg(\frac{\bSigma}{1+s_0}+\Qb\bigg)^{-1}(\alpha_0 \one+\beta_0\rb)\nonumber\\
    &\quad =\Theta \bigg(  (\alpha_0 \one+\beta_0\rb)^\top \bigg(\frac{\bSigma}{1+s_0}+\Qb\bigg)^{-1}(\alpha_0 \one+\beta_0\rb)\bigg).\label{eq:frontier_simplify}
\end{align}
Recall the definition of $\alpha_0$ and $\beta_0$ in \eqref{eq:def_alphabeta0}, we have
\begin{align}
      (\alpha_0 \one+\beta_0\rb)^\top \bigg(\frac{\bSigma}{1+s_0}+\Qb\bigg)^{-1}(\alpha_0 \one+\beta_0\rb)=(\cA_{11}\cA_{rr}-\cA_{r1}^2)\cdot\big(\mu_0(\cA_{11}\mu_0-2\cA_{r1})+\cA_{rr}\big).\label{eq:frontier_simplify_calculate}
\end{align}
By Assumption~\ref{assump:assump5}, we can easily get that $\cA_{11}\cA_{rr}-\cA_{r1}^2=\Theta(\cA_{11}\cA_{rr})$. 
We separate the proof of \eqref{eq:target_frontier} into two cases: $\mu_0\leq C\sqrt{\cA_{rr}/\cA_{11}}$ and $\mu_0\geq C\sqrt{\cA_{rr}/\cA_{11}}$ for some sufficiently large constant $C>0$. Before the proof of \eqref{eq:target_frontier} and \eqref{eq:target_frontier1}, we define several terms.
\begin{equation}
    \label{eq:def_J}
    \begin{split}
        \cJ_0:&=(\alpha_0 \one+\beta_0\rb)^\top (\bSigma/(1+s_0)+\Qb)^{-1}\bSigma(\bSigma/(1+s_0)+\Qb)^{-1}(\alpha_0 \one+\beta_0\rb),\\
    \cJ_1:&=(\alpha_n \one+\beta_n\rb)^\top (\hbSigma+\Qb)^{-1}\bSigma(\hbSigma+\Qb)^{-1}(\alpha_n \one+\beta_n\rb),\\
    \cJ_2:&=(\alpha_0 \one+\beta_0\rb)^\top (\hbSigma+\Qb)^{-1}\bSigma(\hbSigma+\Qb)^{-1}(\alpha_0 \one+\beta_0\rb).
    \end{split}
\end{equation}
It is clear that $\cJ_1/D^2=\sigma_0^2$.
Also from \eqref{eq:frontier_simplify} and \eqref{eq:frontier_simplify_calculate}, we have 
\begin{align}
\label{eq:J0order}
\cJ_0=\Theta(\cA_{11}\cA_{rr})\cdot (\mu_0(\cA_{11}\mu_0-2\cA_{r1})+\cA_{rr}).
\end{align}

\noindent\textbf{Case 1:} If $\mu_0\leq C\sqrt{\cA_{rr}/\cA_{11}}$ holds, we have 
\begin{align*}
    \mu_0(\cA_{11}\mu_0-2\cA_{r1})+\cA_{rr}\leq C^2(\cA_{rr}+\cA_{r1}\cdot\sqrt{\cA_{rr}/\cA_{11}})=O(\cA_{rr}).
\end{align*}
Here, the last equality is by $\cA_{r1}^2/(\cA_{rr}\cdot\cA_{11})\leq1$. For the minimum value of $\mu_0(\cA_{11}\mu_0-2\cA_{r1})+\cA_{rr}$, it takes value when $\mu_0=\cA_{r1}/\cA_{11}$, therefore we have
\begin{align*}
    \mu_0(\cA_{11}\mu_0-2\cA_{r1})+\cA_{rr}\geq \frac{\cA_{11}\cA_{rr}-\cA_{r1}^2}{\cA_{11}}\geq \Omega(\cA_{rr}).
\end{align*}
Hence we conclude that when $\mu_0=O(\sqrt{\cA_{rr}/\cA_{11}})$, $\mu_0(\cA_{11}\mu_0-2\cA_{r1})+\cA_{rr}=\Theta(\cA_{rr})$. 
We have $\big(\mu_0(\cA_{11}\mu_0-2\cA_{r1})+\cA_{rr}\big)=\Theta(\cA_{rr})$. Combined this with \eqref{eq:def_J} and \eqref{eq:J0order} we have
\begin{align}
    \frac{\cJ_0}{\cA_{rr}^2\cA_{11}}=\Theta(1).\label{eq:J0_fact}
\end{align}

We claim that the following facts hold.
\begin{align*}
   1.~\frac{(1+s_0)^2}{(1+s_0)^2-s_{1,\bSigma}}=\Theta(1),\quad 2.~\frac{(1+s_0)^2}{(1+s_0)^2-s_{1,\bSigma}}\cdot \frac{\cJ_2}{\cJ_0}~\cvas~1,
\end{align*}
\begin{align*}
   3.~\frac{\cJ_0}{\cA_{rr}^2\cA_{11}}=\Theta(1),\quad4.~\frac{\cJ_1-\cJ_2}{\cA_{rr}^2\cA_{11}}~\cvas~0.
\end{align*}
Then Equation~\eqref{eq:target_frontier}, which is equivalent to $\frac{(1+s_0)^2}{(1+s_0)^2-s_{1,\bSigma}}\cdot \cJ_1/\cJ_0~\cvas~1$ can be directly concluded from the above facts. Moreover,  by $\cJ_1/D^2=\sigma_0^2$ and $D=\Theta(\cA_{11}\cA_{rr})$ in Lemma~\ref{lemma:analysis_dependence}, we have
\begin{align*}
    \sigma_0^2=\Theta(1/\cA_{11}).
\end{align*}

It is easy to verify that the first fact comes from Lemma~\ref{lemma:constant_level}, the second fact comes from Lemma~\ref{lemma:conclusions_previous} and the third fact comes from \eqref{eq:J0_fact}. Therefore, to prove \eqref{eq:target_frontier}, it only remains to prove the last fact, that is $(\cJ_1-\cJ_2)/\cA_{rr}^2\cA_{11}~\cvas~0$. By Lemma~\ref{lemma:analysis_dependence}, it is not hard to see that 
\begin{align*}
    \frac{\alpha_n}{\cA_{rr}}-\frac{\alpha_0}{\cA_{rr}}~\cvas~0,\quad  \frac{\beta_n}{\sqrt{\cA_{rr}\cdot\cA_{11}}}-\frac{\beta_0}{\sqrt{\cA_{rr}\cdot\cA_{11}}}~\cvas~0.
\end{align*}
We apply perturbation analysis to prove $(\cJ_1-\cJ_2)/\cA_{rr}^2\cA_{11}~\cvas~0$. Define $\cJ^\circ=(\alpha_0 \one+\beta_n\rb)^\top (\hbSigma+\Qb)^{-1}\bSigma(\hbSigma+\Qb)^{-1}(\alpha_0 \one+\beta_n\rb)$, $\cJ_1-\cJ_2$ can be separated into $\cJ_1-\cJ^\circ $ and $\cJ^\circ-\cJ_2$. We only need to prove that  the two terms $(\cJ_1-\cJ^\circ )/(\cA_{rr}^2\cA_{11})$ and $(\cJ^\circ-\cJ_2)/(\cA_{rr}^2\cA_{11})$ both converge to $0$. We need the following equations. Similar to the proof of Lemma~\ref{lemma:same_order_2_to_1}, we have
\begin{align*}
    &\one^\top (\hbSigma+\Qb)^{-1}\bSigma(\hbSigma+\Qb)^{-1}\one=\Theta(\cA_{11}),\\
    &\big|\rb^\top (\hbSigma+\Qb)^{-1}\bSigma(\hbSigma+\Qb)^{-1}\one\big|=O(|\cA_{r1}|),\\
    &\rb^\top (\hbSigma+\Qb)^{-1}\bSigma(\hbSigma+\Qb)^{-1}\rb=\Theta(\cA_{rr}),\\
    &\frac{|\alpha_0|}{\cA_{rr}}=\frac{|\cA_{rr}-\mu_0\cA_{r1}|}{\cA_{rr}}\leq 1+\mu_0\frac{|\cA_{r1}|}{\cA_{rr}}\leq 1+C,\\
    &\frac{|\beta_0|}{\sqrt{\cA_{rr}\cdot\cA_{11}}}=\frac{|\mu_0\cA_{11}-\cA_{r1}|}{\sqrt{\cA_{rr}\cdot\cA_{11}}}\leq 1+\mu_0\frac{|\cA_{11}|}{\sqrt{\cA_{rr}\cdot\cA_{11}}}\leq 1+C.
\end{align*}
Combining the conclusions above with Lemma~\ref{lemma:analysis_dependence}, we can see 
\begin{align*}
    \frac{\cJ_1-\cJ^\circ}{\cA_{rr}^2\cA_{11}}&=O\bigg(\frac{\alpha_n^2-\alpha_0^2}{\cA_{rr}^2}\bigg)+O\bigg(\frac{\alpha_n-\alpha_0}{\cA_{rr}}\bigg)=o(1),\\
    \frac{\cJ^\circ-\cJ_2}{\cA_{rr}^2\cA_{11}}&=O\bigg(\frac{\beta_n^2-\beta_0^2}{\cA_{rr}\cA_{11}}\bigg)+O\bigg(\frac{\beta_n-\beta_0}{\sqrt{\cA_{rr}\cA_{11}}}\bigg)=o(1).
\end{align*}
Here, the convergence also utilize the fact that $|\alpha_n+\alpha_0|/\cA_{rr}\leq 3(1+C)$, $|\beta_n+\beta_0|/\sqrt{\cA_{rr}\cdot\cA_{11}}\leq 3(1+C)$. 
We conclude that $\cJ_1-\cJ_2~\cvas~0$, hence we have $\frac{(1+s_0)^2}{(1+s_0)^2-s_{1,\bSigma}}\cdot \frac{\cJ_1}{\cJ_0}~\cvas~1$ which completes the proof of \eqref{eq:target_frontier}. Similarly we can prove \eqref{eq:target_frontier1}. We hence complete the proof of \eqref{eq:target_frontier} when $\mu_0\leq C\sqrt{\cA_{rr}/\cA_{11}}$.

\noindent\textbf{Case 2:} If $\mu_0\geq C\sqrt{\cA_{rr}/\cA_{11}}$ holds, we  have 
\begin{align*}
    \cA_{11}\mu_0-2\cA_{r1} &= \cA_{11}\mu_0/2+\cA_{11}\mu_0/2-2\cA_{r1}\\
    &\geq \cA_{11}\mu_0/2+\frac{C}{2}\sqrt{\cA_{11}\cA_{rr}}-\sqrt{\cA_{11}\cA_{rr}}\geq \cA_{11}\mu_0/2,
\end{align*}
where the first inequality is by condition $\mu_0\geq C\sqrt{\cA_{rr}/\cA_{11}}$ and the second inequality is by $C$ large enough. Similarly, we can conclude that $ \cA_{11}\mu_0-2\cA_{r1}\leq 2\mu_0\cA_{11}$.  Hence it holds that
\begin{align*}
   \cA_{11}\mu_0^2/2 \leq \mu_0(\cA_{11}\mu_0-2\cA_{r1})+\cA_{rr}\leq 3\cA_{11}\mu_0^2.
\end{align*}
Here, we use the fact $\cA_{rr}\leq \cA_{11}\mu_0^2$. Combined the inequalities above with \eqref{eq:frontier_simplify_calculate}, we have 
\begin{align}
    (\alpha_0 \one+\beta_0\rb)^\top \bigg(\frac{\bSigma}{1+s_0}+\Qb\bigg)^{-1}(\alpha_0 \one+\beta_0\rb)=\Theta\big( \mu_0^2\cA_{rr}\cA_{11}^2\big).\label{eq:frontier_simplify_calculate2}
\end{align}
Combined this with \eqref{eq:def_J} and \eqref{eq:J0order} we have
\begin{align}
    \frac{\cJ_0}{\mu_0^2\cA_{rr}\cA_{11}^2}=\Theta(1).\label{eq:J0_fact1}
\end{align}
Here, the second equality uses \eqref{eq:frontier_simplify} and \eqref{eq:frontier_simplify_calculate2}.
Similar to the case when $\mu_0\leq C\sqrt{\cA_{rr}/\cA_{11}}$, we claim that the following facts hold.
\begin{align*}
   1.~\frac{(1+s_0)^2}{(1+s_0)^2-s_{1,\bSigma}}=\Theta(1),\quad 2.~\frac{(1+s_0)^2}{(1+s_0)^2-s_{1,\bSigma}}\cdot \frac{\cJ_2}{\cJ_0}~\cvas~1,
\end{align*}
\begin{align*}
   3.~\frac{\cJ_0}{\mu_0^2\cA_{rr}\cA_{11}^2}=\Theta(1),\quad4.~\frac{\cJ_1-\cJ_2}{\mu_0^2\cA_{rr}\cA_{11}^2}~\cvas~0.
\end{align*}
Then Equation \eqref{eq:target_frontier}, which is equivalent to $\frac{(1+s_0)^2}{(1+s_0)^2-s_{1,\bSigma}}\cdot \cJ_1/\cJ_0~\cvas~1$, can be directly concluded from the above facts. Moreover,  by $\cJ_1/D^2=\sigma_0^2$ and $D=\Theta(\cA_{11}\cA_{rr})$ in Lemma~\ref{lemma:analysis_dependence}, we have
\begin{align*}
    \sigma_0^2=\Theta(\mu_0^2/\cA_{rr}).
\end{align*}

It is easy to verify that the first fact comes from Lemma~\ref{lemma:constant_level}, the second fact comes from Lemma~\ref{lemma:conclusions_previous} and the third fact comes from \eqref{eq:J0_fact1}. Therefore, to prove \eqref{eq:target_frontier}, it only remains to prove that $(\cJ_1-\cJ_2)/(\mu_0^2\cA_{rr}\cA_{11}^2)~\cvas~0$. By Lemma~\ref{lemma:analysis_dependence}, it is easy to see that 
\begin{align*}
    \frac{\alpha_n}{\mu_0\sqrt{\cA_{rr}\cA_{11}}}-\frac{\alpha_0}{\mu_0\sqrt{\cA_{rr}\cA_{11}}}~\cvas~0,\quad  \frac{\beta_n}{\mu_0\cA_{11}}-\frac{\beta_0}{\mu_0\cA_{11}}~\cvas~0.
\end{align*}
We apply perturbation analysis to prove $(\cJ_1-\cJ_2)/(\mu_0^2\cA_{rr}\cA_{11}^2)~\cvas~0$. Define $\cJ^\circ=(\alpha_0 \one+\beta_n\rb)^\top (\hbSigma+\Qb)^{-1}\bSigma(\hbSigma+\Qb)^{-1}(\alpha_0 \one+\beta_n\rb)$, $\cJ_1-\cJ_2$ can be separated into $\cJ_1-\cJ^\circ $ and $\cJ^\circ-\cJ_2$. We only need to prove that  the two terms $(\cJ_1-\cJ^\circ )/(\mu_0^2\cA_{rr}\cA_{11}^2)$ and $(\cJ^\circ-\cJ_2)/(\mu_0^2\cA_{rr}\cA_{11}^2)$ both converge to $0$. 
We  have the following equations.
\begin{align*}
    &\one^\top (\hbSigma+\Qb)^{-1}\bSigma(\hbSigma+\Qb)^{-1}\one=\Theta(\cA_{11}),\\
    &\big|\rb^\top (\hbSigma+\Qb)^{-1}\bSigma(\hbSigma+\Qb)^{-1}\one\big|=O(|\cA_{r1}|),\\
    &\rb^\top (\hbSigma+\Qb)^{-1}\bSigma(\hbSigma+\Qb)^{-1}\rb=\Theta(\cA_{rr}),\\
      & \frac{|\alpha_0|}{\mu_0\sqrt{\cA_{rr}\cA_{11}}}=\frac{|\cA_{rr}-\mu_0\cA_{r1}|}{\mu_0\sqrt{\cA_{rr}\cA_{11}}}\leq \frac{1}{C}+\mu_0\frac{|\cA_{r1}|}{\mu_0\sqrt{\cA_{rr}\cA_{11}}}\leq \frac{1}{C}+1,\\
    &\frac{|\beta_0|}{\mu_0\cA_{11}}=\frac{|\mu_0\cA_{11}-\cA_{r1}|}{\mu_0\cA_{11}}\leq \frac{1}{C}+\mu_0\frac{|\cA_{11}|}{\mu_0\cA_{11}}\leq \frac{1}{C}+1.
\end{align*}
Combined the conclusions above with Lemma~\ref{lemma:analysis_dependence}, we can see 
\begin{align*}
    \frac{\cJ_1-\cJ^\circ}{\mu_0^2\cA_{rr}\cA_{11}^2}&=O\bigg(\frac{\alpha_n^2-\alpha_0^2}{\mu_0^2\cA_{rr}\cA_{11}}\bigg)+O\bigg(\frac{\alpha_n-\alpha_0}{\mu_0\sqrt{\cA_{rr}\cA_{11}}}\bigg)=o(1),\\
    \frac{\cJ^\circ-\cJ_2}{\mu_0^2\cA_{rr}\cA_{11}^2}&=O\bigg(\frac{\beta_n^2-\beta_0^2}{\mu_0^2\cA_{11}^2}\bigg)+O\bigg(\frac{\beta_n-\beta_0}{\mu_0\cA_{11}}\bigg)=o(1).
\end{align*}
We conclude that $(\cJ_1-\cJ_2)/(\mu_0^2\cA_{rr}\cA_{11}^2)~\cvas~0$, hence we have $\frac{(1+s_0)^2}{(1+s_0)^2-s_{1,\bSigma}}\cdot \frac{\cJ_1}{\cJ_0}~\cvas~1$, which completes the proof of \eqref{eq:target_frontier}. We hence complete the proof of \eqref{eq:target_frontier} when $\mu_0\geq C\sqrt{\cA_{rr}/\cA_{11}}$. 

\medskip
\noindent\textbf{\underline{Proof of \eqref{eq:target_frontier1}}}. We can prove \eqref{eq:target_frontier1} similar to the proof of \eqref{eq:target_frontier}, with a few small differences. We need to prove that the constant $(1+s_0)^2/(1+s_0+s_{1,\Qb})=\Theta(1)$, which can be easily concluded from Lemma~\ref{lemma:constant_level}, and we simply replace the matrix $\bSigma$ by $\hbSigma$ in the expressions of $\cJ_1$ and $\cJ_2$.

\medskip
\noindent\textbf{\underline{Proof of ratio consistency}}.
From \eqref{eq:target_frontier} and \eqref{eq:target_frontier1}, we can conclude 
\begin{align}
  \frac{1+s_0+s_{1,\Qb}}{(1+s_0)^2-s_{1,\bSigma}}\cdot\frac{(\alpha_n \one+\beta_n\rb)^\top (\hbSigma+\Qb)^{-1}\bSigma(\hbSigma+\Qb)^{-1}(\alpha_n \one+\beta_n\rb)}{(\alpha_n \one+\beta_n\rb)^\top (\hbSigma+\Qb)^{-1}\hbSigma(\hbSigma+\Qb)^{-1}(\alpha_n \one+\beta_n\rb)}~\cvas~1\label{eq:ratio_const_frontier}.
\end{align}
Here the convergence comes from \eqref{eq:target_frontier}/\eqref{eq:target_frontier1}. Moreover, Proposition~\ref{prop:converge_s} indicates that
\begin{align*}
    \frac{1+s_0+s_{1,\Qb}}{(1+s_0)^2-s_{1,\bSigma}} /\bigg(1-\frac{c}{p}\tr\hbSigma\big(\hbSigma+\Qb\big)^{-1} \bigg)^2\cvas1.
\end{align*}
Combined this with \eqref{eq:ratio_const_frontier} we have
\begin{align*}
    \frac{\sigma_0^2}{\hat{\sigma}^2}=\bigg(1-\frac{c}{p}\tr\hbSigma\big(\hbSigma+\Qb\big)^{-1} \bigg)^2\cdot\frac{(\alpha_n \one+\beta_n\rb)^\top (\hbSigma+\Qb)^{-1}\bSigma(\hbSigma+\Qb)^{-1}(\alpha_n \one+\beta_n\rb)}{(\alpha_n \one+\beta_n\rb)^\top (\hbSigma+\Qb)^{-1}\hbSigma(\hbSigma+\Qb)^{-1}(\alpha_n \one+\beta_n\rb)}~\cvas~1
\end{align*}
This completes the proof of ratio consistency.

\medskip
\noindent\textbf{\underline{Proof of difference consistency}}.
As for the remaining conclusions on the consistency of Sharpe ratio or variance estimation, from the analysis above we have
\begin{enumerate}
    \item $\sigma_0^2=\Theta(1/\cA_{11})$ when $\mu_0\leq C\sqrt{\cA_{rr}/\cA_{11}}$;
    \item $\sigma_0^2=\Theta(\mu_0^2/{\cA_{rr}})$ when $\mu_0\geq C\sqrt{\cA_{rr}/\cA_{11}}$.
\end{enumerate}
When $\cA_{rr}$ is bounded, if $\mu_0\leq C\sqrt{\cA_{rr}/\cA_{11}}$ we can see that $\mu_0/\sigma_0=\Theta(\mu_0\cdot\sqrt{\cA_{11}})= O(\sqrt{\cA_{rr}})=O(1)$; if $\mu_0\geq C\sqrt{\cA_{rr}/\cA_{11}}$,  we can see that $\mu_0/\sigma_0=\Theta(\mu_0/(\mu_0/\sqrt{\cA_{rr}}))=O(1)$.
Therefore by ratio consistency we have $\frac{\mu_0}{\sigma_0}-\frac{\mu_0}{\hat{\sigma}}~\cvas~0$. When $r_0=O(\mu_0)$, we also have $\frac{r_0}{\sigma_0}-\frac{r_0}{\hat{\sigma}}~\cvas~0$, this completes the proof of the consistency of Sharpe ratio.

To prove $ \sigma_0^2-\hat{\sigma}^2 ~\cvas~0$. By $\|\bSigma/p\|_{\tr}\leq C$, we have $\lambda_{\min}\big(\bSigma/(1+s_0)+\Qb\big)^{-1}\geq \Theta(1/p)$, which indicates $\cA_{11}\geq \Theta(1)$. Therefore, if $\mu_0\leq C\sqrt{\cA_{rr}/\cA_{11}}$ we have $\sigma_0^2=\Theta(1/\cA_{11})=O(1)$; if $\mu_0\geq C\sqrt{\cA_{rr}/\cA_{11}}$ we have $\sigma_0^2=\Theta(\mu_0^2/{\cA_{rr}})$. The condition $\mu_0\leq C\sqrt{\cA_{rr}}$ then indicates that $\sigma_0^2=O(1)$. By ratio consistency we have $ \sigma_0^2-\hat{\sigma}^2 ~\cvas~0$. We complete the proof of Theorem~\ref{thm:frontier}.
\end{proof}

\subsection{Proof of Lemma~\ref{lemma:same_order_2_to_1}}
\label{sec:prooflemma_same_order_2_to_1}
Note that $\Ab$ is semi-positive defined. It is equal for us to prove that for any sequence of $\bxi\in\RR^{p}$, 
\begin{align*}
    \bxi^\top \bigg( \frac{\bSigma}{1+s_0}+\Qb \bigg)^{-1}\bSigma\bigg( \frac{\bSigma}{1+s_0}+\Qb \bigg)^{-1} \bxi=\Theta \bigg(  \bxi^\top \bigg( \frac{\bSigma}{1+s_0}+\Qb \bigg)^{-1}\bxi \bigg).
\end{align*}
We have
\begin{align*}
    &\bxi^\top \bigg( \frac{\bSigma}{1+s_0}+\Qb \bigg)^{-1}\bSigma\bigg( \frac{\bSigma}{1+s_0}+\Qb \bigg)^{-1} \bxi\\
    &\quad =\bxi^\top \bigg( \frac{\bSigma}{1+s_0}+\Qb \bigg)^{-\frac{1}{2}} \bigg( \frac{\bSigma}{1+s_0}+\Qb \bigg)^{-\frac{1}{2}}\bSigma \bigg( \frac{\bSigma}{1+s_0}+\Qb \bigg)^{-\frac{1}{2}} \bigg( \frac{\bSigma}{1+s_0}+\Qb \bigg)^{-\frac{1}{2}} \bxi. 
\end{align*}
The first conclusion in Lemma~\ref{lemma:same_order_2_to_1} holds from the fact that
\begin{align*}
    &\lambda_{\max}\bigg(\bigg( \frac{\bSigma}{1+s_0}+\Qb \bigg)^{-\frac{1}{2}}\bSigma \bigg( \frac{\bSigma}{1+s_0}+\Qb \bigg)^{-\frac{1}{2}}\bigg)\leq 1+s_0,\\
    &\lambda_{\min}\bigg(\bigg( \frac{\bSigma}{1+s_0}+\Qb \bigg)^{-\frac{1}{2}}\bSigma \bigg( \frac{\bSigma}{1+s_0}+\Qb \bigg)^{-\frac{1}{2}}\bigg)\\
    &\qquad =\lambda_{\min}\bigg(\bigg( \frac{\bSigma}{1+s_0}+\Qb \bigg)^{-1}\bSigma \bigg)\geq \frac{1+s_0}{C'(1+s_0)+1}.
\end{align*}
To prove $\bxi^\top \big(\bSigma/(1+s_0)+\Qb\big)^{-1}\bxi=\Theta (\|\bSigma^{-\frac{1}{2}}\bxi\|_2^2)$, we have
\begin{align*}
   \bxi^\top \big(\bSigma/(1+s_0)+\Qb\big)^{-1}\bxi=(\bSigma^{-\frac{1}{2}}\bxi)^\top \big(\Ib/(1+s_0)+\bSigma^{-\frac{1}{2}}\Qb\bSigma^{-\frac{1}{2}}\big)^{-1}\bSigma^{-\frac{1}{2}}\bxi,
\end{align*}
the conclusion directly holds from the fact that the maximal and minimal eigenvalue of $\big(\Ib/(1+s_0)+\bSigma^{-\frac{1}{2}}\Qb\bSigma^{-\frac{1}{2}}\big)^{-1}$ can be upper and lower bounded by a constant, respectively. This completes the proof of Lemma~\ref{lemma:same_order_2_to_1}.

\subsection{Proof of Lemma~\ref{lemma:analysis_dependence}}
\label{sec:prooflemma_analysis_dependence}
We first prove for the case of $\mu_0\leq C\sqrt{\cA_{rr}/\cA_{11}} $. From the results of Proposition~\ref{prop:converge_Tn1Tn2}, we have that
\begin{align*}
    \frac{\rb^\top \big(\hbSigma+\Qb\big)^{-1}\rb}{\cA_{rr}}-1~\cvas~0.
\end{align*}
For the term $\mu_0\rb^\top \big(\hbSigma+\Qb \big)^{-1}\one/\cA_{rr}$, we have 
\begin{align*}
    \frac{\mu_0\rb^\top \big(\hbSigma+\Qb \big)^{-1}\one}{\cA_{rr}}=\frac{\mu_0(\bSigma^{-\frac{1}{2}}\rb)^\top \big(\Zb^\top\Zb/n+\bSigma^{-\frac{1}{2}}\Qb\bSigma^{-\frac{1}{2}} \big)^{-1}\bSigma^{-\frac{1}{2}}\one}{\cA_{rr}}.
\end{align*}
We next prove that 
\begin{align*}
    \mu_0\cdot \bigg\| \frac{\bSigma^{-\frac{1}{2}}\one(\bSigma^{-\frac{1}{2}}\rb)^\top}{\cA_{rr}} \bigg\|_{\tr}
\end{align*}
is bounded. If the results above hold, by applying the results in Proposition~\ref{prop:converge_Tn1Tn2} we will have
\begin{align*}
    \frac{\mu_0(\bSigma^{-\frac{1}{2}}\rb)^\top \big(\Zb^\top\Zb/n+\bSigma^{-\frac{1}{2}}\Qb\bSigma^{-\frac{1}{2}} \big)^{-1}\bSigma^{-\frac{1}{2}}\one}{\cA_{rr}}-\frac{\mu_0(\bSigma^{-\frac{1}{2}}\rb)^\top \big(\Ib/(1+s_0)+\bSigma^{-\frac{1}{2}}\Qb\bSigma^{-\frac{1}{2}} \big)^{-1}\bSigma^{-\frac{1}{2}}\one}{\cA_{rr}}~\cvas~0,
\end{align*}
which directly implies 
\begin{align*}
    \frac{\mu_0\rb^\top \big(\hbSigma+\Qb \big)^{-1}\one}{\cA_{rr}}-\frac{\mu_0\rb^\top \big(\bSigma/(1+s_0)+\Qb \big)^{-1}\one}{\cA_{rr}}~\cvas~0.
\end{align*}
Then the first conclusion in Lemma~\ref{lemma:analysis_dependence} will hold.

It remains only to prove that $\mu_0\cdot \Big\| \frac{\bSigma^{-\frac{1}{2}}\one(\bSigma^{-\frac{1}{2}}\rb)^\top}{\cA_{rr}} \Big\|_{\tr}$ is bounded. Note that $\mu_0=O(\sqrt{\cA_{rr}/\cA_{11}})$, therefore it is equal to prove that
\begin{align*}
\bigg\| \frac{\bSigma^{-\frac{1}{2}}\one(\bSigma^{-\frac{1}{2}}\rb)^\top}{\sqrt{\cA_{rr}\cdot\cA_{11}}} \bigg\|_{\tr}
\end{align*}
is bounded. From Lemma~\ref{lemma:same_order_2_to_1}, we can see that $\sqrt{\cA_{rr}\cdot\cA_{11}}=\Theta(\|\bSigma^{-\frac{1}{2}}\one\|_2\cdot\|\bSigma^{-\frac{1}{2}}\rb\|_2 )$, we conclude that
\begin{align*}
    \mu_0\cdot \bigg\| \frac{\bSigma^{-\frac{1}{2}}\one(\bSigma^{-\frac{1}{2}}\rb)^\top}{\cA_{rr}} \bigg\|_{\tr}
\end{align*}
is bounded, 
We complete the proof of the first conclusion when $\mu_0\leq C\sqrt{\cA_{rr}/\cA_{11}} $. The proof of the second conclusion when $\mu_0\leq C\sqrt{\cA_{rr}/\cA_{11}} $ follows the same procedure and thus we omit it.

\medskip
We then prove the case of $\mu_0\geq C\sqrt{\cA_{rr}/\cA_{11}} $. For the first equation, we have
\begin{align*}
    \frac{\rb^\top \big(\hbSigma+\Qb\big)^{-1}\rb}{\cA_{rr}}-\frac{\cA_{rr}}{\cA_{rr}}~\cvas~0.
\end{align*}
Note that $\mu_0\sqrt{\cA_{11}/\cA_{rr}}\geq C>0$, we divide $\mu_0\sqrt{\cA_{11}/\cA_{rr}}$ into the equation above and get that
\begin{align}
    \frac{\rb^\top \big(\hbSigma+\Qb\big)^{-1}\rb}{\mu_0\sqrt{\cA_{rr}\cA_{11}}}-\frac{\cA_{rr}}{\mu_0\sqrt{\cA_{rr}\cA_{11}}}~\cvas~0.\label{eq:mu0>C_converge1}
\end{align}
We also know  that
\begin{align}
    \frac{\rb^\top \big(\hbSigma+\Qb \big)^{-1}\one}{\sqrt{\cA_{rr}\cA_{11}}}-\frac{\cA_{r1}}{\sqrt{\cA_{rr}\cA_{11}}}~\cvas~0\label{eq:mu0>C_converge2}
\end{align}
since $\Big\| \frac{\bSigma^{-\frac{1}{2}}\one(\bSigma^{-\frac{1}{2}}\rb)^\top}{\sqrt{\cA_{rr}\cdot\cA_{11}}} \Big\|_{\tr}$ is bounded, therefore
\eqref{eq:mu0>C_converge1} and \eqref{eq:mu0>C_converge2} show that
\begin{align*}
    \frac{\alpha_n-\alpha_0}{\mu_0\sqrt{\cA_{rr}\cA_{11}}}=\frac{\rb^\top \big(\hbSigma+\Qb \big)^{-1}\rb-\mu_0\rb^\top \big(\hbSigma+\Qb \big)^{-1}\one }{\mu_0\sqrt{\cA_{rr}\cA_{11}}}-\frac{\cA_{rr}-\mu_0 \cA_{r1}}{\mu_0\sqrt{\cA_{rr}\cA_{11}}}~\cvas~0.
\end{align*}
We can similarly show that 
\begin{align*}
    \frac{\mu_0\one^\top \big(\hbSigma+\Qb \big)^{-1}\one-\rb^\top \big(\hbSigma+\Qb \big)^{-1}\one }{\mu_0\cA_{11}}-\frac{\mu_0\cA_{11}-\cA_{r1}}{\mu_0\cA_{11}}~\cvas~0.
\end{align*}

For the order of $D=BC-A^2$, it is clear from \eqref{eq:mu0>C_converge2} that
\begin{align*}
    \frac{A^2-\cA_{r1}}{\cA_{rr}\cA_{11}}~\cvas~0, \quad \frac{BC}{\cA_{rr}\cA_{11}}-1~\cvas~0.
\end{align*}
Therefore we have $(BC-A^2)/(\cA_{rr}\cA_{11})-(\cA_{rr}\cA_{11}-\cA_{r1}^2)/(\cA_{rr}\cA_{11})~\cvas~0$. By Assumption~\ref{assump:assump5} we have $(\cA_{rr}\cA_{11}-\cA_{r1}^2)/(\cA_{rr}\cA_{11})=\Theta(1)$, therefore $D/(\cA_{rr}\cA_{11})=\Theta(1)$.
This completes the proof of Lemma~\ref{lemma:analysis_dependence}.

\section{Relationship Between Ridge Regularization and Maximal Sharpe Ratio}
\label{sec:discussion_largest_SR}
According to the theory proposed by \cite{markowitz1952harry}, the maximal population Sharpe ratio can be expressed as:
\begin{align*}
    SR_{\max}=\sqrt{\bmu^\top \bSigma^{-1}\bmu}.
\end{align*}
In this section, we discuss the relationship between $SR_{\max}$ and the optimal Sharpe ratio $SR(\Qb)$ we can achieve by choosing $\Qb$. The following proposition demonstrates that there exists $\tilde\Qb$ such that $SR(\tilde\Qb)$ can approximate $SR_{\max}$ very closely.
\begin{proposition}
\label{prop:compare_SRmax_SRQ}
Suppose that Assumptions~\ref{def:X}, \ref{ass:assump1}, \ref{ass:assump2} and \ref{ass:assump3} hold, and  $\|\bSigma^{-1}\|_{\op}$ is bounded. Then for any given $\varepsilon>0$, there exists deterministic sequences of matrices $\tilde\Qb\in\RR^{p\times p}$ such that with probability 1,
    \begin{align*}
    1-\varepsilon\leq\lim_{n\to+\infty} SR(\tilde\Qb)/SR_{\max}\leq 1.
\end{align*}
\end{proposition}
The proof of Proposition~\ref{prop:compare_SRmax_SRQ} is given in Section~\ref{subsec:proofofcompare}. This proposition tells us that a well-designed regularization matrix $\tilde\Qb$ allows $SR(\tilde\Qb)$ to closely approximate the maximal population Sharpe ratio. 
Furthermore, we could use the criteria established in Theorem~\ref{thm:main_theorem} to select $\Qb$ that approximately maximizes $SR(\Qb)$.

Another interesting case arises when $\bSigma$ exhibits a specific factor structure, defined as $\bSigma = \Bb\Bb^\top + \Db$. Here, $\Bb = (\bb_1, \dots, \bb_K) \in \mathbb{R}^{p \times K}$ is knonwn as the loading matrix with each $\|\bb_j\|_2 = \Theta(\sqrt{p})$ for $j \in [K]$, $\operatorname{rank}(\Bb)=K$. $\Db$ is known as the residual covariance matrix, which is a positive definite matrix with eigenvalues bounded away from $0$ and $+\infty$. $\Db$ is typically assumed to be a diagonal matrix, as the off-diagonal correlations have been taken care of mostly by the low-rank factor component $\Bb\Bb^\top$. But for the discussion here, we only need to assume $\Db$ has bounded spectrum from above and below.
Under this model, the following proposition provides further insight into the role of $\Qb$.
\begin{proposition}
\label{prop:compare_SRmax_SRQ_detail}
Assume that Assumptions~\ref{def:X}, \ref{ass:assump1}, \ref{ass:assump2} and \ref{ass:assump3} are satisfied, and  $\bSigma=\Bb\Bb^\top+\Db$ as specified. Let $K$ be fixed and suppose $\bb_j^\top \Db^{-1} \bmu = O(\|\bmu\|_2)$. Then, for any given $\varepsilon>0$, there exists deterministic sequences of matrices $\tilde\Qb$ with bounded operator norm such that with probablity 1,
    \begin{align*}
    1-\varepsilon\leq\lim_{n\to+\infty} SR(\tilde\Qb)/SR_{\max}\leq 1.
\end{align*}
\end{proposition}
The proof of Proposition~\ref{prop:compare_SRmax_SRQ_detail} is given in Section~\ref{subsec:proofofcompare_detail}. This proposition demonstrates that even with the presence of factors, i.e. spiked eigenvalues, it is still possible to identify a regularization matrix $\tilde\Qb$ with a bounded operator norm that can approximate the maximal population Sharpe ratio. Again, we can utilize Theorem~\ref{thm:main_theorem} to evaluate the effectiveness of the design of $\Qb$.

\subsection{Proof of Proposition~\ref{prop:compare_SRmax_SRQ}}
\label{subsec:proofofcompare}
Recall the portfolio allocation vector $\wb\propto (\hbSigma+\Qb)^{-1}\bmu$, where $\hbSigma$ represents the sample covariance matrix and $\Qb$ is a Ridge regularization term. The Sharpe ratio $SR(\Qb) $ can be written as:
\begin{align*}
SR(\Qb)=\frac{\bmu^\top(\hbSigma+\Qb)^{-1}\bmu}{\sqrt{\bmu^\top\tr(\hbSigma+\Qb)^{-1}\bSigma(\hbSigma+\Qb)^{-1}\bmu}}.
\end{align*}
From Proposition~\ref{prop:converge_Tn1Tn2}, it can be seen that:
\begin{align*}
    \frac{1}{\sqrt{1-\frac{s_{1,\bSigma}}{(1+s_0)^2}}}\cdot \frac{\bmu^\top(\bSigma/(1+s_0)+\Qb)^{-1}\bmu}{\sqrt{\bmu^\top(\bSigma/(1+s_0)+\Qb)^{-1}\bSigma(\bSigma/(1+s_0)+\Qb)^{-1}\bmu}} \cdot \frac{1}{SR(\Qb)} ~\cvas~ 1,
\end{align*}
where $s_0$ and $s_{1,\bSigma}$ are defined in Lemma~\ref{lemma:constant_level} by letting $z=0$. If $\Qb$ is carefully designed such that $\Qb= q\cdot\bSigma$, we have
\begin{align*}
&\frac{1}{\sqrt{1-\frac{s_{1,\bSigma}}{(1+s_0)^2}}}\cdot \frac{\bmu^\top(\bSigma/(1+s_0)+q\bSigma)^{-1}\bmu}{\sqrt{\bmu^\top(\bSigma/(1+s_0)+q\bSigma)^{-1}\bSigma(\bSigma/(1+s_0)+q\bSigma)^{-1}\bmu}} \cdot \frac{1}{SR(\Qb)}\\
    &\quad=\frac{1}{\sqrt{1-\frac{s_{1,\bSigma}}{(1+s_0)^2}}}\cdot\sqrt{\bmu^\top \bSigma^{-1}\bmu}/SR(\Qb) ~\cvas~ 1.
\end{align*}
Here, $s_0=\frac{c}{p}\tr\bSigma(\bSigma/(1+s_0)+q\bSigma)^{-1}=\frac{c(1+s_0)}{q(1+s_0)+1}$. It should be noted that $s_0$ is uniformly bounded, and for any $\varepsilon>0$, there exists $M>0$ such that for any given $q>M$ it holds that $0<s_0<\varepsilon/3$. We have $0<-s_{1,\bSigma}/(1+s_0)^2\leq s_0\leq \varepsilon/3$, where the first inequality comes from Lemma~\ref{lemma:constant_level}. By taking $\tilde\Qb=2M\bSigma$, we conclude that with probability 1,
\begin{align*}
    1\geq \lim_{n\to+\infty}SR(\tilde\Qb)/SR_{\max}\geq \lim_{n\to+\infty}\frac{1}{\sqrt{1-\frac{s_{1,\bSigma}}{(1+s_0)^2}}}\geq\lim_{n\to+\infty} \frac{1}{\sqrt{1+\frac{\varepsilon}{3}}}\geq 1-\varepsilon,
\end{align*}
which completes the proof of Proposition~\ref{prop:compare_SRmax_SRQ}.


\subsection{Proof of Proposition~\ref{prop:compare_SRmax_SRQ_detail}}
\label{subsec:proofofcompare_detail}
From Proposition~\ref{prop:converge_Tn1Tn2}, it can be seen that:
\begin{align*}
\frac{1}{\sqrt{1-\frac{s_{1,\bSigma}}{(1+s_0)^2}}}\cdot \frac{\bmu^\top(\bSigma/(1+s_0)+\Qb)^{-1}{\bmu}}{\sqrt{\bmu^\top(\bSigma/(1+s_0)+\Qb)^{-1}\bSigma(\bSigma/(1+s_0)+\Qb)^{-1}\bmu}} \cdot \frac{1}{SR(\Qb)} ~\cvas~ 1,
\end{align*}
where $s_0$ and $s_{1,\bSigma}$ are defined in Lemma~\ref{lemma:constant_level} by letting $z=0$. We first consider the expression of $\bmu^\top(\bSigma/(1+s_0)+\Qb)^{-1}\bmu$ and $\bmu^\top(\bSigma/(1+s_0)+\Qb)^{-1}\bSigma(\bSigma/(1+s_0)+\Qb)^{-1}\bmu$ by specifically choosing $\Qb=q\cdot  \Db$. Letting $\lambda=1+(1+s_0)q$, we have
\begin{align}
    \bigg(\frac{\bSigma}{1+s_0}+q\Db\bigg)^{-1}/(1+s_0)&=\big(\Bb\Bb^\top+\lambda\Db\big)^{-1}\nonumber\\
    &=\frac{1}{\lambda}\cdot\big(\Db^{-1}-\Db^{-1}\Bb(\lambda\Ib_K+\Bb^{\top}\Db^{-1}\Bb)^{-1}\Bb^\top\Db^{-1} \big)\nonumber\\
    &=\frac{1}{\lambda}\cdot\big(\Db^{-1}-\Db^{-1}\Bb\bXi\Bb^\top\Db^{-1} \big),\label{eq:woodbury_compare}
\end{align}
where $\bXi$ is defined as  $\bXi=(\lambda\Ib_K+\Bb^{\top}\Db^{-1}\Bb)^{-1}\in\RR^{K\times K}$. Here, the second equality is by Woodbury formula. We can easily show that $\|\bXi\|_{\op}=O(1/p)$. From \eqref{eq:woodbury_compare}, we have
\begin{align*}
    &\frac{\bmu^\top(\bSigma/(1+s_0)+\Qb)^{-1}\bmu}{\sqrt{\bmu^\top(\bSigma/(1+s_0)+\Qb)^{-1}\bSigma(\bSigma/(1+s_0)+\Qb)^{-1}\bmu}}\\
    &\qquad=\frac{\bmu^\top\big(\Db^{-1}-\Db^{-1}\Bb\bXi\Bb^\top\Db^{-1} \big)\bmu}{\sqrt{\bmu^\top\big(\Db^{-1}-\Db^{-1}\Bb\bXi\Bb^\top\Db^{-1} \big)\bSigma\big(\Db^{-1}-\Db^{-1}\Bb\bXi\Bb^\top\Db^{-1} \big)\bmu}}.
\end{align*}
By the condition that $|\bb_j^\top\Db^{-1}\bmu|=O(\|\bmu\|_2)$, we can see 
\begin{align*}
&\bmu^\top\big(\Db^{-1}-\Db^{-1}\Bb\bXi\Bb^\top\Db^{-1} \big)\bmu=\bmu^\top\Db^{-1}\bmu+O(\|\bmu\|_2^2/p),\\
    &\bmu^\top\big(\Db^{-1}-\Db^{-1}\Bb\bXi\Bb^\top\Db^{-1} \big)\Bb=\lambda\bmu^\top \Db^{-1} \Bb \bXi=O(\|\bmu\|_2^2/p),\\
    &\bmu^\top\big(\Db^{-1}-\Db^{-1}\Bb\bXi\Bb^\top\Db^{-1} \big)\Db\big(\Db^{-1}-\Db^{-1}\Bb\bXi\Bb^\top\Db^{-1} \big)\bmu\\
&\qquad=\bmu^\top\Db^{-1}\bmu-2\bmu^\top\Db^{-1}\Bb\bXi\Bb\Db^{-1}\bmu+\bmu^\top\Db^{-1}\Bb\bXi\Bb^\top\Db^{-1}\Bb\bXi\Bb\Db^{-1}\bmu\\
&\qquad =\bmu^\top\Db^{-1}\bmu+O(\|\bmu\|_2^2/p).
\end{align*}
Here, we use the fact $\Bb^\top\Db^{-1}\Bb\bXi=\Ib_K-\lambda \bXi$ and $\|\bXi\|_{\op}=O(1/p)$. Hence, we conclude that by letting $\Qb=q\cdot \Db$, it holds 
\begin{align*}
    \frac{\bmu^\top(\bSigma/(1+s_0)+\Qb)^{-1}\bmu}{\sqrt{\bmu^\top(\bSigma/(1+s_0)+\Qb)^{-1}\bSigma(\bSigma/(1+s_0)+\Qb)^{-1}\bmu}}&=\frac{\bmu^\top\Db^{-1}\bmu+O(\|\bmu\|_2^2/p)}{\sqrt{\bmu^\top\Db^{-1}\bmu+O(\|\bmu\|_2^2/p)}}\\
&=\sqrt{\bmu^\top\Db^{-1}\bmu}+O(\|\bmu\|_2/p).
\end{align*}
Note also that
\begin{align*}
    SR_{\max}^2=\bmu^\top \bSigma^{-1}\bmu=\bmu^\top\Db^{-1}\bmu-\bmu^\top\Db^{-1}\Bb(\Ib_K+\Bb^\top\Db^{-1}\Bb)^{-1}\Bb^\top\bmu=\bmu^\top\Db^{-1}\bmu+O(\|\bmu\|_2^2/p), 
\end{align*}
we have
\begin{align*}
\lim_{n\to+\infty}\frac{\bmu^\top(\bSigma/(1+s_0)+\Qb)^{-1}\bmu}{\sqrt{\bmu^\top(\bSigma/(1+s_0)+\Qb)^{-1}\bSigma(\bSigma/(1+s_0)+\Qb)^{-1}\bmu}} \cdot \frac{1}{SR_{\max}} =1.
\end{align*}
Therefore, by letting $\Qb=q\cdot\Db$, with probability 1, 
\begin{align*}
1\geq\lim_{n\to+\infty}SR(\Qb)/SR_{\max}=\lim_{n\to+\infty} \frac{1}{\sqrt{1-\frac{s_{1,\bSigma}}{(1+s_0)^2}}}.
\end{align*}
Here, $s_0$ and $s_{1,\bSigma}$ are defined in Lemma~\ref{lemma:constant_level} by letting $z=0$ and $\Qb=q\cdot\Db$. The solution of $s_0$ is given by 
\begin{align*}
    s_0=\frac{c}{p}\tr \bSigma\bigg(\frac{\bSigma}{1+s_0}+q\Db\bigg)^{-1}, 
\end{align*}
rewrite the formula above we can see
\begin{align*}
    s_0=\frac{c(1+s_0)}{p}\tr (\Bb\Bb^\top+\Db)\big( \Bb\Bb^\top+\lambda\Db\big)^{-1}.
\end{align*}
Here, $\lambda=1+(1+s_0)q$. Easy to see that $\|\Bb^\top \big( \Bb\Bb^\top+\lambda\Db\big)^{-1}\Bb\|_{\op}=O(1) $, hence we have
\begin{align*}
    s_0&=\frac{c(1+s_0)}{p}\tr\Db\big( \Bb\Bb^\top+\lambda\Db\big)^{-1}+O(1/p)\\
    &=\frac{c}{p}\tr\bigg(\frac{\Db^{-\frac{1}{2}}\Bb\Bb^\top \Db^{-\frac{1}{2}}+\Ib}{1+s_0} +q\Ib\bigg)^{-1}+O(1/p)\\
    &\leq c/q+O(1/p)
\end{align*}
Here, we use the fact that $s_0$ is uniformly bounded. For any $\varepsilon>0$, there exists constant $M>0$ such that for any $q>M$, $s_0\leq \varepsilon/6+O(1/p)$, hence by letting $\tilde\Qb=2M\Db$, we have $\lim_{n\to+\infty}s_0\leq \varepsilon/6$. Thus when $\Qb=\tilde\Qb=2M\Db$, with probability 1 it holds that $0<\lim_{n\to+\infty}-s_{1,\bSigma}/(1+s_0)^2\leq \lim_{n\to+\infty}s_0\leq \varepsilon/6$, which indicates that 
\begin{align*}
\lim_{n\to+\infty}SR(\tilde\Qb)/SR_{\max}=\lim_{n\to+\infty} \frac{1}{\sqrt{1-\frac{s_{1,\bSigma}}{(1+s_0)^2}}}\geq \frac{1}{\sqrt{1+\varepsilon/6}}\geq 1-\varepsilon
\end{align*}
with probability 1. This completes the proof.

\section{Out-of-sample Sharpe Ratio with Unknown Mean Vector}
\label{sec:unknown_mu}
In this section, we extend to the case where $\bmu$ is unknown. In this scenario, the entire framework must be adjusted. When $\bmu$ is unknown, it is clear that the optimization problem can be rewritten as:
\begin{align*}
    \wb^*=\argmin\limits_{\wb\in\RR^p} \wb^\top \bSigma \wb,\quad \st \wb^\top \hbmu =\mu_0,
\end{align*}
where $\hbmu$ is some estimation of the mean vector $\bmu$. 
Then the out-sample Sharpe ratio can be expressed as 
\begin{align}
SR(\Qb)=\frac{\hbmu^\top(\hbSigma+\Qb)^{-1}\bmu}{\sqrt{\hbmu^\top(\hbSigma+\Qb)^{-1}\bSigma(\hbSigma+\Qb)^{-1}\hbmu}}.\label{eq:unknown_SRQ}
\end{align}
If $\bmu$ is unknown and $\hat{\bmu}$ is not independent of $\hbSigma$, this would introduce significant complexity, as the dependence between $\hat{\bmu}$ and $\hbSigma$ can be highly non-trival, e.g. if we use a machine learning method to fit $\hat{\bmu}$. 
From a technical standpoint, to our knowledge, existing literature commonly assumed independence between  $\hbmu$ and $\hbSigma$, especially in Random Matrix Theory (RMT) analysis, as seen in works such as \citet{li2022spectrally}, \citet{bodnar2022recent}, \citet{bodnar2024reviving}, and \citet{bodnar2024two}.

A natural case we consider here is that the sample data comes from the normal distribution $\cN(\bmu,\bSigma)$, then  the sample mean and  sample covariance are independent.  We give the following additional assumptions when we consider the case where $\bmu$ is unknown: 
\begin{assumption}
\label{ass:assumpadd1}
Given the mean return vector $\rb\in\RR^p$ and the covariance matrix $\bSigma\in\RR^{p\times p}$, we consider the observed data matrix $\Rb\in\RR^{n\times p}$ with the form
\begin{align*}
    \Rb=\one_n \rb^\top+\Xb,
\end{align*}
where $\Xb=\Zb\bSigma^{\frac{1}{2}}\in\RR^{n\times p}$. Here, the elements in the matrix $\Zb\in\RR^{n\times p}$ are i.i.d standard normal distribution. The vector $\bmu$ is given by $\rb-r_0\one$ with known free risk rate $r_0$. 
\end{assumption}
Compared with Assumption~\ref{def:X}, we assumed the normality of the data distribution to ensure the estimated $\hbmu$ and $\hbSigma$ are independent.  When the data is Gaussian, it is easy to get the mean vector $\hbmu=\overline{\Rb}-r_0\one_p$ satisfies that $\bSigma^{-\frac{1}{2}}(\hbmu-\bmu)\sim \cN(0,\Ib_p/n)$, and the sample covariance matrix $\hbSigma$ is given by $\hbSigma=(\Rb-r_0\one_n\one^\top_p-\one^\top_n\hbmu)^\top(\Rb-r_0\one_n\one^\top_p-\one^\top_n\hbmu)/n$. Here, $\overline{\Rb}$ is the sample mean of $\Rb$. The small rank perturbation on matrix $\hbSigma$ will not  affect the large-scale behavior or distribution of the eigenvalues and eigenvectors of the covariance matrix as dimension grows \citep{yao2015sample}, and  the independence of $\hbSigma$ and $\hbmu$ will then maintain all the properties we analyzed above. We have the following theorem with unknown $\bmu$.

\begin{theorem}
\label{thm:unknown_mu}
Suppose Assumptions~\ref{ass:assump1}, \ref{ass:assump2}, \ref{ass:assump3} and \ref{ass:assumpadd1} hold. Additionally,  assume that $\|\bSigma^{-\frac{1}{2}}\bmu\|_2$ is bounded. For any $\Qb\in\cQ$, a good estimator $\hat{SR}(\Qb)$ for $SR(\Qb)$ which is defined in \eqref{eq:unknown_SRQ} is given as follows. 
\begin{align*}
    \hat{SR}(\Qb)=\frac{\hbmu^\top(\hbSigma+\Qb)^{-1}\hbmu-\frac{\tr(\hbSigma+\Qb)^{-1}\hbSigma}{n-\tr(\hbSigma+\Qb)^{-1}\hbSigma}}{\sqrt{\hbmu^\top(\hat{\bSigma}+\Qb)^{-1}\hbSigma(\hat{\bSigma}+\Qb)^{-1}\hbmu}}\cdot\big(1-\frac{c}{p} \tr\hbSigma(\hbSigma+\Qb)^{-1}\big).
\end{align*}
It holds that 
\begin{align*}
    &\hbmu^\top(\hbSigma+\Qb)^{-1}\bmu-\hbmu^\top(\hbSigma+\Qb)^{-1}\hbmu+\frac{\tr(\hbSigma+\Qb)^{-1}\hbSigma}{n-\tr(\hbSigma+\Qb)^{-1}\hbSigma}~\cvas~0;\\
    & \sqrt{\frac{\hbmu^\top(\hat{\bSigma}+\Qb)^{-1}\bSigma(\hat{\bSigma}+\Qb)^{-1}\hbmu}{\hbmu^\top(\hat{\bSigma}+\Qb)^{-1}\hbSigma(\hat{\bSigma}+\Qb)^{-1}\hbmu}}\cdot  \big(1-\frac{c}{p} \tr\hbSigma(\hbSigma+\Qb)^{-1}\big)~\cvas~1.
\end{align*}
If additionally $\bmu^\top\big(\frac{\bSigma}{1+s_0}+\Qb\big)^{-1}\bmu$ is lower bounded, it holds that
\begin{align*}
    \hat{SR}(\Qb)/SR(\Qb)~\cvas~1.
\end{align*}
\end{theorem}
It is clear that $\hbmu^\top(\hbSigma+\Qb)^{-1}\bmu$ and $\hbmu^\top(\hbSigma+\Qb)^{-1}\bSigma(\hbSigma+\Qb)^{-1}\hbmu$ correspond to the numerator and denominator of $SR(\Qb)$, respectively. Under stronger conditions, we have established the convergence of these quantities when $\bmu$ is unknown, providing valuable insights that can also inform practical applications. The condition for a lower bound on $\bmu^\top\big(\frac{\bSigma}{1+s_0}+\Qb\big)^{-1}\bmu$ requires that $\bmu$ does not lie entirely within the space spanned by distant spiked eigenvectors. Otherwise, if $\bmu$ falls solely within this space, $\hbmu$ will perform poorly as a predictor, and the ratio’s convergence will fail in this scenario. 

\begin{remark}
Theorem~\ref{thm:unknown_mu} is established under Gaussian assumption. When data distribution is non-Gaussian, it posts significant challenges.  From a technical perspective, our existing random-matrix framework crucially assumes independence between $\hat{\bmu}$ and $\hbSigma$. Introducing dependence breaks down the foundation of the ``leave-one-out'' technique from Random Matrix Theory \citep{bai2010spectral,rubio2011spectral}, a cornerstone of our analysis. To our best knowledge, no alternative methodologies exist to address the dependency between $\hat{\bmu}$ and $\hbSigma$  in high-dimensional settings. So to handle dependence between $\hat{\bmu}$ and $\hbSigma$, it is necessary to develop new tools in order to study  the asymptotic behavior of the complicated term $\hat{\bmu}(\Rb)^\top (\hbSigma(\Rb) + \Qb)^{-1} \hat{\bmu}(\Rb)$, which depend on our random data $\Rb$ from multiple places. This is beyond the scope of this work, and we view this as an important direction for future research.
\end{remark}

\subsection{Proof of Theorem~\ref{thm:unknown_mu}}
\label{sec:discussion_unknown_mu}

By the definition of $\hbmu$, it is clear that 
\begin{align*}
    \hbmu-\bmu\sim \cN(0,\bSigma/n). 
\end{align*}
Set $\zb=\hbmu-\bmu$, then it is clear that
\begin{align}
\hbmu^\top(\hbSigma+\Qb)^{-1}\bmu=\hbmu^\top(\hbSigma+\Qb)^{-1}\hbmu-(\bmu+\zb)^\top(\hbSigma+\Qb)^{-1}\zb.\label{eq:unknown_decomposition1}
\end{align}
We investigate the term $\bmu^\top(\hbSigma+\Qb)^{-1}\zb$ first. It is easy to see that
\begin{align*}
  \bmu^\top(\hbSigma+\Qb)^{-1}\zb=  \bmu^\top\bSigma^{-\frac{1}{2}}(\bSigma^{-\frac{1}{2}}\hbSigma\bSigma^{-\frac{1}{2}}+\bSigma^{-\frac{1}{2}}\Qb\bSigma^{-\frac{1}{2}})^{-1}\bSigma^{-\frac{1}{2}}\zb
\end{align*}
By the bound of  $\|\bSigma^{-\frac{1}{2}}\bmu\|_2$, it is easy to see that $\|(\bSigma^{-\frac{1}{2}}\hbSigma\bSigma^{-\frac{1}{2}}+\bSigma^{-\frac{1}{2}}\Qb\bSigma^{-\frac{1}{2}})^{-1}\bSigma^{-\frac{1}{2}}\bmu\|_2$ is bounded. Moreover, $\sqrt{n}\cdot\bSigma^{-\frac{1}{2}}\zb$ is standard normal distribution. 
By Borel-Cantelli lemma, it is easy to see that  
\begin{align*}
|\bmu^\top(\hbSigma+\Qb)^{-1}\zb|~\cvas~0.
\end{align*} 
Hence combined the equation above with \eqref{eq:unknown_decomposition1} we have that
\begin{align}
 |\hbmu^\top(\hbSigma+\Qb)^{-1}\bmu-\hbmu^\top(\hbSigma+\Qb)^{-1}\hbmu+\zb^\top  (\hbSigma+\Qb)^{-1}\zb   |~\cvas~0.\label{eq:unknown_decomposition2}
\end{align}
By the concentration inequalities, it is also  easy to see that
\begin{align*}
 \bigg|\zb^\top  (\hbSigma+\Qb)^{-1}\zb-\frac{1}{n}\tr(\hbSigma+\Qb)^{-1}\bSigma  \bigg| ~\cvas~0. 
\end{align*}
\eqref{eq:unknown_decomposition2} can be further transfered to 
\begin{align}
     |\hbmu^\top(\hbSigma+\Qb)^{-1}\bmu-\hbmu^\top(\hbSigma+\Qb)^{-1}\hbmu+\frac{1}{n}\tr(\hbSigma+\Qb)^{-1}\bSigma|~\cvas~0.\label{eq:unknown_decomposition3}
\end{align}
By Proposition~\ref{prop:converge_Tn1Tn2}, it can be easily seen that
\begin{align*}
    &\frac{1}{p}\tr(\hbSigma+\Qb)^{-1}\hbSigma-\frac{1}{p}\tr\bigg(\frac{\bSigma}{1+s_0}+\Qb \bigg)^{-1}\frac{\bSigma}{1+s_0}~\cvas~0,\\
    &\frac{1}{p}\tr(\hbSigma+\Qb)^{-1}\bSigma-\frac{1}{p}\tr\bigg(\frac{\bSigma}{1+s_0}+\Qb \bigg)^{-1}\bSigma~\cvas~0.
\end{align*}
Hence by $1+s_0$ is bounded, we can easily get that
\begin{align*}
    \frac{1+s_0}{p}\tr(\hbSigma+\Qb)^{-1}\hbSigma-\frac{1}{p}\tr(\hbSigma+\Qb)^{-1}\bSigma~\cvas~0.
\end{align*}
Combined the equation above with the results in Proposition~\ref{prop:converge_s} that 
\begin{align*}
    \frac{1}{1-\frac{c}{p}\tr(\hbSigma+\Qb)^{-1}\hbSigma}\cdot\frac{1}{1+s_0}~\cvas~1,
\end{align*}
we have
\begin{align}
    \bigg|\frac{1}{n}\tr(\hbSigma+\Qb)^{-1}\bSigma-\frac{\tr(\hbSigma+\Qb)^{-1}\hbSigma}{n-\tr(\hbSigma+\Qb)^{-1}\hbSigma}\bigg|~\cvas~0\label{eq:unknown_decomposition4}
\end{align}
after  simple algebra calculation. 
Combining \eqref{eq:unknown_decomposition3} and \eqref{eq:unknown_decomposition4} gives the proof of the first conclusion in Theorem~\ref{thm:unknown_mu}. The second conclusion in Theorem~\ref{thm:unknown_mu} is a direct extension of the previous results we analyzed and we hence omit the proof.

As for the convergence of $\hat{SR}(\Qb)/SR(\Qb)$, we can easily see that
\begin{align*}
\hbmu^\top(\hbSigma+\Qb)^{-1}\bmu-\bmu^\top(\hbSigma+\Qb)^{-1}\bmu~\cvas~ 0.
\end{align*}
The lower bound condition of $\bmu^\top\big(\frac{\bSigma}{1+s_0}+\Qb\big)^{-1}\bmu$ indicates that $\hbmu^\top(\hbSigma+\Qb)^{-1}\bmu>0$ and is lower bounded a constant, therefore
\begin{align*}
    \frac{\hbmu^\top(\hbSigma+\Qb)^{-1}\hbmu-\frac{\tr(\hbSigma+\Qb)^{-1}\hbSigma}{n-\tr(\hbSigma+\Qb)^{-1}\hbSigma}}{\hbmu^\top(\hbSigma+\Qb)^{-1}\bmu}~\cvas~1.
\end{align*}
Combing the equation above with the second conclusion we have that
\begin{align*}
    \frac{\hat{SR}(\Qb)}{SR(\Qb)}~\cvas~1,
\end{align*}
which completes the proof of Theorem~\ref{thm:unknown_mu}.

\subsection{Discussion on the maximal Sharpe ratio}
\label{sec:unknown_SRmax}
In this  section, we consider the property of $SR(\Qb)$ with $\bmu$ unknown, and see why $\Qb=C\bSigma$ with a sufficient large $C$ will no longer approach the maximal Sharpe ratio $SR_{\max}=\sqrt{\bmu^\top\bSigma^{-1}\bmu}$. Recall the definition of $SR(\Qb)$ with unknown $\bmu$, 
\begin{align*}
SR(\Qb)=\frac{\hbmu^\top(\hbSigma+\Qb)^{-1}\bmu}{\sqrt{\hbmu^\top(\hbSigma+\Qb)^{-1}\bSigma(\hbSigma+\Qb)^{-1}\hbmu}}. 
\end{align*}
Similar to the proof in Section~\ref{sec:discussion_largest_SR}, with sufficiently large $C>0$, it can be easily seen that the nominator of $SR(\Qb)$ satisfies
\begin{align*}
    \hbmu^\top(\hbSigma+C\bSigma)^{-1}\bmu&=\bmu^\top(\hbSigma+C\bSigma)^{-1}\bmu+\zb^\top(\hbSigma+C\bSigma)^{-1}\bmu\\
    &=\bmu^\top(\hbSigma+C\bSigma)^{-1}\bmu+o_p(1)\\
    &\approx C^{-1}\bmu^\top\bSigma^{-1}\bmu.
\end{align*}
Here, $\zb=\hbmu-\bmu$ and the approximation comes from the fact that when $C$ is large enough, $s_0$ will be small enough. As for the denominator of $SR(\Qb)$, we have
\begin{align*}
\hbmu^\top(\hbSigma+\Qb)^{-1}\bSigma(\hbSigma+\Qb)^{-1}\hbmu&=\bmu^\top(\hbSigma+\Qb)^{-1}\bSigma(\hbSigma+\Qb)^{-1}\bmu+\zb^\top(\hbSigma+\Qb)^{-1}\bSigma(\hbSigma+\Qb)^{-1}\zb+o_p(1).
\end{align*}
When $\Qb=C\bSigma$, it can be easily show that $s_0$ and $s_{1,\bSigma}$ are also sufficiently small, hence from Section~\ref{sec:appendix_proof_thm_bound} we have
\begin{align*}
\hbmu^\top(\hbSigma+\Qb)^{-1}\bSigma(\hbSigma+\Qb)^{-1}\hbmu&\approx \frac{1}{C^2} \cdot \bigg(\bmu^\top\bSigma^{-1}\bmu+\frac{1}{n}\cdot\tr(\Ib_p) \bigg)\\
&\approx\frac{1}{C^2} \cdot \bigg(\bmu^\top\bSigma^{-1}\bmu+c \bigg).
\end{align*}
We can see that when $\bmu$ is unknown, $SR(C\bSigma)$ with sufficient large $C$ will approximate the value $\frac{\bmu^\top\bSigma^{-1}\bmu}{\sqrt{\bmu^\top\bSigma^{-1}\bmu+c}}$.

\section{An Example for Section~\ref{sec:frontier}}
\label{sec:examplecase_frontier}
In this section, we present an  example  with the assumption $\rb=a_1\cdot \one+a_2\cdot \bxi+a_3\cdot \rb_0$ where $\one$, $\bxi$ and $\rb_0$ are orthogonal to each other. Without loss of generality, we assume that $\|\bxi\|_2=\|\rb_0\|_2=1$, and the matrix  $\bSigma/(1+s_0)+\Qb$  has eigenvalue decomposition as 
\begin{align*}
    \frac{\bSigma}{1+s_0}+\Qb=\lambda_1\bxi\bxi^\top+\textbf{Else}.
\end{align*}
Here, the assumption of $\Qb$ that $\lambda_{\min}(\Qb)\geq c'$ for some constant $c'>0$ ensures $\lambda_1=\Omega(1)$, as $\bxi$ is the eigenvector of $\frac{\bSigma}{1+s_0}+\Qb$. It is clear that $\bxi$ represents a distant factor in the matrix if $\lambda_1$ tends to infinity. $\big(\frac{\bSigma}{1+s_0}+\Qb\big)^{-1}$ has the form
\begin{align*}
\bigg(\frac{\bSigma}{1+s_0}+\Qb\bigg)^{-1}=\frac{1}{\lambda_1}\bxi\bxi^\top+\bOmega,
\end{align*}
where $\bOmega$ is the remaining matrix with $\text{rank}(\bOmega)=p-1 $ and $\|\bOmega\|_{\op}=O(1)$. 
By the portfolio optimization \eqref{eq:frontier1}, when $\rb=a_1\one+a_2\cdot \bxi+a_3\cdot \rb_0$, the optimization constraint becomes 
\begin{align*}
    \wb^{\top}(a_2\cdot \bxi+a_3\cdot \rb_0)=\mu_0-a_1,\quad \wb^\top\one=1.
\end{align*}
With a little abuse of notation, we can define $\rb=a_2\cdot \bxi+a_3\cdot \rb_0$ and the constraint constant $\mu_0$ can be changed to $\mu_0-a_1$.

The vector $\one$ could represent a rough market factor, since if the market goes up all stocks will go up accordingly. 
The vector $\bxi$ represents some style factors that influence the returns of the assets. These factors could be things like exposures to Fama-French size or value factors. The vector $\rb_0$  represents the expected residual return, which could reflect asset-specific returns not explained by factors. 
We assume the true expected return is decomposed into the market component, the style factor component and the residual component with coefficients $a_1$, $a_2$ and $a_3$

For this specific example, we calculate the values of $\cA_{r1}$, $\cA_{rr}$, and $\cA_{11}$ in Table~\ref{table:1}. Note that ratio consistency in Theorem~\ref{thm:frontier} holds as long as $\frac{\cA_{r1}^2}{\cA_{11} \cA_{rr}}$ remains smaller than some constant strictly less than $1$. From Table~\ref{table:1}, this condition is equivalent to requiring that $\frac{a_3^2 (\rb_0^\top \bOmega \one)^2}{a_3^2 (\rb_0^\top \bOmega \rb_0)(\one^\top \bOmega \one) + \frac{a_2^2 (\one^\top \bOmega \one)}{\lambda_1}}$ is smaller than a constant less than $1$. 
First note that if we place more weight on $\bxi$ (increasing $a_2$), we improves the chances of satisfying the ratio consistency condition. Even in the worse case of $a_2 = 0$, ratio consistency still holds as long as $\frac{(\rb_0^\top \bOmega \one)^2}{(\rb_0^\top \bOmega \rb_0)(\one^\top \bOmega \one)}$ is smaller than a constant less than $1$, which can be proved under some mild assumptions we discuss next.

\begin{table}[t]
\centering
\caption{ Values of $\cA_{r1}$, $\cA_{rr}$ and $\cA_{11}$, and conditions for the convergence.}
\label{table:1}
\scalebox{0.9}{\begin{tabular}{cccc}
\hline
$\cA_{r 1}$ & $\cA_{1 1}$ & $\cA_{r r}$             & $\frac{\cA_{r1}^2}{\cA_{rr}\cA_{11}}$   \\ 
$a_3\rb_0^\top\bOmega \one$         & $\one^\top\bOmega\one$ & $a_3^2\rb_0^\top\bOmega\rb_0+\frac{a_2^2}{\lambda_1}$ &   $\frac{a_3^2(\rb_0^\top\bOmega \one)^2}{a_3^2(\rb_0^\top\bOmega\rb_0)(\one^\top\bOmega\one)+\frac{a_2^2(\one^\top\bOmega\one)}{\lambda_1}}$  \\ \hline
& $\frac{\hat{\sigma}^2}{\sigma_0^2}~\cvas~1$ & $\frac{\mu_0-r_0}{\sigma_0}-\frac{\mu_0-r_0}{\hat{\sigma}}~\cvas~0$ & $\hat{\sigma}^2-\sigma_0^2~\cvas~0$    \\   
Conditions &    $\frac{\cA_{r1}^2}{\cA_{rr}\cA_{11}}\leq \rho<1$
& $a_3^2\rb_0^\top\bOmega\rb_0+\frac{a_2^2}{\lambda_1}\leq C$                & $\mu_0\leq  a_1+C\sqrt{a_3^2\rb_0^\top\bOmega\rb_0+\frac{a_2^2}{\lambda_1}}$\\ \hline
\end{tabular}}
\end{table}

From the expression of $\frac{\bSigma}{1 + s_0} + \Qb$, we assume $\bxi$ is one of its eigenvectors. So $\bOmega$ must be orthogonal to $\bxi$. We first consider the case when $\rb_0$ or $\one$ is an eigenvector of $\bOmega$, hence we have $\rb_0^\top\bOmega\one=0$ due to orthogonality of $\rb_0$ and $\one$, therefore $\rho=0<1$ satisfies the condition. Next we consider the case where $\bxi,\one,\rb_0$ are not eigenvectors of $\bOmega$ and we further assume that the non-zero eigenvalues of $\bOmega$ are bounded away from $0$. 
Under this mild condition, we can prove that there exists $\rho < 1$ such that $\cA_{r1}^2 / (\cA_{11} \cA_{rr}) \leq \rho < 1$.

\noindent\underline{\textbf{A simple proof:}} When the minimum non-zero eigenvalue of $\bOmega$ is bounded away from $0$, if there does not exist such constant $\rho$, then we can see $1-\frac{(\rb_0^\top \bOmega\one)^2}{\rb_0^\top \bOmega\rb_0 \one^\top\bOmega\one}\to 0$.  Hence, we have $\frac{\bOmega^{\frac{1}{2}}\rb_0}{\|\bOmega^{\frac{1}{2}}\rb_0\|_2}$ and $\frac{\bOmega^{\frac{1}{2}}\one}{\|\bOmega^{\frac{1}{2}}\one\|_2}$ will asymptotically fall into the same direction. Note that $\Big\|\frac{\bOmega^{\frac{1}{2}}\rb_0}{\|\bOmega^{\frac{1}{2}}\rb_0\|_2}\Big\|_2=\Big\|\frac{\bOmega^{\frac{1}{2}}\one}{\|\bOmega^{\frac{1}{2}}\one\|_2}\Big\|_2=1$, we have $\Big\|\frac{\bOmega^{\frac{1}{2}}\rb_0}{\|\bOmega^{\frac{1}{2}}\rb_0\|_2}-\frac{\bOmega^{\frac{1}{2}}\one}{\|\bOmega^{\frac{1}{2}}\one\|_2}\Big\|_2\to 0$. This means $\frac{\rb_0}{\|\bOmega^{\frac{1}{2}}\rb_0\|_2}-\frac{\one}{\|\bOmega^{\frac{1}{2}}\one\|_2}$ will asymptotically fall into the space generated by $\{\bxi\}$. Here, the bounded $\|\bOmega\|_{\op}$ ensures the norm of the vector $\frac{\rb_0}{\|\bOmega^{\frac{1}{2}}\rb_0\|_2}-\frac{\one}{\|\bOmega^{\frac{1}{2}}\one\|_2}$ will not tend to $0$, and by our assumption that the minimum non-zero eigenvalue of $\bOmega$ is bounded away from $0$, this vector will asymptotically fall into the space generated by $\{\bxi\}$ due to $\Big\|\frac{\bOmega^{\frac{1}{2}}\rb_0}{\|\bOmega^{\frac{1}{2}}\rb_0\|_2}-\frac{\bOmega^{\frac{1}{2}}\one}{\|\bOmega^{\frac{1}{2}}\one\|_2}\Big\|_2\to 0$. This violates our assumption that $\rb_0,\one\perp \bxi$. Hence the constant $\rho<1$ exists.



We turn to the condition for Sharpe difference consistency. The Sharpe difference consistency holds  when $\cA_{rr}$ is bounded, which requires $a_3^2\rb_0^\top\bOmega\rb_0+\frac{a_2^2}{\lambda_1}\leq C$. This condition is achievable under reasonable scenarios. Note again $\|\bOmega\|_{\op}$ is bounded.  For instance, we can set $a_2 = O(\sqrt{\lambda_1})$, and the condition will hold provided that $\|\rb_0\|_2$ is bounded and $a_3$ is at a constant level. As for the absolute error of volatility, the target return $\mu_0$ must be set within an appropriate range to ensure accurate estimation. From Theorem~\ref{thm:frontier}, this condition requires $\mu_0 - a_1 \leq C \sqrt{\cA_{rr}}$, which implies $\mu_0\leq  a_1+C\sqrt{a_3^2\rb_0^\top\bOmega\rb_0+\frac{a_2^2}{\lambda_1} }$.

\section{Optimization Over Regularization Matrix}
\label{sec:optimalQ}
In this section, we investigate the optimization over regularization matrix $\Qb$. Section~\ref{sec:comparsion_optimal} shows that despite the fully optimal $\hat{\Qb}$ perform well in $\hat{SR}(\hat{\Qb})$, the performance of $SR(\hat{\Qb})$ is poor. A natural question arises that when will the optimal $\hat{SR}(\hat{\Qb})$ still have good consistency with $SR(\hat{\Qb})$.

To facilitate the analysis, we introduce the following assumptions, which ensures that the empirical and population Sharpe ratios will have good properties. 

\begin{assumption}
\label{assump:optimalQ_1}
There exists universal constants $l,L>0$ such that for all \(\Qb \in \mathcal{Q}\), both \(SR(\Qb)\) and \(\hat{SR}(\Qb)\) satisfy
\(
l \leq SR(\Qb),\ \hat{SR}(\Qb) \leq L
\)
almost surely for all $n$ large enough.
\end{assumption}


Assumption~\ref{assump:optimalQ_1} requires that both \(SR(\Qb)\) and \(\hat{SR}(\Qb)\) are uniformly bounded by universal constants almost surely for all sufficiently large $n$. This is a mild assumption, as the Sharpe ratio of a proper portfolio in empirical studies typically stays at constant level. 
We give another assumption on the candidate set $\cQ$.

\begin{assumption}\label{assump:optimalQ_2}
There exists a sequence of  bijections \( \phi_n: \cB \to \cQ \), where \( \cB \subset \mathbb{R}^k \) is a fixed compact set (independent of $n$) for some constant $ k > 0$.  Furthermore, the sequence $\{ \phi_n \}$ is equicontinuous with respect to the operator norm: for any $\varepsilon > 0$, there exists $\delta > 0$ (independent of $n$) such that for all $n$ and all \( \balpha, \balpha' \in \cB \),
\[
\|\balpha - \balpha'\|_2 \leq \delta \quad \implies \quad \|\phi_n(\balpha) - \phi_n(\balpha')\|_{\op} \leq \varepsilon.
\]
\end{assumption}
Assumption~\ref{assump:optimalQ_2} assumes the candidate set \(\cQ\) to be parameterized by a finite vector \(\balpha=(\alpha_1, \dots, \alpha_k)\) with \(\balpha\) confined to a compact set. This means that the entire family of candidate matrices can be represented as a continuous mapping $\balpha \mapsto \phi_n(\balpha)$, 
and that the operator norm of \(\phi_n(\balpha)\) is equicontinuous with respect to \(\balpha\). Such a condition is common when \(\cQ\) is set to be a well-behaved, finite-dimensional class of matrices.
A specific example is $\phi_n(\balpha)=\alpha_1 \Qb_1 + \alpha_2 \Qb_2 + \cdots + \alpha_k \Qb_k$
where \(\Qb_1, \dots, \Qb_k\) are predetermined matrices  with $\|\Qb_j\|_{\op}$ bounded for all $j$, $\Qb_j$ are linearly independent,  and the coefficients \(\balpha\) vary over a compact set in $\mathbb{R}^k$.

\begin{lemma}
\label{lemma:equi_continuity}
Under the conditions of Theorem~\ref{thm:main_theorem}, given additional Assumption~\ref{assump:optimalQ_1}. For any $\varepsilon > 0$, there exists  $\delta > 0$ only related with $\varepsilon$ (independent of $n$)  such that for any $\Qb_1,\Qb_2\in \cQ$ with $\|\Qb_1 - \Qb_2\|_{\op} \le \delta$, 
\[
\big|\hat{SR}(\Qb_1) - \hat{SR}(\Qb_2)\big| \le \varepsilon \quad \text{and} \quad \big|SR(\Qb_1) - SR(\Qb_2)\big| \le \varepsilon
\]
almost surely for all $n$ large enough.
\end{lemma}

Lemma~\ref{lemma:equi_continuity} proves that Assumption~\ref{assump:optimalQ_1} leads to equicontinuity of the true and estimated Sharpe ratios with respect to $\Qb$.
Based on this lemma, we have the following proposition.

\begin{proposition}
\label{prop:optimal_Q}
Under the conditions of Theorem~\ref{thm:main_theorem}, additionally suppose that Assumptions~\ref{assump:optimalQ_1} and \ref{assump:optimalQ_2} hold. Define 
\begin{align*}
    \hat{\Qb}=\argmax_{\Qb\in \cQ}\hat{SR}(\Qb).
\end{align*}
It holds that
\begin{align*}
\hat{SR}(\hat{\Qb})/SR(\hat{\Qb})~\cvas~1,\quad \hat{SR}(\hat{\Qb})-SR(\hat{\Qb})~\cvas~0.
\end{align*}
\end{proposition}
From Proposition~\ref{prop:optimal_Q}, we theoretically prove that when \(\Qb\) is optimized over a candidate set that is parameterized by a finite number of parameters (i.e., under Assumption~\ref{assump:optimalQ_2}), the in-sample estimator \(\hat{SR}(\hat{\Qb})\) is asymptotically equivalent to the true performance measure \(SR(\hat{\Qb})\). In particular, we show that $\frac{\hat{SR}(\hat{\Qb})}{SR(\hat{\Qb})}~\cvas~1$, 
which implies that any discrepancy between the in-sample and out-of-sample Sharpe ratios vanishes asymptotically. This result confirms that, when the search space for \(\Qb\) is well-behaved and restricted to a finite-dimensional family, the overfitting issue can be controlled, and the optimized candidate achieves consistent performance in the large-sample limit.

\subsection{Proof of Lemma~\ref{lemma:equi_continuity}}
Recall the definition of $SR(\Qb)$ and $\hat{SR}$ that 
\begin{align*}
&\hat{SR}(\Qb)=\left(1-\frac{c}{p}\operatorname{tr}\left[\hbSigma(\hbSigma+\Qb)^{-1}\right]\right)
\cdot \frac{\bmu^\top(\hbSigma+\Qb)^{-1}\bmu}{\sqrt{\bmu^\top(\hbSigma+\Qb)^{-1}\hbSigma(\hbSigma+\Qb)^{-1}\bmu}},\\
&{SR}(\Qb)= \frac{\bmu^\top(\hbSigma+\Qb)^{-1}\bmu}{\sqrt{\bmu^\top(\hbSigma+\Qb)^{-1}\bSigma(\hbSigma+\Qb)^{-1}\bmu}}. 
\end{align*}
From Assumption~\ref{ass:assump2} and \ref{ass:assump3}, we have there exists some constant $c>0$ such that $\Qb\geq c\Ib$ (we allow the specific case $\Qb=\0$ when $c<1$). By the expression of $1-\frac{c}{p}\operatorname{tr}[\hbSigma(\hbSigma+\Qb)^{-1}]$, when $c<1$, it is bounded away from $0$ by universal constant $(1-c)$. When $c\geq 1$, note that the number of spikes of $\hbSigma$ is fixed, we can easily conclude that $1-\frac{c}{p}\operatorname{tr}[\hbSigma(\hbSigma+\Qb)^{-1}]$ is also bounded away from $0$ by some universal constant. Note that $\hat{SR}(\Qb)$ is assumed to be bounded by $L$ almost surely for $n$ sufficiently large,  we conclude that there exists universal constant $L_1,L_2$ such that 
\begin{align}
0<L_1\leq 1-\frac{c}{p}\operatorname{tr}[\hbSigma(\hbSigma+\Qb)^{-1}]\leq 1,\quad 0<L_2\leq \frac{\bmu^\top(\hbSigma+\Qb)^{-1}\bmu}{\sqrt{\bmu^\top(\hbSigma+\Qb)^{-1}\hbSigma(\hbSigma+\Qb)^{-1}\bmu}}\leq L/L_1 \label{eq:bound_condition}
\end{align}
almost surely for all $n$ large enough. We next prove that $\|\rmd \hat{SR}(\Qb)/\rmd \Qb\|_{\op}$ is bounded almost surely for all $n$ large enough, and similar proof can be performed to $\|\rmd SR(\Qb)/\rmd \Qb\|_{\op}$. If the conclusions above hold, then it can be seen that $\hat{SR}(\Qb)$ and $SR(\Qb)$ are Lipchitz with some universal constant almost surely for all $n$ large enough. Then the equicontinuity will hold automatically, which completes the proof. It only remains for us to show $\|\rmd \hat{SR}(\Qb)/\rmd \Qb\|_{\op}$ and $\|\rmd SR(\Qb)/\rmd \Qb\|_{\op}$ are bounded by some universal constant almost surely for all $n$ large enough.

Define the terms
\begin{align*}
&A(\Qb)=1-\frac{c}{p}\operatorname{tr}[\hbSigma(\hbSigma+\Qb)^{-1}],\quad B(\Qb)=\bmu^\top(\hbSigma+\Qb)^{-1}\bmu,\\
&C(\Qb)=\sqrt{D(\Qb)},\quad D(\Qb)=\bmu^\top(\hbSigma+\Qb)^{-1}\hbSigma(\hbSigma+\Qb)^{-1}\bmu.
\end{align*}
It is clear that $\hat{SR}(\Qb)=A(\Qb)\cdot \frac{B(\Qb)}{C(\Qb)}$. By matrix derivatives, we have
\begin{align}
\label{eq:derivative_all_terms}
\begin{aligned}
    &\nabla_{\Qb} A(\Qb)=\frac{c}{p}(\hbSigma+\Qb)^{-1}\hbSigma(\hbSigma+\Qb)^{-1},\\
    &\nabla_{\Qb} B(\Qb)=-(\hbSigma+\Qb)^{-1}\bmu\bmu^\top(\hbSigma+\Qb)^{-1},\\
    &\nabla_{\Qb} C(\Qb)=-\frac1{2C(\Qb)}(\hbSigma+\Qb)^{-1}\Bigl(
\hbSigma (\hbSigma+\Qb)^{-1}\bmu\bmu^\top 
+\bmu\bmu^\top (\hbSigma+\Qb)^{-1}\hbSigma 
\Bigr)(\hbSigma+\Qb)^{-1}.
\end{aligned}
\end{align}
We calculate the gradient $\nabla_{\Qb}\hat{SR}(\Qb)$. We have
\begin{align*}
\nabla_{\Qb}\hat{SR}(\Qb)= \nabla_{\Qb} A(\Qb)\cdot \frac{B(\Qb)}{C(\Qb)}+A(\Qb)\cdot \nabla_{\Qb}\frac{B(\Qb)}{C(\Qb)}.
\end{align*}
From \eqref{eq:bound_condition} and \eqref{eq:derivative_all_terms}, we easily have that $\|\nabla_{\Qb} A(\Qb)\|_{\op}$, $|B(\Qb)/C(\Qb)| $ and $|A(\Qb)|$ is bounded by some universal constant almost surely for all $n$ large enough. To prove $\|\nabla_{\Qb}\hat{SR}(\Qb)\|_{\op}$ is bounded, it remains for us to prove that $\big\|\nabla_{\Qb}\frac{B(\Qb)}{C(\Qb)} \big\|_{\op}$ is bounded. We have 
\begin{align*}
\nabla_{\Qb}\frac{B(\Qb)}{C(\Qb)} =\frac{\nabla_{\Qb} B(\Qb)}{C(\Qb)}-\frac{B(\Qb)\nabla_{\Qb} C(\Qb)}{C^2(\Qb)}.
\end{align*}
By \eqref{eq:derivative_all_terms}, easy to see that 
\begin{align*}
    \bigg\|\frac{\nabla_{\Qb} B(\Qb)}{C(\Qb)} \bigg\|_{\op}&=\bigg\|\frac{(\hbSigma+\Qb)^{-1}\bmu\bmu^\top(\hbSigma+\Qb)^{-1}}{C(\Qb)} \bigg\|_{\op}\\
    &\leq \frac{1}{\lambda_{\min}(\Qb)}\frac{B(\Qb)}{C(\Qb)}\\
    &\leq \frac{L}{L_1\lambda_{\min}(\Qb)},
\end{align*}
where the first inequality $\|(\hbSigma+\Qb)^{-1}\bmu\bmu^\top(\hbSigma+\Qb)^{-1}\|_{\op}=\bmu^\top (\hbSigma+\Qb)^{-2}\bmu$, and the second inequality is by \eqref{eq:bound_condition}. Note also that $\big\|\frac{B(\Qb)\nabla_{\Qb} C(\Qb)}{C^2(\Qb)}\big\|_{\op}\leq \frac{L}{L_1}\big\|\frac{\nabla_{\Qb} C(\Qb)}{C(\Qb)}\big\|_{\op}  $, it only remains for us to investigate $\big\|\frac{\nabla_{\Qb} C(\Qb)}{C(\Qb)}\big\|_{\op}$. From \eqref{eq:derivative_all_terms}, we have
\begin{align*}
   \bigg\|\frac{\nabla_{\Qb} C(\Qb)}{C(\Qb)}\bigg\|_{\op}& \leq \frac{\|(\hbSigma+\Qb)^{-1} \hbSigma (\hbSigma+\Qb)^{-1}\bmu\bmu^\top(\hbSigma+\Qb)^{-1}\|_{\op}}{C^2(\Qb)}\\
   &\leq \bigg\|\frac{(\hbSigma+\Qb)^{-1}\bmu\bmu^\top(\hbSigma+\Qb)^{-1}}{C(\Qb)} \bigg\|_{\op}\cdot \frac{1}{C(\Qb)}.
\end{align*}
Here, the second inequality is by $\|(\hbSigma+\Qb)^{-1} \hbSigma\|_{\op}\leq 1$. Same to the proof above, it suffices to prove that $C(\Qb)$ is bounded away by $0$ by some universal constant. From the assumption that $l\leq \hat{SR}(\Qb)\leq L$, and note that $A(\Qb)$ is lower bounded by $L_1$ and upper bounded by $1$. It suffices to prove that $\bmu^\top (\hbSigma+\Qb)^{-1}\bmu$ is bounded away from $0$. We prove $\bmu^\top(\hbSigma+\Qb)^{-1}\bmu$  bounded away from $0$ from the event  $l\leq SR(\Qb)\leq L$. Note that
\begin{align*}
    (\hbSigma+\Qb)^{-1}\bSigma (\hbSigma+\Qb)^{-1}=\bSigma^{-\frac{1}{2}}\bigg(\frac{\Zb^\top\Zb}{n}+\bSigma^{-\frac{1}{2}}\Qb \bSigma^{-\frac{1}{2}}\bigg)^{-2}\bSigma^{-\frac{1}{2}}\geq c'\bSigma^{-1}
\end{align*}
for some constant $c'>0$. Here, we apply the Bai-Yin Theorem that $\|\frac{\Zb^\top\Zb}{n}\|_{\op}\leq 2(1+\sqrt{c})^2$ almost surely for all $n$ large enough, and $\|\bSigma^{-\frac{1}{2}}\Qb\bSigma^{-\frac{1}{2}}\|_{\op}$ is bounded. We have
\begin{align*}
l\leq SR(\Qb)=\frac{\bmu^\top(\hbSigma+\Qb)^{-1}\bmu}{\sqrt{\bmu^\top(\hbSigma+\Qb)^{-1}\bSigma(\hbSigma+\Qb)^{-1}\bmu}}\leq \bmu^\top(\hbSigma+\Qb)^{-1}\bmu/\sqrt{(c'\bmu^{\top}\bSigma^{-1}\bmu)}.
\end{align*}
Note that $\sqrt{\bmu^{\top}\bSigma^{-1}\bmu}\geq SR(\Qb)\geq l$ is bounded away from $0$, we can easily have that $\bmu^\top (\hbSigma+\Qb)^{-1}\bmu$ is bounded away from $0$. From the discussion above,  we hence prove that $C(\Qb)$ is bounded away from $0$ by some constant, and the constant does not depend on $n$. Therefore,
\begin{align*}
    \bigg\|\frac{\nabla_{\Qb} C(\Qb)}{C(\Qb)}\bigg\|_{\op} \leq C
\end{align*}
for some universal constant $C>0$, which indicates that $\|\rmd \hat{SR}(\Qb)/\rmd \Qb\|_{\op}$ is bounded by some universal constant almost surely for all $n$ large enough. We complete the proof that $\hat{SR}(\Qb)$ is equicontinuous. The same proof can also show the equicontinuity of $SR(\Qb)$, which completes the proof.

\subsection{Proof of Proposition~\ref{prop:optimal_Q}}
We apply $\varepsilon-\delta$ language to prove the ratio and difference consistency. For simplicity, we use $\Qb(\balpha)$ to exactly represent the mapping $\phi_n(\balpha)$. For any given fixed $\varepsilon>0$, from Assumption~ \ref{assump:optimalQ_2} and Lemma~\ref{lemma:equi_continuity}, there exists $\delta>0$ only related with $\varepsilon$, such that as long as $\|\balpha_1-\balpha_2\|_2\leq\delta$, 
\begin{align}
    \big|\hat{SR}(\Qb(\balpha_1)) - \hat{SR}(\Qb(\balpha_2))\big| \le \frac{\varepsilon}{2} \quad \text{and} \quad \big|SR(\Qb(\balpha_1)) - SR(\Qb(\balpha_2))\big| \leq\frac{\varepsilon}{2} \label{eq:optimal_equi-property}
\end{align}
almost surely for all $n$ large enough. 
Since the parameter \(\balpha\) is confined to a fixed compact set, this set has a finite covering \(\mathcal{B}(\delta)=\{B(\balpha^{(1)}), \dots, B(\balpha^{(M)})\}\) where each \(B(\balpha^{(j)})\) is an open ball of radius \(\delta\) centered at \(\balpha^{(j)}\). Here $M(\varepsilon)>0$ is a fixed constant only related with $\varepsilon$ and the compact set $\cB$ due to the assumption that $\cB$ is the fixed compact set. For any candidate \(\balpha\) in the parameter space, there exists at least one ball \(B(\balpha^{(j)})\) such that $\balpha\in B(\balpha^{(j)})$ (\(\|\balpha-\balpha^{(j)}\|_2\le \delta\)). 

For any given \(\varepsilon>0\), the parameter values \(\{\balpha^{(j)}\}_{j=1}^M\) and the corresponding matrices \(\{\Qb(\balpha^{(j)})\}_{j=1}^M\) are fixed. Hence, by Theorem~\ref{thm:main_theorem}, we have almost surely
\begin{equation}
\label{eq:optimal_SR-hatSR}
\sup_{j\in[M]} \lim_{n\to\infty} \big|\hat{SR}(\Qb(\balpha^{(j)})) - SR(\Qb(\balpha^{(j)}))\big| = 0.
\end{equation}
For any optimal \(\hat{\Qb}\), there exists a corresponding \(\hat{\balpha}_n\) such that \(\hat{\Qb} = \Qb(\hat{\balpha}_n)\). By the definition of the finite covering \(\cB(\delta)=\{B(\balpha^{(1)}),\ldots,B(\balpha^{(M)})\}\), there exists a sequence \(\{\balpha_n\}\) with \(\balpha_n\in\{\balpha^{(1)},\ldots,\balpha^{(M)}\}\) satisfying
\[
\|\hat{\balpha}_n - \balpha_n\|_2\le \delta.
\]
Then, by the equicontinuity property in \eqref{eq:optimal_equi-property}, we have that
\begin{align}
\big|\hat{SR}(\Qb(\hat{\balpha}_n)) - \hat{SR}(\Qb(\balpha_n))\big| \le \varepsilon/2,\quad \big|SR(\Qb(\hat{\balpha}_n)) - SR(\Qb(\balpha_n))\big| \le \varepsilon/2. \label{eq:optimal_group_diff}
\end{align} 
\eqref{eq:optimal_SR-hatSR} and \eqref{eq:optimal_group_diff} tell us that for any given fixed $\varepsilon>0$, almost surely
\begin{align*}
\lim_{n\to\infty}  \big|\hat{SR}(\Qb(\hat{\balpha}_n))&-{SR}(\Qb(\hat{\balpha}_n))\big|\leq \lim_{n\to\infty}\big|\hat{SR}(\Qb(\hat{\balpha}_n))-\hat{SR}(\Qb(\balpha_n))\big|\\
&\qquad+\big|SR(\Qb(\hat{\balpha}_n))-SR(\Qb(\balpha_n))\big|+\big|\hat{SR}(\Qb(\balpha_n)) - SR(\Qb(\balpha_n))\big|\\
&\qquad\leq 0 + \varepsilon/2 + \varepsilon/2=\varepsilon.
\end{align*}
Here, the first inequality is by triangle inequality and the second inequality is by \(\balpha_n\in\{\balpha^{(1)},\ldots,\balpha^{(M)}\}\), \eqref{eq:optimal_SR-hatSR} and \eqref{eq:optimal_group_diff}. Since the inequality holds for any given fixed $\varepsilon>0$ almost surely, we conclude that
\begin{align*}
  \big|\hat{SR}(\Qb(\hat{\balpha}_n))-{SR}(\Qb(\hat{\balpha}_n))\big|~\cvas~0.
\end{align*}
The assumption that $SR(\Qb)$ bounded away from $0$ for any $\Qb\in\cQ$ completes the proof of Proposition~\ref{prop:optimal_Q}.

\end{document}